\documentclass[english]{article}
\usepackage[T1]{fontenc}
\usepackage[latin9]{inputenc}
\usepackage{geometry}
\usepackage{verbatim}
\usepackage{float}
\usepackage{bm}
\usepackage{amsmath}
\usepackage{amssymb}
\usepackage{dsfont}
\usepackage{multirow}

\usepackage{grffile}

\DeclareMathOperator{\ind}{\mathds{1}}  

\usepackage[colorlinks=true,allcolors=blue]{hyperref}%

\makeatletter

\floatstyle{ruled}
\newfloat{algorithm}{tbp}{loa}
\providecommand{\algorithmname}{Algorithm}
\floatname{algorithm}{\protect\algorithmname}

\usepackage{babel}
\usepackage{tabularx,mathtools}

\usepackage{cite}\usepackage{amsthm}\usepackage{dsfont}\usepackage{array}\usepackage{mathrsfs}\usepackage{comment}\onecolumn
\usepackage{bbm}
\usepackage{color}

\allowdisplaybreaks

\usepackage{enumitem}
\setlist[itemize]{leftmargin=2em}
\setlist[enumerate]{leftmargin=2em}

\usepackage{babel}
\usepackage[noend]{algpseudocode}

\usepackage{algorithm}
\usepackage{algorithmicx}
\usepackage{arydshln}

\definecolor{yxc}{RGB}{255,0,0}
\definecolor{cc}{RGB}{125,0,0}
\definecolor{ytw}{RGB}{0,125,0}

\newcommand{\ytw}[1]{\textcolor{ytw}{[YTW: #1]}}

\newcommand{\lambdastar}{\lambda^{\star}}
\newcommand{\ur}{\bm{u}}
\newcommand{\ul}{\bm{w}}
\newcommand{\ustar}{\bm{u}^{\star}}
\newcommand{\ltwo}[1]{\|#1\|_2}
\newcommand{\sigmamin}{\sigma_{\mathrm{min}}}
\newcommand{\sigmamax}{\sigma_{\mathrm{max}}}
\newcommand{\var}{\mathsf{Var}}

\newcommand{\defn}{:=}
\newcommand{\Tone}{T_1}

\newcommand{\real}{\ensuremath{\mathbb{R}}}
\newcommand{\mprob}{\ensuremath{\mathbb{P}}}
\newcommand{\Exs}{\ensuremath{\mathbb{E}}}
\newcommand{\vastar}{v_{\bm{a},l}^{\star}}
\newcommand{\vahat}{\widehat{v}_{\bm{a},l}}

\newcommand{\linf}[1]{\|#1\|_\infty}
\newcommand{\Hhatij}{\widehat{H}_{ij}}
\newcommand{\Hij}{H_{ij}}
\newcommand{\vatil}{\widetilde{v}_{\bm{a},l}}

\makeatother

\usepackage{babel}
\begin{document}

\theoremstyle{plain} \newtheorem{lemma}{\textbf{Lemma}} \newtheorem{prop}{\textbf{Proposition}}\newtheorem{theorem}{\textbf{Theorem}}\setcounter{theorem}{0}
\newtheorem{corollary}{\textbf{Corollary}}\newtheorem{claim}{\textbf{Claim}}
\newtheorem{example}{\textbf{Example}}\newtheorem{notation}{\textbf{Notation}}
\newtheorem{definition}{\textbf{Definition}} \newtheorem{fact}{\textbf{Fact}}\newtheorem{remark}{\textbf{Remark}}\newtheorem{assumption}{\textbf{Assumption}}
\theoremstyle{definition}

\theoremstyle{remark}\newtheorem{condition}{\textbf{Condition}}\newtheorem{conjecture}{\textbf{Conjecture}}
\title{Tackling small eigen-gaps: Fine-grained eigenvector estimation and inference under heteroscedastic noise}
\author{Chen Cheng\thanks{Department of Statistics, Stanford University, Stanford, CA 94305, USA; email: \texttt{chencheng@stanford.edu}.} \\
	\and
	Yuting Wei\thanks{Department of Statistics and Data Science, The Wharton School, University of Pennsylvania, Philadelphia, PA 19104, USA; email: \texttt{ytwei@wharton.upenn.edu}.}\\
	\and
	Yuxin Chen\thanks{Department of Electrical and Computer Engineering, Princeton University, Princeton, NJ 08544, USA; email:
	\texttt{yuxin.chen@princeton.edu}.} 
	}
\date{January 2020; ~~ Revised: September 2021}

\maketitle

\begin{abstract}

This paper aims to address two fundamental challenges arising in eigenvector estimation and inference for a low-rank  matrix from noisy observations: (1) how to 
	estimate an unknown eigenvector when the eigen-gap (i.e.~the spacing between the associated eigenvalue and the rest of the spectrum) is particularly small; (2) how to perform estimation and inference on linear functionals of an eigenvector --- a sort of ``fine-grained'' statistical reasoning that goes far beyond the usual $\ell_2$ analysis.  
	
	We investigate how to address these challenges in a setting where the unknown $n\times n$ matrix is symmetric and the additive noise matrix contains independent (and non-symmetric) entries. Based on eigen-decomposition of the asymmetric data matrix, we propose estimation and uncertainty quantification procedures for an unknown eigenvector, which further allow us to reason about  linear functionals of an unknown eigenvector.  The proposed procedures and the accompanying theory enjoy several important features: (1) distribution-free (i.e.~prior knowledge about the noise distributions is not needed); (2) adaptive to heteroscedastic noise; (3) minimax optimal under Gaussian noise. Along the way, we establish valid procedures to construct confidence intervals for the unknown eigenvalues. All this is guaranteed even in the presence of a small eigen-gap (up to $O(\sqrt{n/\mathrm{poly}\log (n)}\,)$ times smaller than the requirement in prior theory), which goes significantly beyond what generic matrix perturbation theory has to offer.  

\end{abstract}


\tableofcontents


\section{Introduction}

A variety of science and engineering applications ask for spectral analysis of low-rank matrices in high dimension \cite{wainwright2019high}. Setting the stage, imagine that we are interested in a large low-rank matrix
\begin{equation}
	\bm{M}^{\star}= \sum_{l=1}^r \lambda^{\star}_l \bm{u}_l^{\star}\bm{u}_l^{\star\top}\in\mathbb{R}^{n\times n},
	\label{eq:model-rankr}
\end{equation}
where $r$ ($r\ll n$) represents the rank, and $\lambda_l^{\star}$ stands for the $l$th largest eigenvalue of $\bm{M}^{\star}$, with \mbox{$\bm{u}_l^{\star}\in \mathbb{R}^n$} the associated eigenvector. Suppose, however, that we do not have access to perfect measurements about the entries of this matrix; rather, the observations we have available, represented by a data matrix $\bm{M} = \bm{M}^{\star} + \bm{H}$,
are contaminated by a substantial amount of random noise (reflected by the noise matrix $\bm{H}$). The aim is to perform reliable estimation and inference on the unseen eigenvectors of $\bm{M}^{\star}$ on the basis of noisy data.  
Motivated by the abundance of applications (e.g.~collaborative filtering, harmonic retrieval, sensor network localization, joint shape matching \cite{CanPla10,hua1990matrix,chen2014robust,javanmard2013localization,chen2014near}), research on eigenspace estimation in this context has flourished in the past several years, typically built upon proper exploitation of low-rank structures. We have now been equipped with a rich suite of modern statistical theory that delivers statistical performance guarantees  for a number of spectral estimators (e.g.~\cite{peche2006largest,cai2018rate,abbe2017entrywise,o2018random,vu2011singular,yu2015useful,keshavan2010matrix,chi2018nonconvex,wang2015singular,zhong2017eigenvector,chen2018asymmetry,koltchinskii2020efficient,cho2017asymptotic,cai2019subspace,zhong2017eigenvector,zhang2018heteroskedastic, koltchinskii2016asymptotics}); see \cite{chen2020spectral} for a contemporary overview of spectral methods.

\subsection{Motivation and challenges}

Despite a large body of work tackling the above problem, there are several fundamental yet unaddressed challenges that deserve further attention.

\paragraph{Stringent requirements on eigen-gaps.} A crucial identifiability issue stands out when estimating individual eigenvectors. In general, one cannot possibly disentangle the eigenvectors $\bm{u}_1^{\star},\cdots,\bm{u}_r^{\star}$ unless there is sufficient spacing between adjacent eigenvalues. After all, even in the noiseless setting, one can only hope to recover the subspace spanned by $\{\bm{u}_i^{\star}\}$, rather than individual eigenvectors, if all non-zero eigenvalues are identical.

In principle, a ``sufficient'' eigen-gap criterion in the presence of noise is dictated by the noise levels, or more precisely, the signal-to-noise ratios (SNRs). However, generic linear algebra theory typically imposes fairly stringent, and hence pessimistic, eigen-gap requirements for both eigenvector and eigenvalue estimation. More concretely, imagine we wish to estimate the $l$th eigenvector $\bm{u}_l^{\star}$: 
generic matrix perturbation theory (e.g.~the Davis-Kahan $\sin\bm{\Theta}$ theorem or the Wedin theorem \cite{davis1970rotation,wedin1972perturbation,yu2015useful}) typically cannot guarantee meaningful estimation of $\bm{u}_l^{\star}$ unless 
\begin{align}
	(\text{classical theory}) \qquad \underset{\text{eigen-gap}}{\underbrace{ \min_{k:k\neq l} \, | \lambda_l^{\star} - \lambda_k^{\star} | }} > \| \bm{H} \| .
	\label{eq:eigen-separation}
\end{align}
This eigen-gap requirement can be problematic and hard to satisfy as the noise size grows, casting doubts on our ability to perform informative inference on individual eigenvectors.

\paragraph{Fine-grained estimation and inference (beyond $\ell_2$ guarantees).}
In many applications, it is often the case that the ultimate goal is not to characterize the $\ell_2$ or ``bulk'' behavior (e.g.~the mean squared estimation error) of an eigenvector estimator, but rather to reason about the eigenvectors along a few preconceived yet important directions. Take community recovery for instance: the eigenvector of a certain adjacency matrix  encodes community membership of a set of users \cite{abbe2017community}; if we wish to infer the community memberships of a few important users, or the similarities between a few pairs of users, it boils down to assessing the entrywise behavior or certain pairwise linear functional of an eigenvector estimator. Another example is concerned with harmonic retrieval: the leading eigenvector of a properly arranged Toeplitz matrix represents the time-domain response of the sinusoidal signals of interest \cite{chi2014compressive}; thus, retrieving the underlying frequency  involves inferring the Fourier coefficients of this eigenvector at some given frequencies. These problems can be formulated as estimation and inference for {\em linear functionals of an eigenvector}, namely,  quantities of the form $\bm{a}^{\top}\bm{u}_l^{\star}$ $(1\leq l\leq r)$ with $\bm{a}$ some prescribed vector.
In principle, this task can be viewed as a sort of \textit{fine-grained} statistical analysis, given that it pursues highly ``local'' and ``delicate'' information of interest. 

Towards estimating $\bm{a}^{\top}\bm{u}_l^{\star}$, a natural starting point is a ``plug-in'' estimator, which computes a reasonable estimate $\widehat{\bm{u}}_l$ of $\bm{u}_l^{\star}$  and outputs $\bm{a}^{\top}\widehat{\bm{u}}_l$. There are several critical challenges that we shall bear in mind. To begin with, a dominant fraction of prior theory focuses on $\ell_2$ risk analysis of an eigenvector estimator, which is often too coarse to deliver tight uncertainty assessment for the plug-in estimator. In fact, the $\ell_2$ risk bounds alone often lead to highly conservative estimates for the risk under consideration, making it hard to assess the performance of the plug-in estimator. To further complicate matters, there is often a severe bias issue surrounding the plug-in estimator. Even when an estimator $\widehat{\bm{u}}_l$ is nearly unbiased in a strong entrywise sense (meaning that the estimation bias is dominated by the variability in every single coordinate), this property alone does not preclude the possibility of bias accumulation along the  direction $\bm{a}$. Addressing these issues calls for refined risk analysis as well as  careful examination into the bias-variance trade-off.

\subsection{A glimpse of our approach and our contributions}

\paragraph{Eigen-decomposition meets statistical asymmetry.} 
The current paper makes progress in a setting where the noise matrix $\bm{H}$ consists of independent (but possibly heterogeneous and heteroscedastic) zero-mean components. Our approach is inspired by the findings of \cite{chen2018asymmetry}. Consider, for example, the case when $\bm{H}$ is a random and asymmetric matrix and when $\bm{M}^{\star}$ is a rank-1 symmetric matrix. The results in \cite{chen2018asymmetry} reveal that an eigen-decomposition approach applied to $\bm{M}$ (without symmetrization) achieves appealing statistical accuracy when estimating the leading eigenvalue of $\bm{M}^{\star}$. The key enabler is an implicit bias reduction feature when computing vanilla eigen-decomposition of an asymmetric data matrix. 
While \cite{chen2018asymmetry} only provides highly partial results when going beyond the rank-1 case, it hints at the potential benefits of eigen-decomposition in super-resolving the spectrum.

\paragraph{Our contributions.}
The main contributions of this paper are summarized below, all of which are built upon an eigen-decomposition approach applied to the asymmetric data matrix $\bm{M}$ (without symmetrization). As an important advantage, the proposed procedures and their accompanying theory are distribution-free, meaning that they do not require prior knowledge about the distribution of the noise and thus are fully adaptive to heteroscedasticity of data.

\begin{itemize}
	\item We demonstrate that the $l$th eigenvector $\bm{u}^{\star}_l$ and eigenvalue $\lambda_l^{\star}$ can be estimated with near-optimal accuracy even when the eigen-gap \eqref{eq:eigen-separation} is extremely small. More concretely, for various noise distributions our results only require 
\begin{align}
	(\text{our theory}) \qquad \min_{k:k\neq l} \, | \lambda_l^{\star} - \lambda_k^{\star} |  \gtrsim \| \bm{H} \| \frac{  \mathrm{poly}\log (n)}{\sqrt{n}}  ,
\label{eq:eigen-separation-new}
\end{align}
		which is about $O(\sqrt{n})$ times less stringent (up to log factor) than generic linear algebra theory (cf.~\eqref{eq:eigen-separation}). 

	\item We propose a new estimator for the linear functional $\bm{a}^{\top}\bm{u}^{\star}_l$ (obtained via proper de-biasing of certain plug-in estimators) that achieves minimax-optimal statistical accuracy. In addition, we demonstrate how to construct valid confidence intervals for  $\bm{a}^{\top}\bm{u}^{\star}_l$. 
		

	\item Additionally, we demonstrate how to compute valid confidence intervals for the eigenvalues of interest. Interestingly, de-biasing is not needed at all for performing inference on eigenvalues, as the eigen-decomposition approach implicitly alleviates the estimation bias.  
\end{itemize}
Our findings unveil new insights into the capability of spectral methods in a statistical context, going far beyond what conventional matrix perturbation theory has to offer.

\subsection{Notation}

For any vector $\bm{x}$, denote by $\|\bm{x}\|_2$ and $\|\bm{x}\|_{\infty}$ the $\ell_2$ norm and the $\ell_{\infty}$ norm of $\bm{x}$, respectively.
For any matrix $\bm{M}$, we let $\|\bm{M}\|$, $\|\bm{M}\|_{\mathrm{F}}$,  $\|\bm{M}\|_{\infty}$ and $\|\bm{M}\|_{2,\infty}$ represent the spectral norm, the Frobenius norm, the entrywise $\ell_{\infty}$ norm (i.e.~$\|\bm{M}\|_{\infty}:=\max_{i,j}|M_{ij}|$), and the two-to-infinity norm (i.e.~$\|\bm{M}\|_{2,\infty}:=\sup_{\|\bm{x}\|_{2}=1}\|\bm{M}\bm{x}\|_{\infty}$) of $\bm{M}$, respectively.
The  notation $f(n)=O\left(g(n)\right)$ or
$f(n)\lesssim g(n)$ means that there is a universal constant $c>0$ such
that $\left|f(n)\right|\leq c|g(n)|$, $f(n)\gtrsim g(n)$ means that there is a universal constant $c>0$ such
that $|f(n)|\geq c\left|g(n)\right|$, and $f(n)\asymp g(n)$ means that there exist constants $c_{1},c_{2}>0$
such that $c_{1}|g(n)|\leq|f(n)|\leq c_{2}|g(n)|$. The notation $f(n)\gg g(n)$ (resp.~$f(n)\ll g(n)$) means that there exists a sufficiently large (resp.~small) constant $c>0$ such that $f(n) \geq cg(n)$ (resp.~$f(n)\leq cg(n)$). 

In addition, we denote by $u_{l,j}$, $w_{l,j}$ and $u_{l,j}^{\star}$ the $j$th entry of $\bm{u}_l$, $\bm{w}_l$ and $\bm{u}_l^{\star}$, respectively. 
Let $\bm{e}_1,\cdots,\bm{e}_n$ represent the standard basis vectors in $\mathbb{R}^n$,  $\bm{1}_n\in \mathbb{R}^n$ the all-one vector, and $\bm{I}_n\in \mathbb{R}^{n\times n}$ the identity matrix. We shall also abbreviate the interval $[b-c, b+c]$ to $[b\pm c]$, and denote $\min|b\pm c|=\min\{|b+c|,|b-c|\}$, $\max|b\pm c|=\max\{|b+c|,|b-c|\}$, and $\min |\!|\!| \bm{b} \pm \bm{c} |\!|\!| = \min\{ |\!|\!| \bm{b}-\bm{c}|\!|\!|, |\!|\!|\bm{b}+\bm{c}|\!|\!| \}$ for any norm $|\!|\!|\cdot|\!|\!|$. Moreover, denote by $\Phi(\cdot)$ the cumulative distribution function (CDF) of a standard Gaussian random variable.  

\section{Problem formulation}

\paragraph{Model.} Imagine we seek to estimate an unknown rank-$r$
matrix
\begin{equation}
	\bm{M}^{\star}= \sum_{l=1}^r \lambda^{\star}_l \bm{u}_l^{\star} \bm{u}_l^{\star\top} 
	=: \bm{U}^{\star}\bm{\Sigma}^{\star}\bm{U}^{\star\top} ~\in~\mathbb{R}^{n\times n},
	\label{eq:model-rankr}
\end{equation}
where $\lambda_1^{\star} \geq \lambda_2^{\star} \geq \cdots \geq \lambda_r^{\star}$ denote the $r$ nonzero eigenvalues of $\bm{M}^{\star}$, and $\bm{u}_1^{\star}$, $\cdots$, $\bm{u}_r^{\star}$ represent the associated (normalized) eigenvectors. Here, for notational convenience we let
\begin{align}
	\bm{U}^\star := \begin{bmatrix}
	\bm{u}_1^\star, \cdots, \bm{u}_r^\star
	\end{bmatrix} \in \mathbb{R}^{n \times r} 
	\qquad \text{and} \qquad
	\bm{\Sigma}^{\star} := {\footnotesize \left[\begin{array}{ccc}
		\lambda_{1}^{\star}\\
 		& \ddots\\
 		&  & \lambda_{r}^{\star}
	\end{array}\right] \in \mathbb{R}^{r \times r} } .
\end{align}
We shall also define
\begin{align}
	\lambda_{\max}^{\star} := \max_{1\leq l\leq r}|\lambda_l^{\star}|, \qquad \lambda_{\min}^{\star} := \min_{1\leq l\leq r} |\lambda_l^{\star}| \qquad \text{and} \qquad \kappa := \lambda_{\max}^{\star} \,/\, \lambda_{\min}^{\star} .
\end{align}

What we have observed is a corrupted
copy of $\bm{M}^{\star}$, namely, 
\begin{equation}
	\bm{M}=\bm{M}^{\star}+\bm{H}, 
	\label{eq:model-rank1-noisy}
\end{equation}
where $\bm{H}\in\mathbb{R}^{n\times n}$ stands for a random noise
matrix. This paper focuses on the family of noise matrices satisfying the following assumptions:
\begin{assumption}
\label{assumption-H}
The noise matrix $\bm{H} \in \mathbb{R}^{n\times n}$
satisfies the following properties.
\begin{enumerate}
\item \textbf{Independence and zero mean}. The entries $\{H_{ij}\}_{1\leq i,j\leq n}$
are independent zero-mean random variables; this indicates that the
matrices $\bm{H}$ and $\bm{M}$ are, in general, {\em asymmetric}.
\item \textbf{Heteroscedasticity and unknown variances}. 
Let $\sigma_{ij}^{2} := \mathsf{Var}(H_{ij})$
denote the variance of $H_{ij}$. To accommodate more realistic scenarios,
we allow $\sigma_{ij}^{2}$ to vary across entries --- commonly referred
to as \emph{heteroscedastic noise}. In addition, we do not {a
priori} know $\{\sigma_{ij}^{2}\}_{1\leq i,j\leq n}$. Throughout
this paper, we assume
\begin{equation}
	0\leq \sigma_{\min}^{2}\leq\sigma_{ij}^{2}\leq\sigma_{\max}^{2},\qquad1\leq i,j\leq n.
	\label{eq:noise-var-bounds}
\end{equation}
%
%
\item \textbf{Magnitudes}. Each $H_{ij}$ ($1\leq i,j\leq n$) satisfies
either of the following conditions: 
\begin{enumerate} 
\item[(a)] $|H_{ij}|\leq B$; 
\item[(b)] $H_{ij}$ has a symmetric distribution
obeying $\mathbb{P}\{|H_{ij}|>B\}\leq c_{\mathrm{b}}n^{-12}$ for
		some constant $c_{\mathrm{b}}>0$, and $\mathbb{E}[H_{ij}^2\ind_{\{|H_{ij}|>B\}}]=o(\sigma_{ij}^2)$. 
\end{enumerate} 
\end{enumerate}
\end{assumption}
\begin{remark}
	Here, the quantities $B$, $\{\sigma_{i,j}\}$, $\sigma_{\min}$ and $\sigma_{\max}$ may all depend on $n$. 
\end{remark}
\begin{remark}
	As we shall see momentarily, the lower bound $\sigma_{\min}^{2}$ on the noise variance is imposed only for our statistical inference theory (or more precisely, it is imposed in order to ensure the plausibility to estimate the variance of the estimation error in an accurate manner).  This lower bound $\sigma_{\min}^{2}$ is not needed at all in our eigenvector estimation guarantees (e.g., Theorems~\ref{thm:rankr-bounds-simple} and \ref{thm:rankr-evalue-bounds-simple}). 
\end{remark}

\noindent The careful reader would naturally ask when we would have an asymmetric noise matrix $\bm{H}$.  This might happen when, for example, one has collected two independent samples about each entry of $\bm{M}^{\star}$ and chooses to arrange the observed data in an asymmetric manner.  Moreover, when the noise matrix is a symmetric Gaussian random matrix and possibly contains missing data,  \cite[Appendix J]{chen2018asymmetry} points out some simple asymmetrization tricks that allow one to convert a symmetric data matrix $\bm{M}$ to an asymmetric counterpart with independent components; see Appendix~\ref{sec:example-homoscedastic-Gaussian} for an example.

In addition to the above assumptions on the noise, we assume that
the unknown matrix $\bm{M}^{\star}$ satisfies a certain incoherence
condition, as commonly seen in the low-rank matrix estimation literature.

\begin{definition}[\textbf{Incoherence}]
	\label{defn:incoherence}
	The matrix $\bm{M}^{\star}$ with eigen-decomposition $\bm{M}^{\star}=\bm{U}^{\star}\bm{\Sigma}^{\star}\bm{U}^{\star\top}$
	 is said to be $\mu$-incoherent
	if $\|\bm{U}^{\star}\|_{\infty}\leq\sqrt{\mu/n}$.

\end{definition}
%

\paragraph{Goal.} The primary goal of the current paper is to perform certain ``fine-grained'' statistical
estimation and inference on the unknown eigenvectors $\{\bm{u}^{\star}_1, \cdots, \bm{u}^{\star}_r\}$. 
To be more specific,
we aim at developing statistically efficient methods to estimate and construct confidence
intervals for linear functionals taking the form $\bm{a}^{\top}\bm{u}^{\star}_l$ ($1\leq l\leq r$), where $\bm{a}\in\mathbb{R}^{n}$ is a pre-determined fixed vector. Along the way, 
we shall also demonstrate how to perform estimation and inference on the unknown eigenvalues $\{\lambda^{\star}_1, \cdots, \lambda_r^{\star}\}$.  Ideally, all these tasks can be accomplished even when the associated eigen-gaps are very small, without requiring prior knowledge about the noise distributions and noise levels.

\section{Estimation}
\label{sec:estimation}

This section presents our algorithms and the accompanying theory for estimating an eigenvector $\bm{u}^{\star}_l$, a linear functional $\bm{a}^{\top}\bm{u}_l^{\star}$ of this eigenvector,  as well as the associated eigenvalue $\lambda_l^{\star}$.  

Two popular estimation schemes immediately come into mind: (1) performing eigen-decomposition after 
symmetrizing the data matrix (e.g.~replacing $\bm{M}$ by $\frac{1}{2}(\bm{M}+\bm{M}^{\top})$); (2) computing singular value decomposition (SVD) of the asymmetric data matrix $\bm{M}$.  Contrastingly, this paper adopts a far less widely used, and in fact far less widely studied, strategy based on eigen-decomposition without symmetrization; namely, we attempt estimation of $\{\bm{u}^{\star}_l\}$ and $\{\lambda_l^{\star}\}$ via the eigenvectors and the eigenvalues of $\bm{M}$, respectively, despite the asymmetric nature of $\bm{M}$  in general. While conventional wisdom does not advocate eigen-decomposition of asymmetric matrices (due to, say, potential numerical instability), our investigation uncovers remarkable advantages of this approach under the statistical context considered in the present paper.

\subsection{Notation: eigen-decomposition without symmetrization}
\label{sec:notation-eigen}

Before continuing, we introduce several additional notation that shall be adopted throughout. Owing to the asymmetry of
$\bm{M}$, the left and the right eigenvectors of $\bm{M}$ do not coincide. 
%
\begin{notation}
\label{notation:leading-evectors}
Let $\lambda_1,\cdots,\lambda_r$ denote the top-$r$ leading eigenvalues of $\bm{M}$ (so that $\min_{1\leq l\leq r}|\lambda_l|$ is larger in magnitude than any other eigenvalue of $\bm{M}$). Assume that they are sorted by their real parts, namely, 
$\mathrm{Re}(\lambda_1) \geq \mathrm{Re}(\lambda_2) \geq \cdots \geq \mathrm{Re}(\lambda_r)$, and denote by $\ur_l$
(resp.~$\ul_l$) the right (resp.~left) eigenvector
of $\bm{M}$ associated with $\lambda_l$; this means 
\begin{align}
	\bm{M}\ur_l = \lambda_l \ur_l \qquad 
	\text{and} \qquad
	\bm{M}^{\top}\ul_l = \lambda_l \ul_l.
\end{align}
In addition, if $\ur_l$ and $\ul_l$ are both real-valued, then we assume without loss of generality that $\langle\ur_l,\ul_l\rangle \geq 0$.
\end{notation}

\begin{remark}
	As we shall justify shortly in Theorem~\ref{thm:rankr-bounds-simple}, even though $\bm{M}$ is in general asymmetric, the eigenvalue $\lambda_l$ and the eigenvectors $\bm{u}_l$ and $\bm{w}_l$ are, with high probability, real-valued under the assumptions imposed in this paper. 
\end{remark}

\subsection{Estimation algorithms}

We are now ready to present our procedures for estimating the $l$th eigenvector of $\bm{M}^{\star}$, on the basis of eigen-decomposition of $\bm{M}$ 
(see Notation~\ref{notation:leading-evectors}).
\begin{itemize}
	\item {\bf Estimator for $\bm{u}_l^{\star}$ (with the aim of achieving low $\ell_2$ risk):}  
	\begin{equation}
		\widehat{\bm{u}}_{l} := \frac{1}{\|\bm{u}_{l}+\bm{w}_{l}\|_{2}}\big(\bm{u}_{l}+\bm{w}_{l}\big); 
		\label{eq:defn-hat-ul}
	\end{equation}
\item {\bf Estimator for $\bm{a}^{\top}\bm{u}_l^{\star}$ for a preconceived direction $\bm{a}$: }
	\begin{align}
		\label{eq:ua-estimator-large}
		\widehat{u}_{\bm{a},l} := \min\Bigg\{ \sqrt{\Big|\tfrac{ (\bm{a}^{\top}\ur_l) (\bm{a}^{\top}\ul_l)}{\ul^{\top}_l \ur_l}\Big|}, \, \|\bm{a}\|_2 \Bigg\} .
	\end{align}
\end{itemize}
\noindent 
Here, the rationale behind $\widehat{u}_{\bm{a},l}$ is this:  both $\bm{a}^{\top}\ur_l$ and $\bm{a}^{\top}\ul_l$ might systematically {\em under-estimate} the quantity of interest $\bm{a}^{\top}\bm{u}_l^{\star}$. As a result, one is advised to first alleviate the bias via proper de-shrinking, which is the role played by the factor $\frac{1}{\ul^{\top}_l \ur_l}$. The choice of this factor is based on in-depth understanding of the behavior of $\bm{a}^{\top}\ur_l$ and $\bm{a}^{\top}\ul_l$, and 
will be better understood after we delve into technical details. As an important feature, the proposed estimation procedures do not rely on prior knowledge about the noise distributions.  
\begin{remark}
	The estimator \eqref{eq:ua-estimator-large} 
	involves the term $\|\bm{a}\|_2$ due to the trivial upper bound $|\bm{a}^{\top}\bm{u}_l^{\star}| \leq \|\bm{a}\|_2 \|\bm{u}_l^{\star}\|_2=\|\bm{a}\|_2 $.
\end{remark}

\subsection{Theoretical guarantees}
\label{sec:estimation-theory}

We now embark on theoretical development for the proposed estimators. As alluded to previously, whether we can reliably estimate and infer an eigenvector $\bm{u}_l^{\star}$ depends largely upon the spacing between the $l$th eigenvalue and its adjacent eigenvalues. In light of this, we define formally the eigen-gap metric w.r.t.~the $l$th eigenvalue of $\bm{M}^{\star}$ as follows
\begin{equation}
	\Delta_l^\star \defn \begin{cases} \min\limits_{1 \leq k \leq r, k \neq l} | \lambda_l^\star - \lambda_k^\star|, \quad & \text{if }r>1; \\ \infty, & \text{otherwise}. \end{cases}
	\label{eq:eigengap}
\end{equation}
Most of our theoretical guarantees rely on this crucial metric. Moreover, our theory in this section is built upon a set of assumptions on the noise levels:  
\begin{assumption}
	\label{assumption:noise-size-rankr-revised}
	Suppose that the noise parameters defined in Assumption \ref{assumption-H} satisfy
\begin{subequations}
	\begin{align}
		\Delta_{l}^{\star} & > 2c_4 \kappa^2 r^2 \sigma_{\max} \sqrt{\mu\log n}  \label{eq:condition-Delta-noise-revised} \\
		B\log n & \leq \sigma_{\mathrm{max}} \sqrt{n \log n}  \leq c_5 \lambda_{\min}^{\star} / \kappa^3 \label{eq:condition-lambdamin-rankr-revised}
	\end{align}
\end{subequations}
for some sufficiently large (resp.~small) universal constant $c_4>0$ (resp.~$c_5>0$). 
\end{assumption}
\begin{remark} 
Consider, for example, the case where $r,\kappa,\mu\asymp 1$: the lower bound on  $\Delta_{l}^{\star}$ in \eqref{eq:condition-Delta-noise-revised} is $O(\sqrt{n})$ times smaller than the lower bound on $\lambda_{\min}^{\star}$ in \eqref{eq:condition-lambdamin-rankr-revised}, meaning that the eigen-gap can be considerably smaller than the minimum eigenvalue of the truth.  
\end{remark}

The following theorem delivers statistical guarantees for the proposed eigenvector estimators, with the proof postponed to Appendix~\ref{sec:proof-rank-r-all}. 
We recall the notation $\min  |\!|\!| \bm{z}\pm\bm{u}_{l}^{\star} |\!|\!|  := \min\{ |\!|\!| \bm{z}-\bm{u}_{l}^{\star} |\!|\!|, |\!|\!| \bm{z}+\bm{u}_{l}^{\star} |\!|\!|  \}   $  for any norm $|\!|\!|\cdot|\!|\!|$. 

\begin{theorem}[eigenvector estimation] 
\label{thm:rankr-bounds-simple}
Suppose that $\mu\kappa^2 r^4 \lesssim n$, and that Assumptions~\ref{assumption-H}-\ref{assumption:noise-size-rankr-revised} hold.  With probability at least $1-O(n^{-6})$,  the eigenvalue $\lambda_l$ and the associated eigenvectors $\ur_l$ and $\ul_l$ (see Notation \ref{notation:leading-evectors}) are all real-valued, and one has the following: 
\begin{subequations}
\begin{enumerate}
	\item ($\ell_2$ and $\ell_{\infty}$ guarantees)
	\begin{align} 
		\min \left\Vert \bm{u}_{l}\pm\bm{u}_{l}^{\star}\right\Vert _{2}  
		& \lesssim
		\frac{\sigma_{\max}\sqrt{\kappa^{6}n\log n}}{\lambda_{\min}^{\star}} + \frac{\sigma_{\max}}{\Delta_{l}^{\star}}\sqrt{\mu\kappa^{4}r^{2}\log n} , 
	 	\label{eq:rankr-l2-bound-2-simple} \\
		\min  \left\Vert \bm{u}_{l}\pm\bm{u}_{l}^{\star}\right\Vert _{\infty}  & \lesssim\,  \frac{\sigma_{\mathrm{max}}}{\lambda_{\mathrm{min}}^{\star}}\sqrt{\mu\kappa^{4}r\log n}+\frac{\sigma_{\mathrm{max}}}{\Delta_{l}^{\star}}\sqrt{\frac{\mu^{2}\kappa^{4}r^{3}\log n}{n}} ;
		\label{eq:rankr-l-infty-bound-simple}
	\end{align}
		these hold unchanged if ${\bm{u}}_l$ is replaced by either ${\bm{w}}_l$ or $\widehat{\bm{u}}_l$ (cf.~\eqref{eq:defn-hat-ul}); 
	\item (statistical guarantees for linear forms of an eigenvector) for any fixed vector $\bm a$ with $\ltwo{\bm a} = 1$, the proposed estimator \eqref{eq:ua-estimator-large} obeys\footnote{If the rank $r=1$, then the 2nd and the 3rd terms on the right-hand side of \eqref{eq:rankr-linear-form-bound-simple} are set to be zero.}
	\begin{align} 
		\min\left|\widehat{u}_{\bm{a},l}\pm\bm{a}^{\top}\bm{u}_{l}^{\star}\right|
		\lesssim\,
		\tfrac{\sigma_{\mathrm{max}}r^2\sqrt{\mu\kappa^{4}\log n}}{\lambda_{\mathrm{min}}^{\star}}
		+\tfrac{\sigma_{\max}^{2}\mu r^{2}\kappa^{4}\log n}{(\Delta_{l}^{\star})^{2}}|\bm{a}^{\top}\bm{u}_{l}^{\star}|
		+\sigma_{\mathrm{max}}\sqrt{\mu\kappa^{4}r^{3}\log n}\max_{k\neq l}\tfrac{|\bm{a}^{\top}\bm{u}_{k}^{\star}|}{|\lambda_{l}^{\star}-\lambda_{k}^{\star}|}	. 
		\label{eq:rankr-linear-form-bound-simple}
	\end{align}
	%
\end{enumerate}
\end{subequations}
\end{theorem}
Interestingly, even though we work with eigen-decomposition of asymmetric matrices, the obtained eigenvalue and eigenvector are provably real-valued as long as a fairly mild eigen-gap condition is satisfied. This presents an important feature that cannot be derived from generic matrix perturbation theory.

To interpret the effectiveness of this theorem, we find it convenient to concentrate on the case with $r,\kappa, \mu \asymp 1$ under i.i.d.~Gaussian noise. The implications in this case are summarized below.   
\begin{itemize}

	\item {\em $\ell_2$ and $\ell_{\infty}$ guarantees.} The $\ell_2$ and $\ell_{\infty}$ statistical guarantees derived in Theorem~\ref{thm:rankr-bounds-simple} read
	\begin{subequations}
	\label{eq:L2-Linf-simple-r1}
	\begin{align}
		\min\left\Vert \bm{u}_{l}\pm\bm{u}_{l}^{\star}\right\Vert _{2} & \lesssim\frac{\sigma_{\max}\sqrt{n\log n}}{\lambda_{\min}^{\star}}+\frac{\sigma_{\max}\sqrt{\log n}}{\Delta_{l}^{\star}}=:\mathcal{E}_{2} ;
		\label{eq:L2-simple-r1}\\
		\min\left\Vert \bm{u}_{l}\pm\bm{u}_{l}^{\star}\right\Vert _{\infty} & \lesssim\frac{1}{\sqrt{n}}\mathcal{E}_{2}. \label{eq:Linf-simple-r1}
	\end{align}
	\end{subequations}
	 Our $\ell_{\infty}$ error bound is about $O(\sqrt{n})$ times smaller than the $\ell_2$ risk bound, implying that the energy of the estimation error is more or less dispersed across all entries. 

	\item {\em Improved eigen-gap requirements.} Consistent estimation of $\bm{u}_l^{\star}$ --- in the sense of $\min\|\bm{u}_{l}\pm\bm{u}_{l}^{\star}\|_{2}=o(\|\bm{u}_{l}^{\star}\|_{2})$
		and $\min\|\bm{u}_{l}\pm\bm{u}_{l}^{\star}\|_{\infty}=o(\|\bm{u}_{l}^{\star}\|_{\infty})$ --- is guaranteed as long as\footnote{Note that when $\{H_{ij}\}$ are i.i.d.~Gaussian, we have $\|\bm{H}\|\asymp \sigma_{\max}\sqrt{n}$ with high probability. } 
		\begin{subequations}
			\begin{align}
				\|\bm{H}\|\log n & \lesssim\|\bm{M}^{\star}\|;\label{eq:H-norm-constraint-r1}\\
				\Delta_{l}^{\star} & \gtrsim\frac{\|\bm{H}\|\log n}{\sqrt{n}}.\label{eq:Delta-constraint-r1}
			\end{align} 
		\end{subequations}
		While the condition \eqref{eq:H-norm-constraint-r1} is commonly seen in prior literature (up to some log factor), the eigen-gap requirement \eqref{eq:Delta-constraint-r1} is in stark contrast to classical matrix perturbation theory (e.g.~the Davis-Kahan $\sin\bm{\Theta}$ theorem or the Wedin theorem \cite{davis1970rotation, wedin1972perturbation, chen2020spectral}). In fact, prior theory typically requires the spacing between $\lambda_l^{\star}$ and the rest of the eigenvalues to at least exceed
		\begin{align}
			\Delta_l^{\star} \gtrsim  \|\bm{H}\|  \qquad (\text{prior theory}), 
		\end{align}
		which is about $O(\sqrt{n}/\log n)$ times more stringent than our requirement \eqref{eq:Delta-constraint-r1}.

	\item {\em The influence of eigen-gaps and $\{\bm{a}^{\top}\bm{u}_i^{\star}\}$ upon estimation accuracy.} For any fixed unit vector $\bm{a}$, our theoretical guarantees for estimating $\bm{a}^{\top}\bm{u}_l^{\star}$ read
\begin{align}
\min\left|\widehat{u}_{\bm{a},l}\pm\bm{a}^{\top}\bm{u}_{l}^{\star}\right| 
& \lesssim\, \frac{\sigma_{\mathrm{max}}\sqrt{\log n}}{\lambda_{\mathrm{min}}^{\star}}
+
\frac{\sigma_{\max}^{2}\log n}{(\Delta_{l}^{\star})^{2}}|\bm{a}^{\top}\bm{u}_{l}^{\star}|
+
\underset{\text{``interferers''}}{\underbrace{\sigma_{\mathrm{max}}\sqrt{\log n}\max_{k\neq l}\frac{\big|\bm{a}^{\top}\bm{u}_{k}^{\star}\big|}{\big|\lambda_{l}^{\star}-\lambda_{k}^{\star}\big|}}}.\label{eq:au-smiple-r1}
\end{align}
		The first term on the right-hand side of \eqref{eq:au-smiple-r1} is a universal term that is no larger than $O(1/\sqrt{n})$ times the $\ell_{2}$ bound \eqref{eq:L2-simple-r1}. The other two terms are more complicated, which depend on not only the spacing of the eigenvalues but also the sizes of $\{\bm{a}^{\top}\bm{u}_i^{\star}\}_{1\leq i\leq r}$. More concretely, (1) the influence of the target quantity $\bm{a}^{\top}\bm{u}_l^{\star}$ upon the estimation error is captured by the multiplicative factor $\frac{\sigma_{\max}^{2}}{(\Delta_{l}^{\star})^{2}}$, which scales inverse quadratically in $\Delta_l^{\star}$; (2) the linear functionals of other eigenvectors (namely,  $\{\bm{a}^{\top}\bm{u}_k^{\star}\}_{k\neq l}$) essentially behave as ``interferers'' that might degrade estimation fidelity; in particular, the influence of $\bm{a}^{\top} \bm{u}_k^{\star}$ ($k\neq l$) upon estimation loss can be understood through the coefficient $\sigma_{\max} / |\lambda_{l}^{\star}-\lambda_{k}^{\star}|$, which is inversely proportional to the associated eigen-gap $|\lambda_{l}^{\star}-\lambda_{k}^{\star}|$. Intuitively, if $\bm{a}^{\top}\bm{u}_k^{\star}$ ($k\neq l$) becomes large, it results in stronger ``interference'' when estimating $\bm{a}^{\top}\bm{u}_l^{\star}$; this adverse effect would be easier to mitigate if the gap $|\lambda_l^{\star} - \lambda_k^{\star} |$ increases (so that it becomes easier to differentiate the $l$th and the $k$th eigenvectors). 
	%


	\item {\em The rank-1 case.} When $r=1$, the preceding theory can be significantly simplified as follows
		\begin{align*}
			\min\left\Vert \bm{u}_{l}\pm\bm{u}_{l}^{\star}\right\Vert _{2} & \lesssim\frac{\sigma_{\max}\sqrt{n\log n}}{\lambda_{\min}^{\star}},\qquad
			\min\left\Vert \bm{u}_{l}\pm\bm{u}_{l}^{\star}\right\Vert _{\infty}\lesssim\frac{\sigma_{\max}\sqrt{\log n}}{\lambda_{\min}^{\star}} , \\
 			& \quad\quad  \min\left|\widehat{u}_{\bm{a},l}\pm\bm{a}^{\top}\bm{u}_{l}^{\star}\right|\lesssim\frac{\sigma_{\mathrm{max}}\sqrt{\log n}}{\lambda_{\mathrm{min}}^{\star}}.
		\end{align*}
		In particular, the same estimation error bound --- which is about $\sqrt{n}$ times smaller than  the $\ell_2$ loss --- holds for any arbitrary direction as specified by $\bm{a}$. In other words, the estimation error is more or less identical over any  pre-determined direction.

\end{itemize}
Encouragingly, the above performance guarantees are minimax optimal up to some logarithmic factor, as we shall elucidate in Section~\ref{sec:minimax-lower-bounds}.

As it turns out, the feasibility of faithful eigenvector estimation is largely dictated by our ability to locate the $l$th eigenvalue $\lambda^{\star}_l$ and to disentangle it from the rest of the spectrum, which becomes particularly challenging if the eigen-gap $\Delta_l^{\star}$ is very small. In light of this, we develop the following eigenvalue perturbation theory that significantly improves upon generic linear algebra theory.

%
\begin{theorem}[eigenvalue estimation] 
\label{thm:rankr-evalue-bounds-simple}
	Suppose that $\mu\kappa^2 r^4 \lesssim n$, and that Assumptions~\ref{assumption-H}-\ref{assumption:noise-size-rankr-revised} hold. With probability at least $1-O(n^{-6})$, one has
\begin{equation} 
	\label{eq:rankr-eigenvalue-bound-simple}
	|\lambda_l - \lambda_l^\star| \leq c_4 \sigma_{\mathrm{max}}  \sqrt{\mu \kappa^2 r^4 \log n} . 
\end{equation}
\end{theorem}
In words, as long as the eigen-gap $\Delta_l^\star$ exceeds \eqref{eq:condition-Delta-noise-revised}, the vanilla eigen-decomposition approach (without symmetrization) produces an estimate of $\lambda^{\star}_l$ with an estimation error at most $\Delta_l^\star/2$, meaning that we have managed to locate $\lambda^{\star}_l$ with a high resolution. As mentioned previously, this eigen-gap requirement can be interpreted as $\Delta_l^{\star} \gtrsim \|\bm{H}\| \frac{\log n}{\sqrt{n}}$ in the Gaussian noise case with $r,\mu,\kappa \asymp 1$.   
For the sake of comparison, we remind  the reader that classical linear algebra theory (e.g.~Weyl's inequality and the Bauer-Fike inequality) only provides perturbation bounds with fairly low resolution --- namely, bounds at best on the order $|\lambda_l - \lambda^{\star}_l| \leq \|\bm{H}\|$  --- which are highly insufficient and in fact suboptimal for our purpose.

It is worth emphasizing that the ability to obtain highly accurate eigenvalue estimates is the key to tackling small eigen-gaps. 
Given that we can estimate $\lambda_l^\star$ with a precision  $\sigma_{\mathrm{max}}  \sqrt{\mu \kappa^2 r^4 \log n}$, 
we can distinguish $\lambda_l^\star$ from the rest of the spectrum even when the associated eigen-gap is on the same order. 
This also helps explain why a large body of prior theory fell short of accommodating small eigen-gaps. 
As discussed in \cite[Section~4.3]{chen2018asymmetry}, the standard eigen-decomposition approach when applied to symmetric matrices suffers from a non-negligible bias issue, 
which fails to yield a high-precision eigenvalue estimate. While the recent work \cite{li2021minimax} developed a de-biasing scheme that allows one to cope with small eigen-gaps in the presence of an i.i.d.~symmetric Gaussian noise matrix,  the current paper covers a remarkably broader class of noise distributions by exploiting the special property of statistical asymmetry.

Finally,  as mentioned before, one might be tempted to first symmetrize the data matrix (i.e.~replacing $\bm{M}$ with $\frac{1}{2}(\bm{M}+\bm{M}^{\top})$) before computing eigen-decomposition. While \cite{chen2018asymmetry} has discussed the non-negligible bias of this approach in eigenvalue estimation, it remains unclear whether or not this natural approach suffers from a bias issue as well when estimating the eigenvectors, and, more crucially, whether or not we can still expect reliable eigenvector estimation after symmetrizing the data matrix.   Addressing these questions is instrumental in understanding whether asymmetry plays an important role or it merely provides theoretical convenience. While a theory for this is yet to be developed, 
we shall compare these two approaches numerically and demonstrate the gains of our approach in Section~\ref{SecNumericals}.

\subsection{Minimax lower bounds}
\label{sec:minimax-lower-bounds}

To assess the tightness of our statistical guarantees, we develop localized minimax lower bounds on eigenvector estimation, 
aimed at providing in-depth understanding about how the estimation difficulty changes as a function of certain salient parameters of the problem instance. 
To derive these lower bounds, we consider the non-asymptotic local minimax framework (see, e.g.~ \cite{cai2004adaptation,cai2015framework,cai2018adaptive}), an approach built upon the concept of the hardest local alternatives that finds its roots in \cite{stein1956efficient}.

Before stating our main results, let us first define several sets that are properly localized around the truth.
Let $\mathbb{S}^n \subseteq \mathbb{R}^n$ represent the set of $n\times n$ symmetric matrices. 
Denoting by $\lambda_{l}(\bm{A})$ the $l$th largest eigenvalue of $\bm{A}$ and  $\bm{u}_{l}(\bm{A})$  the associated eigenvector of a symmetric matrix $\bm{A}$, we define 
\begin{subequations}
\begin{align}
\mathcal{\mathcal{M}}_{0}(\bm{M}^{\star}) &:=\left\{ \bm{A}\in \mathbb{S}^n \mid\mathsf{rank}(\bm{A})=r,\,\lambda_{i}(\bm{A})=\lambda_{i}^{\star}\text{ }(1\leq i\leq r),\,\|\bm{A}-\bm{M}^{\star}\|_{\mathrm{F}}\leq\frac{\sigma_{\min}}{2}\right\} \\
\mathcal{\mathcal{M}}_{1}(\bm{M}^{\star}) &:=\left\{ \bm{A} \in \mathbb{S}^n \mid\mathsf{rank}(\bm{A})=r,\,\lambda_{i}(\bm{A})=\lambda_{i}^{\star}\text{ }(1\leq i\leq r),\,\|\bm{u}_{l}(\bm{A})-\bm{u}_{l}^{\star}\|_{2}\leq\frac{\sigma_{\min}}{4|\lambda_{l}^{\star}|}\right\} \\
\mathcal{\mathcal{M}}_{2}(\bm{M}^{\star}) &:=\left\{ \bm{A} \in \mathbb{S}^n \mid\mathsf{rank}(\bm{A})=r,\,\lambda_{i}(\bm{A})=\lambda_{i}^{\star}\text{ }(1\leq i\leq r),\,\|\bm{u}_{l}(\bm{A})-\bm{u}_{l}^{\star}\|_{2}\leq c_{4}\frac{\sigma_{\min}\sqrt{n}}{|\lambda_{l}^{\star}|}\right\}
\end{align}
\end{subequations}
for some sufficiently large constant $c_{4}>0$. 

Now we are ready to state our results, the proof of which is deferred to Appendix~\ref{sec:proof-minimax-lower-bounds}.

\begin{theorem}[minimax lower bounds]
\label{thm:minimax-l2-au-eigengap}
Consider any $1\leq l\leq r\le n/2$. Suppose that $H_{ij}\overset{\mathrm{ind.}}{\sim}\mathcal{N}(0,\sigma_{ij}^{2})$ with $\sigma_{ij}^2\geq \sigma_{\min}^2$,
and assume that $4\sigma_{\min}\sqrt{n}\leq|\lambda_{l}^{\star}|$ and
$\sigma_{\min}\leq\Delta_{l}^{\star}$. 
\begin{enumerate}
\item 
There exist some constants $c_{0},c_{1}>0$ such that\begin{subequations}
\begin{align}
\inf_{\widehat{\bm{u}}_{l}}\sup_{\bm{A}\in{\mathcal{M}}_{0}(\bm{M}^{\star})}\mathbb{E}\Big[\min\|\widehat{\bm{u}}_{l}\pm\bm{u}_{l}(\bm{A})\|_{2}\Big] & \geq c_{0}\frac{\sigma_{\min}}{\Delta_{l}^{\star}};\label{eq:minimax-L2-u}\\
\inf_{\widehat{u}_{\bm{a},l}}\sup_{\bm{A}\in{\mathcal{M}}_{0}(\bm{M}^{\star})}\mathbb{E}\Big[\min\big|\widehat{u}_{\bm{a},l}\pm\bm{a}^{\top}\bm{u}_{l}(\bm{A})\big|\Big] & \geq c_{1}\Bigg\{\left|\bm{a}^{\top}\bm{u}_{l}^{\star}\right|\frac{\sigma_{\min}^{2}}{\Delta_{l}^{\star2}}+\sigma_{\min}\max_{k:k\neq l}\frac{\left|\bm{a}^{\top}\bm{u}_{k}^{\star}\right|}{|\lambda_{l}^{\star}-\lambda_{k}^{\star}|}\Bigg\}.\label{eq:minimax-au-eigen-gap}
\end{align}
\item 
There exists some constant $c_{2}>0$ such that
\begin{align}
\inf_{\widehat{u}_{\bm{a},l}}\sup_{\bm{A}\in{\mathcal{M}}_{1}(\bm{M}^{\star})}\mathbb{E}\Big[\min\big|\widehat{u}_{\bm{a},l}\pm\bm{a}^{\top}\bm{u}_{l}(\bm{A})\big|\Big] & \geq c_{2}\frac{\sigma_{\min}\big\|\bm{a}-(\bm{a}^{\top}\bm{u}_{l}^{\star})\bm{a}\big\|_{2}}{\big|\lambda_{l}^{\star}\big|}.\label{eq:minimax-au-sigma}
\end{align}
\item 
Then there exists
some  constant $c_{3}>0$ such that
\begin{align}
\inf_{\widehat{\bm{u}}_{l}}\sup_{\bm{A}\in {\mathcal{M}}_{2}(\bm{M}^{\star})}\mathbb{E}\Big[\min\|\widehat{\bm{u}}_{l}\pm\bm{u}_{l}(\bm{A})\|_{2}\Big] & \geq c_{3}\frac{\sigma_{\min}\sqrt{n}}{\big|\lambda_{l}^{\star}\big|}.\label{eq:minimax-L2-sigma}
\end{align}
\end{subequations}
\end{enumerate}
Here, the infimum in \eqref{eq:minimax-L2-u} and \eqref{eq:minimax-L2-sigma} is over all eigenvector
estimators, while the infimum in \eqref{eq:minimax-au-eigen-gap}
and \eqref{eq:minimax-au-sigma} is over all estimators for the linear
form of the $l$th eigenvector.
\end{theorem}
\begin{remark}
	Note that each part of Theorem~\ref{thm:minimax-l2-au-eigengap} captures one bottleneck for the estimation task. Given that enlarging the set of matrices under consideration can only lead to increased minimax lower bound,  we can take the supremum over the union $\mathcal{M}_{0}(\bm{M}^{\star}) \cup \mathcal{M}_{1}(\bm{M}^{\star}) \cup \mathcal{M}_{2}(\bm{M}^{\star})$  in order to yield a minimax lower bound that reflects all these bottlenecks at once. 
\end{remark}


Consider again the case with $r,\kappa,\mu \asymp 1$ under i.i.d.~Gaussian noise, and the case when $|\bm{a}^{\top}\bm{u}_l^{\star}|\leq (1-\epsilon)\|\bm{a}\|_2$ for an arbitrarily small constant $\epsilon>0$ (so that $\bm{a}$ is not perfectly aligned with $\bm{u}_l^{\star}$). The theoretical guarantees derived for our estimators (see Theorem~\ref{thm:rankr-bounds-simple}, or \eqref{eq:L2-Linf-simple-r1} and \eqref{eq:au-smiple-r1}) match the minimax lower bounds in Theorem~\ref{thm:minimax-l2-au-eigengap} up to some logarithmic factor, including both $\ell_2$ loss and the risk for estimating an arbitrary linear form $\bm{a}^{\top}\bm{u}_l^{\star}$. In particular, the dependencies of the estimation risk on $\Delta_l^{\star}$,  $|\lambda_l^{\star}-\lambda_k^{\star}|$, and $\{\bm{a}^{\top}\bm{u}_k^{\star}\}$ in Theorem~\ref{thm:rankr-bounds-simple} --- which might seem complicated at first glance --- are all optimal modulo some log factor. All this confirms the effectiveness and optimality of the proposed estimator in fine-grained eigenvector estimation.

\section{Inference and uncertainty quantification}
\label{sec:inference}

This section moves one step further to the task of statistical inference, with the aim of constructing valid and efficient confidence
intervals for the linear functionals $\bm{a}^{\top}\bm{u}^{\star}_l$ and the eigenvalues $\lambda^\star_{l}$  ($1\leq l\leq r$). More precisely, for any
target coverage level $0<1-\alpha<1$, we seek to compute intervals such that 
\begin{align}
	\mathbb{P}\left\{ \bm{a}^{\top}\bm{u}^{\star}_l \in[c_{\mathrm{lb}}^{\bm{a}},c_{\mathrm{ub}}^{\bm{a}}]\right\} \approx 1-\alpha
\qquad\text{or}\qquad
	\mathbb{P}\left\{ -\bm{a}^{\top}\bm{u}^{\star}_l\in[c_{\mathrm{lb}}^{\bm{a}},c_{\mathrm{ub}}^{\bm{a}}]\right\} \approx 1-\alpha.
\end{align}
and
\begin{align}
	\mathbb{P}\left\{ \lambda^\star_l \in[c_{\mathrm{lb}}^{\lambda},c_{\mathrm{ub}}^{\lambda}]\right\} \approx 1-\alpha
\end{align}
Here, we have taken into account the difficulty in distinguishing $\bm{u}^{\star}_l$ 
from $-\bm{u}^{\star}_l$ using only the data $\bm{M}$. 
As we shall see, providing precise confidence intervals for the linear functional is much more challenging than the estimation task, and our theory requires additional assumptions beyond Assumption~\ref{assumption:noise-size-rankr-revised} in order for our procedures to have exact coverage.

\subsection{Algorithms}
\label{sec:construction-confidence-intervals}

\subsubsection{Confidence intervals for linear forms of eigenvectors}
\label{sec:CI-linear-form}

We start by constructing confidence intervals for the linear form $\bm{a}^{\top}\bm{u}^{\star}_l$. 
Towards this, two ingredients are needed:  (1) a nearly unbiased
estimate of $\bm{a}^{\top}\ustar_l$, and (2) a valid length
of the interval (which is typically built upon a certain variance
estimate).  We describe these ingredients as follows.

\paragraph{A modified nearly unbiased estimator $\widehat{u}_{\bm{a},l}^{\mathsf{modified}}$. }  While the estimator $\widehat{u}_{\bm{a},l}$ (cf.~\eqref{eq:ua-estimator-large}) discussed previously enjoys minimax optimal statistical accuracy,  we find it more convenient to work with a modified estimator when conducting uncertainty quantification, particularly in the regime when $\bm{a}^{\top}\bm{u}_l^{\star}$ is very small.\footnote{More precisely, in the regime where $\bm{a}^{\top}\bm{u}_l^{\star}$ is very small, the uncertainty of $\widehat{u}_{\bm{a},l}^{\mathsf{modified}}$ is nearly Gaussian, while that of $\widehat{u}_{\bm{a},l}$ is non-Gaussian and more complicated to describe.}   More specifically, we introduce
	\begin{equation}
		\widehat{u}_{\bm{a},l}^{\mathsf{modified}} := 
		\begin{cases}
			\frac{1}{2}\bm{a}^{\top}\big(\ur_l + \ul_l\big), 
			& \text{if }\big|\bm{a}^{\top}\bm{u}_l^{\star}\big|\text{ is ``small''},\\
			\min\left\{ \sqrt{\left|\frac{1}{\ul^{\top}_l \ur_l}\left(\bm{a}^{\top}\ur_l\right)\left(\bm{a}^{\top}\ul_l\right)\right|}, \, \|\bm{a}\|_2 \right\}, \quad & \text{otherwise.}
		\end{cases}
		\label{eq:proposed-debiased-estimate}
	\end{equation}
\noindent We shall make precise in Algorithm \ref{alg:CI-au} what it means by ``small'' in a practical and data-dependent manner. 
As before, the proposed procedures are fully data-driven, without requiring prior knowledge about the noise distributions.  A few informal yet important remarks are in order. 
\begin{itemize}
	\item When $\bm{a}^{\top}\bm{u}_l^{\star}$ is very ``small'' in magnitude,  both $\bm{a}^{\top}\ur_l$ and $\bm{a}^{\top}\ul_l$ serve as unbiased estimates of $\bm{a}^{\top}\bm{u}_l^{\star}$, and the averaging operation further reduces the uncertainty. Here, the quantity $\sqrt{\widehat{v}_{a,l}}$ in Algorithm \ref{alg:CI-au} captures the level of the smallest possible uncertainty when estimating $\bm{a}^{\top}\bm{u}_l^{\star}$. 

	\item When $\bm{a}^{\top}\bm{u}_l^{\star}$ is not very ``small'', the procedure is identical to estimator $\widehat{u}_{\bm{a},l}$ proposed previously.

\end{itemize}


\paragraph{Construction of confidence intervals.}
As it turns out, the estimation error of the modified estimator $\widehat{u}_{\bm{a},l}^{\mathsf{modified}}$ is well approximated by a zero-mean Gaussian random variable with  tractable variance $v_{\bm{a},l}^{\star}$ (to be formalized shortly). Motivated by this observation, we propose to first obtain an estimate of the variance of 
$\widehat{u}_{\bm{a},l}^{\mathsf{modified}}$ --- denoted by $\widehat{v}_{\bm{a},l}$. 
The proposed confidence interval for any given coverage level $0<1-\alpha<1$ is then given by
\begin{equation}
\mathsf{CI}_{1-\alpha}^{\bm{a}} := \left[\,\widehat{u}_{\bm{a},l}^{\mathsf{modified}} \pm\Phi^{-1}(1-\alpha/2)\sqrt{\widehat{v}_{\bm{a},l}}\,\right],
\label{eq:proposed-CI-2cases}
\end{equation}
where we abbreviate $[b\pm c] := [b-c,b+c]$, and $\Phi(\cdot)$
is the CDF of a standard Gaussian.

\paragraph{Variance estimates.}
 We shall take a moment to discuss how to obtain the variance estimate $\widehat{v}_{\bm{a},l}$. 
As will be seen momentarily, the proposed estimator $\widehat{u}_{\bm{a},l}^{\mathsf{modified}}$ admits --- modulo
some global sign --- the following first-order approximation for a broad range of settings:
\begin{align}
	\widehat{u}_{\bm{a},l}^{\mathsf{modified}} \approx\bm{a}^{\top}\bm{u}^{\star}_l + \underset{\text{uncertainty term}}{\underbrace{ \frac{1}{2\lambda^{\star}_l}\bm{a}_l^{\perp \top}\big(\bm{H}+\bm{H}^{\top}\big)\bm{u}^{\star}_l }}
	\qquad(\text{up to global sign}), 
\end{align}
where $\bm{a}^{\perp}_l := \bm{a} - (\bm{a}^{\top}\bm{u}_l^{\star}) \bm{u}_l^{\star}$.
For a broad family of noise distributions, the uncertainty term is approximately zero-mean Gaussian with variance 
\begin{equation}
v_{\bm{a},l}^{\star}:=\mathsf{Var}\left[\frac{1}{2\lambda_{l}^{\star}}(\bm{a}_{l}^{\perp})^{\top}\big(\bm{H}+\bm{H}^{\top}\big)\bm{u}_{l}^{\star}\right]=\frac{1}{4\lambda_{l}^{\star2}}\sum_{1\leq i,j\leq n}\big(a_{l,i}^{\perp}u_{l,j}^{\star}+a_{l,j}^{\perp}u_{l,i}^{\star}\big)^{2}\sigma_{ij}^{2},
	\label{eq:defn-va-star}
\end{equation}
where $a_{l,j}^{\perp}$ denotes the $j$th entry of $\bm{a}_l^{\perp}$.  
At first glance, evaluating this variance quantity precisely requires prior knowledge about
the noise level in addition to the truth $(\bm{u}^{\star}_l,\lambda^{\star}_l)$. To enable a model-agnostic and data-driven estimate of $v_{\bm{a},l}^{\star}$, we make the observation that $\frac{1}{4\lambda^{\star2}_l}\sum_{i, j}\big(a_{l,i}^{\perp}u_{l,j}^{\star}+a^{\perp}_{l,j}u_{i}^{\star}\big)^{2} H_{ij}^{2} $ is very close to $v_{\bm{a},l}^{\star}$ based on the concentration of measure phenomenon. This in turn motivates us to estimate $v_{\bm{a},l}^{\star}$ via a plug-in approach: (1) replacing
$\bm{u}_l^{\star}$ with an estimate $\widehat{\bm{u}}_l$, (2) replacing $\bm{a}_{l}^{\perp}$ with a plug-in estimate $\widehat{\bm{a}}_{l}^{\perp}$,  (3) using $\lambda_l$ in place of $\lambda^{\star}_l$, and (4) replacing $H_{ij}$ with an estimate $\widehat{H}_{ij}$, where
%
\begin{align}
	\label{eq:defn-uhat-Hhat-main}
	\widehat{\bm{u}} _l :=\tfrac{1}{\|\ur_l+\ul_l\|_2}(\ur_l+\ul_l) \qquad \text{and} \qquad 
	\widehat{\bm{H}}=[\widehat{H}_{ij}]_{1\leq i,j\leq n} := \bm{M}- \bm{M}_{\mathsf{svd},r},
\end{align}
with $\bm{M}_{\mathsf{svd},r} := \arg \min_{\text{rank}(\bm{Z})\leq r} \|\bm{M}-\bm{Z}\|_{\mathrm{F}}$ the best rank-$r$ approximation of $\bm{M}$. 
As a byproduct, this variance estimate in turn allows us to specify whether $\bm{a}^{\top}\bm{u}^{\star}_l$ is ``small'' (the case where $\bm{a}^{\top}\bm{u}^{\star}_l$ is comparable to or smaller than the typical size of the uncertainty component).  
\begin{remark}
For estimating $\bm{H}$ (or equivalently, $\bm{M}^{\star}$), it has been shown that the SVD-based approach achieves appealing entrywise accuracy (e.g. \cite{abbe2017entrywise, ma2017implicit}). 
\end{remark}
For ease of reference, the proposed procedure is summarized in Algorithm \ref{alg:CI-au}. The computational cost of this procedure mostly lies in computing the eigen-decomposition and the SVD of $\bm{M}$.

\begin{algorithm}[H]
\caption{Constructing a confidence interval for the linear form $\bm{a}^{\top}\bm{u}^{\star}_l$.}
\label{alg:CI-au}
\begin{algorithmic}[1]

\State Compute the $l$th eigenvalue $\lambda_l$ of $\bm{M}$, and the associated right eigenvector $\ur_l$ and left eigenvector $\ul_l$
	such that $\mathrm{Re}(\ul_l^{\top}\ur_l) \geq 0$ (see Notation~\ref{notation:leading-evectors}).

\State Compute the following estimator
\begin{equation}
	\widehat{u}_{\bm{a},l}^{\mathsf{modified}}:=\begin{cases}
	\frac{1}{2}\bm{a}^{\top}\big(\ur_{l}+\ul_{l}\big), & \text{if }\big|\bm{a}^{\top}\ur_{l}\big|\leq c_{\mathrm{b}}\sqrt{\widehat{v}_{\bm{a},l}}\log^{1.5} n,\\
	\min\Bigg\{ \sqrt{\left|\frac{1}{\ul_{l}^{\top}\ur_{l}}\left(\bm{a}^{\top}\ur_{l}\right)\left(\bm{a}^{\top}\ul_{l}\right)\right|},\, \|\bm{a}\|_2 \Bigg\}, & \text{else},
\end{cases}
\label{eq:defn-ua-hat}
\end{equation}
where $c_{\mathrm{b}}>0$ is some sufficiently large constant. Here, $\widehat{v}_{\bm{a},l}=\Call{Estimate-variance}{\widehat{\bm{a}}_{l}^{\perp},\widehat{\bm{u}}_l,\lambda_l}$ with 
 $\widehat{\bm{a}}_{l}^{\perp} := \bm{a}-(\bm{a}^{\top}\widehat{\bm{u}}_{l})\widehat{\bm{u}}_{l}$ and $\widehat{\bm{u}}_l \defn  \frac{1}{\|\ur_l+\ul_l\|_2} (\ur_l+\ul_l)$. 

\State For any prescribed coverage level $1-\alpha$, compute the confidence interval as 
\begin{align}
	\mathsf{CI}_{1-\alpha}^{\bm{a}}  :=\left[\, \widehat{u}_{\bm{a},l}^{\mathsf{modified}} \pm\Phi^{-1}(1-\alpha/2)\sqrt{\widehat{v}_{\bm{a},l}} \,\right].  
\end{align}

\end{algorithmic}
\end{algorithm}

\begin{algorithm}[H]
\label{alg:one-factor}
\begin{algorithmic}[1]
\Function{Estimate-variance}{$\bm{a},\bm{u},\lambda$}
\State Compute  $\widehat{\bm{H}} \defn \bm{M}- \bm{M}_{\mathsf{svd},r}$ with $\bm{M}_{\mathsf{svd},r} \defn \arg\min_{\text{rank}(\bm{Z})\leq r} \|\bm{M}-\bm{Z}\|_{\mathrm{F}}$.  
\State \Return $v := \frac{1}{4\lambda^{2}}\sum_{1\leq i,j\leq n}\big(a_{i}u_{j}+a_{j}u_{i}\big)^{2}\widehat{H}_{ij}^{2}$. 
\EndFunction
\end{algorithmic}
\end{algorithm}


\subsubsection{Confidence intervals for eigenvalues}
\label{sec:inference-eigenvalues}

Moving beyond linear forms of eigenvectors, one might also be interested in performing inference on the eigenvalues of interest.  As it turns out, this task is simpler than inferring linear functionals of eigenvectors; one can simply estimate $\lambda_l^{\star}$ via the $l$th eigenvalue $\lambda_l$ (see Notation~\ref{notation:leading-evectors}) and compute a confidence interval based on the distributional characterization of $\lambda_l$.  There is absolutely no need for careful de-biasing, since the eigen-decomposition approach automatically exploits the asymmetry structure of noise to suppress bias. Similar to Section~\ref{sec:CI-linear-form}, our procedure for performing inference on $\lambda_l^{\star}$ is distribution-free (i.e.~it does not require prior knowledge about the noise variance)
and adaptive to heteroscedasticity of noise.

To be more specific, a crucial observation, which we will make precise shortly, is as follows 
\begin{align}
	\lambda_l ~\approx~ \lambda_l^{\star} + \bm{u}_l^{\star\top} \bm{H} \bm{u}_l^{\star},
\end{align}
where the uncertainty term $\bm{u}_l^{\star\top} \bm{H} \bm{u}_l^{\star}$ is approximately Gaussian with variance 
\begin{align}
	\label{eq:vstar-lambda-l}
	v_{\lambda,l}^{\star} := \mathsf{Var}\left[\bm{u}_{l}^{\star\top}\bm{H}\bm{u}_{l}^{\star}\right]=\sum_{1\leq i,j\leq n}\big(u_{l,i}^{\star}u_{l,j}^{\star}\big)^{2}\sigma_{ij}^{2}.	 
\end{align}
Similar to the estimator $\widehat{v}_{\bm{a},l}$, we proposee to estimate $v_{\lambda,l}^{\star}$ via the following estimator
\begin{align}
	\label{eq:v-hat-lambda-defn}
	\widehat{v}_{\lambda,l}:= \sum_{1\leq i,j\leq n}\big(\widehat{u}_{l,i} \widehat{u}_{l,j} \big)^{2}\widehat{H}_{ij}^{2} ,
\end{align}
with $\widehat{\bm{u}}_l$ and $\widehat{\bm{H}}$ defined in \eqref{eq:defn-uhat-Hhat-main}. All this suggests the following confidence interval for $\lambda_l^{\star}$: 
\begin{equation}
\mathsf{CI}_{1-\alpha}^{\lambda} := \left[\,\lambda_{l} \pm \Phi^{-1}(1-\alpha/2)\sqrt{\widehat{v}_{\lambda,l}}\,\right].
\label{eq:proposed-CI-2cases}
\end{equation}
See Algorithm~\ref{alg:CI-lambda} for the complete procedure for performing inference on $\lambda_l^{\star}$.

\begin{algorithm}[H]
	\caption{Constructing a confidence interval for $\lambda_l^{\star}$.}
\label{alg:CI-lambda}
\begin{algorithmic}[1]

\State Compute the $l$th eigenvalue $\lambda_l$ of $\bm{M}$, and the associated right eigenvector $\ur_l$ and left eigenvector $\ul_l$
	such that $\ul_l^{\top}\ur_l \geq 0$ (see Notation~\ref{notation:leading-evectors}).

\State For any prescribed coverage level $1-\alpha$, compute the confidence interval as 
\begin{align}
	\mathsf{CI}_{1-\alpha}^{\lambda} := \left[\,\lambda_{l} \pm \Phi^{-1}(1-\alpha/2)\sqrt{\widehat{v}_{\lambda,l}}\,\right].
\end{align}
Here, $\widehat{v}_{\lambda,l}=\Call{Estimate-variance}{\widehat{\bm{u}}_l,\widehat{\bm{u}}_l,1}$ with 
$\widehat{\bm{u}}_l \defn  \frac{1}{\|\ur_l+\ul_l\|_2} (\ur_l+\ul_l)$.

\end{algorithmic}
\end{algorithm}

\subsection{Theoretical guarantees}
\label{sec:theory-inference}

We now set out to justify the validity of the proposed confidence intervals. Before proceeding, we need to impose another set of assumptions on the noise levels:  
\begin{assumption}
	\label{assumption:noise-size-rankr}
	Suppose that the noise parameters defined in Assumption \ref{assumption-H} satisfy
\begin{subequations}
	\begin{align}
		\sigma_{\max}\kappa^{3}\sqrt{\mu^3r^{3}} \log^{1.5} n& =o(\Delta_{l}^{\star}),  \label{eq:condition-Delta-noise} \\
		\sigma_{\max}\mu\kappa^2 r^2\sqrt{n}\log^{1.5} n & =o(\lambda_{\min}^{\star}),  \label{eq:condition-lambdamin-rankr}\\
\frac{\sigma_{\max}}{\sigma_{\min}}\asymp1, & \quad B=o\left(\sigma_{\min}\sqrt{\frac{n}{\mu\log n}}\right) .
	\label{eq:assumption-sigma-B-new}
	\end{align}
\end{subequations}
\end{assumption}
%

\noindent Armed with the above assumptions, we are ready to present theoretical guarantees. 
\begin{theorem}[\textbf{Validity of confidence intervals~(rank-$r$)}]
\label{thm:confidence-interval-validity-rankr}
	Assume $\bm{M}^{\star}$ is rank-$r$ and $\mu$-incoherent with $\mu^{4} \kappa^8 r^2 \log^{2}n=o(n)$. Given any integer $1\leq l\leq r$, 
	under Assumptions \ref{assumption-H} and \ref{assumption:noise-size-rankr}, the confidence interval returned by Algorithm~\ref{alg:CI-lambda} obeys
	\begin{align}
  		\mathbb{P}\big\{ \lambda^{\star}_l \in\mathsf{CI}_{1-\alpha}^{\lambda} \big\} = 1-\alpha+o(1) 
	\end{align}
	uniformly over all $0<\alpha<1$.  
	In addition, fix an arbitrarily small constant $0<\epsilon<1$, and consider any fixed vector $\bm{a}$ with $\ltwo{\bm a} = 1$ obeying  
	\begin{subequations}
	\label{eq:au-constraint-rankr-general}
	\begin{align}
		\big|\bm{a}^{\top}\bm{u}_{l}^{\star}\big| &\leq 1-\epsilon, \quad
		\big|\bm{a}^{\top}\bm{u}_{l}^{\star}\big|  =o\left(\frac{\Delta_l^{\star 2} }{|\lambda_l^\star|  \sigma_{\mathrm{max}} \kappa^4 r^2 \mu \log n} \right),  \label{eq:aul-target}  \\
		\big|\bm{a}^{\top}\bm{u}_{k}^{\star}\big| & =o\left(\frac{|\lambda_{l}^{\star}-\lambda_{k}^{\star}|}{|\lambda_{l}^{\star}|\sqrt{\mu\kappa^{4}r^{3}\log n}}\right),\qquad\forall k\neq l.
		\label{eq:auk-interfere}
	\end{align}
	\end{subequations}
	Then the confidence interval returned by Algorithms~\ref{alg:CI-au} obeys 
\begin{align}
	\mathbb{P}\left\{ \bm{a}^{\top}\bm{u}^{\star}_l \in\mathsf{CI}_{1-\alpha}^{\bm{a}}\right\} =  1-\alpha  +o(1)\quad\text{or}\quad\mathbb{P}\left\{ -\bm{a}^{\top}\bm{u}^{\star}_l\in\mathsf{CI}_{1-\alpha}^{\bm{a}}\right\} =1-\alpha+o(1) 
\end{align}
	uniformly over all $0<\alpha<1$. 
\end{theorem}
%

We immediately make note of a direct consequence of Theorem~\ref{thm:confidence-interval-validity-rankr} (by taking $\bm{a}$ to be the $k$th standard basis vector $\bm{e}_k$), which concerns entrywise confidence intervals for $\bm{u}^{\star}_l$. 
\begin{corollary}[\textbf{Validity of entrywise confidence intervals~(rank-$r$)}]
\label{thm:confidence-interval-validity-rankr-entrywise}
	Assume that $\bm{M}^{\star}$ is rank-$r$ and $\mu$-incoherent with $\mu^{4}=o(n/\log n)$. Suppose that Assumptions \ref{assumption-H} and \ref{assumption:noise-size-rankr} hold and that $\sqrt{\frac{\mu^{2}\kappa^{6}r^{6}\log n}{n}}\big|\lambda_{l}^{\star}\big|=o\big(\Delta_{l}^{\star}\big)$. 
	Then the confidence interval constructed in Algorithms~\ref{alg:CI-au} with $\bm{a}=\bm{e}_k$ satisfies
\begin{align}
	\mathbb{P}\left\{ {u}^{\star}_{l,k} \in\mathsf{CI}_{1-\alpha}^{\bm{\bm{e}_k}}\right\} =  1-\alpha  +o(1)\quad\text{or}\quad\mathbb{P}\left\{ -{u}^{\star}_{l,k} \in\mathsf{CI}_{1-\alpha}^{\bm{e}_k}\right\} =1-\alpha+o(1); 
\end{align}
this holds true uniformly over all $0<\alpha<1$. 
\end{corollary}

The above results confirm the validity of our inferential procedure in high dimension. 
Our theory applies to a very broad family of noise distributions, without the need of any prior knowledge about detailed noise distributions or noise levels $\{\sigma_{ij}\}$. The results are fully adaptive to heteroscedasticity of noise,  making them appealing for practical scenarios. 
In the sequel, we discuss several important implications of our results in a more quantitative manner. For simplicity of discussion, we shall concentrate on the case where $r,\kappa,\mu\asymp 1$ and consider Gaussian noise with $\sigma_{\max}\asymp \sigma_{\min}$.  We shall also assume $\|\bm{a}\|_2=1$ without loss of generality. 

\begin{itemize}

	\item We start with the rank-1 case (i.e.~$r=1$), in which we have $\Delta_l^{\star}=\infty$ according to the definition \eqref{eq:eigengap}. 
		In this case, the conditions \eqref{eq:au-constraint-rankr-general} admit considerable simplification:   
	\begin{align}
		|\bm{a}^{\top}\bm{u}_1^{\star}| \leq 1-\epsilon
	\end{align}
	for any small constant $\epsilon>0$. 
	This means that our theory covers a very wide range of linear functionals $\bm{a}^{\top}\bm{u}^{\star}_1$; basically, the proposed confidence interval finds theoretical support unless the preconceived direction $\bm{a}$ is already highly aligned with the truth $\bm{u}^{\star}_1$.

	\item Going beyond the rank-1 case, the eigen-gap condition imposed in Assumption~\ref{assumption:noise-size-rankr} and Corollary~\ref{thm:confidence-interval-validity-rankr-entrywise} reads
	\begin{align}
		\Delta_{l}^{\star}\gtrsim\sigma_{\max}\mathrm{poly}\log n\asymp\frac{\|\bm{H}\|\mathrm{poly}\log n}{\sqrt{n}} ,		
	\end{align}
	which is allowed to be substantially smaller than the spectral norm $\|\bm{H}\|$ of the noise matrix. 
		
\item Under such eigen-gap requirements, we develop an informative distributional theory underlying the proposed estimators $\widehat{u}_{\bm{a},l}$, provided that
	$\{\bm{a}^{\top}\bm{u}_k^{\star}\}_{1\leq k\leq r}$ are not too large in the sense that 
	%
	\begin{align}
		|\bm{a}^\top \bm{u}_l^\star|  \lesssim  \frac{{\Delta_l^\star}^2}{|\lambda_l^\star| \sigma_{\mathrm{max}}  \log^{1.5} n}; \qquad |\bm{a}^\top \bm{u}_k^\star| \lesssim \frac{|\lambda_l^\star -\lambda_k^\star|}{|\lambda_l^\star| \log n}, \quad k \neq l.
		\label{eq:au-simplified}	
	\end{align}
	%
	In words, the condition \eqref{eq:au-simplified} requires that both the target quantity $\bm{a}^{\top}\bm{u}_{l}^{\star}$ and the ``interferers'' are dominated by the respective (normalized) eigen-gaps.  An important instance that automatically satisfies such conditions has been singled out in Corollary~\ref{thm:confidence-interval-validity-rankr-entrywise}, leading to appealing entrywise distributional characterizations and inferential procedures.

	\item While our theorems focus on asymmetric noise matrices, one can combine them with a simple asymmetrization trick to yield valid confidence intervals when $\bm{H}$ is a symmetric matrix with homoscedastic Gaussian noise.  See Appendix~\ref{sec:example-homoscedastic-Gaussian} for detailed discussion. 

\end{itemize}

We note, however, that the additional requirement \eqref{eq:au-simplified} makes our inference results less general than our estimation guarantees in Section~\ref{sec:estimation-theory}, except for the rank-1 case. In fact,  if \eqref{eq:au-simplified} is violated, then the influence of the eigen-gaps might become the dominant factor in the uncertainty term and needs to be quantified in a precise fashion. Achieving this calls for more refined theoretical analysis, and we leave it to future investigation.

\subsection{Key ingredients behind Theorem~\ref{thm:confidence-interval-validity-rankr}} 
Next, we single out two key ingredients that shed light on not only the validity of, but also the statistical efficiency of, our inferential procedures. To be specific, 
Theorem~\ref{thm:confidence-interval-validity-rankr} is mainly built upon distributional guarantees developed for the proposed estimator $\widehat{u}_{\bm{a},l}^{\mathsf{modified}}$ and the $l$th eigenvalue $\lambda_l$ of $\bm{M}$. In a nutshell,  the quantity $\widehat{u}_{\bm{a},l}^{\mathsf{modified}}$ (resp.~$\lambda_l$) is a nearly unbiased estimator of $\bm{a}^{\top}\bm{u}_l^{\star}$ (resp.~$\lambda_l^{\star}$) and is approximately Gaussian. We formalize this distributional theory as follows, whose proof is deferred to Appendix \ref{sec:proof-thm-distribution-validity-rankr}.

\begin{theorem}[\textbf{Distributional characterization (rank-$r$)}]
\label{thm:distribution-validity-rankr}
Instate the assumptions of Theorem~\ref{thm:confidence-interval-validity-rankr}.
Let $\widehat{u}_{\bm{a},l}^{\mathsf{modified}}$ (cf.~Algorithm \ref{alg:CI-au})
be the proposed estimate for $\bm{a}^{\top}\bm{u}^{\star}_l$. 
Then one can write
\begin{subequations}
\begin{align}
	\widehat{u}_{\bm{a},l}^{\mathsf{modified}} = \bm{a}^{\top}\bm{u}^{\star}_l +  \sqrt{{v}_{\bm{a},l}^{\star}}\,W_{\bm{a},l}+  \zeta
	\quad\text{or}\quad
	- \widehat{u}_{\bm{a},l}^{\mathsf{modified}} =\bm{a}^{\top}\bm{u}^{\star}_l +\sqrt{{v}_{\bm{a},l}^{\star}}\,W_{\bm{a},l}+ \zeta
\end{align}
and
\begin{align}
 	\lambda_l =\lambda^{\star}_l  + \sqrt{v_{\lambda,l}^{\star}}\,W_{\lambda,l}+\xi, 
	\label{eq:ua-hat-expansion-rankr}
\end{align}
\end{subequations}
where ${v}_{\bm{a},l}^{\star}$ and $v_{\lambda,l}^{\star}$ are defined respectively in  \eqref{eq:defn-va-star} and \eqref{eq:vstar-lambda-l}, and 
\begin{align}
	W_{\bm{a},l}:= \frac{(\bm{a}_l^{\perp})^\top\left(\bm{H}+\bm{H}^{\top}\right)\bm{u}^{\star}_l}{2\lambda^{\star}_l\sqrt{{v}_{\bm{a},l}^{\star}}}
	\qquad \text{and}\qquad
	W_{\lambda,l}:=\frac{\bm{u}_l^{\star\top}\bm{H}\bm{u}_l^{\star}}{\sqrt{v_{\lambda,l}^{\star}}}
\end{align}
with $\bm{a}_l^{\perp}:=\bm{a}-(\bm{a}^{\top}\bm{u}_l^{\star})\bm{u}_l^{\star}$.  
The residual terms  obey $|\zeta|=o\big(\sqrt{{v}_{\bm{a},l}^{\star}}\big)$ and $|\xi|=o\big(\sqrt{{v}_{\lambda,l}^{\star}}\big)$ with probability at least $1-O(n^{-5})$. In addition, one has
	\begin{align}
		\label{eq:property-two-inprod-rankr}
		\sup_{z\in\real}\big|\,\mprob(W_{\bm{a},l}\leq z)-\Phi(z)\,\big|=o(1)\qquad\text{and}\qquad\sup_{z\in\real}\big|\,\mprob(W_{\lambda,l}\leq z)-\Phi(z)\,\big|=o(1). 		
	\end{align}
%
\end{theorem}

Informally, this theorem reveals the tightness of the following first-order approximation
\begin{align}
	\widehat{u}_{\bm{a},l}^{\mathsf{modified}} \approx \bm{a}^{\top} \bm{u}^{\star}_l+\frac{1}{2\lambda_l^{\star}}(\bm{a}_l^{\perp})^\top(\bm{H}+\bm{H}^{\top})\bm{u}^{\star}_l ~~(\text{up to global sign})
	\qquad\text{and}\qquad
	\lambda_l\approx\lambda_l^{\star}+\bm{u}_l^{\star\top}\bm{H}\bm{u}^{\star}_l	
\end{align}
for a wide range of directions $\bm{a}$. In addition, the first-order approximations are close in distribution to Gaussian random variables.  
Such distributional characterizations would immediately lead to $(1-\alpha)$-confidence intervals for any $\alpha$, if an oracle had informed us of the quantities $v_{\bm{a},l}^{\star}$ and $v_{\lambda,l}^{\star}$. Fortunately, the proposed 
$\widehat{v}_{\bm{a},l}$ and $\widehat{v}_{\lambda,l}$ (see Algorithms~\ref{alg:CI-au} and \ref{alg:CI-lambda}) serve as highly accurate estimates of  $v_{\bm{a},l}^{\star}$ and $v_{\lambda,l}^{\star}$, respectively, and can be employed in place of $v_{\bm{a},l}^{\star}$ and $v_{\lambda,l}^{\star}$. This observation is asserted by the following theorem; the proof is postponed to Appendix \ref{sec:proof-lemma-var-control}.

\begin{theorem}[\textbf{Accuracy of variance estimates~(rank-$r$)}]
\label{Lem-var-control-rankr}	
Instate the assumptions of Theorem~\ref{thm:confidence-interval-validity-rankr}. 
With probability at least $1 - O(n^{-10})$, the variance estimators proposed in Algorithms~\ref{alg:CI-au}-\ref{alg:CI-lambda} satisfy
\begin{align}
	\widehat{v}_{\bm{a},l} = (1+ o(1)) v_{\bm{a},l}^{\star} 
	\qquad \text{and} \qquad 
	\widehat{v}_{\lambda,l} = (1+ o(1)) v_{\lambda,l}^{\star}.
\end{align}
\end{theorem}
Clearly, putting Theorems~\ref{thm:distribution-validity-rankr}-\ref{Lem-var-control-rankr} together immediately establishes Theorem~\ref{thm:confidence-interval-validity-rankr}.

\section{Numerical experiments}
\label{SecNumericals}

This section consists of numerical experiments in various settings, in order to demonstrate the performance
of our estimation and inference procedures and their accompanying theory.

\subsection{Eigen-decomposition after symmetrization?}

As mentioned in the discussions after Theorem~\ref{thm:rankr-evalue-bounds-simple}, one may consider first symmetrizing the data matrix 
with $\frac{1}{2}(\bm{M}+\bm{M}^{\top})$ before computing the eigen-decomposition. 
It is unclear whether this procedure provides sensible eigenvector estimators in the presence of small eigen-gaps and heteroscedastic noise. 
This subsection is devoted to understanding the potential sub-optimality of this approach via a simple example. Specifically,  consider the rank-$2$ model where 
\begin{align}
\label{eq:Mstar-rank2}
\bm{M}^\star = \lambda_1^{\star}\bm{u}_1^\star \bm{u}_1^{\star \top} + \lambda_2^{\star}\bm{u}_2^\star \bm{u}_2^{\star \top},
\end{align}
with $\bm{H}$ satisfying our usual assumptions. To simplify presentation, we define 
\begin{enumerate}
	\item[1.] $\mathsf{Spectral}\text{-}\mathsf{asym}$: eigen-decomposition applied to the observed asymmetric data matrix $\bm{M}$;
	\item[2.] $\mathsf{Spectral}\text{-}\mathsf{sym}$: eigen-decomposition applied to the symmetrized data matrix $\frac{1}{2}\left(\bm{M} + \bm{M}^{\top}\right)$.
\end{enumerate}

\begin{figure}
	\centering
	\begin{tabular}{cc}
	\includegraphics[width=0.4\textwidth]{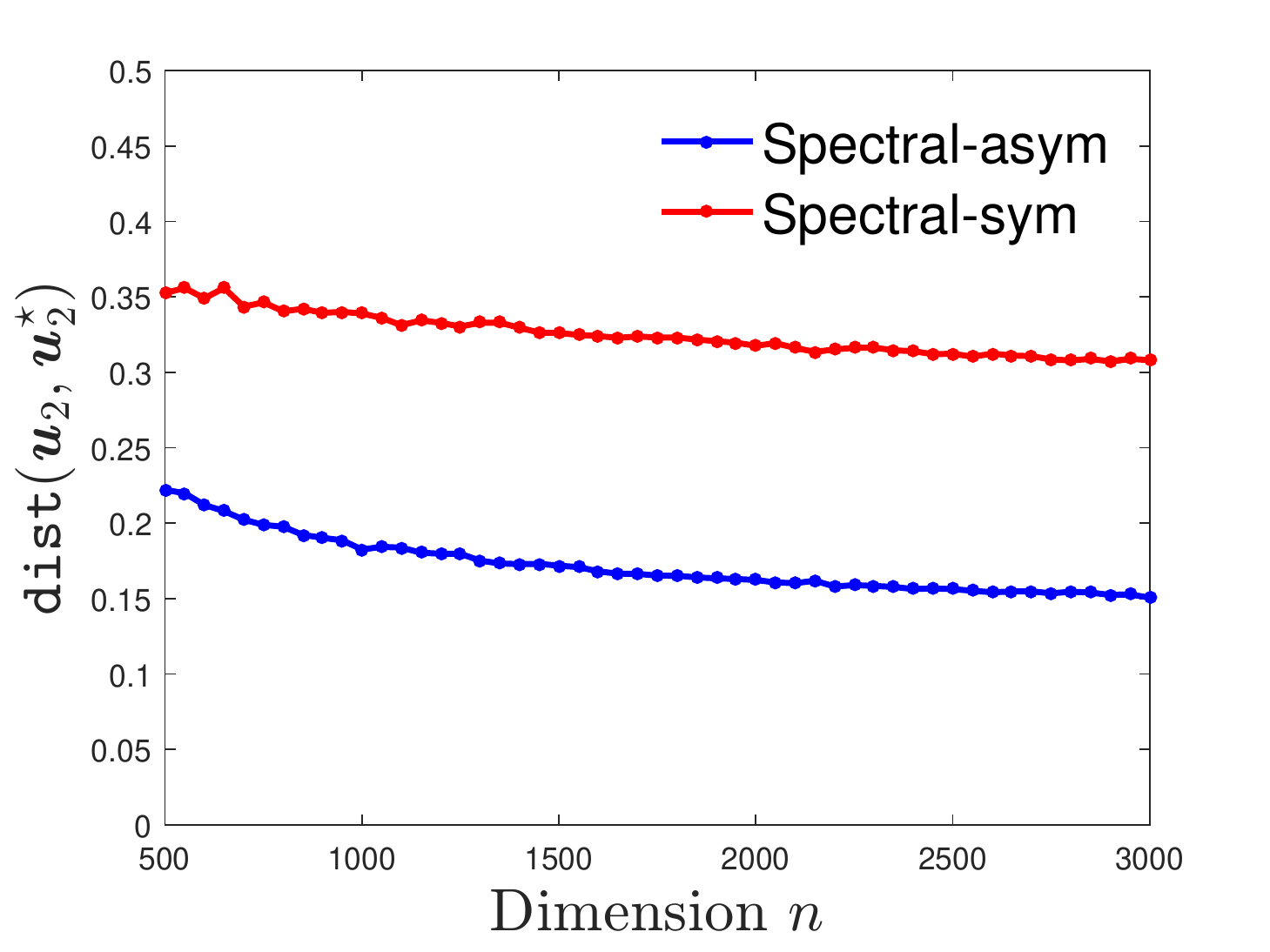} &
	\includegraphics[width=0.4\textwidth]{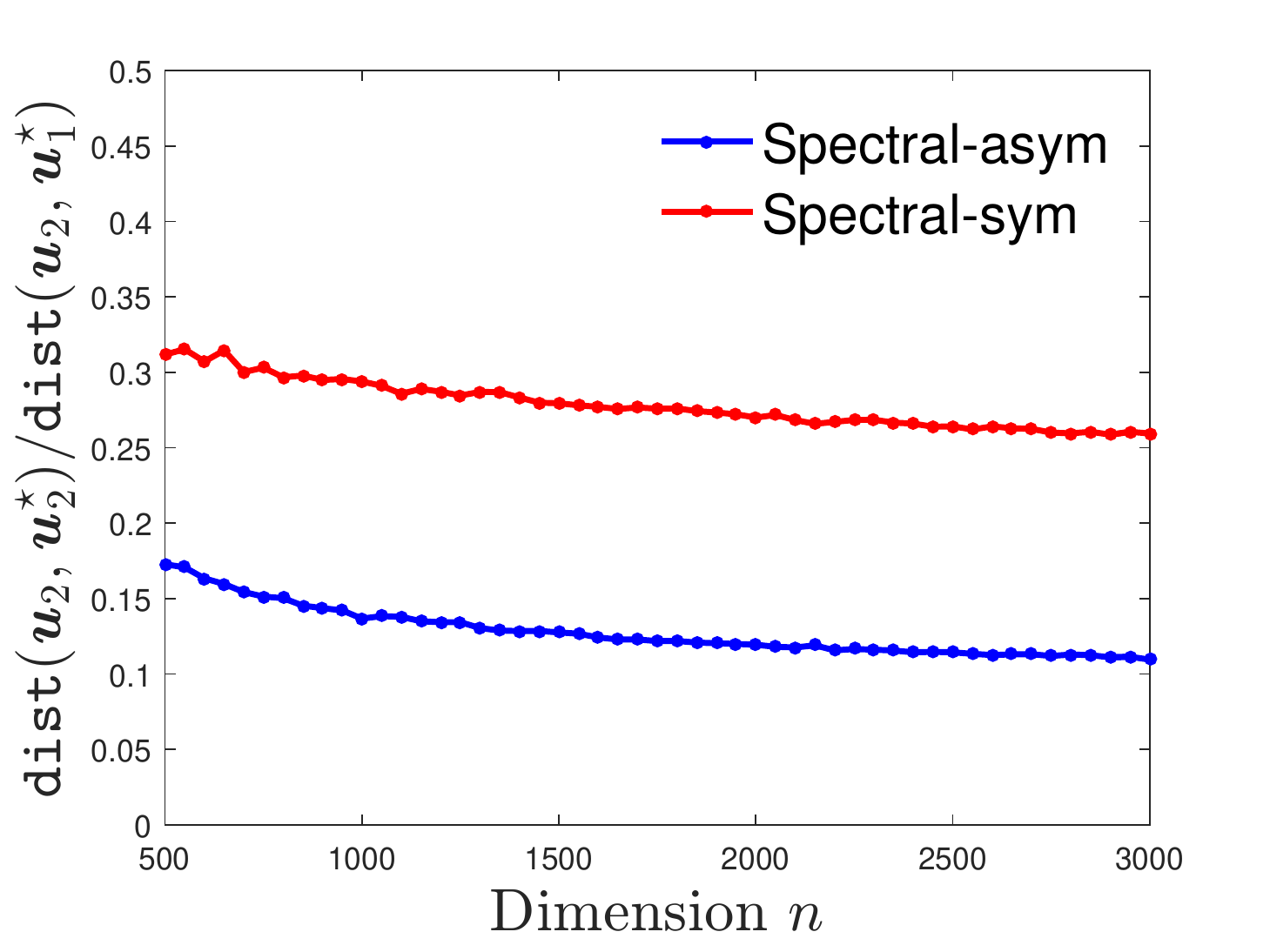} \tabularnewline
		(a) $\ell_2$ estimation error $\mathsf{dist}(\bm{u}_2 , \bm{u}_2^\star)$ & (b) relative $\ell_2$ error $\mathsf{dist}(\bm{u}_2, \bm{u}_2^\star) / \mathsf{dist}( \bm{u}_2, \bm{u}_1^\star )$ \tabularnewline
	\end{tabular}
	\caption{Numerical performance of $\mathsf{Spectral}\text{-}\mathsf{asym}$ vs.~$\mathsf{Spectral}\text{-}\mathsf{sym}$ when $\lambda_1^\star = 1$ and $\lambda_2^\star = 0.95$. Here, we define $\mathsf{dist}(\bm{u},\bm{v}):=\min\{ \|\bm{u}-\bm{v}\|_2, \|\bm{u}+\bm{v}\|_2 \}$.  (a) plots the $\ell_2$ estimation error vs.~the dimension $n$, whereas (b) plots the relative $\ell_2$ error vs.~$n$. For each $n$, the results are averaged over $1000$ trials, with $\sigma_1 = {1}/{\sqrt{n \log n}}$ and $\sigma_2 = {0.1}/{\sqrt{n \log n}}$ (see~the expression \eqref{eq:noise-variance-example-rank2}). }
	\label{fig:symvsasy}	
\end{figure}

\paragraph{A numerical example with heteroscedastic Gaussian noise.} We start by looking at an example with  
\begin{equation}
	\bm{u}_1^\star = \frac{1}{\sqrt{n}} \bm{1}_n \qquad \text{and} \qquad \bm{u}_2^\star = \frac{1}{\sqrt{n}} \begin{bmatrix}
\bm{1}_{n/2} \\
-\bm{1}_{n/2}
\end{bmatrix}.
	\label{eq:u1-u2-example-rank2}
\end{equation}
The noise matrix $\bm{H}$ is assumed to have independent zero-mean Gaussian entries with variance
\begin{equation}
	\label{eq:noise-variance-example-rank2}
	\mathsf{Var}(\bm{H}) := \big[ \mathsf{Var}(H_{ij}) \big]_{1\leq i,j\leq n} = \sigma_1^2 \begin{bmatrix}
		\bm{1}_{n/2 }\bm{1}_{n/2 }^{\top} - \frac{1}{2} \bm{I}_{n/2} & \bm{0} \\
\bm{0} & \bm{0}
	\end{bmatrix} + \sigma_2^2 \left(\bm{1}_{n }\bm{1}_{n}^{\top} - \frac 1 2 \bm{I}_n\right),
\end{equation}
and hence the variance of the symmetrized data satisfies
\begin{equation}
	\mathsf{Var}\left(\frac{\bm{H} + \bm{H}^\top}{2}\right) := \big[\mathsf{Var}\big(\tfrac{1}{2}(H_{ij}+H_{ji})\big)\big]_{1\leq i,j\leq n} 
	=  \frac{\sigma_1^2}{2} \begin{bmatrix}
		\bm{1}_{n/2 }\bm{1}_{n/2 }^{\top} & \bm{0} \\
\bm{0} & \bm{0}
	\end{bmatrix} + \frac{\sigma_2^2}{2} \bm{1}_{n }\bm{1}_{n }^{\top}.
\end{equation}
We plot in Fig.~\ref{fig:symvsasy} the numerical performance of both $\mathsf{Spectral}\text{-}\mathsf{asym}$ and $\mathsf{Spectral}\text{-}\mathsf{sym}$ in estimating $\bm{u}_2^{\star}$. For $\mathsf{Spectral}\text{-}\mathsf{asym}$, the estimator $\widehat{\bm{u}}_{2}$ is constructed as in the expression \eqref{eq:defn-hat-ul}, whereas the second eigenvector of $\frac{1}{2}\left(\bm{M} + \bm{M}^{\top}\right)$ is used for $\mathsf{Spectral}\text{-}\mathsf{sym}$.
As can be seen, $\mathsf{Spectral}\text{-}\mathsf{asym}$ strictly outperforms $\mathsf{Spectral}\text{-}\mathsf{sym}$ in all cases, thus unveiling the real benefits of exploiting asymmetry in eigen-decomposition. The interested reader is deferred to Appendix~\ref{sec:interpretation-symmetrization} for some high-level interpretation of this phenomenon.  

\subsection{Estimation}

This section is devoted to numerically studying the efficiency of our estimators for linear functionals of the eigenvectors. Given any fixed $\bm a$ with $\ltwo{\bm a} = 1$, our estimator is constructed as in the expression \eqref{eq:ua-estimator-large}.
Theorem~\ref{thm:rankr-bounds-simple} (in particular, the upper bound \eqref{eq:rankr-linear-form-bound-simple}) together with Theorem~\ref{thm:minimax-l2-au-eigengap} implies that the estimation error is controlled by the eigen-gap, the true signal strength $|\bm a^\top \bm u^\star_l|$, and the magnitude of the ``interferers''. 
In the following, we examine qualitatively the effects of these quantities upon the estimation errors through some simple yet representative examples. 

\paragraph{Settings.} Consider the rank-2 case as in Eq.~\eqref{eq:Mstar-rank2} again, where 
the leading eigenvalue is set to be $1$ and the orthonormal pair $\bm{u}_1^\star$ and $\bm{u}_2^\star$ are randomly generated. 
We focus on estimating the linear functionals of the form $\bm a^{\top} \bm{u}_2^\star$. 
In the following, several quantities that are selected to examine how they affect the estimation error $|\widehat{u}_{\bm{a},2} - \bm{a}^\top \bm{u}_2^\star|$ include 
	the ground-truth $|\bm{a}^\top \bm{u}_2^\star|$, 
	the interferer $|\bm{a}^\top \bm{u}_1^\star|$, 
	and the eigen-gap $\Delta_2^\star = \lambda_1^\star - \lambda_2^\star$.
In the following, we consider two scenarios: the case where the observation matrix is generated from a rank-2 underlying matrix $\bm M^{\star}$ plus heteroscedastic Gaussian noise; and the case where entries are missing at random with sampling rate $0.1$ on top of the above-mentioned observation model (more details are given below). 

\paragraph{Heteroscedastic Gaussian noise.} 
Consider a heteroscedastic Gaussian noise matrix $\bm{H}$ with independent entries $H_{ij} \sim \mathcal{N}(0, \sigma_{ij}^2)$ obeying
\begin{align}
\mathsf{Var}(\bm{H}):=\big[\mathsf{Var}(H_{ij})\big]_{1\leq i,j\leq n}= 
{\footnotesize \begin{bmatrix}\sigma_{1}^{2}\\
	(\sigma_{1}+\delta_{\sigma})^{2}\\
\vdots\\
(\sigma_{1}+(n-1)\delta_{\sigma})^{2}
\end{bmatrix}}
\cdot  \bm{1}_{n}^{\top}, 
	\label{eq:variance-Gaussian-numerics}
\end{align}
where $\delta_{\sigma}$ determines the spacing between adjacent standard deviation. 
The parameters chosen to be
$$
\sigma_1 = 0.1/\sqrt{n \log n}, \quad \delta_{\sigma} = 0.9/((n-1)\sqrt{n \log n}).
$$

\paragraph{Missing data model.} Suppose that we only get to observe a fraction of the entries of $\bm{M}^{\star}$; more precisely, each entry of $\bm{M}^{\star}$ is observed independently with probability $p$. In this setting, we can take the observed data matrix via zero-padding and rescaling as follows
\begin{align}
\label{eq:missing-data}
M_{ij}=\begin{cases}
	\frac{1}{p}(M_{ij}^{\star}+ \sqrt p \widetilde{H}_{ij}), & \text{with probability }p,\\
0, & \text{otherwise}.
\end{cases}
\end{align}
This way we guarantee that $\mathbb{E}[\bm{M}]=\bm{M}^{\star}$, and $\widetilde{\bm{H}} = [\widetilde{H}_{ij}]_{1\leq i,j\leq n} $ is a heteroscedastic Gaussian noise with variance 
$\mathsf{Var}(\widetilde{\bm{H}}):=\big[\mathsf{Var}(\widetilde{H}_{ij})\big]_{1\leq i,j\leq n}$ equal to \eqref{eq:variance-Gaussian-numerics}. 
In this case, the matrix $\bm{H} := \bm{M} - \bm{M}^{\star} $ obeys
\begin{align*}
	|H_{ij}| \leq \left|\frac{1}{p}M_{ij}^\star\right| + \left|\frac{1}{\sqrt p} \widetilde{H}_{ij}\right|  \lesssim \frac{\mu}{np} + \frac{\widetilde{\sigma}_{\max} \sqrt{\log n}}{\sqrt p}, \qquad \mathbb{E} [H_{ij}^2] \lesssim \frac{\mu}{n^2 p} + \widetilde{\sigma}_{\max}^2 
\end{align*}
 with high probability, 
where $\widetilde{\sigma}_{\max} := \sigma_{1}+(n-1)\delta_{\sigma}$. When the sampling rate exceeds $p \geq \frac{c_{\mathrm{p}}  \mu \log^2 n}{n}$ and the noise size is below $ \widetilde{\sigma}_{\max} \leq \frac{1}{\sqrt{c_{\mathrm{p}}  n \log n}}$ for some constant $c_{\mathrm{p}}\geq 1$, one has
\begin{align*}
	|H_{ij}| \lesssim \frac{\mu}{np} +\frac{\widetilde{\sigma}_{\max} \sqrt{\log n}}{\sqrt p} \lesssim  \frac{1}{ {c_{\mathrm{p}}} \log n}, \qquad \sqrt{\mathbb{E} [H_{ij}^2] } \lesssim \frac{\sqrt \mu}{n \sqrt p} + \widetilde{\sigma}_{\max} \lesssim \frac{1}{ \sqrt{c_{\mathrm{p}} n \log n}}, 
\end{align*}
thus satisfying Assumption~\ref{assumption:noise-size-rankr-revised} for $c_{\mathrm{p}}$ sufficiently large. In the numerical experiments here, we shall choose $p = 0.1$ and
$$
\sigma_1 = 1/\sqrt{10 n \log n}, \quad \delta_{\sigma} = 9/((n-1)\sqrt{10 n \log n}).
$$

\begin{figure}[htbp!]
\centering
	\includegraphics[width=0.4\textwidth]{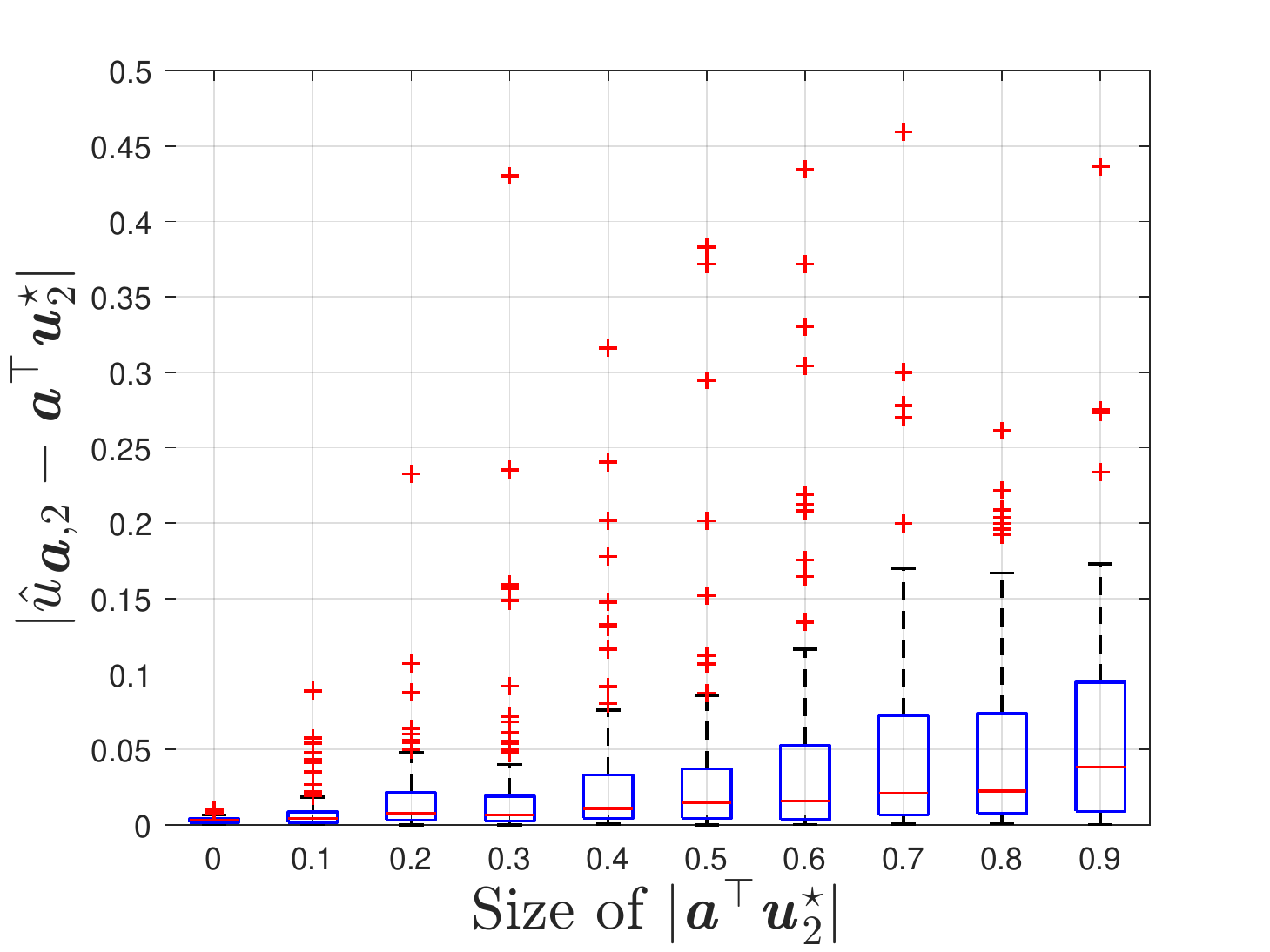}
	\includegraphics[width=0.4\textwidth]{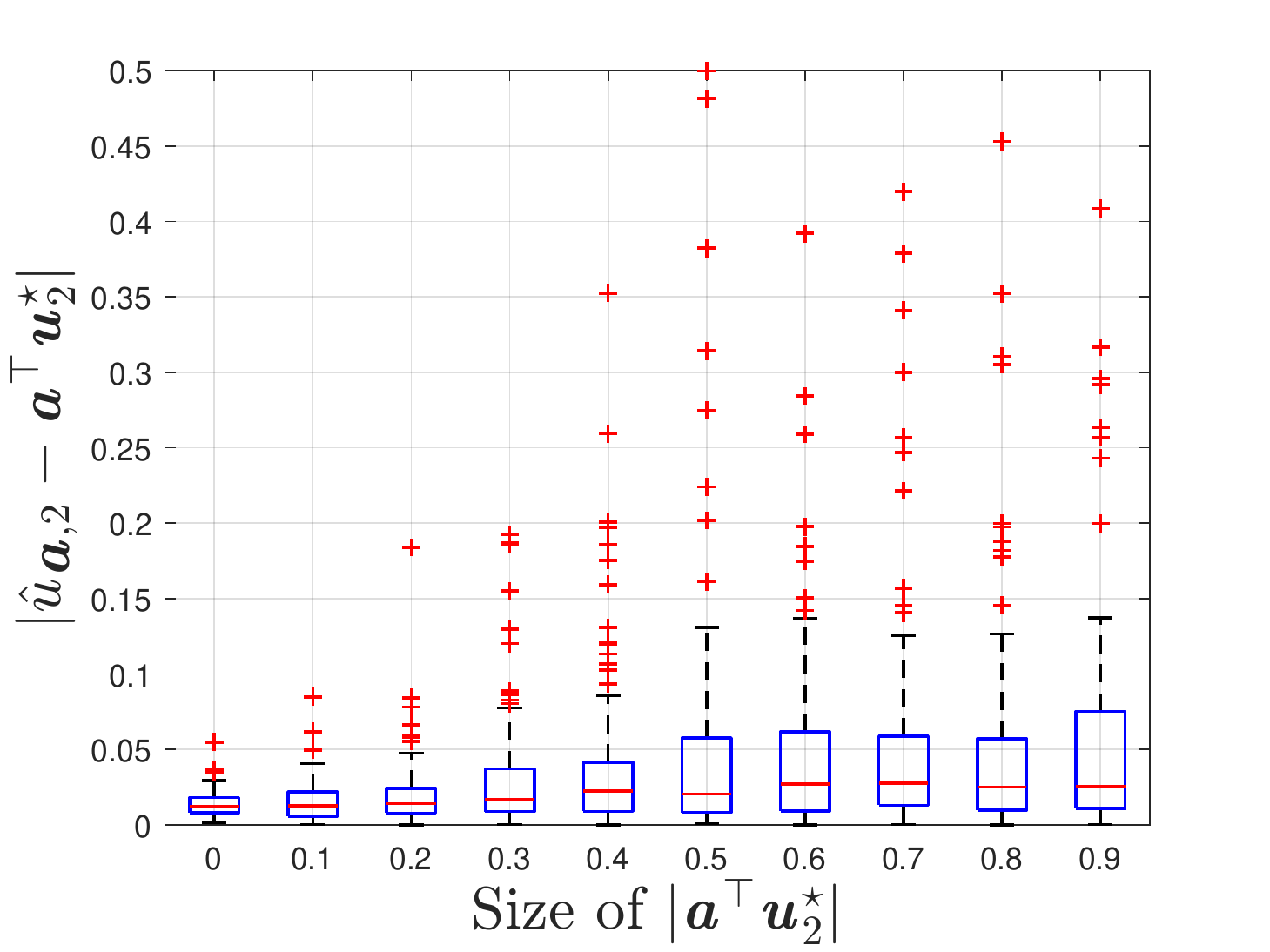}
	\caption{Estimation error vs.~size of the ground-truth $|\bm{a}^\top \bm{u}_2^\star|$. 
	The left figure corresponds to the noisy observation case where $|\Delta_2^\star| = 0.01$ and the right figure corresponds to the missing observation case where $|\Delta_2^\star| = 0.05$.
	In both cases, $\bm a$ is chosen such that $\bm a^{\top} \bm{u}_1^\star = 0$.} 
	\label{fig:estimation-vs-ground-truth}
\end{figure}

\paragraph{Estimation error vs.~size of the ground-truth.} 
To study the qualitative effect of the magnitude of the group-truth, various values of $|\bm{a}^\top \bm{u}_2^\star|$ are considered. 
For each configuration, the dimension $n$ is set to be $500$ and we run $100$ trials for the box plot. 
The experiments are conducted for a variety of eigen-gaps and directions $\bm a$. A clear positive correspondence is seen between the size $|\bm{a}^\top \bm{u}_2^\star|$ and the estimation error; see Fig.~\ref{fig:estimation-vs-ground-truth}.

\begin{figure}[H]
		\centering
		\includegraphics[width=0.4\textwidth]{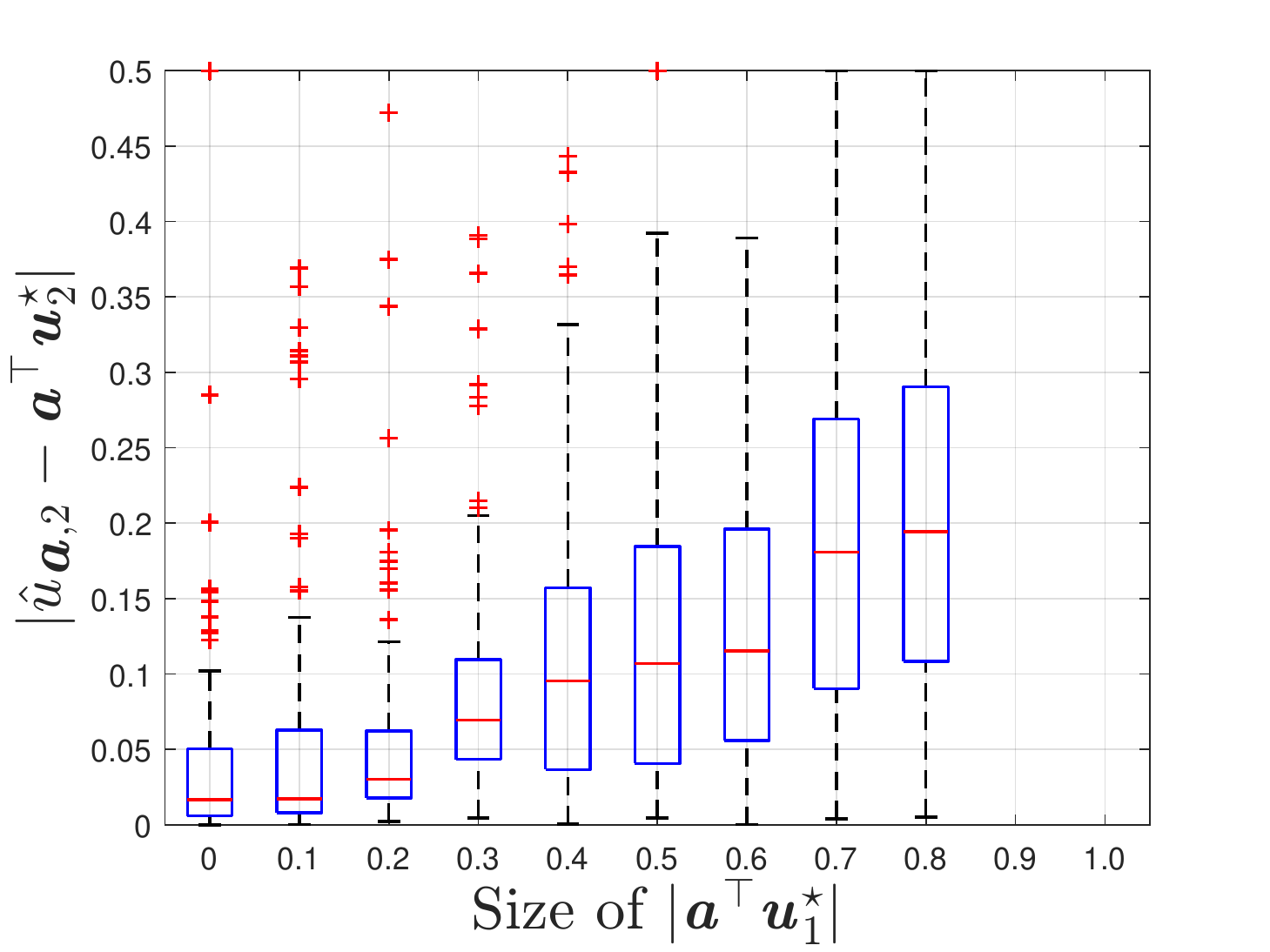}
		\includegraphics[width=0.4\textwidth]{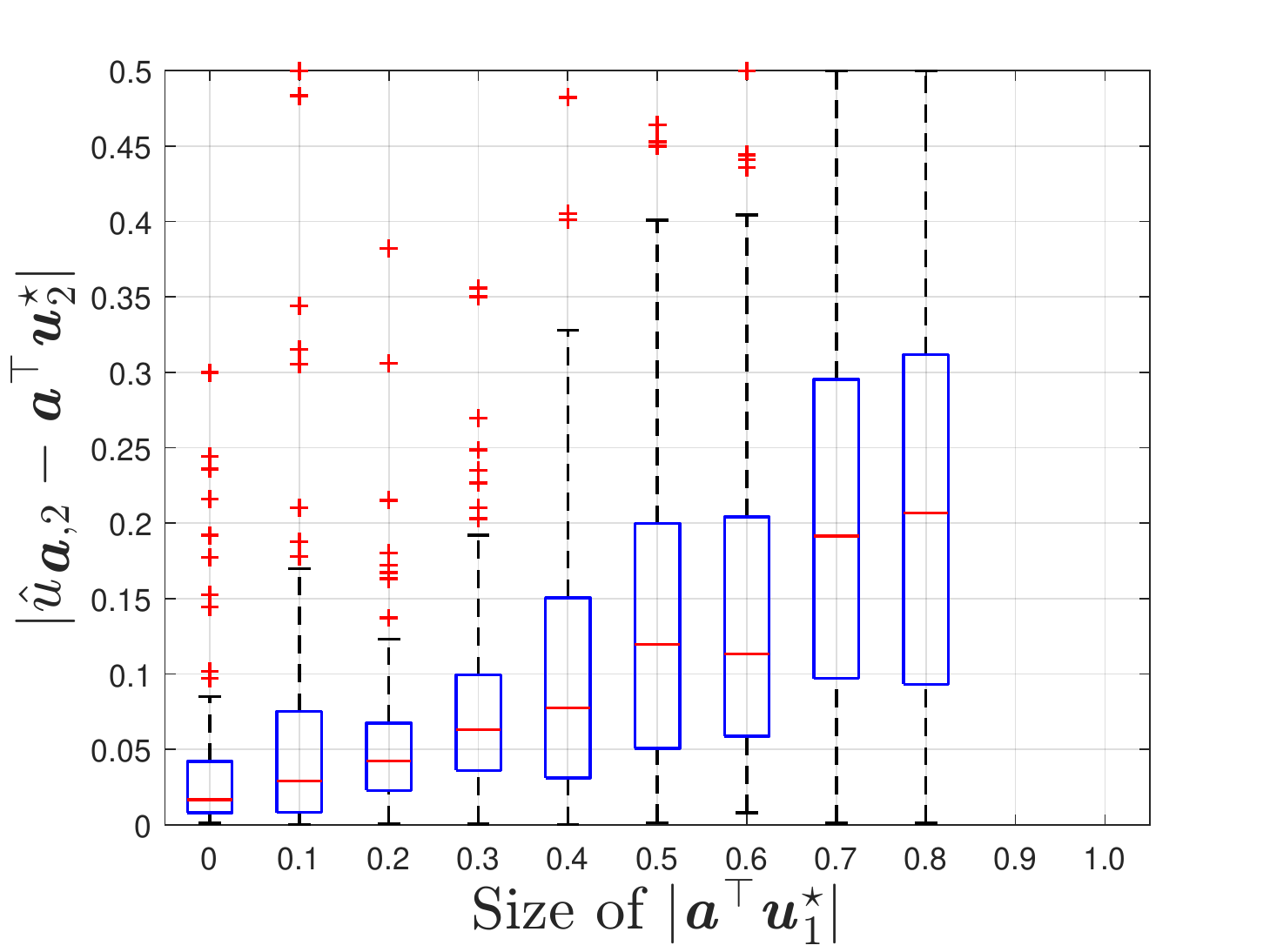}
		\caption{Estimation error vs.~size of the interferer $|\bm{a}^\top \bm{u}_1^\star|$.
		The left figure corresponds to the full observation case where $|\Delta_2^\star| = 0.01$, while the right figure corresponds to the missing observation case where $|\Delta_2^\star| = 0.05$.
		} \label{fig:estimation-vs-interferer}
	\end{figure}

\paragraph{Estimation error vs size of the interferer.} 
To study the qualitative effect of the magnitude of the interferers, we consider a range of values for $|\bm{a}^\top \bm{u}_1^\star|$ while holding $|\bm{a}^\top \bm{u}_2^\star| = 0.5$ unchanged. 
Again, for each configuration, the dimension $n$ is set to be $500$ and we run $100$ trials for the box plot. 
The experiments are run with various values of eigen-gaps. A negative dependency of the estimation error on the interferer is observed. In particular, Fig.~\ref{fig:estimation-vs-interferer} illustrates how the  estimation errors grow as the interferer gets stronger.
	
\begin{figure}[H]
		\centering
		\includegraphics[width=0.4\textwidth]{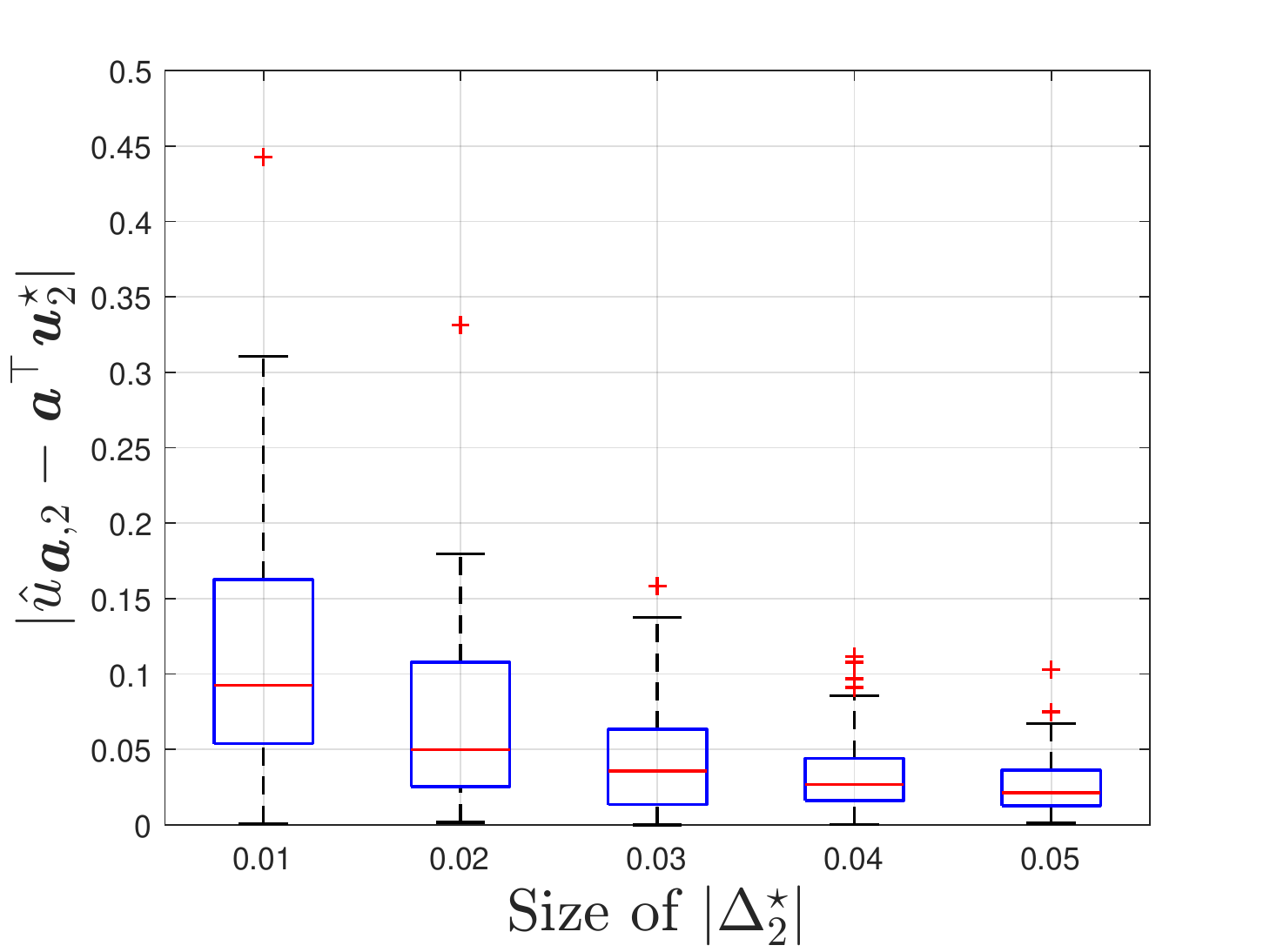} 
		\includegraphics[width=0.4\textwidth]{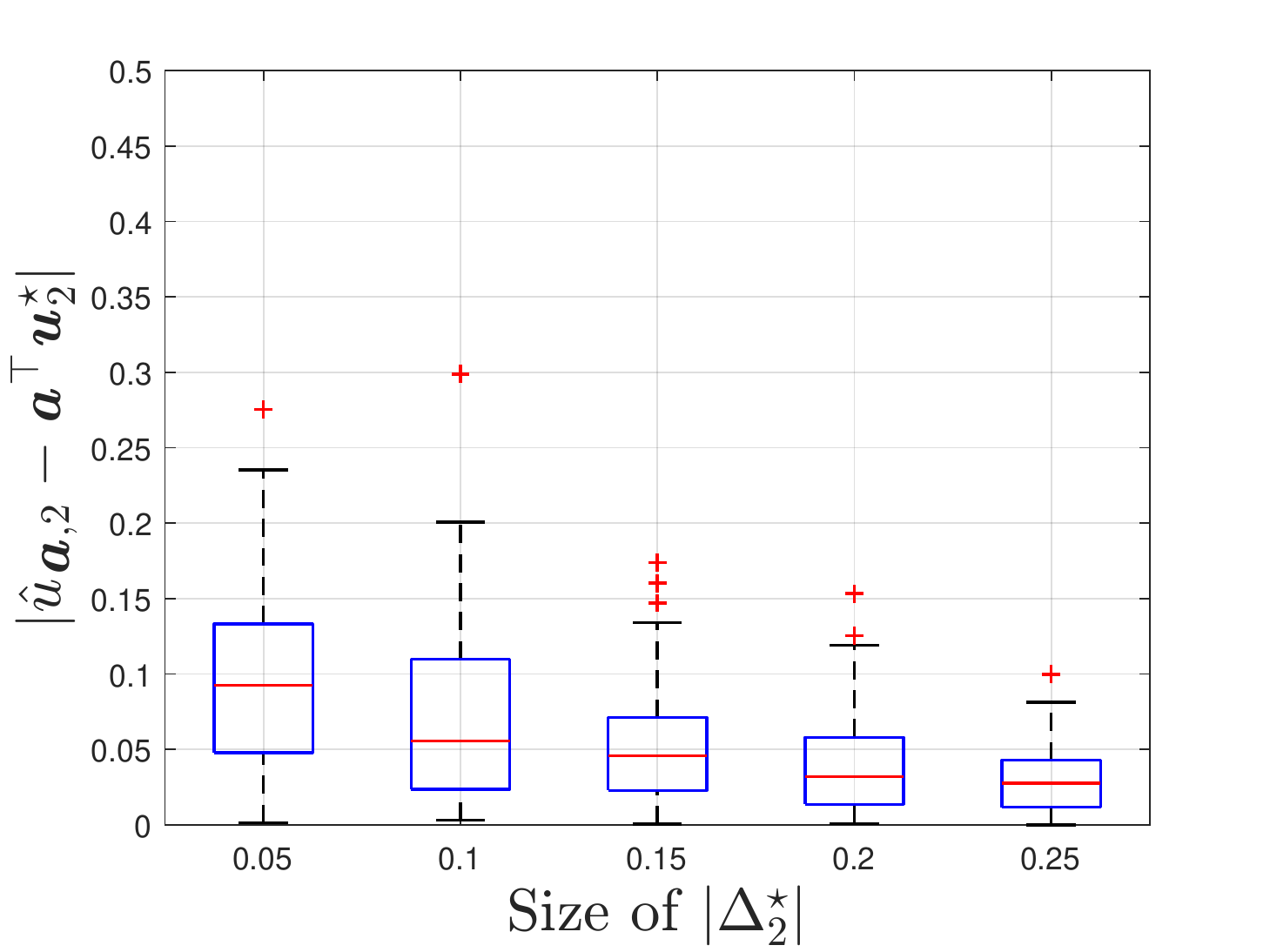} 
		\caption{Estimation error vs eigen-gap $|\Delta_2^\star|$.
		The left figure corresponds to the full observation case, while the right figure corresponds to the missing observation case.
		}\label{fig:estimation-vs-eigengap}
	\end{figure}

\paragraph{Estimation error vs eigen-gap.} 
In this part, we consider various magnitudes of the eigen-gap and study how the estimation error is influenced. 
Similar to what Theorem~\ref{thm:rankr-bounds-simple} and Theorem~\ref{thm:minimax-l2-au-eigengap} predict, 
the estimation task becomes easier when the eigen-gap gets larger. 
Fig.~\ref{fig:estimation-vs-eigengap} manifests this relationship in the case when $|\bm{a}^\top \bm{u}_1^\star| = 0.5$ and $|\bm{a}^\top \bm{u}_2^\star| = 0.2$.



\subsection{Uncertainty quantification}

This subsection performs numerical experiments to validate our inferential procedures as well as the accompanying theory. Throughout this subsection, we consider the rank-$2$ case again as in the expression \eqref{eq:ua-estimator-large}, where $\bm{M}^\star \defn \lambda_1^{\star}\bm{u}_1^\star \bm{u}_1^{\star \top} + \lambda_2^{\star}\bm{u}_2^\star \bm{u}_2^{\star \top}$.  The ground truth $\bm{M}^\star$ is produced by setting $\lambda_1^\star = 1$ with the corresponding eigenvectors $\bm{u}_1^\star$ and $\bm{u}_2^\star$ randomly generated.
We aim to construct 95\% confidence intervals (namely, $\alpha=0.05$) for 
\begin{enumerate}
	\item[1.] the linear form $\bm{a}^\top \bm{u}^\star_2$ for a given vector $\bm{a}$; 
	\item[2.] the eigenvalue $\lambda_2^\star$. 
\end{enumerate}
The unit vector $\bm{a}$ is chosen such that $\bm{a}^\top \bm{u}_1^\star = 0$ and $\bm{a}^\top \bm{u}_2^\star = 0.5$.
Given that $|\bm{a}^\top \bm{u}_2^\star|$ is quite large, we are expected to have $|\bm{a}^{\top} \bm{u}_2| \geq c_{\mathrm{b}}\sqrt{\widehat{v}_{\bm{a},2}}\log^{1.5} n$ with high probability, and hence we can simply take
$\widehat{u}_{\bm{a},2}^{\mathsf{modified}} = \widehat{u}_{\bm{a},2} $ in these numerical experiments  
(according to \eqref{eq:defn-ua-hat} and \eqref{eq:ua-estimator-large}).  

\begin{figure}[htbp!]
	\begin{tabular}{ccc}
		\includegraphics[width=0.315\textwidth]{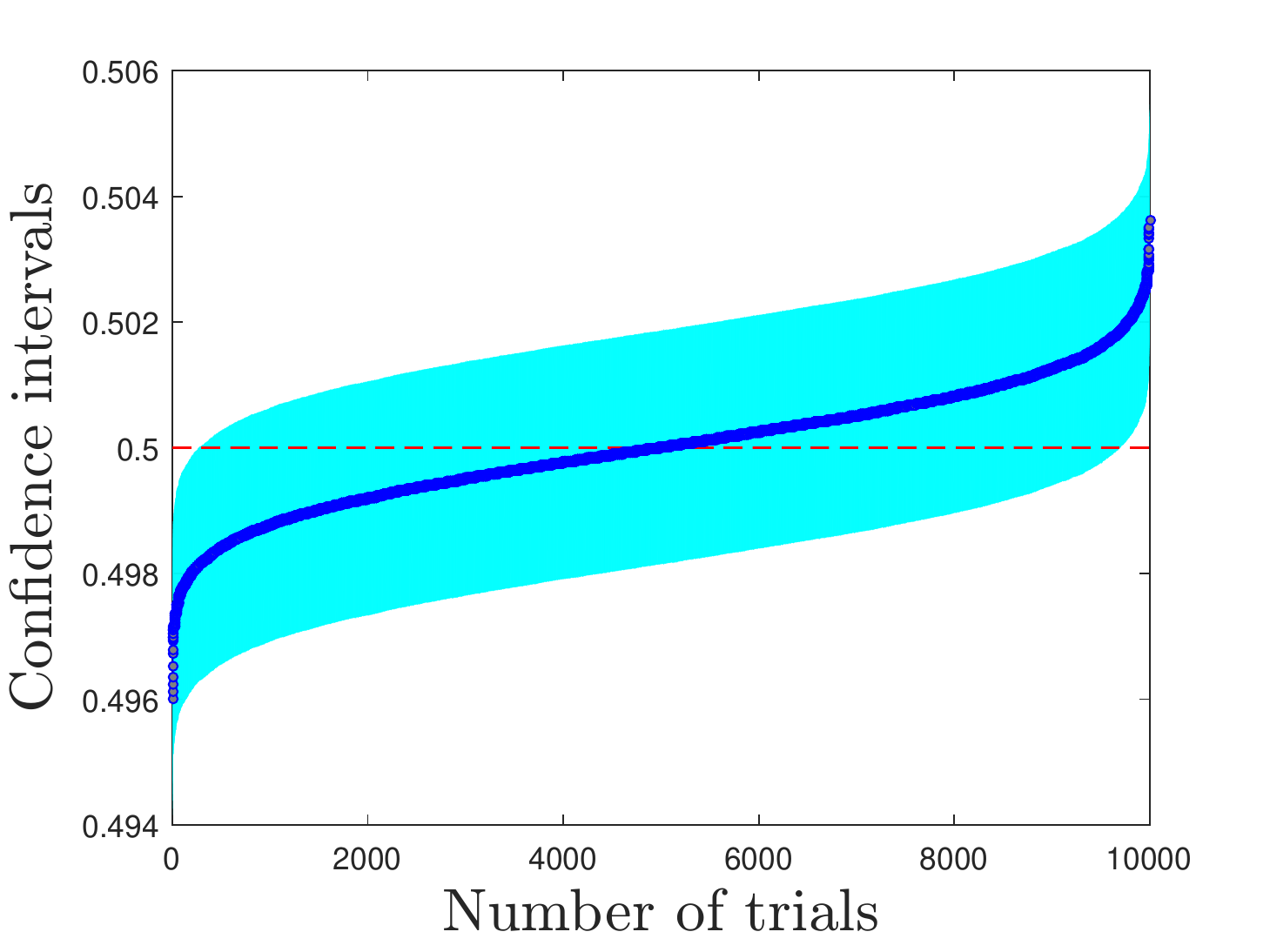} 
		& \includegraphics[width=0.315\textwidth]{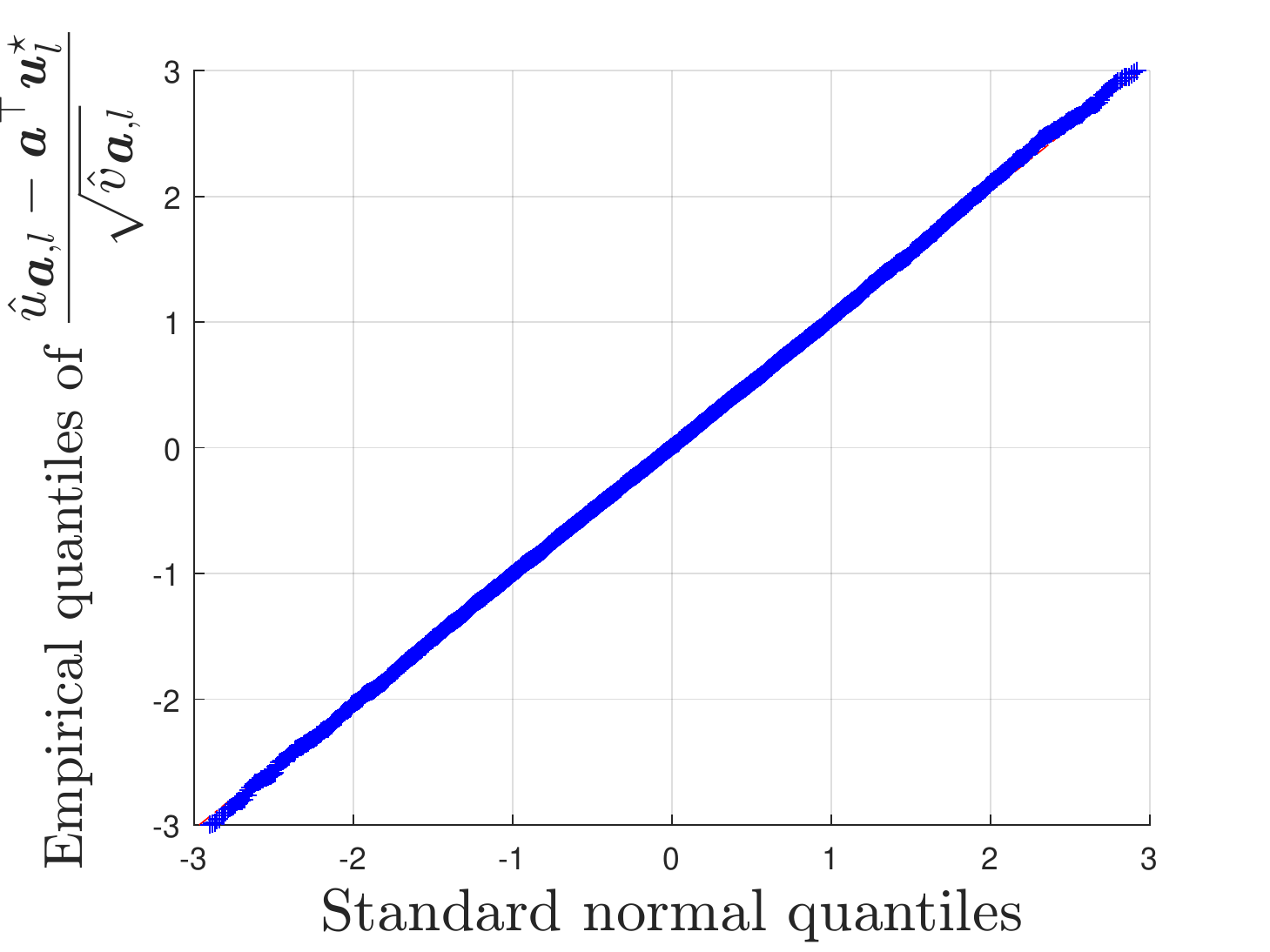} 
		& \includegraphics[width=0.315\textwidth]{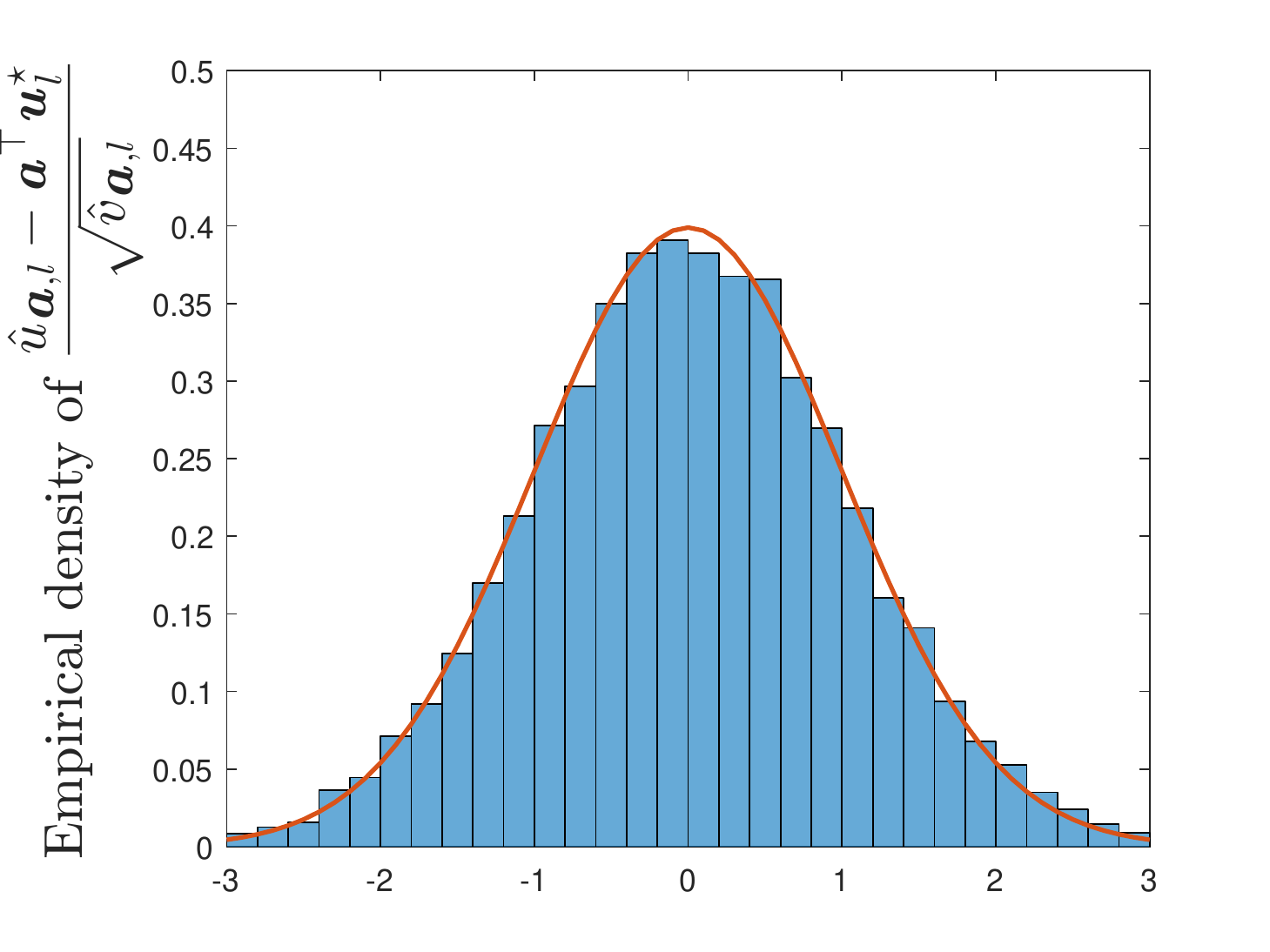} 
		\tabularnewline
		(a) CIs for $\bm{a}^\top \bm{u}_l^\star$ 
		& (b) Q-Q plot for $\frac{\widehat{u}_{\bm{a},l} - \bm{a}^\top \bm{u}_l^\star}{\sqrt{\widehat{v}_{\bm{a}, l}}}$ 
		& (c) Histogram for $\frac{\widehat{u}_{\bm{a},l} - \bm{a}^\top \bm{u}_l^\star}{\sqrt{\widehat{v}_{\bm{a}, l}}}$ 
		\tabularnewline 
		\includegraphics[width=0.315\textwidth]{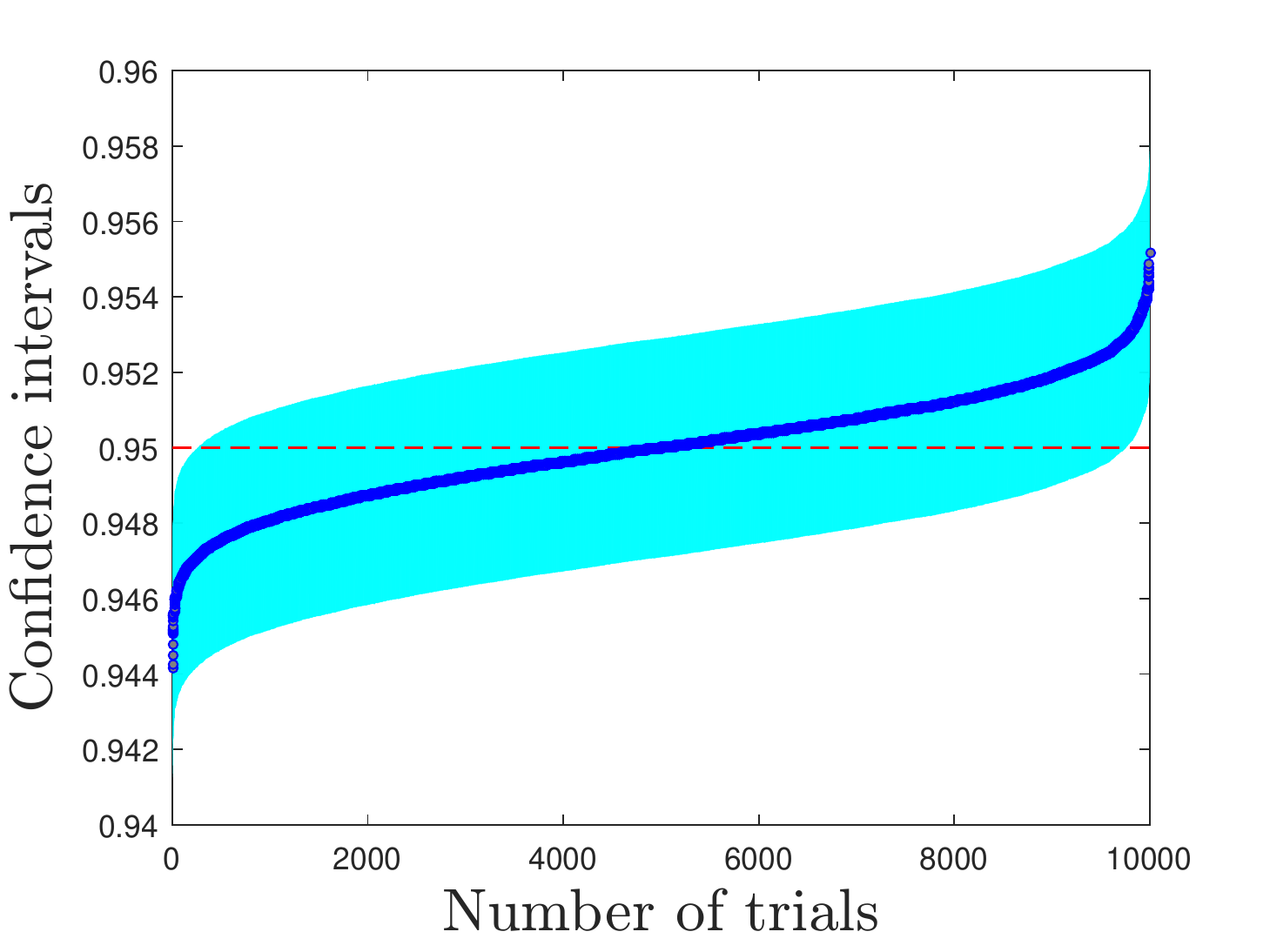} 
		&\includegraphics[width=0.315\textwidth]{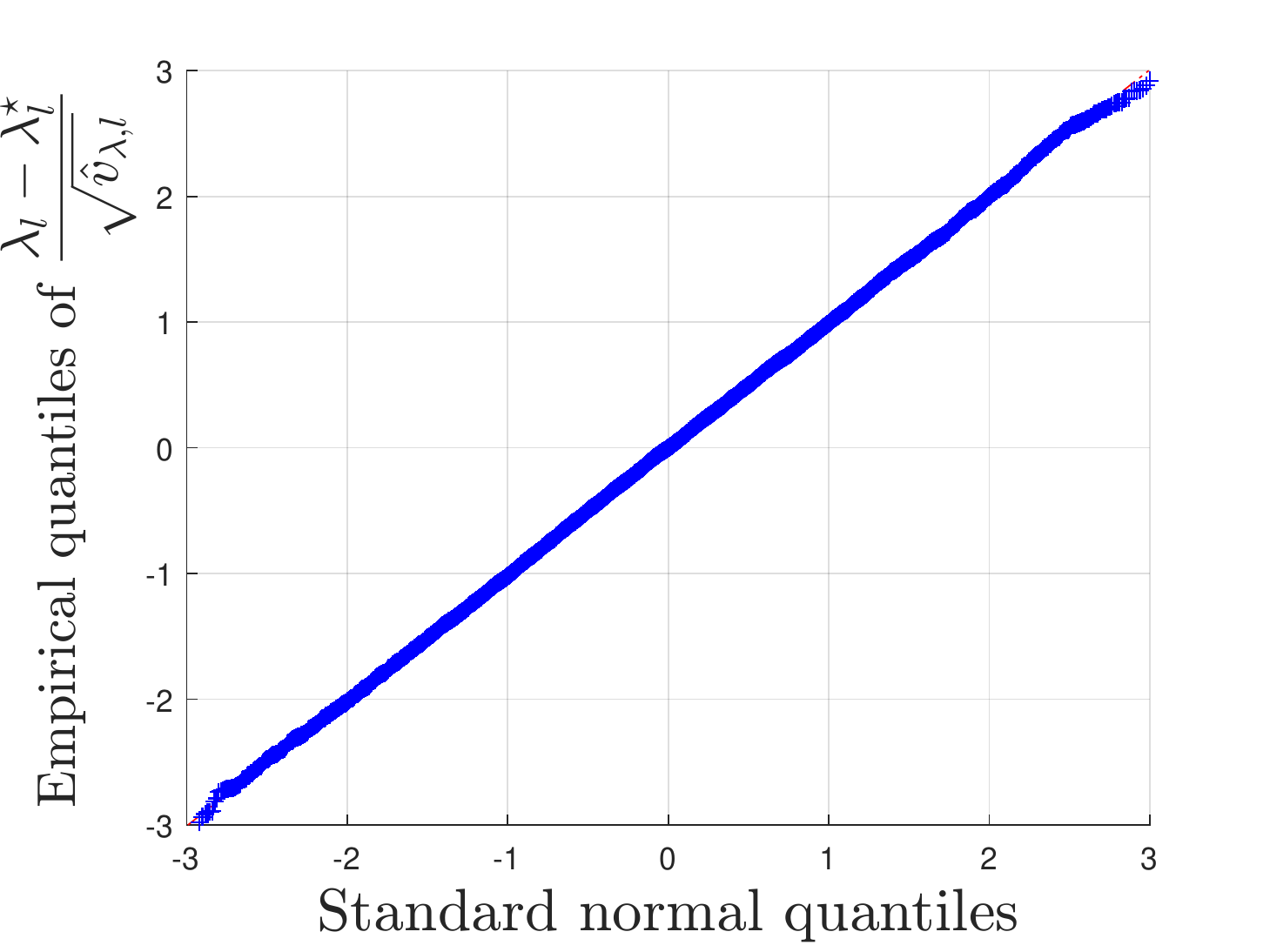} 
		& \includegraphics[width=0.315\textwidth]{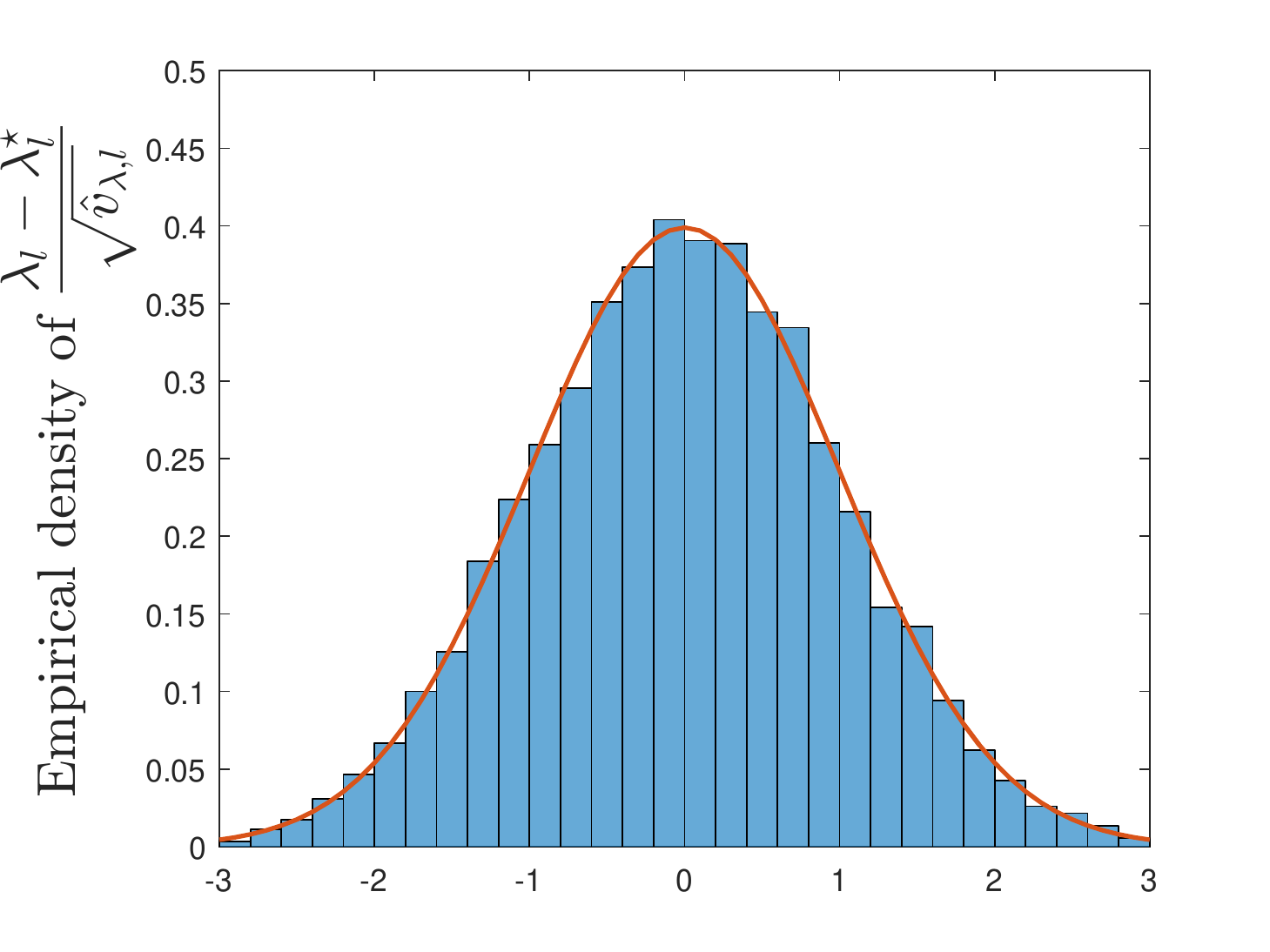} 
		\tabularnewline
		(d) CIs for the eigenvalue $\lambda_l^\star$ 
		& (e) Q-Q plot for $\frac{\lambda_l - \lambda_l^\star}{\sqrt{\widehat{v}_{\lambda, l}}}$ 
		& (f) Histogram for $\frac{\lambda_l - \lambda_l^\star}{\sqrt{\widehat{v}_{\lambda, l}}}$ 
		\tabularnewline 
		\includegraphics[width=0.315\textwidth]{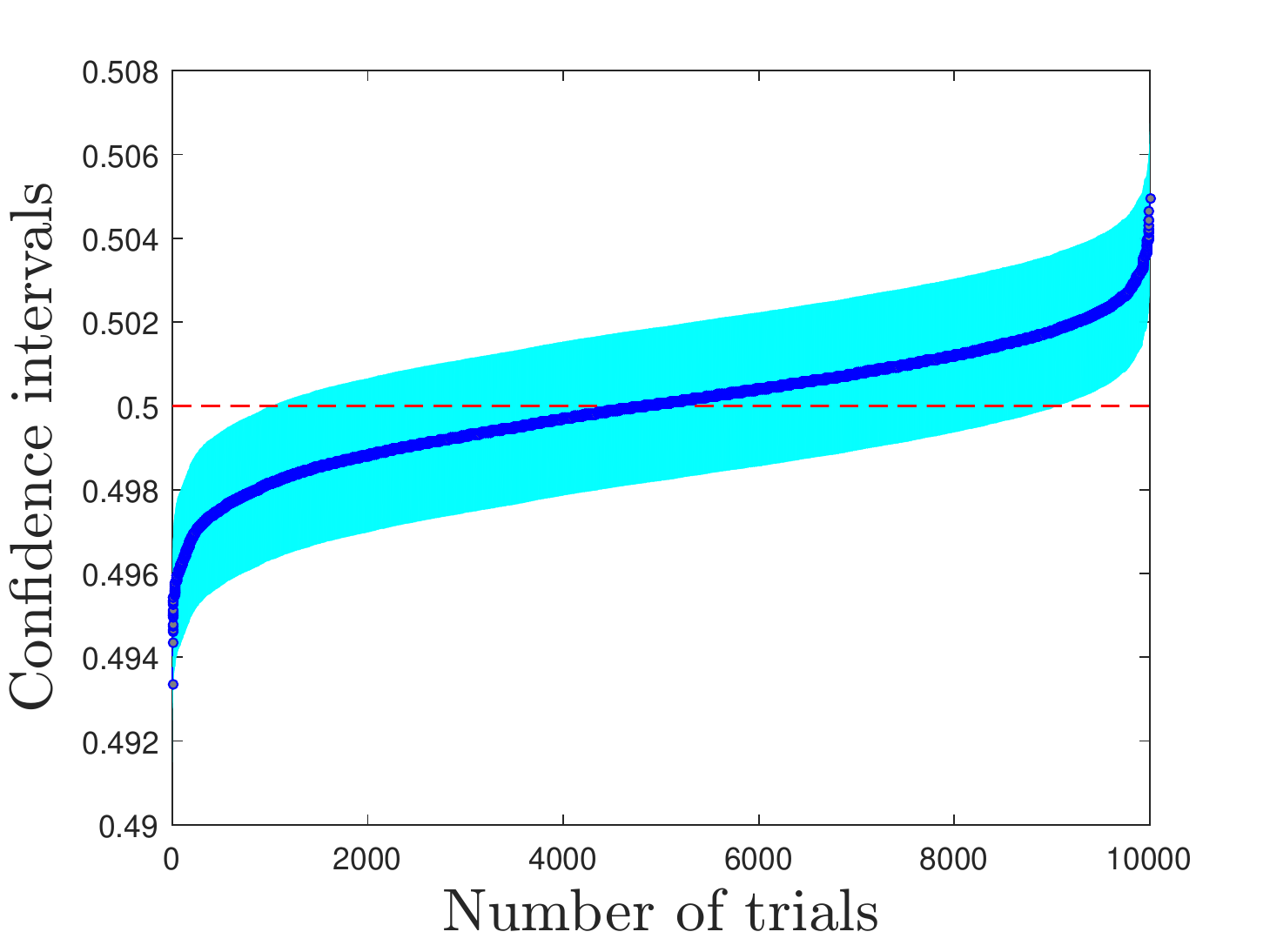} 
		&\includegraphics[width=0.315\textwidth]{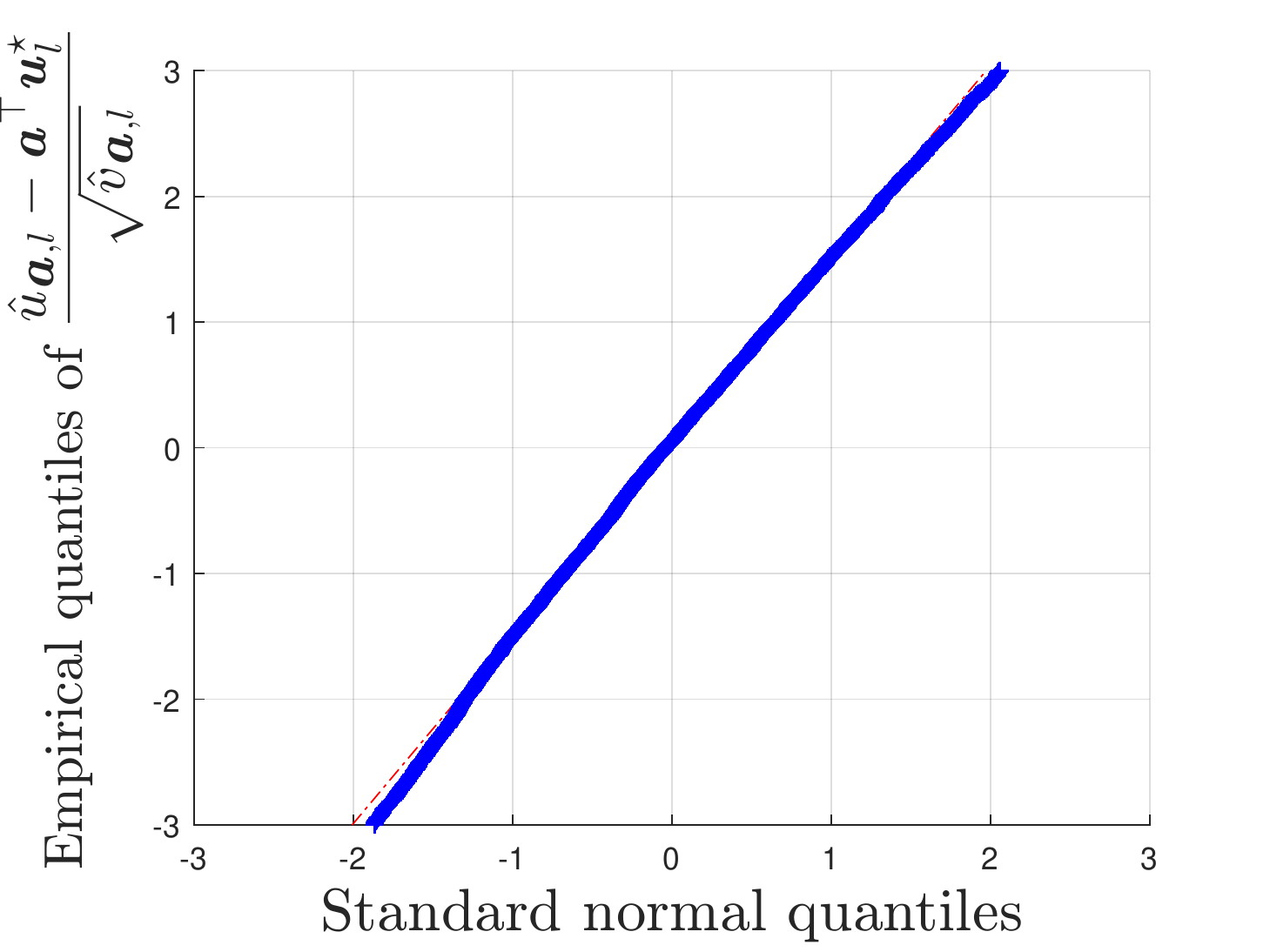} 
		& \includegraphics[width=0.315\textwidth]{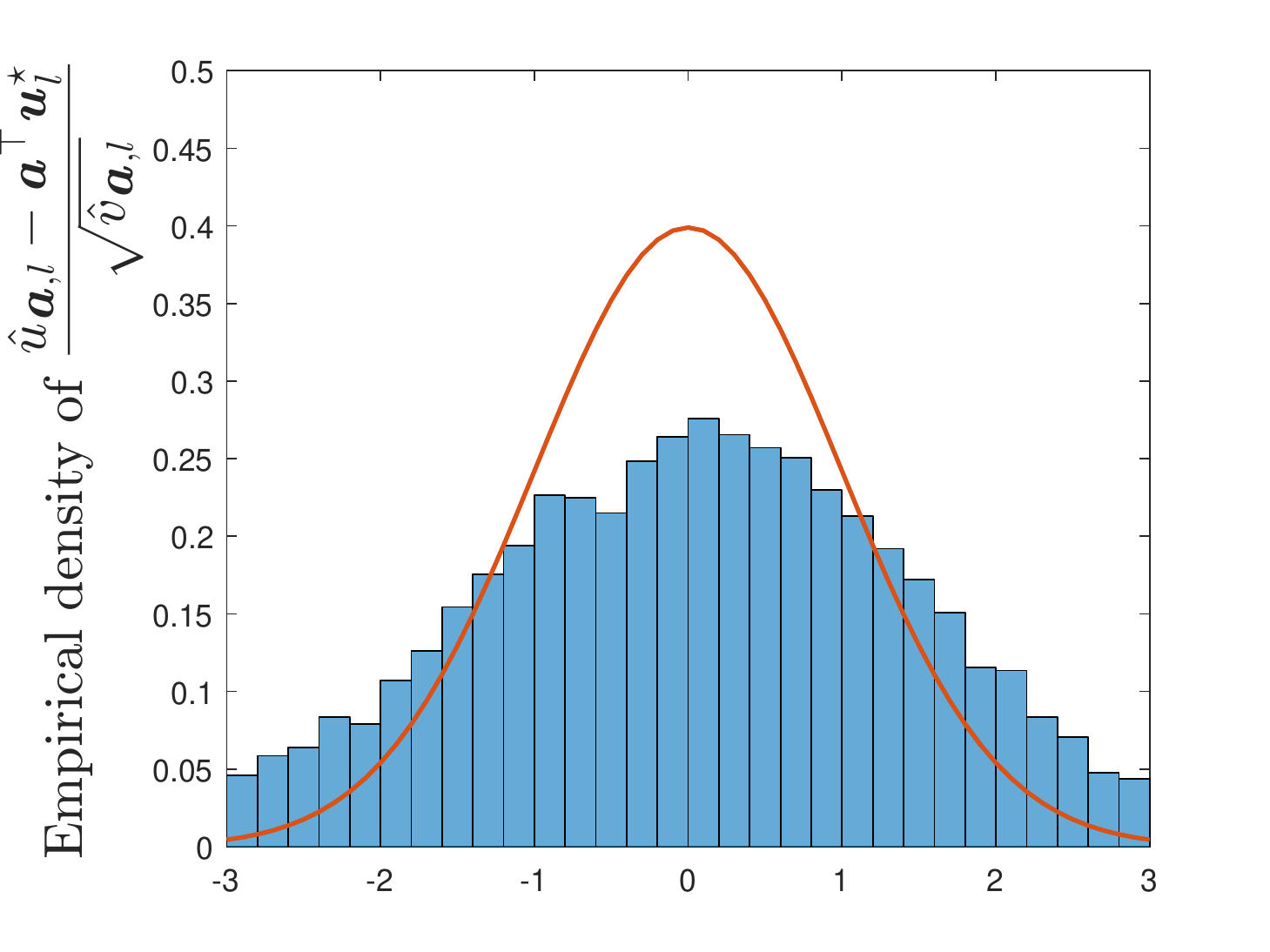}
		\tabularnewline
		(g) CIs for $\bm{a}^\top \bm{u}_l^\star$ 
		& (h) Q-Q plot for $\frac{\widehat{u}_{\bm{a},l} - \bm{a}^\top \bm{u}_l^\star}{\sqrt{\widehat{v}_{\bm{a}, l}}}$ 
		& (i) Histogram for $\frac{\widehat{u}_{\bm{a},l} - \bm{a}^\top \bm{u}_l^\star}{\sqrt{\widehat{v}_{\bm{a}, l}}}$.
	\end{tabular}
	\caption{Numerical results for inference  for the linear form $\bm{a}^\top \bm{u}_l^\star$ and the eigenvalue $\lambda_l^\star = 0.95$ ($l=2$) under heteroscedastic Gaussian noise. In (a)-(f), we take $\bm{a}^\top \bm{u}_1^\star =0$ (no interferer), while in (g), (h) and (i) $\bm{a}^\top \bm{u}_1^\star = 0.05$ (with interferer). In both two settings, we take $n = 1000$, set $\sigma_1 = 0.1 / \sqrt{n \log n}$ and $\delta_{\sigma} = 0.4/((n-1)\sqrt{n\log n})$ (cf.~\eqref{eq:variance-Gaussian-numerics}), and run independent $10000$ trials.  In (a), (d) and (g), the confidence intervals are sorted respectively by the magnitudes of the estimators $\widehat{u}_{\bm{a},l}$ and $\lambda_l$ in these trials.    Here, $\bm{u}_2$ is chosen such that $\bm{u}_2^{\top}\bm{u}_2^{\star}\geq 0$. In (c), (f) and (i), the empirical densities  are compared to the pdf of the standard normal (red curve).}
	\label{fig:rank2-CI-Gaussian}	
\end{figure}

\paragraph{Heteroscedastic Gaussian noise.}
The noise setting is the same as in expression \eqref{eq:variance-Gaussian-numerics}. 
The numerical results are displayed in Fig.~\ref{fig:rank2-CI-Gaussian} and Tab.~\ref{tab:conf-interval-all}.   We also examine the necessity of our requirement on the ``interferers'' (i.e.~$|\bm{a}^\top \bm{u}_k^\star| \lesssim |\lambda_l^\star - \lambda_k^\star|/(|\lambda_l^\star| \log n), k \neq l$) as stated in Theorem \ref{thm:confidence-interval-validity-rankr}. Specifically, we make comparisons of the following two settings: (1) the ``no interferer'' case where $\bm{a}^\top \bm{u}_1^\star = 0$ and $\bm{a}^\top \bm{u}_2^\star = 0.5$, as plotted in  Fig.~\ref{fig:rank2-CI-Gaussian}(a)-\ref{fig:rank2-CI-Gaussian}(f); and (2) the case with a strong interferer where $\bm{a}^\top \bm{u}_1^\star = 0.05$ and $\bm{a}^\top \bm{u}_2^\star = 0.5$, as plotted in  Fig.~\ref{fig:rank2-CI-Gaussian}(g)-\ref{fig:rank2-CI-Gaussian}(i). Numerically, our distributional guarantees are inaccurate when there exists a strong interferer, which is consistent with what our theory predicts.

\begin{figure}[htbp!]
	\centering
	\begin{tabular}{ccc}
		\includegraphics[width=0.315\textwidth]{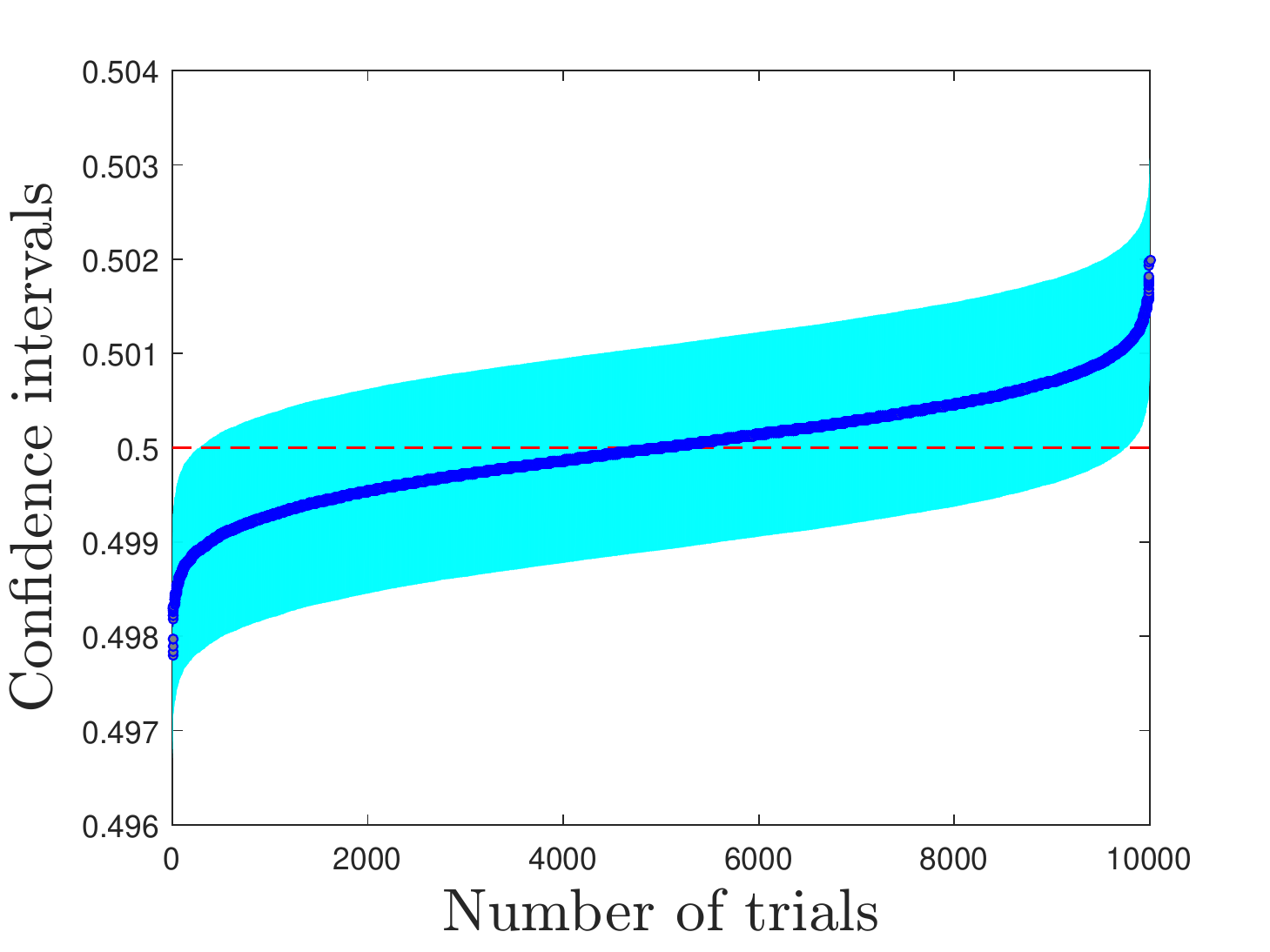} &
		\includegraphics[width=0.315\textwidth]{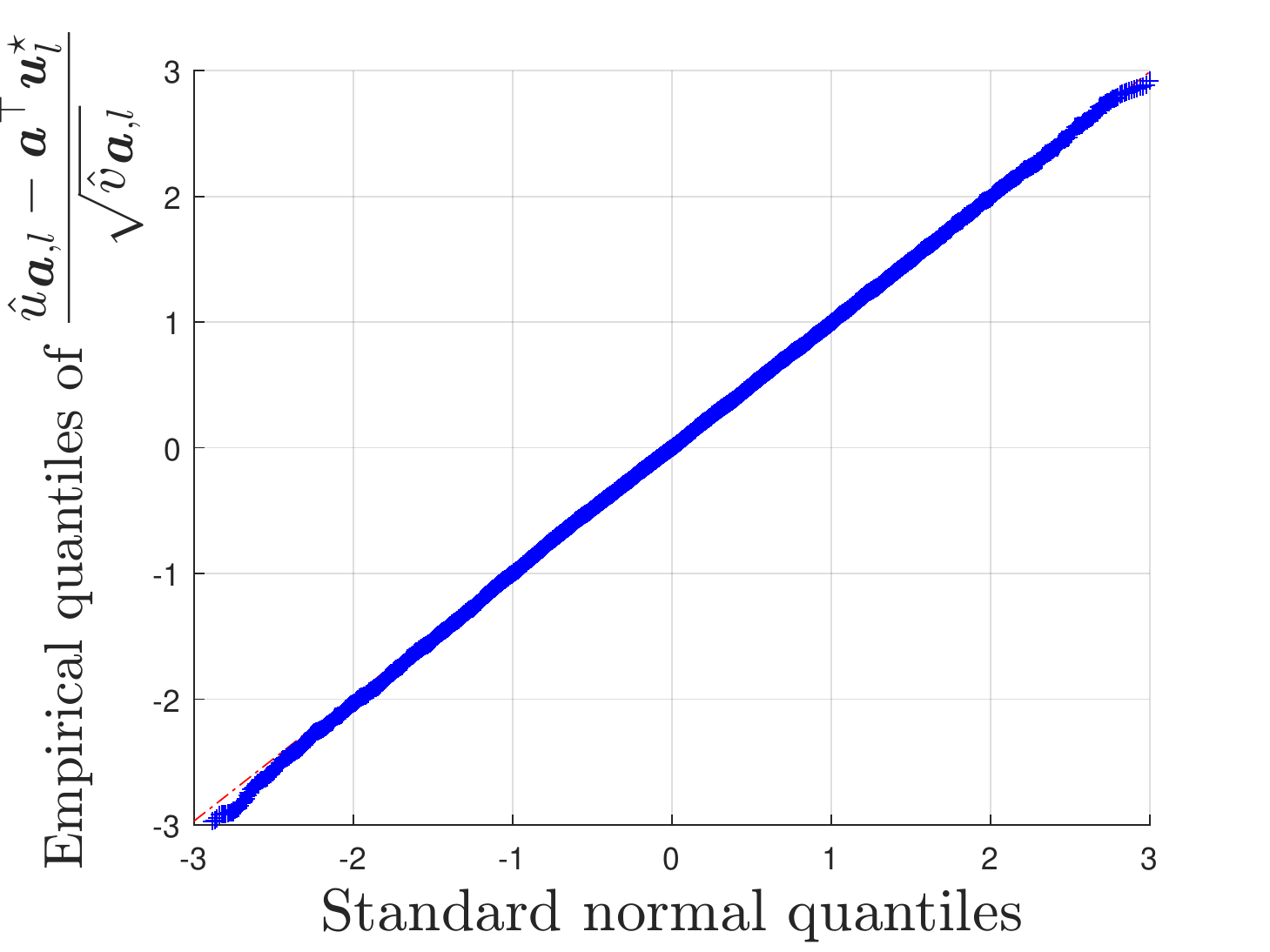} & \includegraphics[width=0.315\textwidth]{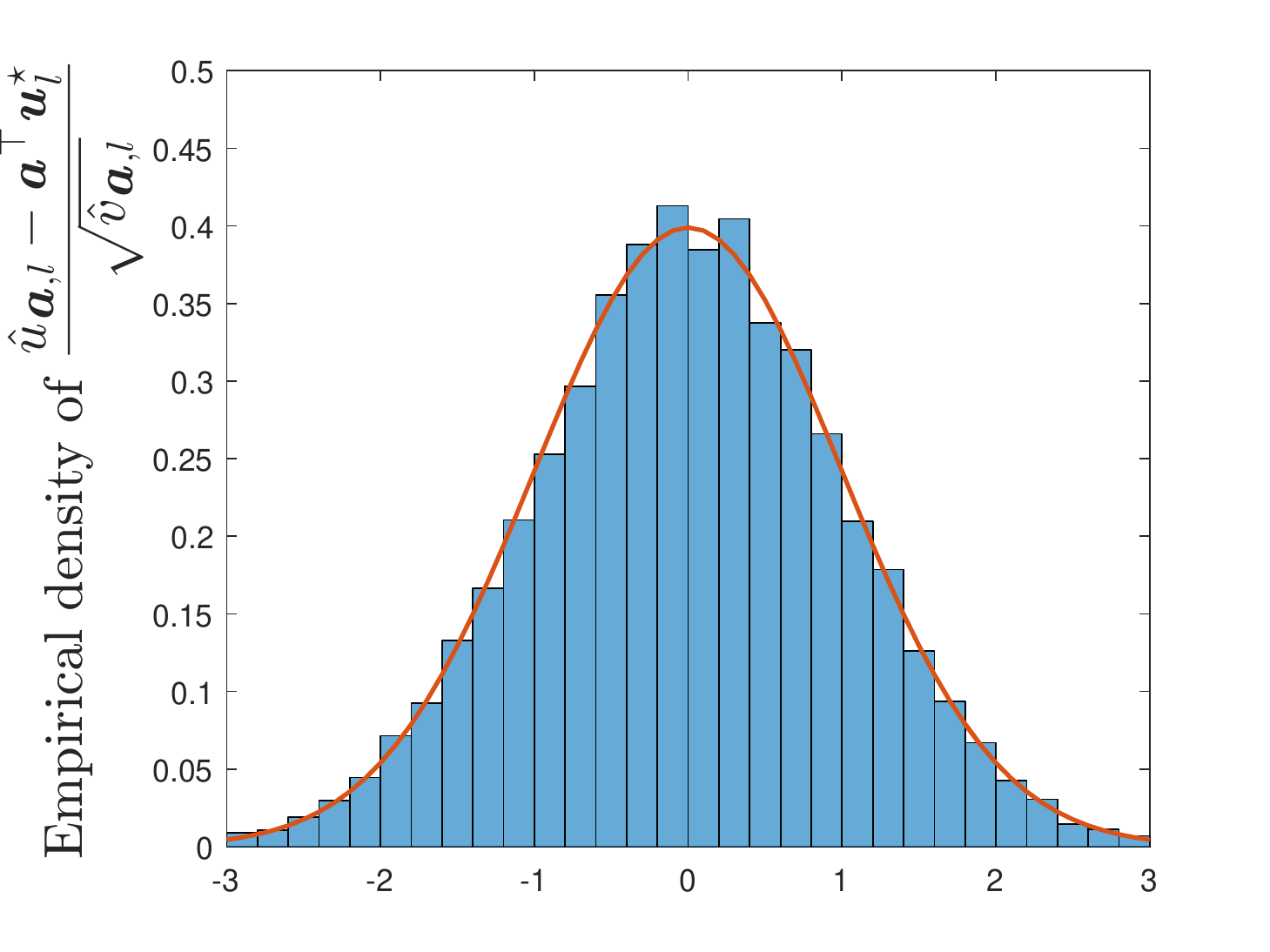} 
		\tabularnewline
		(a) CIs for $\bm{a}^\top \bm{u}_l^\star$ 
		& (b) Q-Q plot for $\frac{\widehat{u}_{\bm{a},l} - \bm{a}^\top \bm{u}_l^\star}{\sqrt{\widehat{v}_{\bm{a}, l}}}$ 
		& (c) Histogram for $\frac{\widehat{u}_{\bm{a},l} - \bm{a}^\top \bm{u}_l^\star}{\sqrt{\widehat{v}_{\bm{a}, l}}}$ 
		\tabularnewline 
		\includegraphics[width=0.315\textwidth]{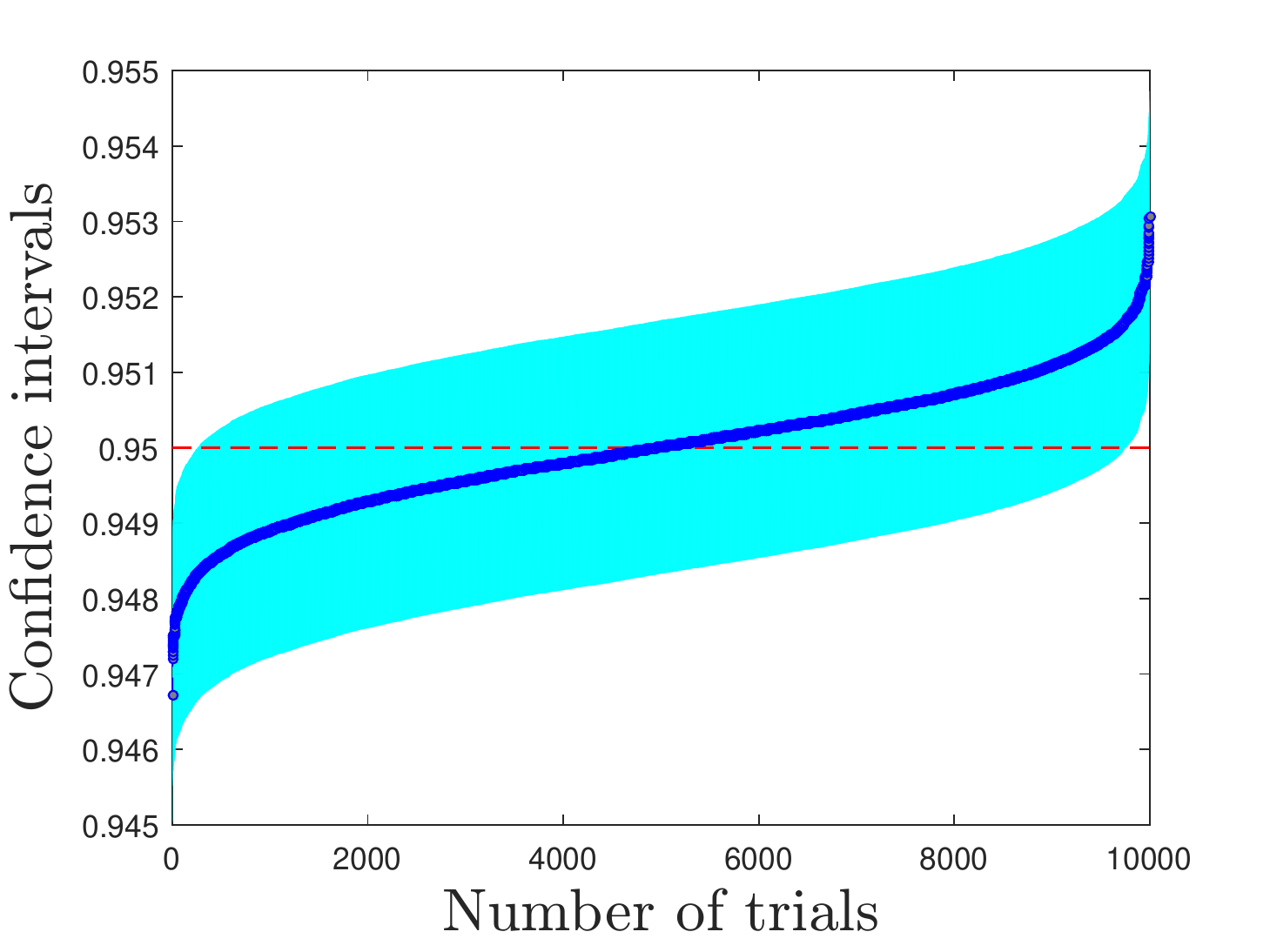} 
		&\includegraphics[width=0.315\textwidth]{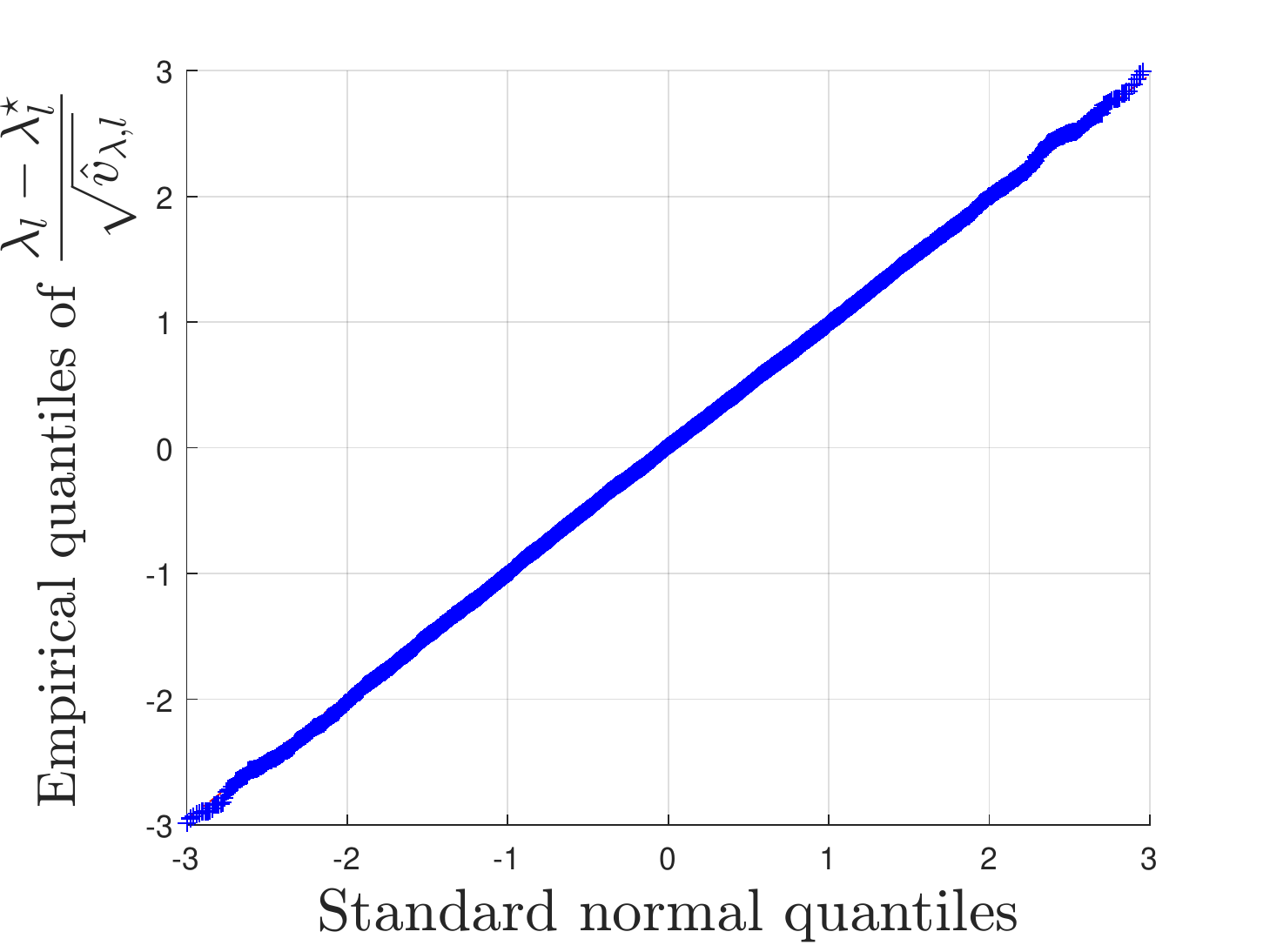} 
		& \includegraphics[width=0.315\textwidth]{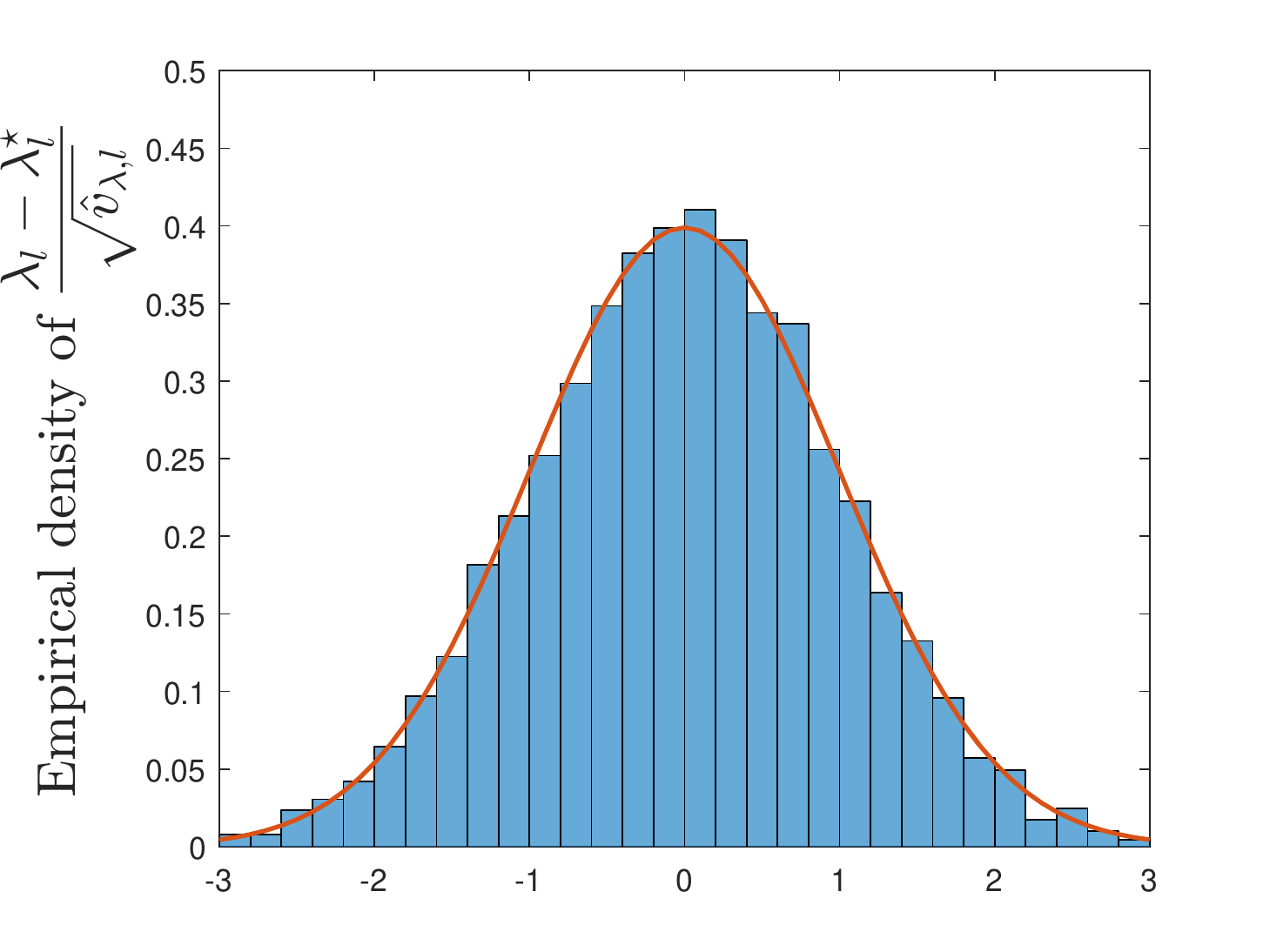} 
		\tabularnewline
		(d) CIs for the eigenvalue $\lambda_l^\star$ 
		& (e) Q-Q plot for $\frac{\lambda_l - \lambda_l^\star}{\sqrt{\widehat{v}_{\lambda, l}}}$ 
		& (f) Histogram for $\frac{\lambda_l - \lambda_l^\star}{\sqrt{\widehat{v}_{\lambda, l}}}$ 
	\end{tabular}
	\caption{Numerical results for inference for the linear form $\bm{a}^\top \bm{u}_l^\star$ and the eigenvalue $\lambda_l^\star=0.95$ ($l=2$) under heteroscedastic Bernoulli noise. We take $n = 1000$, set $\sigma_1 = 0.1 / \sqrt{n \log n}$ and $\delta_{\sigma} = 0.4/((n-1)\sqrt{n\log n})$ (cf.~\eqref{eq:variance-Gaussian-numerics}), and run independent $10000$ trials.  In (a) and (d), the confidence intervals are sorted respectively by the magnitudes of the estimators $\widehat{u}_{\bm{a},l}$ and $\lambda_l$ in these trials. Here, $\bm{u}_2$ is chosen such that $\bm{u}_2^{\top}\bm{u}_2^{\star}\geq 0$. In (c) and (f), the empirical densities  are compared to the pdf of the standard normal (red curve).}
	\label{fig:rank2-CI-Bernoulli}	
\end{figure}

\paragraph{Heteroscedastic Bernoulli noise.}
Consider a heteroscedastic Bernoulli noise matrix $\bm{H}$ with independent entries such that $H_{ij} = -\sigma_{ij}$ with probability $1/2$ and $H_{ij} = \sigma_{ij}$ otherwise. The variance matrix is also chosen to satisfy \eqref{eq:variance-Gaussian-numerics}. The numerical results are plotted in Fig.~\ref{fig:rank2-CI-Bernoulli} and Tab.~\ref{tab:conf-interval-all}.

\begin{figure}[htbp!]
	\centering
	\begin{tabular}{ccc}
		\includegraphics[width=0.315\textwidth]{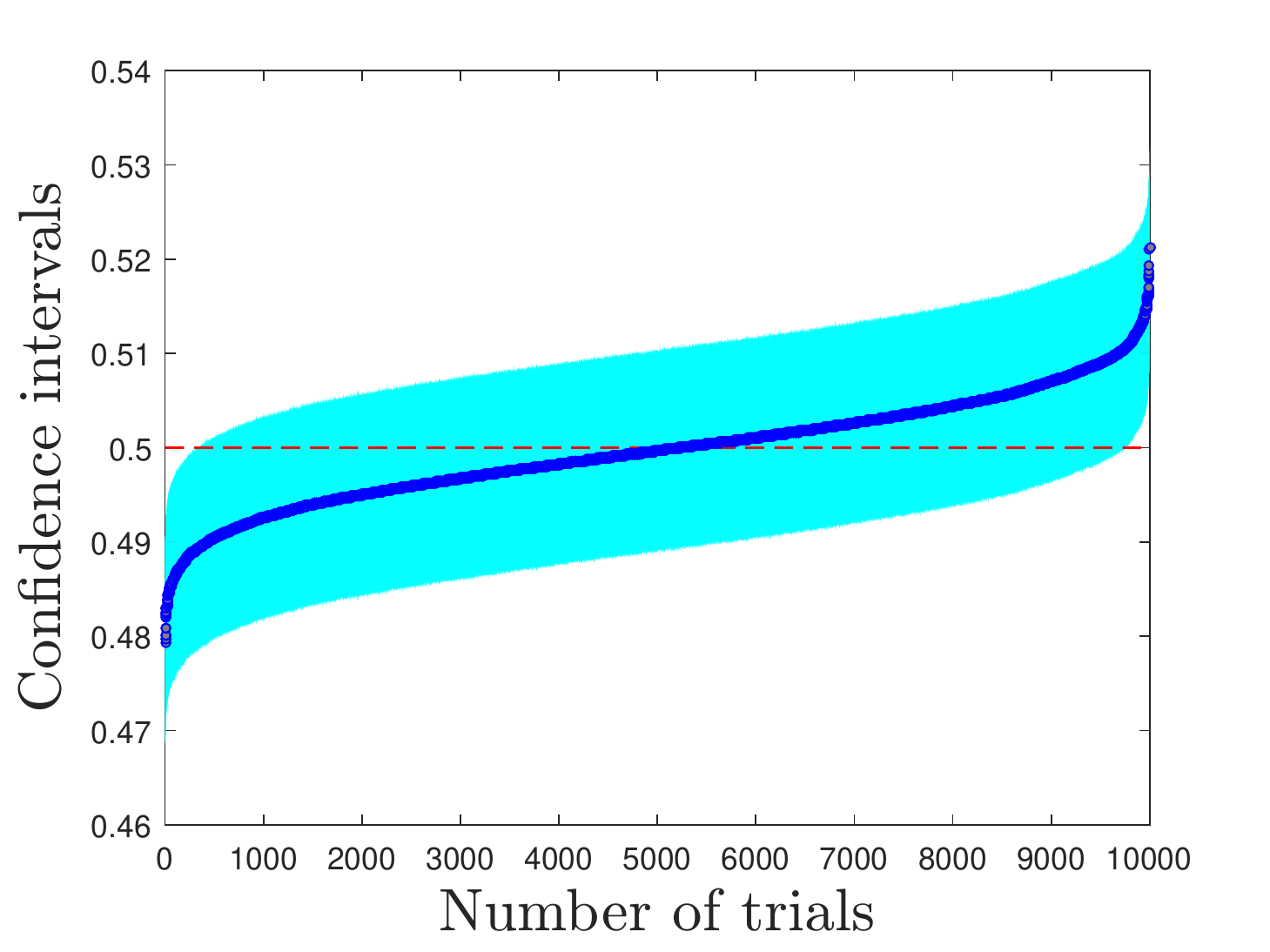} &
		\includegraphics[width=0.315\textwidth]{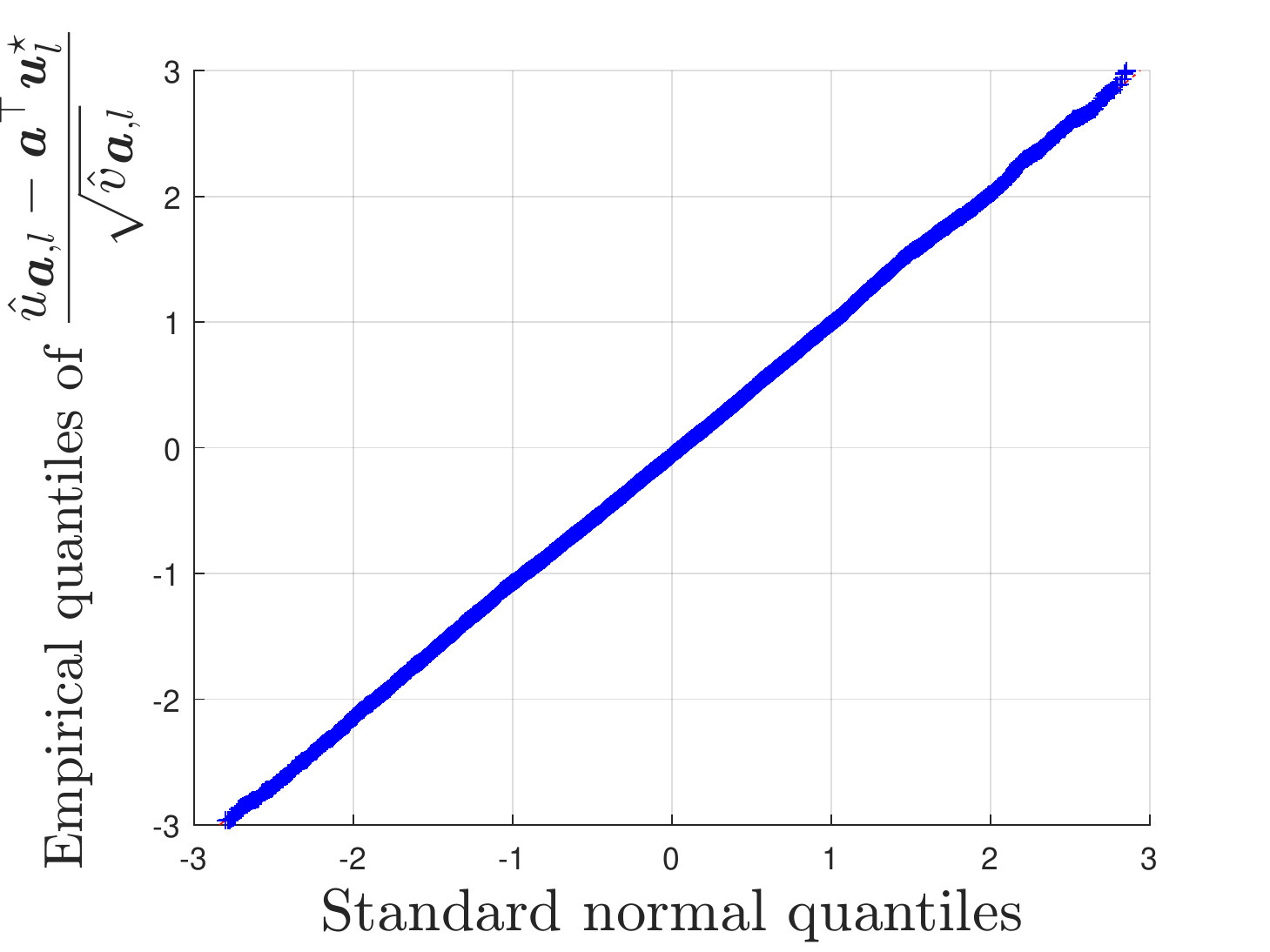} & \includegraphics[width=0.315\textwidth]{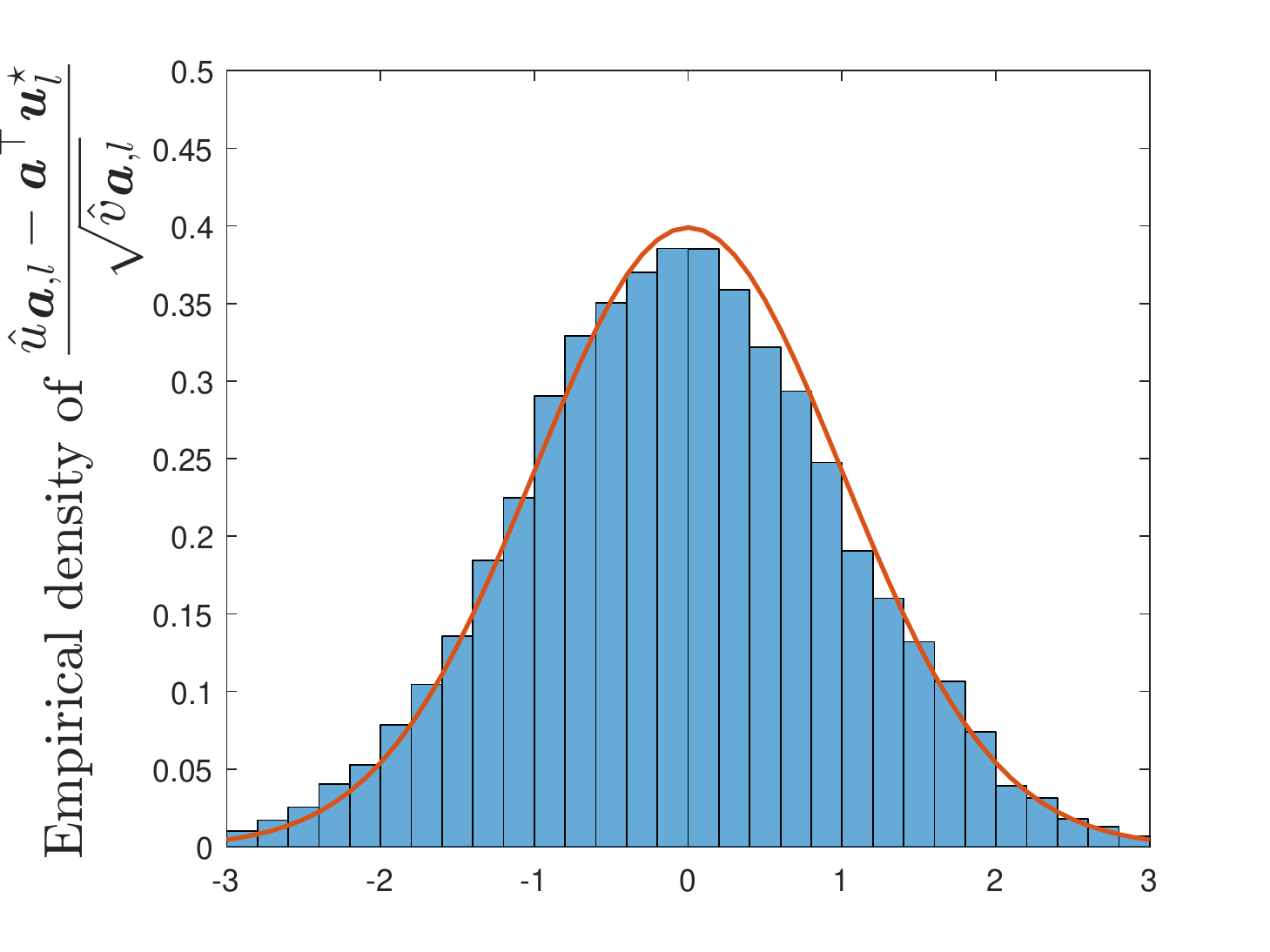} 
		\tabularnewline
		(a) CIs for $\bm{a}^\top \bm{u}_l^\star$ 
		& (b) Q-Q plot for $\frac{\widehat{u}_{\bm{a},l} - \bm{a}^\top \bm{u}_l^\star}{\sqrt{\widehat{v}_{\bm{a}, l}}}$ 
		& (c) Histogram for $\frac{\widehat{u}_{\bm{a},l} - \bm{a}^\top \bm{u}_l^\star}{\sqrt{\widehat{v}_{\bm{a}, l}}}$ 
		\tabularnewline 
		\includegraphics[width=0.315\textwidth]{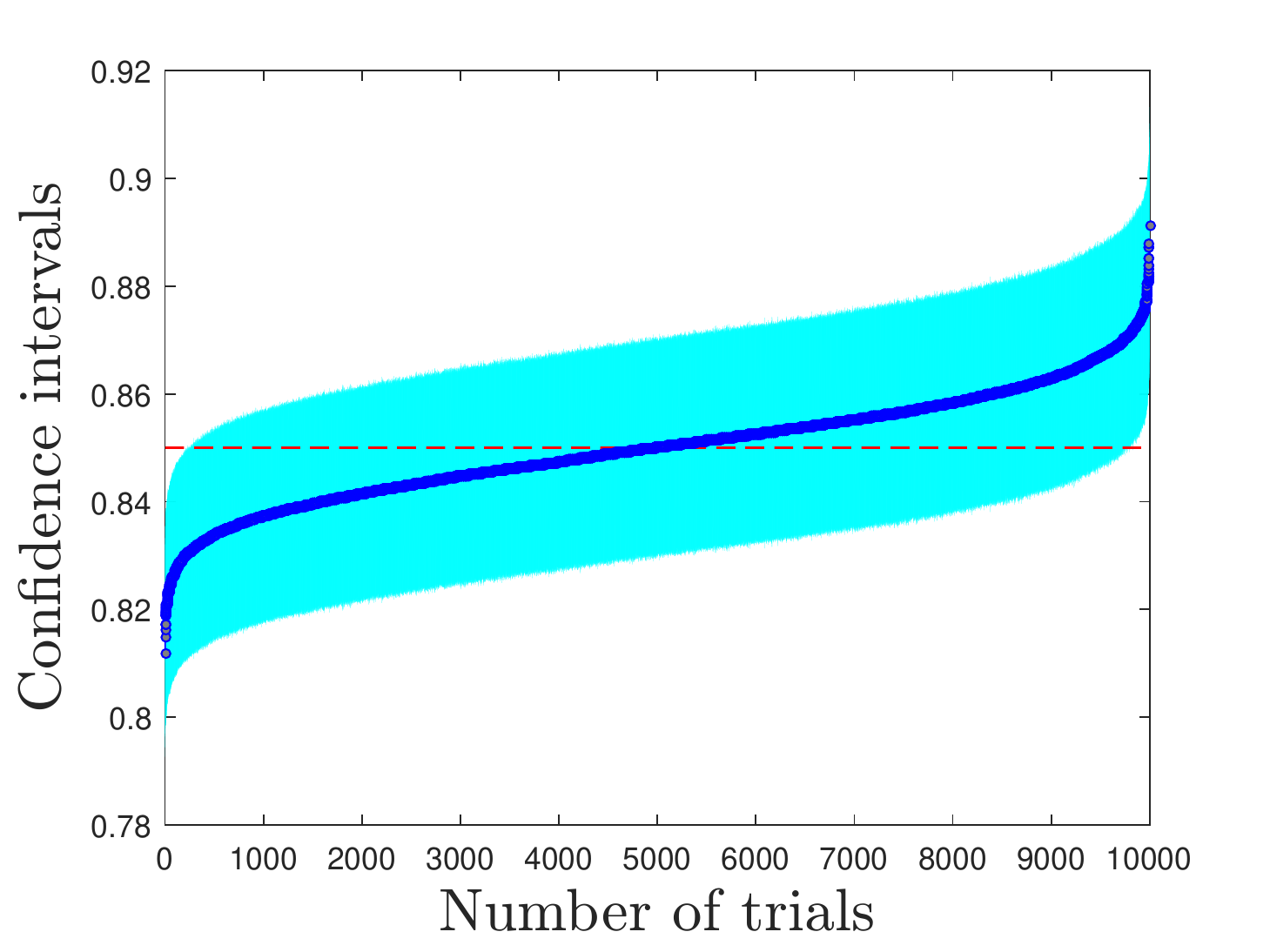} 
		&\includegraphics[width=0.315\textwidth]{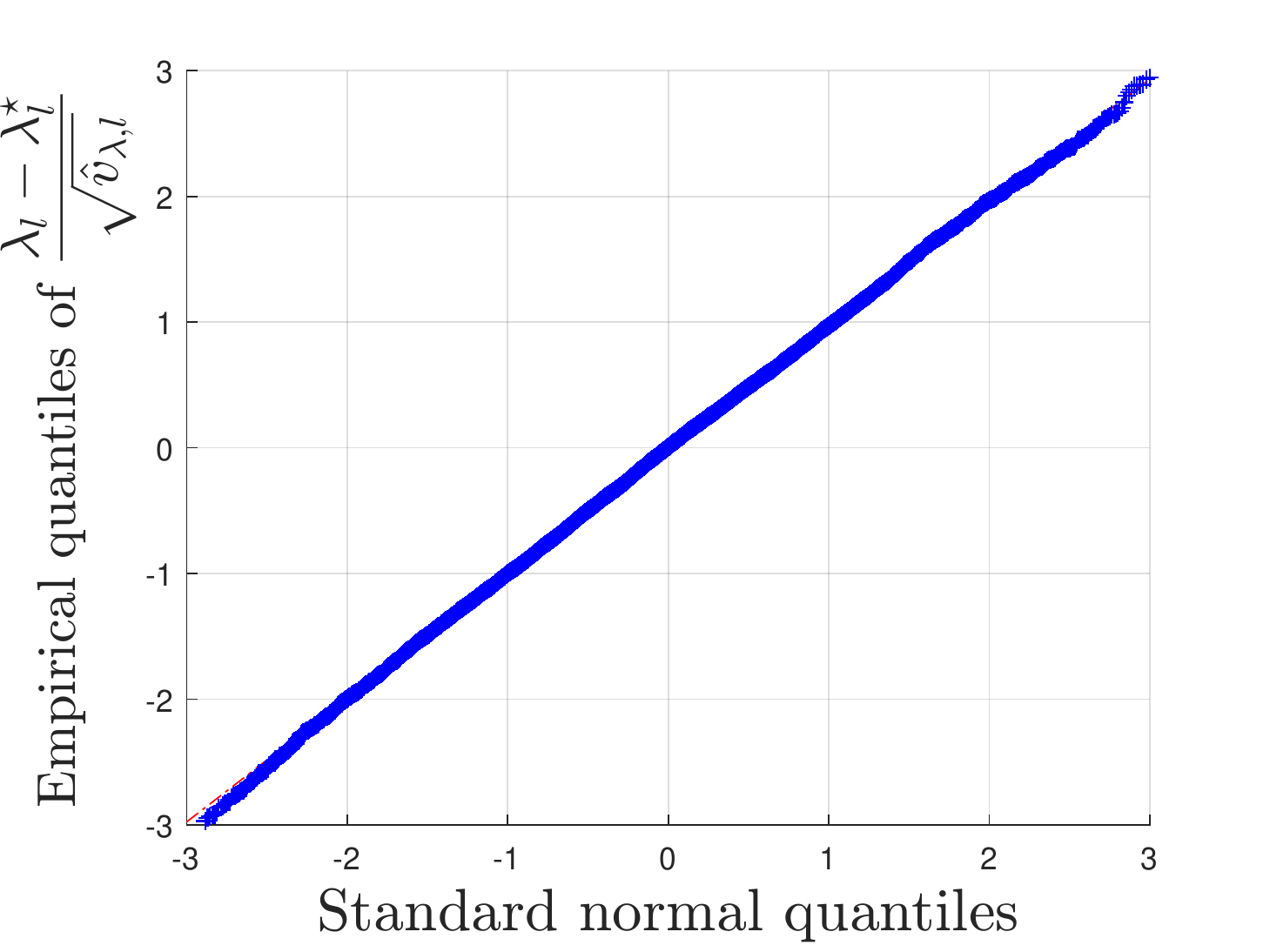} 
		& \includegraphics[width=0.315\textwidth]{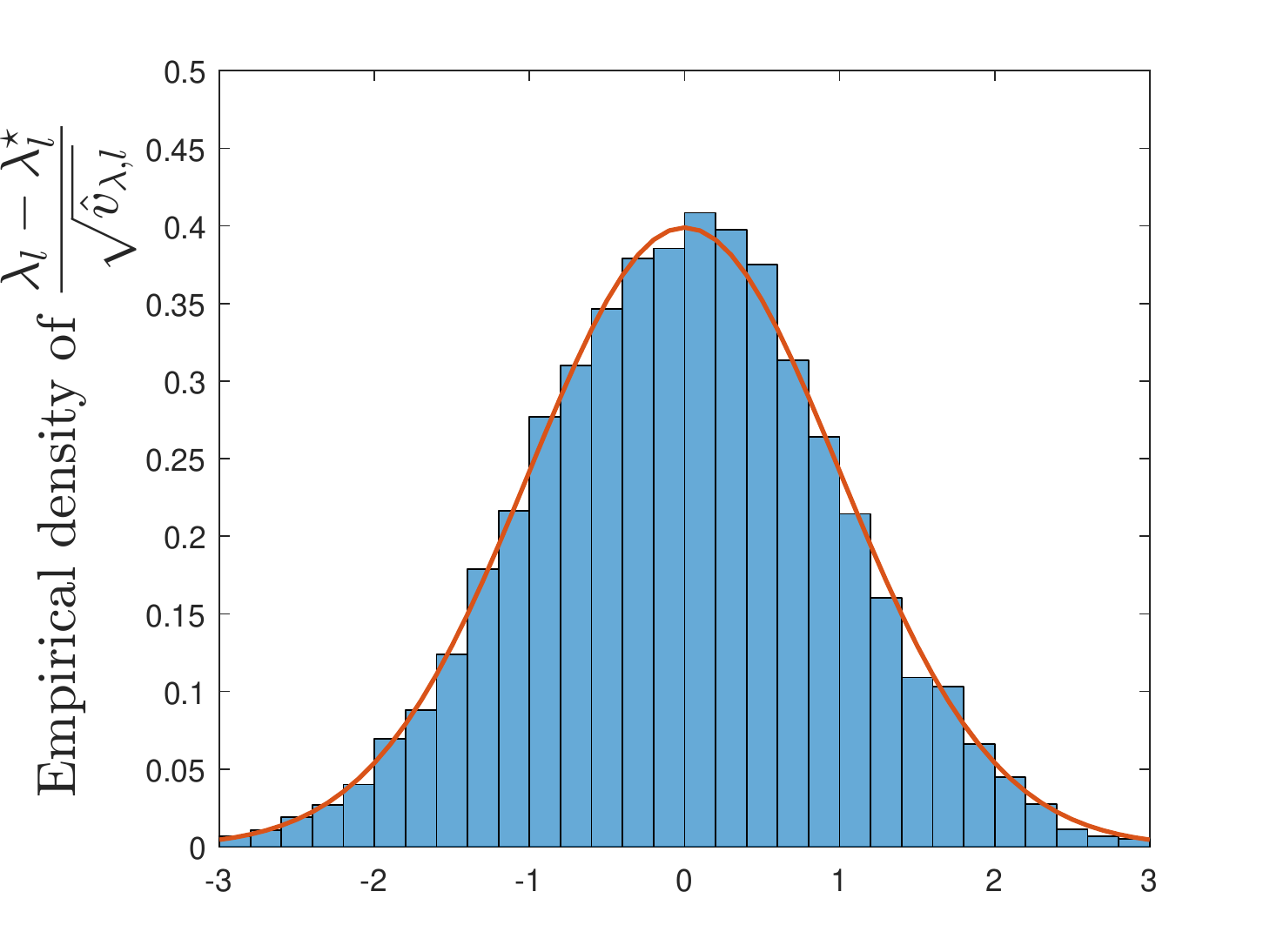} 
		\tabularnewline
		(d) CIs for the eigenvalue $\lambda_l^\star$ 
		& (e) Q-Q plot for $\frac{\lambda_l - \lambda_l^\star}{\sqrt{\widehat{v}_{\lambda, l}}}$ 
		& (f) Histogram for $\frac{\lambda_l - \lambda_l^\star}{\sqrt{\widehat{v}_{\lambda, l}}}$ 
	\end{tabular}
	\caption{Numerical results for inference for the linear form $\bm{a}^\top \bm{u}_l^\star$ and the eigenvalue $\lambda_l^\star = 0.85$ ($l=2$) in the presence of missing data. We take $n = 1000$, set $p = 0.1$, $\sigma_1 = 1 / \sqrt{10 n \log n}$ and $\delta_{\sigma} = 0.4/((n-1)\sqrt{10 n\log n})$ (cf.~\eqref{eq:missing-data}), and run $10000$ independent trials. In (a) and (g), the confidence intervals are sorted respectively by the magnitudes of the estimators $\widehat{u}_{\bm{a},l}$ and $\lambda_l$ in these trials. Here, $\bm{u}_2$ is chosen such that $\bm{u}_2^{\top}\bm{u}_2^{\star}\geq 0$. In (c) and (f), the empirical densities  are compared to the pdf of the standard normal (red curve).}
	\label{fig:rank2-CI-MC}	
\end{figure}

\paragraph{Missing data model.}
In the case when we only get to observe a fraction of the entries of $\bm{M}^{\star}$ as in model~\eqref{eq:missing-data}, 
in the same way, we aim to provide confidence intervals with precise coverages for both the linear functionals 
and the eigenvalues. 
Our numerical results are shown in Fig.~\ref{fig:rank2-CI-MC} and Tab.~\ref{tab:conf-interval-all}.

\paragraph{Conclusion.}  
In all of the above numerical experiments, the confidence intervals and the Q-Q (quantile-quantile) plots we produce match the theoretical predictions in a reasonably well manner, thus corroborating the validity and practicability of our theoretical results. We have numerically verified the need of controlling the ``interferers'' in Fig.~\ref{fig:rank2-CI-Gaussian}(g)-\ref{fig:rank2-CI-Gaussian}(i).

\begin{table}
	\centering
	\begin{tabular}{c|c|c}
		\hline
		Settings & Target & Numerical coverage rates \tabularnewline
		\hline
		\multirow{2}{*}{heteroscedastic Gaussian noise} & linear form $\bm{a}^{\top}\bm{u}_2^{\star}$ & 0.9422 \tabularnewline
		& eigenvalue $\lambda_2^{\star}$ & 0.9496 \tabularnewline
		\hline
		\multirow{2}{*}{heteroscedastic Bernoulli noise} & linear form $\bm{a}^{\top}\bm{u}_2^{\star}$ & 0.9470 \tabularnewline
		& eigenvalue $\lambda_2^{\star}$ & 0.9494 \tabularnewline
		\hline
		\multirow{2}{*}{missing data} & linear form $\bm{a}^{\top}\bm{u}_2^{\star}$ & 0.9398 \tabularnewline
		& eigenvalue $\lambda_2^{\star}$ & 0.9519 \tabularnewline
		\hline
	\end{tabular}
	\caption{Numerical coverage rates for our $95\%$ confidence intervals over $10000$ independent trials. \label{tab:conf-interval-all}}
\end{table}

\section{Prior art}
\label{sec:prior-art}

Recent years have witnessed a flurry of activity in noisy low-rank matrix estimation  \cite{keshavan2010matrix,CanPla10,Negahban2012rscMC,MR2906869,montanari2018adapting,donoho2014minimax,chen2015fast,chen2018projected,ma2017implicit,donoho2018optimal,liu2018pca,johnstone2018pca,zhang2018heteroskedastic,chen2020bridging,yan2021inference}. 
Despite the fundamental importance of estimating linear functionals of eigenvectors (or singular vectors), how to accomplish this task remains largely unknown. 
Only until very recently, researchers started to understand estimation errors for a special type of linear functionals, namely, the entrywise estimation error for the leading eigenvector or the $\ell_{2,\infty}$ error for the rank-$r$ eigenspace. 
Examples include \cite{fan2016eigenvector,chen2017spectral,abbe2017entrywise,cape2019two,pmlr-v83-eldridge18a,ma2017implicit,chen2018asymmetry,chen2019noisy,pananjady2019value,chen2019nonconvex,cai2019subspace,lei2019unified}, which have been motivated by various applications including matrix completion, community detection, top-$K$ ranking, and so on. 
While tight $\ell_{\infty}$ eigenvector (resp.~$\ell_{2,\infty}$ eigenspace) perturbation bounds in the presence of symmetric noise matrices have been derived in the prior works \cite{abbe2017entrywise,cai2019subspace,lei2019unified}, all of these papers required the associated eigen-gap to exceed the spectral norm of the noise matrix, thereby falling short of addressing the scenarios with small eigen-gaps.  
 The work \cite{ma2017implicit} also considered controlling the perturbation error of certain Fourier coefficients of the leading singular vector in a blind deconvolution problem, which, however, does not generalize to other linear functionals. \cite{koltchinskii2016perturbation} developed sharp concentration bounds for estimating linear forms of singular vectors under i.i.d.~Gaussian noise,  which, however, required the true singular values to be sufficiently separated (i.e.~with a spacing much larger than the minimal eigen-gap studied herein). 
The recent work \cite{li2021minimax} is perhaps the only one that studied finite-grained eigenvector perturbation in the face of a small eigen-gap, which, however, is restricted to the case with i.i.d.~Gaussian noise.

Given that  estimation of linear functionals of eigenvectors is already largely under-explored, it is perhaps not surprising to see the lack of investigation about inference and uncertainty quantification for these quantities. This is in stark constrat to sparse estimation and learning problems, for which the construction of confidence regions has been extensively studied \cite{zhang2014confidence,van2014asymptotically,javanmard2014confidence,cai2017confidence,belloni2011inference,ren2015asymptotic,ning2017general,ma2017inter,jankova2018biased,miolane2018distribution,celentano2020lasso}. A few exceptions are worth mentioning: (1) \cite{carpentier2015uncertainty,carpentier2015signal,carpentier2016constructing} identified $\ell_2$ confidence regions that are likely to cover the low-rank matrix of interest, which, however, might be loose in terms of the pre-constant; (2) focusing on low-rank matrix completion, the recent work \cite{chen2019inference} developed a de-biasing strategy that constructs both confidence regions for low-rank factors and entrywise confidence intervals for the unknown matrix, attaining statistical optimality in terms of both the pre-constant and the rate; an independent work by Xia et al.~\cite{xia2019statistical} analyzed a similar de-biasing strategy with the aid of double sample splitting, and shows asymptotic normality of linear forms of the matrix estimator; 
(3) \cite{xia2018confidence,xia2019data} developed a spectral projector to construct confidence regions for singular subspaces in the presence of i.i.d.~additive noise; (4) \cite{koltchinskii2020efficient} considered estimating  linear forms of eigenvectors in a different covariance estimation model, whose analysis relies on the Gaussianity assumption; 
(5) \cite{fan2019asymptotic} characterized the asymptotic normality of bilinear forms of eigenvectors, which accommodates heterogeneous noise;
and (6) \cite{yan2021inference} established the $\ell_{2,\infty}$ distributional guarantees for two spectral estimators (i.e., plain SVD and heteroskedastic PCA) tailored to PCA with heteroskedastic and missing data.

Additionally, the  bulk distribution of the eigenvalues of i.i.d.~random matrices has been studied in the physics literature (e.g.~\cite{sommers1988spectrum, khoruzhenko1996large, brezin1998non, chalker1998eigenvector,feinberg1997non,lytova2018delocalization}), which falls short of the characterizing distributions of extreme eigenvalues. Several more recent papers started to consider a super-position of a low-rank matrix and an i.i.d.~noise matrix, and studied the locations and distributions of the eigenvalues beyond the bulk \cite{tao2013outliers,rajagopalan2015outlier,benaych2016outliers}. These results did not cover a general class of heteroscadestic noise and did not allow the rank $r$ to grow with $n$. What is more, none of these papers studied how to construct valid confidence intervals for extreme eigenvalues, let alone inference for individual eigenvectors.


\section{Discussion}
\label{sec:discussion}

The present paper contributes towards ``fine-grained'' statistical analysis, by developing guaranteed estimation and inference algorithms for linear functionals of the unknown eigenvectors and eigenvalues.  The proposed procedures are model agnostic and are able to accommodate heteroscedastic noise,  without the need of prior knowledge about the noise levels. The validity of our procedures is guaranteed even when the eigen-gap is extremely small, a condition that goes significantly beyond what we have learned from generic matrix perturbation theory. 
The key enabler of our findings lies in an appealing bias reduction feature of eigen-decomposition when coping with asymmetric noise matrices. 

Our studies leave open several interesting questions worthy of future investigation. For instance, our current theory for confidence intervals falls short of accommodating the scenarios when quantity $|\bm{a}^{\top}\bm{u}^{\star}_l|$ far exceeds the associated eigen-gap --- how to determine the fundamental inference limits for such scenarios and, perhaps more importantly, how to attain the limits efficiently?  In addition, our theory is likely suboptimal in terms of the dependency on the rank $r$ and the condition number $\kappa$.  Can we further improve the theoretical support in these regards? 
Furthermore, our analysis framework shed some light on how to perform inference on functions of the eigenvectors.  It would be interesting to develop a unified framework that leads to valid confidence intervals for a broader class of functions (e.g.~quadratic functionals, or more general polynomials) of the eigenvectors. 
Moving beyond estimation and inference for individual eigenvectors, 
we remark that our current theory falls short of delivering useful eigenspace perturbation guarantees unless there exists a sufficient eigen-gap between any adjacent pair of eigenvalues. 
The challenges are at least two-fold when dealing with an asymmetric data matrix: (1) the eigenvectors are, in general, not orthogonal to each other, and (2) the eigenvectors might be complex-valued even when the observed data are real-valued,  both of which are in stark contrast to what happens for a symmetric data matrix.    
How to extend our theory to accommodate more general eigenspace perturbation is left for future investigation. 



\section*{Acknowledgment}

C.~Cheng is supported by the William R.~Hewlett Stanford graduate fellowship. 
Y.~Wei is supported in part by the grants NSF  DMS-2015447/2147546 and CCF-2007911. 
Y.~Chen is supported in part by the AFOSR YIP award FA9550-19-1-0030, by the ONR grant N00014-19-1-2120, by the ARO grants W911NF-20-1-0097 and W911NF-18-1-0303, and by the NSF grants CCF-1907661, IIS-1900140, IIS-2100158 and DMS-2014279. 
C.~Cheng thanks Lihua Lei for helpful discussions.

\appendix

\section{Preliminaries}
\label{sec:preliminaries}

Denote by $\ur_l$ (resp.~$\lambda_l$) the $l$th leading right eigenvector (resp.~eigenvalue) of $\bm{M}=\sum_{j=1}^r \lambda_j^{\star}\bm{u}_j^{\star}\bm{u}_j^{\star\top}+\bm{H}$. We make note of several facts about $\ur_l$ and $\lambda_l$ --- previously established
in \cite{chen2018asymmetry} --- that will prove useful throughout.

To begin with, a simple application of the Neumann series yields the
following expansion \cite[Theorem 2]{chen2018asymmetry}.

\begin{lemma}[\textbf{Neumann expansion}]
\label{lem:Neumann-expansion}
Let $\ur_l$ and $\lambda_l$ be the $l$th leading right eigenvector and the $l$th leading eigenvalue of $\bm{M}$ (cf.~\eqref{eq:model-rank1-noisy}),
	respectively. If $\left\|\bm{H}\right\| < |\lambda_l^{\star}|$, then one has
\begin{equation}
	\label{eq:u-expansion-rankr}
	\bm{u}_l = \sum_{j=1}^r \frac{\lambda_j^\star}{\lambda_l} \big(\bm{u}_j^{\star \top} \bm{u}_l \big) \left\{  \sum_{s=0}^{\infty} \frac{1}{\lambda_l^s} \bm{H}^s\bm{u}_j^\star\right\} .
\end{equation}
\end{lemma}
\noindent As an immediate consequence, we have the following expansion for $\bm{a}^{\top}\bm{u}_l$, which forms the basis of our estimators: 
\begin{align}
\frac{\lambda_{l}}{\lambda_{l}^{\star}\big(\bm{u}_{l}^{\star\top}\bm{u}_{l}\big)}\bm{a}^{\top}\bm{u}_{l} & =\sum_{j=1}^{r}\frac{\lambda_{j}^{\star}}{\lambda_{l}^{\star}}\cdot\frac{\bm{u}_{j}^{\star\top}\bm{u}_{l}}{\bm{u}_{l}^{\star\top}\bm{u}_{l}}\left\{ \sum_{s=0}^{\infty}\frac{1}{\lambda_{l}^{s}}\bm{a}^{\top}\bm{H}^{s}\bm{u}_{j}^{\star}\right\} \nonumber\\
 & =\bm{a}^{\top}\bm{u}_{l}^{\star}+\frac{1}{\lambda_{l}}\bm{a}^{\top}\bm{H}\bm{u}_{l}^{\star}+\sum_{s=2}^{\infty}\frac{1}{\lambda_{l}^{s}}\bm{a}^{\top}\bm{H}^{s}\bm{u}_{l}^{\star}+\sum_{j:j\neq l}\frac{\lambda_{j}^{\star}}{\lambda_{l}^{\star}}\frac{\bm{u}_{j}^{\star\top}\bm{u}_{l}}{\bm{u}_{l}^{\star\top}\bm{u}_{l}}\left\{ \sum_{s=0}^{\infty}\frac{1}{\lambda_{l}^{s}}\bm{a}^{\top}\bm{H}^{s}\bm{u}_{j}^{\star}\right\} .
\label{eq:au-expansion-neumann}
\end{align}

As we shall see shortly, our theory relies heavily on approximating \eqref{eq:au-expansion-neumann} by the lower-order terms (w.r.t.~$\bm{H}$). This requires controlling the influence of the high-order terms. Towards this end, we make note of several useful bounds about $\bm{H}$ that have been established in \cite{chen2018asymmetry}.
\begin{lemma}
\label{lemma:perturbation-H}
Fix any vector $\bm{a}\in\mathbb{R}^{n}$, and suppose that $\|\bm{u}_l^{\star}\|_{\infty}\leq \sqrt{\mu/n}$. Suppose the noise matrix $\bm{H}$ obeys Assumption \ref{assumption-H},
and assume the existence of some sufficiently small constant $c_{1}>0$
such that 
\begin{equation}
 \max\big\{\sigma_{\mathrm{max}} \sqrt{n \log n}, B \log n \big\} \leq c_{1} \lambda_{\mathrm{min}}^{\star} . 
	\label{eq:Bsigma_cond_1}
\end{equation}
Then there exist some universal constants $c_2,c_3>0$ such that with probability at least $1-O(n^{-10})$,  
\begin{subequations}
\begin{align}
	\left|\bm{a}^{\top}\bm{H}^{s}\bm{u}_i^{\star}\right| & \leq\sqrt{\frac{\mu}{n}}\left(c_{2} \max\big\{ \sigma_{\max}\sqrt{n\log n}, B\log n \big\} \right)^{s}\|\bm{a}\|_{2},
\label{eq:linear-form-Hs-general}\\
	\|\bm{H}\| & \leq c_{3} \max\big\{ \sigma_{\max}\sqrt{n\log n}, B\log n \big\} . 
\label{eq:H-operator-norm-general}
\end{align}
\end{subequations}

\end{lemma}
\begin{remark}
The work \cite{chen2018asymmetry} established the bound \eqref{eq:linear-form-Hs-general}  only for the case with $s\leq 20\log n$.  Fortunately, the case with $s> 20\log n$ follows immediately by combining the crude bound  $\left|\bm{a}^{\top}\bm{H}^{s}\bm{u}^{\star}\right| \leq \|\bm{H}\|^s \|\bm{a}\|_2$ and the inequality \eqref{eq:H-operator-norm-general} (by choosing $c_2 = 2c_3$ and using the fact that $(1/2)^s\ll \sqrt{1/n}$). 
\end{remark}

Combining Lemmas~\ref{lem:Neumann-expansion}-\ref{lemma:perturbation-H}, the paper \cite{chen2018asymmetry} establishes the following result.
\begin{lemma}
	\label{lemma:perturbation-linear-form-rankr}
	Consider a rank-r symmetric matrix $\bm{M}^{\star}=\sum_{i=1}^r\lambda_i^{\star}\bm{u}_i^{\star}\bm{u}_i^{\star\top}\in\mathbb{R}^{n\times n}$
	with incoherence parameter $\mu$ (cf.~Definition \ref{defn:incoherence}). Suppose that the noise matrix $\bm{H}$ obeys Assumption \ref{assumption-H} and Condition \eqref{eq:Bsigma_cond_1}.  
	Then for any fixed vector $\bm{a}\in\mathbb{R}^{n}$ and any $1\leq i\leq r$, with probability at least $1-O(n^{-10})$ one has 
	\begin{align}
	\left|\bm{a}^{\top}\left(\ur_{i}-\sum_{k=1}^r\frac{\bm{u}_k^{\star\top}\ur_{i}}{\lambda_i/\lambda_k^{\star}}\bm{u}_k^{\star}\right)\right| & \leq c_2 \frac{ \max\left\{\sigma_{\mathrm{max}} \sqrt{n \log n}, B \log n \right\} }{\lambda_{\mathrm{min}}^{\star}}\sqrt\frac{\mu \kappa^2 r}{n}\|\bm{a}\|_{2}
	\label{eq:linear-form-perturbation-rankr}
	\end{align}
	for some universal constant $c_2>0$. 
	This result holds unchanged if the $i$th right eigenvector $\ur_{i}$ is replaced by the $i$th left eigenvector $\ul_{i}$. 
\end{lemma}

\begin{remark}	
Under the additional condition that $B\log n \leq  \sigma_{\max}\sqrt{n\log n}$, the preceding bounds simplify to
\begin{subequations}
\begin{align}
\left|\bm{a}^{\top}\bm{H}^{s}\bm{u}_i^{\star}\right| & \leq\sqrt{\frac{\mu}{n}}\left(c_{2} \sigma_{\max}\sqrt{n\log n} \right)^{s}\|\bm{a}\|_{2},
\label{eq:linear-form-Hs} \\
\|\bm{H}\| & \leq c_{3} \sigma_{\max}\sqrt{n\log n} , 
\label{eq:H-operator-norm} \\
\left|\bm{a}^{\top}\left(\ur_{i}-\sum_{k=1}^{r}\frac{\bm{u}_{k}^{\star\top}\ur_{i}}{\lambda_{i}/\lambda_{k}^{\star}}\bm{u}_{k}^{\star}\right)\right|
	&\leq c_{2}\frac{\sigma_{\mathrm{max}}\sqrt{\mu\kappa^{2}r\log n}}{\lambda_{\mathrm{min}}^{\star}}\|\bm{a}\|_{2} .
\label{eq:linear-form-perturbation-rankr-simple}
\end{align}
\end{subequations}
Another immediate consequence under our assumptions is that, with probability at least $1-O(n^{-10})$, 
\begin{align}
	\label{eq:H-norm-simple-UB}
	\|\bm{H}\|\leq \lambda_{\min}^{\star} / 10.
\end{align}
\end{remark}

\section{Proof for eigenvector and eigenvalue estimation}
\label{sec:proof-rank-r-all}

In this section, we shall start by proving the eigenvalue perturbation bound stated in Theorem~\ref{thm:rankr-evalue-bounds-simple}, which plays a pivotal role in establishing the eigenvector estimation guarantees in Theorem~\ref{thm:rankr-bounds-simple}.

\subsection{Proof of Theorem~\ref{thm:rankr-evalue-bounds-simple}}
\label{sec:proof-outline-thm:eigengap-condition}

This section aims to establish a slightly stronger version of Theorem~\ref{thm:rankr-evalue-bounds-simple}, stated as follows. 
\begin{theorem}
\label{thm:eigengap-condition} 
Assume that $\mu\kappa^{2}r^{4}\leq c_3 n$ for some sufficiently small constant $c_3>0$. Suppose that the noise parameters defined in Assumption \ref{assumption-H} satisfy
\begin{subequations}
\begin{align}
	\Delta_l^\star > & 2c_4 \kappa^2 r^2  \max\left\{\sigma_{\mathrm{max}} \sqrt{n \log n}, B \log n \right\} \sqrt{\tfrac{\mu}{n}}  \label{eq:deltai-separation-condition} \\
	& \max\left\{\sigma_{\mathrm{max}} \sqrt{n \log n}, B \log n \right\}  \leq c_5 \lambda_{\min}^{\star} 
\end{align}
\end{subequations}
for some sufficiently large (resp.~small) constant $c_4>0$ (resp.~$c_5>0$). Then given any integer $1\leq l\leq r$, with probability $1 - O(n^{-8})$, the eigenvalue $\lambda_l$ and the associated eigenvectors $\ur_l$ and $\ul_l$ (see Notation \ref{notation:leading-evectors}) are all real-valued, and one has
\begin{equation} 
	\label{eq:rankr-eigenvalue-bound}
	|\lambda_l - \lambda_l^\star| \leq c_4 \kappa r^2 \max\left\{\sigma_{\mathrm{max}} \sqrt{n \log n}, B \log n \right\} \sqrt{\tfrac{\mu }{n}}.
\end{equation}
\end{theorem}
A few immediate consequences of this theorem and the Bauer-Fike theorem are summarized as follows.
\begin{corollary}
\label{cor:eigengap-condition-bounds}
Suppose that $\mu\kappa^2 r^4 \lesssim n$, and that Assumptions~\ref{assumption-H}-\ref{assumption:noise-size-rankr-revised} hold.  With probability at least $1-O(n^{-6})$, 
\begin{align}
	|\lambda_{l}-\lambda_{l}^{\star}|\leq\min\Big\{\frac{\Delta_{l}^{\star}}{2},\, c_3\sigma_{\max}\sqrt{n\log n}\Big\}
	\quad \text{and} \quad 
	\Big|\frac{\lambda_{l}-\lambda_{l}^{\star}}{\lambda_{l}^{\star}}\Big|\leq \frac{1}{100}
	\label{eq:lambdal-coarse-bounds}
\end{align}
for some sufficiently small constant $c_3>0$. 
In addition, if either $\sigma_{\max}\sqrt{n\log n} =o( \lambda_{\min}^{\star} )$ or $\Delta_{l}^{\star}=o( \lambda_{\min}^{\star} )$ holds, then with probability at least $1-O(n^{-6})$, 
\begin{align}
	\lambda_l = (1+o(1)) \lambda_l^{\star}.
	\label{eq:lambdal-coarse-bounds2}
\end{align}

\end{corollary}

\subsubsection{Proof outline}
\label{sec:proof-outline-1}

To establish this theorem, we borrow a powerful idea from \cite{tao2013outliers} that converts eigenvalue analysis to  zero counting of certain complex-valued functions.  Specifically, consider the following functions
\begin{subequations}
\begin{align}
	f(z) & \,\defn\, \det \left(\bm{I} +\bm{U}^{\star \top} (\bm{H}- z\bm{I})^{-1} \bm{U}^{\star} \bm{\Sigma}^\star\right), \label{eq:f(z)}\\
	g(z) & \,\defn\, \det \left(\bm{I} + \bm{U}^{\star \top} (-z)^{-1} \bm{U}^\star \bm{\Sigma}^\star\right). \label{eq:g(z)}  
\end{align}
\end{subequations}
%
The intimate connection between these two functions and the eigenvalues of $\bm{M}$ and $\bm{M}^{\star}$ is formalized in the following observation made by \cite{tao2013outliers}. 
\begin{claim}
	\label{lemma:eigenvalue-criterion}
	If $\lambda_{\mathrm{min}} > 2 \left\|\bm{H}\right\|$, then the zeros of $f(\cdot)$ (resp.~$g(\cdot)$) on the region $\mathcal{K} := \left\{z \in \mathbb{C}: |z| > \left\|\bm{H}\right\| \right\} \cup \left\{\infty\right\}$ are exactly the $r$ leading eigenvalues of $\bm{M}$ (resp.~$\bm{M}^\star$). 
\end{claim}

With this claim in mind, we turn attention to studying the zeros of $f(\cdot)$ and $g(\cdot)$. Given that $f(\cdot)$ can be viewed as a perturbed version of $g(\cdot)$ (since $f(\cdot)$ can be obtained by adding a perturbation $\bm{H}$ to $g(\cdot)$ in a certain way), we hope that the zeros of $f(\cdot)$ do not deviate by much from the zeros of $g(\cdot)$. Towards justifying this, we look at the following $\gamma$-neighborhood of $\{\lambda_1^\star, \lambda_2^\star, \cdots, \lambda_r^\star\}$:  
\begin{equation}
	\mathcal{D}(\gamma) \defn 
	\bigcup_{k=1}^r \mathcal{B}(\lambda_k^\star, \gamma),  
	\label{eq:D-gamma}
\end{equation}
where $\mathcal{B}(\lambda, \gamma)$ is a ball of radius $\gamma>0$ centered at $\lambda$, namely, $\mathcal{B}(\lambda, \gamma) \defn \{z \in \mathbb{C}: |z - \lambda| \leq \gamma\}$. A crucial part of the proof is to demonstrate that all zeros of $f(\cdot)$ lie in $\mathcal{D}(\gamma)$ --- or equivalently, in the $\gamma$-neighborhood of the zeros of $g(\cdot)$ --- for sufficiently small $\gamma$. 

\begin{remark} Somewhat surprisingly, $\gamma$ is allowed to be as small as $\frac{\mathrm{poly}\log (n)} {\sqrt{n}}\|\bm{H}\|$ when $r,\kappa,\mu \asymp 1$.   
\end{remark}

In what follows, we shall assume that 
\begin{align}
	\gamma < \lambda_{\min}^{\star}/4 \quad \text{and} \quad \gamma <\Delta_l^{\star}/2.
\end{align}
In view of \eqref{eq:H-norm-simple-UB}, one has $\|\bm{H}\|< \lambda_{\min}^{\star}/4$, thus indicating that 
\begin{align}
	\label{eq:z-lower-bound}
	|z| \geq \lambda_{\min}^{\star} - \gamma > \lambda_{\min}^{\star} / 2 > \|\bm{H}\| \qquad \text{for all }z\in \mathcal{D}(\gamma).
\end{align}
Our proof is based on the following observations:
\begin{itemize}
	\item[(i)] Given that $\gamma < \Delta_l^{\star} / 2$, one has  $\mathcal{B}(\lambda_l^\star, \gamma) \cap  \mathcal{B}(\lambda_k^\star, \gamma) = \emptyset$ for any $k\neq l$, and hence $g(\cdot)$ has exactly 1 zero in $\mathcal{B}(\lambda_l^\star, \gamma)$; in addition, it is clear that $g(\cdot)$ has $l-1$ (resp.~$r-l$) zeros in $\cup_{k: k < l}\mathcal{B}(\lambda_k^\star, \gamma)$ (resp.~$\cup_{k: k > l}\mathcal{B}(\lambda_k^\star, \gamma)$). 

	\item[(ii)] Suppose that in each connected component of $\mathcal{D}(\gamma)$, the functions $f(\cdot)$ and $g(\cdot)$ always have the same number of zeros. If this were true, then one would have (1) a unique zero of $f(\cdot)$ in $\mathcal{B}(\lambda_l^\star, \gamma)$; (2) $l-1$ zeros of $f(\cdot)$ in $\cup_{1\leq k< l}\mathcal{B}(\lambda_k^\star, \gamma)$; (3) $r-l$ zeros of $f(\cdot)$ in $\cup_{k > l}\mathcal{B}(\lambda_k^\star, \gamma)$. 

	\item[(iii)] Recalling how we sort $\{\lambda_k^{\star}\}$ and $\{\lambda_k\}$, we see that the real part of any $z\in \cup_{1\leq k< l}\mathcal{B}(\lambda_k^\star, \gamma)$ must exceed $\lambda_{l-1}^{\star} - \gamma > \lambda_{l}^{\star} + \gamma$ and, similarly, the real part of any $z\in \cup_{k > l}\mathcal{B}(\lambda_k^\star, \gamma)$ is at most $\lambda_{l+1}^{\star} + \gamma < \lambda_{l}^{\star} - \gamma$. Thus, the above arguments reveal that the unique zero of $f(\cdot)$ in $\mathcal{B}(\lambda_l^\star, \gamma)$ is $\lambda_l$, which obeys  
 $|\lambda_l-\lambda_l^{\star}|\leq \gamma$. 

	\item[(iv)] In addition, note that the complex conjugate $\overline{ \lambda_l}$ of $\lambda_l$ is also an eigenvalue of $\bm{M}$, given that $\bm{M}$ is a real-valued matrix.  However, since there is only one zero of $f(\cdot)$ residing in $\mathcal{B}(\lambda_l^\star, \gamma)$, we necessarily have $\overline{ \lambda_l } = \lambda_l$, meaning that $\lambda_l$ is real-valued. A similar argument justifies that the associated right eigenvector $\ur_l$ and left eigenvector $\ul_l$ are real-valued as well. 

\end{itemize}
Theorem \ref{thm:eigengap-condition} is thus established if we are allowed to pick 
\begin{align}
	\label{eq:choice-gamma}
	\gamma = c_4 \kappa r^2  \max\left\{B \log n, \sigma_{\mathrm{max}} \sqrt{n \log n} \right\} \sqrt\frac{\mu}{n}.
\end{align}
Clearly, this choice satisfies $\gamma< \lambda_{\min}^{\star}/4$ and $\gamma < \Delta_l^{\star}/2$ under the assumptions of this theorem.  

The remaining proof then boils down to justifying that $f(\cdot)$ and $g(\cdot)$ have the same number of zeros in each connected component of $\mathcal{D}(\gamma)$. Towards this end, we resort to Rouch\'e's theorem in complex analysis \cite{conway2012functions}.
\begin{theorem}[\textbf{Rouch\'e's theorem}] Let $f(\cdot)$ and $g(\cdot)$ be two complex-valued functions that are holomorphic inside a region $\mathcal{R}$ with closed contour $\partial \mathcal{R}$. If $| f(z) -  g(z)| < |g(z)|$ for all $z\in \partial \mathcal{R}$, then $ f(\cdot)$ and $ g(\cdot)$ have the same number of zeros inside $\mathcal{R}$. 
\end{theorem}
In order to invoke Rouch\'e's theorem to justify the claim in (ii), it suffices to fulfill the requirement $| f(z) -  g(z)| < |g(z)|$ on $\mathcal{D}(\gamma)$. 
Given that $|\bm{z}| > \|\bm{H}\|$ for all $z\in \mathcal{D}(\gamma)$ (see \eqref{eq:z-lower-bound}), we can apply the Neumann series $(\bm{H}-z\bm{I})^{-1}=-\sum_{s=0}^{\infty}z^{-s-1}\bm{H}^{s}$ to reach
\begin{align*}
f(z) & =\det\left(\bm{I}-\bm{U}^{\star\top}\left(\sum_{s=0}^{\infty}z^{-s-1}\bm{H}^{s}\right)\bm{U}^{\star}\bm{\Sigma}^{\star}\right)\nonumber\\
 & =\det\left(\bm{I}+\bm{U}^{\star\top}(-z)^{-1}\bm{U}^{\star}\bm{\Sigma}^{\star}-\sum_{s=1}^{\infty}z^{-s-1}\bm{U}^{\star\top}\bm{H}^{s}\bm{U}^{\star}\bm{\Sigma}^{\star}\right)\nonumber\\
 & =\det\left(\left(\bm{I}+\bm{U}^{\star\top}(-z)^{-1}\bm{U}^{\star}\bm{\Sigma}^{\star}\right)\left(\bm{I}-\left(\bm{I}+\bm{U}^{\star\top}(-z)^{-1}\bm{U}^{\star}\bm{\Sigma}^{\star}\right)^{-1}\sum_{s=1}^{\infty}z^{-s-1}\bm{U}^{\star\top}\bm{H}^{s}\bm{U}^{\star}\bm{\Sigma}^{\star}\right)\right)\nonumber\\
 & =g(z)\det\left(\bm{I}-\left(\bm{I}+\bm{U}^{\star\top}(-z)^{-1}\bm{U}^{\star}\bm{\Sigma}^{\star}\right)^{-1}\sum_{s=1}^{\infty}z^{-s-1}\bm{U}^{\star\top}\bm{H}^{s}\bm{U}^{\star}\bm{\Sigma}^{\star}\right)\nonumber\\
 & = g(z)\det\Bigg(\bm{I}-\underset{=:\, \bm{\Delta}}{\underbrace{\left(\bm{I}-\bm{\Sigma}^{\star}/z\right)^{-1}\sum_{s=1}^{\infty}z^{-s-1}\bm{U}^{\star\top}\bm{H}^{s}\bm{U}^{\star}\bm{\Sigma}^{\star}}}\Bigg),\nonumber
\end{align*}
where the last line holds since $\bm{U}^{\star\top}\bm{U}^\star=\bm{I}$. In addition, observe that the zeros of $g(\cdot)$ are in the interior of $\mathcal{D}$ and that  $g(z) \neq 0$ for all $z\in \partial \mathcal{D}(\gamma)$. Hence, on $\partial \mathcal{D}(\gamma)$ we have
\begin{align}
 & |f(z)-g(z)|<|g(z)| \nonumber\\
\Longleftrightarrow\qquad & |g(z)|\cdot |\det(\bm{I} - \bm{\Delta})-1|<|g(z)| \nonumber\\
	\Longleftrightarrow\qquad & |\det(\bm{I} - \bm{\Delta})-1|<1 .  \label{eq:equivalent-Rouches}
\end{align}
To justify \eqref{eq:equivalent-Rouches} on $\partial \mathcal{D}$, we make the following observations: 
\begin{claim}  
\label{claim:Delta-condition}
	The condition \eqref{eq:equivalent-Rouches} can be guaranteed as long as $\|\bm{\Delta}\|  < 1/(2r)$.
\end{claim}
\begin{claim}  
\label{claim:delta-norm-upper-bound}
There exists a universal constant $c_3>0$ such that with probability $1 - O(n^{-8})$, one has
\begin{equation} \label{eq:delta-norm-upper-bound}
\|\bm{\Delta}\| \leq c_3 \frac{ \max\left\{B \log n, \sigma_{\mathrm{max}} \sqrt{n \log n} \right\}}{\gamma} \sqrt\frac{\mu \kappa^2 r^2}{n} .
\end{equation}
\end{claim}
%

In summary, in order to guarantee $|f(z)-g(z)|<|g(z)|$ on $\partial \mathcal{D}(\gamma)$, it suffices to ensure that  $\|\bm{\Delta}\|\leq 1/(2r)$ (by \eqref{eq:equivalent-Rouches} and Claim~\ref{claim:Delta-condition}), which would hold as long as we take $\gamma=2c_3 \max\left\{B \log n, \sigma_{\mathrm{max}} \sqrt{n \log n} \right\} \sqrt\frac{\mu \kappa^2 r^4}{n}$ (by Claim \ref{claim:delta-norm-upper-bound}). This in turn establishes Theorem \ref{thm:eigengap-condition} (and hence Theorem~\ref{thm:rankr-evalue-bounds-simple}) as long as $c_4\geq 2c_3$.  

Finally, the proofs of the auxiliary claims are postponed to Appendix \ref{sec:proof-claims-theorem-eigenvap-condition}. 
%


\subsubsection{Proofs for auxiliary claims in Appendix~\ref{sec:proof-outline-1}}
\label{sec:proof-claims-theorem-eigenvap-condition}

\paragraph{Proof of Claim \ref{lemma:eigenvalue-criterion}.} Note that both $f(\cdot)$ and $g(\cdot)$ are holomorphic on $\mathcal{K}$. 
Making use of the identity $\det \left(\bm{I} + \bm{A} \bm{B} \right) = \det \left(\bm{I} + \bm{B}\bm{A} \right)$, we have
\begin{align}
f(z) & = \det\left(\bm{I} + \bm{U}^\star \bm{\Sigma}^\star \bm{U}^{\star \top} ( \bm{H}- \bm{z} \bm{I})^{-1}\right) 
= \det\left(\bm{I} + \bm{M}^{\star } ( \bm{H}- \bm{z} \bm{I})^{-1}\right) \nonumber \\
& = \det \left((\bm{H}-z\bm{I}+\bm{M}^\star)(\bm{H}-z\bm{I})^{-1}\right) \nonumber \\
& = \frac{\det (\bm{M}-z\bm{I})}{\det(\bm{H}-z\bm{I})} \\
\text{and} \qquad
g(z) & = \det\left(\bm{I} + \bm{U}^\star \bm{\Sigma}^\star \bm{U}^{\star \top}/(-z)\right) 
  = \prod_{i=1}^r \left(1 - \frac{\lambda_i^\star}{z}\right).
\end{align}
These identities make clear that the zeros of $g(\cdot)$ (resp.~$f(\cdot)$) on $\mathcal{K}$ are all eigenvalues of $\bm{M}^{\star}$ (resp.~$\bm{M}$). In particular, the zeros of $g(\cdot)$ are precisely $\{\lambda_i^\star\}_{1\leq i\leq r}$. 

It remains to show that there are exactly $r$ zeros of $f(\cdot)$ lying in $\mathcal{K}$, and that they are exactly $\lambda_1, \cdots, \lambda_r$. This is equivalent to showing that the set of eigenvalues of $\bm{M}$ contained in the region $\mathcal{K}$ is $\{\lambda_i\}_{1\leq i\leq r}$. Towards this, define $\bm{M}(t) = \bm{M}^\star + t \bm{H}$. From the Bauer-Fike theorem, all eigenvalues of $\bm{M}(t)$ lie in the set
\begin{equation}
\mathcal{B}(0,t\|\bm{H}\|)\,\cup\, \Big\{  \underset{=: \,\mathcal{D}(t\|\bm{H}\|)}{\underbrace{\cup_{k=1}^{r}\mathcal{B}(\lambda_{k}^{\star},t\|\bm{H}\|)}} \Big\}.
\end{equation}
Given that $\|\bm{H}\| < \lambda_{\mathrm{min}} / 4$ (according to \eqref{eq:H-norm-simple-UB} under our assumptions), it is easily seen that $\mathcal{B}(0, t\|\bm{H}\|)$ does not intersect with $\mathcal{D}(t \|\bm{H}\|)$ for all $0\leq t \leq 1$. Additionally, the set of the eigenvalues of $\bm{M}(t)$ depends continuously on $t$ \cite[Theorem 6]{embree2001generalizing}, requiring $\bm{M}(t)$ to have the same number of zeros in $\mathcal{B}(0, t\|\bm{H}\|)$ for all $0\leq t\leq 1$. 
Hence,  $\bm{M}=\bm{M}(1)$ has exactly $r$ eigenvalues in $\mathcal{D}(\|\bm{H}\|)$ (since $\bm{M}^{\star}=\bm{M}(0)$ has $r$ eigenvalues in this region).
These eigenvalues necessarily have magnitudes larger than any point in $\mathcal{B}(0, \|\bm{H}\|)$ (since $\|\bm{H}\| < \lambda_{\mathrm{min}} / 4$), thus 
 indicating that the eigenvalues of $\bm{M}$ in $\mathcal{K}$ are precisely $\lambda_1,\cdots,\lambda_r$.

\paragraph{Proof of Claim \ref{claim:Delta-condition}.} 
Denoting by $\mu_{1},\cdots,\mu_{r}\in\mathbb{C}$ the $r$ eigenvalues
of $\bm{I}- \bm{\Delta}\in \mathbb{R}^{r\times r}$, one has
\[
\big|\det(\bm{I}-\bm{\Delta})-1\big|=\Big|\prod_{i=1}^{r}\mu_{i}-1\Big|.
\]
By virtue of the elementary inequality\footnote{This holds since 
\begin{align*}
	\left|\prod_{i=1}^{r}(1+a_{i})-1\right|
	=\Big|\sum_{k=1}^{r}\sum_{1\leq i_{1},\cdots,i_{k}\leq r}a_{i_{1}}\cdots a_{i_{k}}\Big| 
	 \leq\sum_{k=1}^{r}\sum_{1\leq i_{1},\cdots,i_{k}\leq r}|a_{i_{1}}| \cdots |a_{i_{k}}|=\prod_{i=1}^{r}(1+|a_{i}|)-1.
\end{align*}
} $\left|\prod_{i=1}^{r}(1+a_{i})-1\right|\leq\prod_{i=1}^{r}(1+|a_{i}|)-1$, we arrive at
\[
\big|\det(\bm{I}-\bm{\Delta})-1\big|\leq\prod_{i=1}^{r}\left(1+|\mu_{i}-1|\right)-1\leq(1+\|\bm{\Delta}\|)^{r}-1,
\]
where the last inequality follows from the Bauer-Fike theorem (which
forces that $|\mu_{i}-1|\leq\|\bm{\Delta}\|$). This means that Condition \eqref{eq:equivalent-Rouches} holds if $(1+\|\bm{\Delta}\|)^r<2$; the latter condition is clearly guaranteed if $\|\bm{\Delta}\|\leq 1/(2r)$.

\paragraph{Proof of Claim \ref{claim:delta-norm-upper-bound}.}

Since $\gamma < \lambda_{\mathrm{min}} / 4$, it always holds that $|z| > 3\lambda_{\mathrm{min}} / 4$ for all $z\in \partial \mathcal{D}(\gamma)$. Hence
\begin{align}
\|\bm{\Delta}\| = \Big\|\left(z\bm{I} - \bm{\Sigma}^\star \right)^{-1}\sum_{s=1}^{\infty} z^{-s}\bm{U}^{\star \top} \bm{H}^s \bm{U}^\star \bm{\Sigma}^\star\Big\| \nonumber 
	& \leq \big\|\left(z\bm{I} - \bm{\Sigma}^\star\right)^{-1} \big\| \cdot \sum_{s=1}^{\infty} \frac{\big\|\bm{U}^{\star \top} \bm{H}^s \bm{U}^\star \bm{\Sigma}^\star\big\|}{|z|^{s}} \nonumber \\ 
	& \leq \|\bm{\Sigma}^\star \| \cdot \big\|\left(z\bm{I} - \bm{\Sigma}^\star\right)^{-1} \big\| \cdot \sum_{s=1}^{\infty}  \frac{\big\|\bm{U}^{\star \top} \bm{H}^s \bm{U}^\star \big\|}{|z|^{s}} . 
	\label{eq:delta-norm-bound-1} 
\end{align}
Recall that $\|\bm{\Sigma}^\star\| = \lambda_{\mathrm{max}}$. Also,  on $\partial \mathcal{D}(\gamma)$ we have
\begin{equation}
	\big\|\left(z\bm{I} - \bm{\Sigma}^\star \right)^{-1} \big\|  
	\leq \max_{1 \leq k \leq r, \,z \in \partial \mathcal{D}(\gamma) } \left|\frac{1}{z-\lambda_k^\star}\right| 
	\leq \frac{1}{\gamma}.
\end{equation}
Taken collectively, the above results yield
\begin{align}
\|\bm{\Delta}\| 
& \leq \frac{\lambda_{\mathrm{max}}}{\gamma} \sum_{s=1}^{\infty} \frac{\left\|\bm{U}^{\star \top} \bm{H}^s \bm{U}^\star\right\|}{\left(\frac{3}{4} \lambda_{\mathrm{min}}\right)^s} 
\stackrel{\mathrm{(i)}}{\leq} \frac{\lambda_{\mathrm{max}}r}{\gamma}  \sum_{s=1}^{\infty}\frac{ \max_{1 \leq i, j \leq r} |\bm{u}_i^{\star \top} \bm{H}^s \bm{u}_j^\star |}{\left(\frac{3}{4} \lambda_{\mathrm{min}}\right)^s},
\end{align}
where $\mathrm{(i)}$ makes use of the elementary inequality $\|\bm{A}\| \leq r \|\bm{A}\|_\infty$ for any 
$\bm{A} \in \mathbb{R}^{r \times r}$. Invoke Lemma \ref{lemma:perturbation-H} and a union bound to show that: with probability $1-O(n^{-10}r^2)$, 
\begin{align*}
\|\bm{\Delta}\| & \leq\frac{\lambda_{\mathrm{max}}r}{\gamma}\sum_{s=1}^{\infty}\left(\frac{c_{2}\max\left\{ \sigma_{\mathrm{max}}\sqrt{n\log n},  B\log n\right\} }{\frac{3}{4}\lambda_{\mathrm{min}}}\right)^{s}\sqrt{\frac{\mu}{n}}\nonumber\\
& \overset{\mathrm{(ii)}}{\leq}\frac{\lambda_{\mathrm{max}}r}{\gamma}\cdot\frac{c_{2}\max\left\{ \sigma_{\mathrm{max}}\sqrt{n\log n}, B\log n\right\} }{\frac{3}{4}\lambda_{\mathrm{min}}}\cdot\sum_{s=1}^{\infty}\left(\frac{1}{2}\right)^{s}\sqrt{\frac{\mu}{n}}\nonumber\\
 & \leq \frac{4c_2 \max\left\{ \sigma_{\mathrm{max}}\sqrt{n\log n}, B\log n\right\} }{3\gamma}\sqrt{\frac{\mu\kappa^{2}r^{2}}{n}},
\end{align*}
where (ii) holds as long as $\frac{c_{2}\max\left\{ B\log n,\sigma_{\mathrm{max}}\sqrt{n\log n}\right\} }{\frac{3}{4}\lambda_{\mathrm{min}}}\leq\frac{1}{2}$.

\subsection{Proof of Theorem~\ref{thm:rankr-bounds-simple}} \label{sec:proof-of-rankr-bounds}

As mentioned previously, Theorem~\ref{thm:rankr-evalue-bounds-simple} enables us to prove the statistical guarantees  stated in Theorem \ref{thm:rankr-bounds-simple}. In the sequel, we shall first establish perturbation bounds for the vanilla eigenvector estimator (i.e.~$\bm{u}_l$ and $\bm{w}_l$) as well as the estimator $\widehat{\bm{u}}_{l}=\frac{1}{\|\bm{u}_{l}+\bm{w}_{l}\|_{2}}(\bm{u}_{l}+\bm{w}_{l})$. 
\begin{theorem} 
\label{thm:rankr-bounds}
Consider any $1\leq l\leq r$. Suppose that $\mu\kappa^2 r^4 \lesssim n$ and that Assumptions~\ref{assumption-H}-\ref{assumption:noise-size-rankr-revised} hold.
Then with probability at least $1-O(n^{-6})$, we have: 
\begin{subequations}
\begin{enumerate}
	\item ($\ell_2$ perturbation of the $l$th eigenvector)
	\begin{align} 
		\min \left\Vert \bm{u}_{l}\pm\bm{u}_{l}^{\star}\right\Vert _{2}  
		& \lesssim
		\frac{\kappa^{2}\sigma_{\max}\sqrt{\mu r^{2}\log n}}{\Delta_{l}^{\star}}
		+
		\frac{\kappa^{2}\sigma_{\max}\sqrt{n\log n}}{\lambda_{\max}^{\star}}
	 	\label{eq:rankr-l2-bound-2} \\
		\left|\bm{u}_{k}^{\star\top}\bm{u}_{l}\right|	
		& \lesssim
		\tfrac{\sigma_{\mathrm{max}}}{\Delta_{l}^{\star}}\sqrt{\mu\kappa^{4}r\log n}, \qquad k\neq l, \label{eq:ui-uk-ub-thm} \\
		|\bm{u}_{l}^{\star\top}\bm{u}_{l}| 
		&\geq
		1-O\left((\kappa^{4}\sigma_{\max}^{2}n\log n)\left\{ \tfrac{1}{(\Delta_{l}^{\star})^{2}}\tfrac{\mu r^{2}}{n}+\tfrac{1}{(\lambda_{\mathrm{max}}^{\star})^{2}}\right\} \right); \label{eq:ui-ui-ub-thm}
	\end{align}
	\item (perturbation of linear forms of the $l$th eigenvector) for any fixed vector $\bm{a}$ with $\|\bm{a}\|_2=1$,
	\begin{align} 
		\min\left|\bm{a}^{\top}\left(\bm{u}_{l}\pm\bm{u}_{l}^{\star}\right)\right| & \lesssim\,\tfrac{\sigma_{\mathrm{max}}\sqrt{\mu\kappa^{2}r^{4}\log n}}{\lambda_{\mathrm{min}}^{\star}}+\sigma_{\mathrm{max}}\sqrt{\mu\kappa^{4}r^{3}\log n}\cdot\max_{k\neq l}\tfrac{|\bm{a}^{\top}\bm{u}_{k}^{\star}|}{|\lambda_{l}^{\star}-\lambda_{k}^{\star}|}\nonumber\\
		+ & \tfrac{\sigma_{\max}^{2}\mu\kappa^{4}r^{2}\log n}{(\Delta_{l}^{\star})^{2}}|\bm{a}^{\top}\bm{u}_{l}^{\star}|+\tfrac{\sigma_{\max}^{2}\kappa^{4}n\log n}{(\lambda_{\max}^{\star})^{2}}|\bm{a}^{\top}\bm{u}_{l}^{\star}|;
		\label{eq:rankr-linear-form-bound}
	\end{align}
	\item ($\ell_{\infty}$ perturbation of the $l$th eigenvector) 
	\begin{align} 
		\min  \left\Vert \bm{u}_{l}\pm\bm{u}_{l}^{\star}\right\Vert _{\infty}  \lesssim\, & \tfrac{\sigma_{\mathrm{max}}}{\lambda_{\mathrm{min}}^{\star}}\sqrt{\mu\kappa^{4}r\log n}+\tfrac{\sigma_{\mathrm{max}}}{\Delta_{l}^{\star}}\sqrt{\tfrac{\mu^{2}\kappa^{4}r^{3}\log n}{n}}. 
		\label{eq:rankr-l-infty-bound}
	\end{align}
	%
	%
\end{enumerate}
\end{subequations}
In addition, the above results hold unchanged if $\bm{u}_{l}$ is replaced by either $\bm{w}_{l}$ or $\widehat{\bm{u}}_{l}=\frac{1}{\|\bm{u}_{l}+\bm{w}_{l}\|_{2}}(\bm{u}_{l}+\bm{w}_{l})$.  
\end{theorem}
\begin{proof} See Appendix~\ref{sec:proof-theorem-rankr-bounds}.\end{proof}
%
%

In words, Theorem~\ref{thm:rankr-bounds} reveals appealing statistical performance for the estimators $\bm{u}_l$, $\bm{w}_l$ and $\widehat{\bm{u}}_l$. However, when estimating linear functionals of eigenvectors via, say, the plug-in estimator  $\bm{a}^{\top}\bm{u}_l$, the bound \eqref{eq:rankr-linear-form-bound} suffers from an additional term $\tfrac{\sigma_{\max}^{2}\kappa^{4}n\log n}{(\lambda_{\max}^{\star})^{2}}|\bm{a}^{\top}\bm{u}_{l}^{\star}|$
compared to the desired bound \eqref{eq:rankr-linear-form-bound-simple} in Theorem~\ref{thm:rankr-bounds-simple}.  As it turns out, this extra term arises due to a systematic bias of the plug-in estimates. To compensate for the bias, one needs to properly enlarge the plug-in estimator, thus leading us to the proposed estimators $\widehat{u}_{\bm{a},l}$. We now establish the claimed performance guarantees for these properly corrected estimators for $\bm{a}^{\top}\bm{u}_l^{\star}$; the proof is deferred to Appendix \ref{sec:thm-linear-form-ua-coarse}, which builds heavily upon the analysis of Theorem~\ref{thm:rankr-bounds}.

\begin{theorem}
	\label{thm:linear-form-ua-coarse}
	Instate the assumptions of Theorem~\ref{thm:rankr-bounds}. Fix any vector $\bm{a}$ with $\|\bm{a}\|_2=1$. With probability exceeding $1-O(n^{-6})$, the estimator $\widehat{u}_{\bm{a},l} := \min \left\{\sqrt{\Big|\frac{(\bm{a}^{\top}\bm{u}_{l})(\bm{a}^{\top}\bm{w}_{l})}{\bm{u}_{l}^{\top}\bm{w}_{l}}\Big|}, 1 \right\}$ satisfies
	\begin{align}
		\min\left|\widehat{u}_{\bm{a},l}\pm\bm{a}^{\top}\bm{u}_{l}^{\star}\right| 
		& \lesssim\,\tfrac{\sigma_{\mathrm{max}}r^{2}\sqrt{\mu\kappa^{2}\log n}}{\lambda_{\mathrm{min}}^{\star}}+\tfrac{\sigma_{\max}^{2}\mu r^{2}\kappa^{4}\log n}{(\Delta_{l}^{\star})^{2}}|\bm{a}^{\top}\bm{u}_{l}^{\star}|+
		\sigma_{\mathrm{max}}\sqrt{\mu\kappa^{4}r^{3}\log n}\max_{k\neq l}\tfrac{|\bm{a}^{\top}\bm{u}_{k}^{\star}|}{|\lambda_{l}^{\star}-\lambda_{k}^{\star}|}. 
		\label{eq:au-coarse-estimate-UB}
	\end{align}	
	%
\end{theorem}

Theorems~\ref{thm:rankr-bounds}-\ref{thm:linear-form-ua-coarse} taken together establish Theorem~\ref{thm:rankr-bounds-simple}.

\subsubsection{Proof of Theorem \ref{thm:rankr-bounds}}
\label{sec:proof-theorem-rankr-bounds}

We shall start by proving the results for $\bm{u}_l$; the proofs for the results w.r.t.~$\bm{w}_l$ are clearly identical.  

\paragraph{Proofs for $\ell_2$ perturbation bounds.}
To derive the $\ell_2$ bound for $\bm{u}_l$, we start by considering the distance between $\bm{u}_l$ and the subspace spanned by the true eigenvectors $\{\bm{u}_k^\star\}_{1 \leq k \leq r}$. The Neumann series (cf.~Lemma~\ref{lem:Neumann-expansion}) tells us that: when $\|\bm{H}\| \leq \lambda_{\mathrm{min}}^{\star} / 4$,
\begin{align}
\left\Vert \bm{u}_{l}-\sum_{k=1}^{r}\frac{\lambda_{k}^{\star}\bm{u}_{k}^{\star\top}\bm{u}_{l}}{\lambda_{l}}\bm{u}_{k}^{\star}\right\Vert _{2} & =\left\Vert \sum_{k=1}^{r}\frac{\lambda_{k}^{\star}\left(\bm{u}_{k}^{\star\top}\bm{u}_{l}\right)}{\lambda_{l}}\left\{ \sum_{s=1}^{\infty}\frac{1}{\lambda_{l}^{s}}\bm{H}^{s}\bm{u}_{k}^{\star}\right\} \right\Vert _{2} \nonumber\\
 & \leq\sqrt{\sum_{k=1}^{r}\left|\frac{\lambda_{k}^{\star}\left(\bm{u}_{k}^{\star\top}\bm{u}_{l}\right)}{\lambda_{l}}\right|^{2}}\left\Vert \Bigg(\sum_{s=1}^{\infty}\frac{1}{\lambda_{l}^{s}}\bm{H}^{s}\Bigg)\left[\bm{u}_{1}^{\star},\cdots,\bm{u}_{r}^{\star}\right]\right\Vert \nonumber\\
 & \leq \sqrt{\frac{(\lambda_{\max}^{\star})^{2}}{|\lambda_{l}|^{2}}\sum_{k=1}^{r}\left|\bm{u}_{k}^{\star\top}\bm{u}_{l}\right|^{2}}\left\Vert \Bigg(\sum_{s=1}^{\infty}\frac{1}{\lambda_{l}^{s}}\bm{H}^{s}\Bigg)\right\Vert \Big\|\left[\bm{u}_{1}^{\star},\cdots,\bm{u}_{r}^{\star}\right]\Big\| \nonumber\\
 & \leq\frac{\lambda_{\max}^{\star}}{|\lambda_{l}|}\sum_{s=1}^{\infty}\left(\frac{\left\Vert \bm{H}\right\Vert }{|\lambda_{l}|}\right)^{s},  \label{eq:UB15}
\end{align}
where the last inequality holds since $\|\left[\bm{u}_{1}^{\star},\cdots,\bm{u}_{r}^{\star}\right]\|=1$
and $\sum_{k=1}^{r}\left|\bm{u}_{k}^{\star\top}\bm{u}_{l}\right|^{2}\leq\|\bm{u}_{l}\|_{2}^{2}=1$
(given that the $\bm{u}_{k}^{\star}$'s are orthonormal).
In view of the Bauer-Fike theorem and the bound $\|\bm{H}\|\leq \lambda_{\min}^{\star}/4$,  we have the crude lower bound
\begin{align}
	|\lambda_l| \geq  \lambda_{\min}^{\star} - \|\bm{H}\| \geq  3\lambda_{\min}^{\star}/4.
\end{align}
This taken collectively with the inequality \eqref{eq:UB15} yields
\begin{align}
\left\|\bm{u}_l - \sum_{k=1}^r \frac{\lambda_k^\star \bm{u}_k^{\star \top}\bm{u}_l }{\lambda_l}\bm{u}_k^\star\right\|_2  
& \leq \frac{\lambda_{\max}^{\star}}{|\lambda_l|} \cdot  \frac{\|\bm{H}\|}{|\lambda_l| - \left\|\bm{H}\right\|} \leq \frac{\lambda_{\max}^{\star}}{\frac{3}{4} \lambda_{\min}^{\star}} \cdot \frac{ \|\bm{H}\| }{\frac{1}{2} \lambda_{\min}^{\star}} = \frac{8\kappa^2}{3} \cdot \frac{\|\bm{H}\|}{\lambda_{\max}^{\star}}.  
\label{eq:ui-bound-14}
\end{align}
In addition, note that the Euclidean projection of $\bm{u}_l$ onto the subspace spanned by $\{\bm{u}_k^\star\}_{1 \leq k \leq r}$ is given by $\sum_{k=1}^r (\bm{u}_k^{\star \top} \bm{u}_l) \bm{u}_k^\star$. This together with \eqref{eq:ui-bound-14} implies that 
\begin{align}
\frac{8\kappa^{2}}{3}\cdot\frac{\|\bm{H}\|}{\lambda_{\max}^{\star}} & \geq\left\Vert \bm{u}_{l}-\sum_{k=1}^{r}\frac{\lambda_{k}^{\star}\bm{u}_{k}^{\star\top}\bm{u}_{l}}{\lambda_{l}}\bm{u}_{k}^{\star}\right\Vert _{2}\geq\min_{\bm{z}\in\mathsf{span}\{\bm{u}_{1}^{\star},\cdots,\bm{u}_{r}^{\star}\}}\left\Vert \bm{u}_{l}-\bm{z}\right\Vert _{2}\nonumber\\
 & =\left\Vert \bm{u}_{l}-\sum_{k=1}^{r}(\bm{u}_{k}^{\star\top}\bm{u}_{l})\bm{u}_{k}^{\star}\right\Vert _{2}
 =\sqrt{1-\sum_{k=1}^{r}\left|\bm{u}_{k}^{\star\top}\bm{u}_{l}\right|^{2}},
\label{eq:l2-subspace-bound}
\end{align}
where the last identity arises from the Pythagorean theorem.

The above inequality \eqref{eq:l2-subspace-bound} indicates that most of the energy of $\bm{u}_l$ lies in $\mathsf{span}\{\bm{u}_1^{\star},\cdots,\bm{u}_r^{\star}\}$. In order to show that $|\bm{u}_{l}^{\star\top}\bm{u}_{l}| \approx 1$, it suffices to justify that $|\bm{u}_{k}^{\star\top}\bm{u}_{l}|$ is very small for any $k\neq l$. To this end, taking $\bm{a} = \bm{u}_k^\star$ in Lemma \ref{lemma:perturbation-linear-form-rankr} and making use of Condition \eqref{eq:condition-lambdamin-rankr-revised}, we arrive at
\begin{align}
\left|\bm{u}_{k}^{\star\top}\left(\bm{u}_{l}-\sum_{i=1}^{r}\frac{\bm{u}_{i}^{\star\top}\bm{u}_{l}}{\lambda_{l}/\lambda_{i}^{\star}}\bm{u}_{i}^{\star}\right)\right| & \leq c_{2}\frac{\sigma_{\mathrm{max}}}{\lambda_{\mathrm{min}}^{\star}}\sqrt{\mu\kappa^{2}r\log n}.
\label{eq:uk-ui-expansion-bound}
\end{align}
In addition, the orthonormality of $\{\bm{u}_i^{\star}\}_{1\leq i\leq r}$ indicates that
\begin{align}
\bm{u}_{k}^{\star\top}\left(\bm{u}_{l}-\sum_{i=1}^{r}\frac{\bm{u}_{i}^{\star\top}\bm{u}_{l}}{\lambda_{l}/\lambda_{i}^{\star}}\bm{u}_{i}^{\star}\right)=\bm{u}_{k}^{\star\top}\bm{u}_{l}-\frac{\lambda_{k}^{\star}}{\lambda_{l}}\bm{u}_{k}^{\star\top}\bm{u}_{l}=\Big(1-\frac{\lambda_{k}^{\star}}{\lambda_{l}}\Big)\bm{u}_{k}^{\star\top}\bm{u}_{l}
\label{eq:uk-ui-expansion}
\end{align}
and, as a result, 
\begin{align}
\left|\bm{u}_{k}^{\star\top} \bm{u}_{l} \right| & \leq  \left|1-\frac{\lambda_{k}^{\star}}{\lambda_{l}}\right|^{-1} c_{2}\frac{\sigma_{\mathrm{max}}}{\lambda_{\mathrm{min}}^{\star}}\sqrt{\mu\kappa^{2}r\log n}.
\label{eq:uk-ui-expansion-bound2}
\end{align} 
In order to control $\left|1-\frac{\lambda_{k}^{\star}}{\lambda_{l}}\right|^{-1}$, we note that: for all $k\neq l$, it follows from the triangle inequality and Theorem~\ref{thm:eigengap-condition} that 
\begin{align}
	& |\lambda_{k}^{\star}-\lambda_{l}|\geq|\lambda_{k}^{\star}-\lambda_{l}^{\star}|-|\lambda_{l}^{\star}-\lambda_{l}|\geq|\lambda_{k}^{\star}-\lambda_{l}^{\star}|-\Delta_{l}^{\star}/2\geq|\lambda_{k}^{\star}-\lambda_{l}^{\star}|/2\geq\Delta_{l}^{\star}/2,	 
	\label{eq:lambdak-lambdai-gap-LB} \\
	& \Longrightarrow \qquad 
\left|1-\frac{\lambda_{k}^{\star}}{\lambda_{l}}\right|^{-1}  =\frac{|\lambda_{l}|}{|\lambda_{k}^{\star}-\lambda_{l}|}\leq\frac{\lambda_{\max}^{\star}+\|\bm{H}\|}{\Delta_{l}^{\star}/2}\leq\frac{4\lambda_{\max}^{\star}}{\Delta_{l}^{\star}}, \nonumber
\end{align}
where we have also used the condition $\|\bm{H}\|\leq \lambda^{\star}_{\max}$. Substitution into \eqref{eq:uk-ui-expansion-bound2} yields
\begin{align}
\left|\bm{u}_{k}^{\star\top}\bm{u}_{l}\right| & \leq\frac{4\lambda_{\max}^{\star}}{\Delta_{l}^{\star}}\cdot c_{2}\frac{\sigma_{\mathrm{max}}}{\lambda_{\min}^{\star}}\sqrt{\mu\kappa^{2}r\log n}=\frac{4c_{2}\sigma_{\mathrm{max}}}{\Delta_{l}^{\star}}\sqrt{\mu\kappa^{4}r\log n}
\label{eq:uk-ui-UB2}
\end{align}
for all $k\neq l$. 
Putting the above results together, we arrive at
\begin{align}
\min\left\{ \left\Vert \bm{u}_{l}\pm\bm{u}_{l}^{\star}\right\Vert _{2}^{2}\right\}  & =2-2|\bm{u}_{l}^{\star\top}\bm{u}_{l}|\leq2-2|\bm{u}_{l}^{\star\top}\bm{u}_{l}|^{2}\nonumber\\
 & \leq2\left\{ 1-\sum_{k=1}^{r}|\bm{u}_{k}^{\star\top}\bm{u}_{l}|^{2}\right\} +2\left\{ \sum_{k\neq l,k=1}^{r}|\bm{u}_{k}^{\star\top}\bm{u}_{l}|^{2}\right\} \nonumber \nonumber\\
 & \lesssim\left(\kappa^{2}\frac{\left\Vert \bm{H}\right\Vert }{\lambda_{\max}^{\star}}\right)^{2}+\left(\frac{\sigma_{\mathrm{max}}}{\Delta_{l}^{\star}}\right)^{2}\mu\kappa^{4}r^{2}\log n \nonumber\\
 & \lesssim\left(\kappa^{2}\frac{\sigma_{\max}\sqrt{n\log n}}{\lambda_{\max}^{\star}}\right)^{2}+\left(\frac{\sigma_{\mathrm{max}}}{\Delta_{l}^{\star}}\right)^{2}\mu\kappa^{4}r^{2}\log n,  \label{eq:ui-uistar-l2-bound}
\end{align}
where the penultimate inequality relies on \eqref{eq:l2-subspace-bound} and \eqref{eq:uk-ui-UB2}, and the last line follows from \eqref{eq:H-operator-norm}. 

Finally, the advertised bound on $|\bm{u}_l^{\star\top}\bm{u}_l|$ can be immediately established by combining the identity $|\bm{u}_{l}^{\star\top}\bm{u}_{l}|=1-\frac{1}{2}\min \{ \Vert \bm{u}_{l}\pm\bm{u}_{l}^{\star} \Vert _{2}^{2} \} $ and the above bound \eqref{eq:ui-uistar-l2-bound}.

\paragraph{Proof for perturbation of linear forms of eigenvectors.}
Invoke the triangle inequality to obtain
\begin{align}
\left|\bm{a}^{\top}\left(\bm{u}_{l}-\frac{\bm{u}_{l}^{\star\top}\bm{u}_{l}}{\lambda_{l}/\lambda_{l}^{\star}}\bm{u}_{l}^{\star}\right)\right| & \leq\left|\bm{a}^{\top}\left(\bm{u}_{l}-\sum_{k=1}^{r}\frac{\bm{u}_{k}^{\star\top}\bm{u}_{l}}{\lambda_{l}/\lambda_{k}^{\star}}\bm{u}_{k}^{\star}\right)\right|+\sum_{k\neq l,k=1}^{r}\left|\frac{\bm{u}_{k}^{\star\top}\bm{u}_{l}}{\lambda_{l}/\lambda_{k}^{\star}}\left(\bm{a}^{\top}\bm{u}_{k}^{\star}\right)\right|. 
\label{eq:linear-form-decompose-1}
\end{align}
The first term on the right-hand side of \eqref{eq:linear-form-decompose-1} can be controlled via  Lemma \ref{lemma:perturbation-linear-form-rankr} and Condition \eqref{eq:condition-lambdamin-rankr-revised}: 
\begin{align*}
\left|\bm{a}^{\top}\left(\bm{u}_{l}-\sum_{k=1}^{r}\frac{\bm{u}_{k}^{\star\top}\bm{u}_{l}}{\lambda_{l}/\lambda_{k}^{\star}}\bm{u}_{k}^{\star}\right)\right| & \lesssim\frac{\sigma_{\mathrm{max}}}{\lambda_{\mathrm{min}}^{\star}}\sqrt{\mu\kappa^{2}r\log n}\cdot\|\bm{a}\|_{2}.
\end{align*} 
Regarding the second term on the right-hand side of \eqref{eq:linear-form-decompose-1}, we can invoke \eqref{eq:uk-ui-expansion-bound2} to derive
\begin{equation}
\left|\frac{\bm{u}_{k}^{\star\top}\bm{u}_{l}}{\lambda_{l}/\lambda_{k}^{\star}}\right|\lesssim\frac{|\lambda_{k}^{\star}|}{|\lambda_{l}-\lambda_{k}^{\star}|}\frac{\sigma_{\max}}{\lambda_{\min}^{\star}}\sqrt{\mu\kappa^{2}r\log n}\leq\frac{\lambda_{\max}^{\star}}{|\lambda_{l}^{\star}-\lambda_{k}^{\star}|/2}\frac{\sigma_{\max}}{\lambda_{\min}^{\star}}\sqrt{\mu\kappa^{2}r\log n}\lesssim\frac{\sigma_{\max}}{|\lambda_{l}^{\star}-\lambda_{k}^{\star}|}\sqrt{\mu\kappa^{4}r\log n}, \label{eq:rank-r-inner-product-bound}
\end{equation}
where the penultimate inequality results from \eqref{eq:lambdak-lambdai-gap-LB}. As a consequence, one has
\begin{align}
\sum_{k:k\neq l,1\leq k\leq r}\left|\frac{\bm{u}_{k}^{\star\top}\bm{u}_{l}}{\lambda_{l}/\lambda_{k}^{\star}}\left(\bm{a}^{\top}\bm{u}_{k}^{\star}\right)\right|\lesssim r\cdot\sigma_{\max}\sqrt{\mu\kappa^{4}r\log n}\cdot\left\{ \max_{k:k\neq l}\frac{|\bm{a}^{\top}\bm{u}_{k}^{\star}|}{\big|\lambda_{l}^{\star}-\lambda_{k}^{\star}\big|}\right\} .
\end{align}
Substituting the above bounds into \eqref{eq:linear-form-decompose-1} and rearranging terms, we obtain
\begin{equation}
\left|\bm{a}^{\top}\left(\bm{u}_{l}-\frac{\bm{u}_{l}^{\star\top}\bm{u}_{l}}{\lambda_{l}/\lambda_{l}^{\star}}\bm{u}_{l}^{\star}\right)\right|\lesssim\sigma_{\mathrm{max}}\sqrt{\mu\kappa^{2}r\log n}\left\{ \frac{1}{\lambda_{\mathrm{min}}^{\star}}\|\bm{a}\|_{2}+\kappa r\max_{k:k\neq l}\frac{|\bm{a}^{\top}\bm{u}_{k}^{\star}|}{\big|\lambda_{l}^{\star}-\lambda_{k}^{\star}\big|}\right\} . \label{eq:linear-form-decompose-1-bound}
\end{equation}
%

%

%


Next, recall from Theorem~\ref{thm:eigengap-condition} that
\begin{align*}
|\lambda_{l}-\lambda_{l}^{\star}| & \leq
c_4 \kappa r^2  \sigma_{\mathrm{max}} \sqrt{\mu \log n} 
<\frac{1}{2}|\lambda_{l}^{\star}|\\
\Longrightarrow\qquad\frac{\lambda_{l}}{\lambda_{l}^{\star}} & \in\left[1-\frac{|\lambda_{l}-\lambda_{l}^{\star}|}{|\lambda_{l}^{\star}|},\,1+\frac{|\lambda_{l}-\lambda_{l}^{\star}|}{|\lambda_{l}^{\star}|}\right] \subset \left[ 0.5, 1.5 \right],
\end{align*}
provided that $\mu\kappa^2 r^4 \lesssim n$ and that Condition~\eqref{eq:condition-lambdamin-rankr-revised} holds. 
We then make use of the bounds \eqref{eq:rankr-eigenvalue-bound} and \eqref{eq:rankr-l2-bound-2} to deduce that
\begin{align*}
\min\left|1\pm\frac{\bm{u}_{l}^{\star\top}\bm{u}_{l}}{\lambda_{l}/\lambda_{l}^{\star}}\right|= & \left|1-\frac{\lambda_{l}^{\star}}{\lambda_{l}}+\frac{\lambda_{l}^{\star}}{\lambda_{l}}-\frac{|\bm{u}_{l}^{\star\top}\bm{u}_{l}|}{\lambda_{l}/\lambda_{l}^{\star}}\right|\leq\left|1-\frac{\lambda_{l}^{\star}}{\lambda_{l}}\right|+\left|\frac{\lambda_{l}^{\star}}{\lambda_{l}}\right|\cdot\left|1-\left|\bm{u}_{l}^{\star\top}\bm{u}_{l}\right|\right|\\
\lesssim & \frac{|\lambda_{l}-\lambda_{l}^{\star}|}{\lambda_{\min}^{\star}}+\left|1-\left|\bm{u}_{l}^{\star\top}\bm{u}_{l}\right|\right|\nonumber\\
\lesssim & \frac{\kappa r^{2}\sigma_{\max}\sqrt{\mu\log n}}{\lambda_{\min}^{\star}}+\left(\kappa^{4}\sigma_{\max}^{2}n\log n\right)\left\{ \frac{1}{(\Delta_{l}^{\star})^{2}}\frac{\mu r^{2}}{n}+\frac{1}{(\lambda_{\max}^{\star})^{2}}\right\} .
\end{align*}
The above bounds taken collectively demonstrate that
\begin{align}
\min\left|\bm{a}^{\top}\left(\bm{u}_{l}\pm\bm{u}_{l}^{\star}\right)\right|\leq & \left|\bm{a}^{\top}\left(\bm{u}_{l}-\mathsf{sign}\Big(\frac{\bm{u}_{l}^{\star\top}\bm{u}_{l}}{\lambda_{l}/\lambda_{l}^{\star}}\Big)\bm{u}_{l}^{\star}\right)\right| \nonumber\\
\leq & \left|\bm{a}^{\top}\left(\bm{u}_{l}-\frac{\bm{u}_{l}^{\star\top}\bm{u}_{l}}{\lambda_{l}/\lambda_{l}^{\star}}\bm{u}_{l}^{\star}\right)\right|+\left|\mathsf{sign}\Big(\frac{\bm{u}_{l}^{\star\top}\bm{u}_{l}}{\lambda_{l}/\lambda_{l}^{\star}}\Big)\bm{a}^{\top}\bm{u}_{l}^{\star}-\frac{\bm{u}_{l}^{\star\top}\bm{u}_{l}}{\lambda_{l}/\lambda_{l}^{\star}}\bm{a}^{\top}\bm{u}_{l}^{\star}\right| \nonumber\\
= & \left|\bm{a}^{\top}\left(\bm{u}_{l}-\frac{\bm{u}_{l}^{\star\top}\bm{u}_{l}}{\lambda_{l}/\lambda_{l}^{\star}}\bm{u}_{l}^{\star}\right)\right|+\min\left|1\pm\frac{\bm{u}_{l}^{\star\top}\bm{u}_{l}}{\lambda_{l}/\lambda_{l}^{\star}}\right|\cdot|\bm{a}^{\top}\bm{u}_{l}^{\star}|\nonumber\\
\lesssim & ~\mathcal{E}_{\mathsf{au}}, 
	\label{eq:aul-UB123}
\end{align}
where 
\begin{align*}
\mathcal{E}_{\mathsf{au}} & :=\,\sigma_{\mathrm{max}}\sqrt{\mu\kappa^{2}r\log n}\left\{ \frac{1}{\lambda_{\mathrm{min}}^{\star}}\|\bm{a}\|_{2}+\kappa r\max_{k:k\neq l}\frac{|\bm{a}^{\top}\bm{u}_{k}^{\star}|}{\big|\lambda_{l}^{\star}-\lambda_{k}^{\star}\big|}\right\} \\
 & \quad+\left\{ \frac{\kappa r^{2}\sigma_{\max}\sqrt{\mu\log n}}{\lambda_{\min}^{\star}}+\left(\kappa^{4}\sigma_{\max}^{2}n\log n\right)\left\{ \frac{1}{(\Delta_{l}^{\star})^{2}}\frac{\mu r^{2}}{n}+\frac{1}{(\lambda_{\max}^{\star})^{2}}\right\} \right\} |\bm{a}^{\top}\bm{u}_{l}^{\star}|.
\end{align*}

Finally, when it comes to the $\ell_{\infty}$ perturbation bound, we shall simply take $\bm{a} = \bm{e}_k$ $(1 \leq k \leq r)$ in the above inequality and use the incoherence condition $|\bm{e}_k^{\top}\bm{u}_l^{\star}|\leq \sqrt{\mu/n}$ ($1\leq k\leq r$) to obtain 
\begin{align*}
\left|\bm{e}_{k}^{\top}\left(\bm{u}_{l}-\mathsf{sign}\Big(\frac{\bm{u}_{l}^{\star\top}\bm{u}_{l}}{\lambda_{l}/\lambda_{l}^{\star}}\Big)\bm{u}_{l}^{\star}\right)\right| & \lesssim\,\sigma_{\mathrm{max}}\sqrt{\mu\kappa^{2}r\log n}\left\{ \frac{1}{\lambda_{\mathrm{min}}^{\star}}+\frac{\kappa r}{\Delta_{l}^{\star}}\sqrt{\frac{\mu}{n}}\right\} \\
 & \qquad+\left\{ \frac{\kappa r^{2}\sigma_{\max}\sqrt{\mu\log n}}{\lambda_{\min}^{\star}}+\left(\kappa^{4}\sigma_{\max}^{2}n\log n\right)\left\{ \frac{1}{(\Delta_{l}^{\star})^{2}}\frac{\mu r^{2}}{n}+\frac{1}{(\lambda_{\max}^{\star})^{2}}\right\} \right\} \sqrt{\frac{\mu}{n}}\\
 & \lesssim\,\frac{\sigma_{\mathrm{max}}\sqrt{\mu\kappa^{4}r\log n}}{\lambda_{\min}^{\star}}+\frac{\sigma_{\mathrm{max}}}{\Delta_{l}^{\star}}\sqrt{\frac{\mu^{2}\kappa^{4}r^{3}\log n}{n}},
\end{align*}
where the last line results from Condition \eqref{eq:condition-Delta-noise-revised} and $\mu \kappa^2 r^4 \lesssim n$.  This together with a union bound yields that: with high probability, one has
\[
\min\|\bm{u}_{l}\pm\bm{u}_{l}^{\star}\|_{\infty}\leq\Big\|\bm{u}_{l}-\mathsf{sign}\Big(\frac{\bm{u}_{l}^{\star\top}\bm{u}_{l}}{\lambda_{l}/\lambda_{l}^{\star}}\Big)\bm{u}_{l}^{\star}\Big\|_{\infty}\lesssim\,\frac{\sigma_{\mathrm{max}}\sqrt{\mu\kappa^{4}r\log n}}{\lambda_{\min}^{\star}}+\frac{\sigma_{\mathrm{max}}}{\Delta_{l}^{\star}}\sqrt{\frac{\mu^{2}\kappa^{4}r^{3}\log n}{n}}
\]
as claimed.

\paragraph{Proof of the claims w.r.t.~$\widehat{\bm{u}}_l:=\frac{\bm{u}_l+\bm{w}_l}{\|\bm{u}_l+\bm{w}_l\|_2}$.}  With the above results in place, we can easily establish these claims w.r.t.~$\widehat{\bm{u}}_l$ as well. In what follows, we demonstrate how to establish the bound \eqref{eq:rankr-linear-form-bound} regarding the linear form; the proofs for $\ell_2$ and $\ell_{\infty}$ bounds follow from nearly identical arguments and are omitted for brevity. 

To begin with, repeating the analysis \eqref{eq:aul-UB123} for $\bm{w}_l$ yields
\begin{align}
\left|\bm{a}^{\top}\left(\bm{u}_{l}-\mathsf{sign}\Big(\frac{\bm{u}_{l}^{\star\top}\bm{w}_{l}}{\lambda_{l}/\lambda_{l}^{\star}}\Big)\bm{u}_{l}^{\star}\right)\right| & \lesssim\mathcal{E}_{\mathsf{au}},
	\label{eq:aul-UB124}
\end{align}
One can therefore combine \eqref{eq:aul-UB123} and \eqref{eq:aul-UB124} to obtain
\begin{align}
\left|\bm{a}^{\top}\Big(\frac{\bm{u}_{l}+\bm{w}_{l}}{2}\Big)-\mathsf{sign}\Big(\frac{\bm{u}_{l}^{\star\top}\bm{w}_{l}}{\lambda_{l}/\lambda_{l}^{\star}}\Big)\bm{a}^{\top}\bm{u}_{l}^{\star}\right|
 & \leq\frac{1}{2}\left|\bm{a}^{\top}\left(\bm{u}_{l}-\mathsf{sign}\Big(\frac{\bm{u}_{l}^{\star\top}\bm{u}_{l}}{\lambda_{l}/\lambda_{l}^{\star}}\Big)\bm{u}_{l}^{\star}\right)\right|+\frac{1}{2}\left|\bm{a}^{\top}\left(\bm{w}_{l}-\mathsf{sign}\Big(\frac{\bm{u}_{l}^{\star\top}\bm{w}_{l}}{\lambda_{l}/\lambda_{l}^{\star}}\Big)\bm{u}_{l}^{\star}\right)\right| \nonumber\\
 & \lesssim\mathcal{E}_{\mathsf{au}}, 
	\label{eq:aul-UB126}
\end{align}
where the first inequality arises from the triangle inequality as well as the fact $\mathsf{sign}(\bm{u}_{l}^{\star\top}\bm{w}_{l})=\mathsf{sign}(\bm{u}_{l}^{\star\top}\bm{u}_{l})$ (see \eqref{eq:ur-ul-uhat-rankr} in Lemma~\ref{lem:basic-uR-uL-rankr}).  This in turn allows us to deduce that
\begin{align}
 & \left|\bm{a}^{\top}\Big(\frac{\bm{u}_{l}+\bm{w}_{l}}{\|\bm{u}_{l}+\bm{w}_{l}\|_{2}}\Big)-\mathsf{sign}\Big(\frac{\bm{u}_{l}^{\star\top}\bm{w}_{l}}{\lambda_{l}/\lambda_{l}^{\star}}\Big)\bm{a}^{\top}\bm{u}_{l}^{\star}\right|\nonumber\\
 & \quad=\left|\frac{2}{\|\bm{u}_{l}+\bm{w}_{l}\|_{2}}\bm{a}^{\top}\Big(\frac{\bm{u}_{l}+\bm{w}_{l}}{2}\Big)-\frac{2}{\|\bm{u}_{l}+\bm{w}_{l}\|_{2}}\mathsf{sign}\Big(\frac{\bm{u}_{l}^{\star\top}\bm{w}_{l}}{\lambda_{l}/\lambda_{l}^{\star}}\Big)\bm{a}^{\top}\bm{u}_{l}^{\star}+\mathsf{sign}\Big(\frac{\bm{u}_{l}^{\star\top}\bm{w}_{l}}{\lambda_{l}/\lambda_{l}^{\star}}\Big)\left(\frac{2}{\|\bm{u}_{l}+\bm{w}_{l}\|_{2}}-1\right)\bm{a}^{\top}\bm{u}_{l}^{\star}\right|\nonumber\\
 & \quad\leq\frac{2}{\|\bm{u}_{l}+\bm{w}_{l}\|_{2}}\left|\bm{a}^{\top}\Big(\frac{\bm{u}_{l}+\bm{w}_{l}}{2}\Big)-\mathsf{sign}\Big(\frac{\bm{u}_{l}^{\star\top}\bm{w}_{l}}{\lambda_{l}/\lambda_{l}^{\star}}\Big)\bm{a}^{\top}\bm{u}_{l}^{\star}\right|+\left(\frac{2}{\|\bm{u}_{l}+\bm{w}_{l}\|_{2}}-1\right)\left|\bm{a}^{\top}\bm{u}_{l}^{\star}\right|\nonumber\\
 & \quad\lesssim\mathcal{E}_{\mathsf{au}}+(2-\|\bm{u}_{l}+\bm{w}_{l}\|_{2})\left|\bm{a}^{\top}\bm{u}_{l}^{\star}\right|\lesssim\mathcal{E}_{\mathsf{au}},\label{eq:au-sign-UB-125} 
\end{align}
where the first inequality follows from the triangle inequality,  the penultimate inequality relies on \eqref{eq:aul-UB126} and the fact $\|\bm{u}_l+\bm{w}_l\|_2\asymp 1$ (see Lemma~\ref{lem:basic-uR-uL-rankr}), and the last inequality follows from the bound $\big|2-\|\bm{u}_{l}+\bm{w}_{l}\|_{2}\big|=O\left(\frac{\kappa^{4}\sigma_{\max}^{2}n\log n}{(\lambda_{\mathrm{max}}^{\star})^{2}}+\frac{\mu\kappa^{4}r^{2}\sigma_{\max}^{2}\log n}{(\Delta_{l}^{\star})^{2}}\right)$ (cf.~Lemma~\ref{lem:basic-uR-uL-rankr}) in addition to the expression of $\mathcal{E}_{\mathsf{au}}$. This establishes the bound \eqref{eq:rankr-linear-form-bound} when $\bm{u}_l$ is replaced by $\widehat{\bm{u}}_l$.

\subsubsection{Proof of Theorem~\ref{thm:linear-form-ua-coarse}}
\label{sec:thm-linear-form-ua-coarse}


We start by considering the coarse estimator $\widehat{u}_{\bm{a},l}$. Recall that $\bm{a}$ is a fixed unit vector obeying $\|\bm{a}\|_2=1$. The following trivial bound holds true:
\begin{equation}
	\min\left|\widehat{u}_{\bm{a},l} \pm\bm{a}^{\top}\bm{u}_{l}^{\star}\right| \leq 1 .
\end{equation}
Consequently, it suffices to focus on the scenario when the upper bound on the right-hand side of \eqref{eq:au-coarse-estimate-UB} is bounded above by some small constant; that is, the case where
\begin{subequations}
	\begin{align}
	\tfrac{\sigma_{\mathrm{max}}r^2\sqrt{\mu\kappa^{4}\log n}}{\lambda_{\mathrm{min}}^{\star}} & \leq c_7   \label{eq:condition-proof-rank-r-linear-form-bounds-1} \\
	\tfrac{\sigma_{\max}^{2}\mu r^{2}\kappa^{4}\log n}{(\Delta_{l}^{\star})^{2}}|\bm{a}^{\top}\bm{u}_{l}^{\star}| & \leq c_8  \label{eq:condition-proof-rank-r-linear-form-bounds-2}\\
	\sigma_{\mathrm{max}}\sqrt{\mu\kappa^{4}r^{3}\log n}\max_{k\neq l}\tfrac{|\bm{a}^{\top}\bm{u}_{k}^{\star}|}{|\lambda_{l}^{\star}-\lambda_{k}^{\star}|} & \leq c_9 
	\label{eq:condition-proof-rank-r-linear-form-bounds-3}
	\end{align}
\end{subequations}
for some sufficiently small constants $c_7,c_8,c_9>0$.

The key step is to invoke the Neumann series (cf.~Lemma~\ref{lem:Neumann-expansion}) for both $\bm{u}_l$ and $\bm{w}_l$, which yields
\begin{align}
\big(\widehat{u}_{\bm{a},l}\big)^{2} & =\Bigg|\frac{\big(\bm{a}^{\top}\bm{u}_{l}\big)\big(\bm{a}^{\top}\bm{w}_{l}\big)}{\left(\sum_{k=1}^{r}\frac{\bm{u}_{k}^{\star\top}\bm{w}_{l}}{\lambda_{l}/\lambda_{k}^{\star}}\sum_{s=0}^{\infty}\frac{(\bm{H}^{\top})^{s}\bm{u}_{k}^{\star}}{\lambda_{l}^{s}}\right)^{\top}\left(\sum_{k=1}^{r}\frac{\bm{u}_{k}^{\star\top}\bm{u}_{l}}{\lambda_{l}/\lambda_{k}^{\star}}\sum_{s=0}^{\infty}\frac{\bm{H}^{s}\bm{u}_{k}^{\star}}{\lambda_{l}^{s}} \right)} \Bigg|\nonumber\\
 & =\Bigg|\frac{\left(\bm{a}^{\top}\bm{u}_{l}^{\star}\cdot\frac{\bm{u}_{l}^{\star\top}\bm{u}_{l}}{\lambda_{l}/\lambda_{l}^{\star}}+\delta_{1}\right)\left(\bm{a}^{\top}\bm{u}_{l}^{\star}\cdot\frac{\bm{u}_{l}^{\star\top}\bm{w}_{l}}{\lambda_{l}/\lambda_{l}^{\star}}+\delta_{2}\right)}{\frac{\bm{u}_{l}^{\star\top}\bm{w}_{l}}{\lambda_{l}/\lambda_{l}^{\star}}\cdot\frac{\bm{u}_{l}^{\star\top}\bm{u}_{l}}{\lambda_{l}/\lambda_{l}^{\star}}+\delta_{3}}\Bigg|.  
	\label{eq:estimator-squared-decompose}
\end{align}
Here, $\delta_1,\delta_2,\delta_3$ are defined as follows
\begin{subequations}
	\begin{align}
	\delta_1 & := \bm{a}^{\top}\left(\bm{u}_{l}-\tfrac{\bm{u}_{l}^{\star\top}\bm{u}_{l}}{\lambda_{l}/\lambda_{l}^{\star}}\bm{u}_{l}^{\star}\right),   \label{eq:estimator-squared-decompose-1} \\
	\delta_2 & :=\bm{a}^{\top}\left(\bm{w}_{l}-\tfrac{\bm{u}_{l}^{\star\top}\bm{w}_{l}}{\lambda_{l}/\lambda_{l}^{\star}}\bm{u}_{l}^{\star}\right),  \label{eq:estimator-squared-decompose-2}\\
	\delta_3 & := \mathcal{U}_1 + \mathcal{U}_2, 
	\label{eq:estimator-squared-decompose-3}
	\end{align}
\end{subequations}
where the quantities $\mathcal{U}_1$ and $\mathcal{U}_2$ are given by
\begin{subequations}
	\begin{align}
		\mathcal{U}_1 & := \sum_{k: k \neq l} \frac{\bm{u}_{k}^{\star \top} \bm{w}_l}{\lambda_l / \lambda_{k}^\star} \cdot \frac{\bm{u}_{k}^{\star \top} \bm{u}_l}{\lambda_l / \lambda_{k}^\star}, \label{eq:estimator-squared-decompose-3-1} \\
		\mathcal{U}_{2} &:= \sum_{k_{1}=1}^{r}\sum_{k_{2}=1}^{r}\sum_{s_{1},s_{2}:s_{1}+s_{2}\geq 1}\frac{\bm{u}_{k_{1}}^{\star\top}\bm{w}_{l}}{\lambda_{l}/\lambda_{k_{1}}^{\star}}\cdot\frac{\bm{u}_{k_{2}}^{\star\top}\bm{u}_{l}}{\lambda_{l}/\lambda_{k_{2}}^{\star}}\cdot\frac{\bm{u}_{k_{1}}^{\star\top}\bm{H}^{s_{1}+s_{2}}\bm{u}_{k_{2}}^{\star}}{\lambda_{l}^{s_{1}+s_{2}}}.
		\label{eq:estimator-squared-decompose-3-2}
	\end{align}
\end{subequations}
This motivates us to control each term on the right-hand side of \eqref{eq:estimator-squared-decompose} separately. 
\begin{itemize}
	\item To begin with, it comes directly from the inequalities \eqref{eq:condition-lambdamin-rankr-revised} and \eqref{eq:condition-proof-rank-r-linear-form-bounds-3} that
\begin{align*}
\frac{\sigma_{\max}\sqrt{n\kappa^{6}\log n}}{\lambda_{\min}^{\star}} + \frac{\sigma_{\max}}{\Delta_{l}^{\star}}\sqrt{\mu\kappa^{4}r^{2}\log n} \leq c_5 + \frac{c_9}{\sqrt r},
\end{align*}
thus allowing us to invoke  Theorem \ref{thm:rankr-bounds} to obtain, with probability at least $1-O(n^{-6})$, that
\begin{equation}
\min \left\{|\bm{u}_l^{\star \top} \bm{w}_l|, |\bm{u}_l^{\star \top} \bm{u}_l| \right\} \geq 3/4. \label{eq:rank-r-l2-lower-bound}
\end{equation}

\item
In view of \eqref{eq:linear-form-decompose-1-bound}, with probability $1-O(n^{-6})$ one has
\begin{equation}
	\max \left\{|\delta_1|, |\delta_2| \right\} \lesssim \frac{\sigma_{\mathrm{max}}\sqrt{\mu\kappa^{2}r\log n}}{\lambda_{\mathrm{min}}^{\star}} + \sigma_{\mathrm{max}}\sqrt{\mu\kappa^{4}r^{3}\log n}\max_{k\neq l}\frac{\big|\bm{a}^{\top}\bm{u}_{k}^{\star}\big|}{\big|\lambda_{l}^{\star}-\lambda_{k}^{\star}\big|}. \label{eq:rank-r-delta12-bound}
\end{equation}

\item
Next, substituting \eqref{eq:rank-r-inner-product-bound} into \eqref{eq:estimator-squared-decompose-3-1} yields
\begin{align}
|\mathcal{U}_1| & \leq r  \max_{k \neq l} \Bigg| \frac{\bm{u}_{k}^{\star \top} \bm{w}_l}{\lambda_l / \lambda_{k}^\star} \cdot \frac{\bm{u}_{k}^{\star \top} \bm{u}_l}{\lambda_l / \lambda_{k}^\star}  \Bigg|
\lesssim r  \left(\frac{ \sigma_{\max}}{\Delta_l^\star}\sqrt{\mu\kappa^{4}r\log n}\right)^2 
	= \frac{\sigma_{\mathrm{max}}^2 \mu r^2\kappa^{4} \log n}{\left(\Delta_l^{\star}\right)^2} .
\end{align}
In addition, one has
\begin{align*}
|\mathcal{U}_{2}| & \overset{(\mathrm{i})}{\lesssim}\sum_{k_{1}=1}^{r}\sum_{k_{2}=1}^{r}\sum_{s_{1},s_{2}:s_{1}+s_{2}\geq1}\frac{1}{|\lambda_{l}^{\star}/\lambda_{k_{1}}^{\star}\cdot\lambda_{l}^{\star}/\lambda_{k_{2}}^{\star}|}\cdot\Bigg|\frac{\bm{u}_{k_{1}}^{\star\top}\bm{H}^{s_{1}+s_{2}}\bm{u}_{k_{2}}^{\star}}{(\lambda_{l}^{\star}/2)^{s_{1}+s_{2}}}\Bigg|\\
 & \overset{(\mathrm{ii})}{\leq}\kappa^{2}r^{2}\sum_{s_{1},s_{2}:s_{1}+s_{2}\geq1}\Big(\frac{2}{\lambda_{\min}^{\star}}\Big)^{s_{1}+s_{2}}\big|\bm{u}_{k_{1}}^{\star\top}\bm{H}^{s_{1}+s_{2}}\bm{u}_{k_{2}}^{\star}\big|\\
 & \overset{(\mathrm{iii})}{\leq}\kappa^{2}r^{2}\sum_{s=1}^{\infty}\sum_{s_{1},s_{2}:s_{1}+s_{2}=s}\left(\frac{2c_{2}\sigma_{\max}\sqrt{n\log n}}{\lambda_{\min}^{\star}}\right)^{s}\sqrt{\frac{\mu}{n}}\\
 & \leq\sqrt{\frac{\mu}{n}}\kappa^{2}r^{2}\sum_{s=1}^{\infty}2^{s}\left(\frac{2c_{2}\sigma_{\max}\sqrt{n\log n}}{\lambda_{\min}^{\star}}\right)^{s}\\
 & =\sqrt{\frac{\mu}{n}}\kappa^{2}r^{2}\frac{\frac{4c_{2}\sigma_{\max}\sqrt{n\log n}}{\lambda_{\min}^{\star}}}{1-\frac{4c_{2}\sigma_{\max}\sqrt{n\log n}}{\lambda_{\min}^{\star}}}\\
 & \overset{(\mathsf{iv})}{\asymp}\frac{\sigma_{\mathrm{max}}\sqrt{\mu\kappa^{4}r^{4}\log n}}{\lambda_{\mathrm{min}}^{\star}},
\end{align*}
where (i) follows from Corollary~\ref{cor:eigengap-condition-bounds}
(so that $|\lambda_{l}|\geq|\lambda_{l}^{\star}|/2$), (ii) makes
use of the fact $|\lambda_{k}^{\star}/\lambda_{l}^{\star}|\leq\kappa$,
(iii) arises from Lemma \ref{lemma:perturbation-H}, and the last
line relies on the assumption \eqref{eq:condition-lambdamin-rankr-revised}
(so that $1-\frac{4c_{2}\sigma_{\max}\sqrt{n\log n}}{\lambda_{\min}^{\star}}\ge\frac{1}{2}$). 
The above two bounds taken together yield
\begin{equation}
|\delta_3| \leq |\mathcal{U}_{1}| + |\mathcal{U}_{2}| \lesssim \frac{\sigma_{\mathrm{max}}^2 \mu r^2\kappa^{4} \log n}{\left(\Delta_l^{\star}\right)^2} + \frac{\sigma_{\mathrm{max}} r^2 \sqrt{\mu \kappa^4 \log n}}{\lambda_{\mathrm{min}}^\star}
	\leq 1/8, \label{eq:rank-r-delta3-bound}
\end{equation}
where the last inequality holds as long as the constants $c_7$ and $c_8$ in \eqref{eq:condition-proof-rank-r-linear-form-bounds-2} and \eqref{eq:condition-proof-rank-r-linear-form-bounds-3} are sufficiently small.

\end{itemize}

With the above bounds in mind, we are ready to control the estimation error of $\widehat{u}_{\bm{a},l}$. We divide the proof into two separate cases as follows.

\paragraph{Case 1: when \texorpdfstring{$|\bm{a}^{\top}\bm{u}_l^{\star}|$}{the truth}  is ``small''.} Consider the case where
$$
|\bm{a}^\top \bm{u}_l^\star| \leq  \max \left\{|\delta_1|, |\delta_2| \right\}.
$$
Continue the bound in \eqref{eq:estimator-squared-decompose} to deduce that
\begin{align}
\big(\widehat{u}_{\bm{a},l}\big)^{2} & =\left|\frac{\left(\bm{a}^{\top}\bm{u}_{l}^{\star}\cdot\frac{\bm{u}_{l}^{\star\top}\bm{u}_{l}}{\lambda_{l}/\lambda_{l}^{\star}}+\delta_{1}\right)\left(\bm{a}^{\top}\bm{u}_{l}^{\star}\cdot\frac{\bm{u}_{l}^{\star\top}\bm{w}_{l}}{\lambda_{l}/\lambda_{l}^{\star}}+\delta_{2}\right)}{\frac{\bm{u}_{l}^{\star\top}\bm{w}_{l}}{\lambda_{l}/\lambda_{l}^{\star}}\cdot\frac{\bm{u}_{l}^{\star\top}\bm{u}_{l}}{\lambda_{l}/\lambda_{l}^{\star}}+\delta_{3}}\right|\nonumber\\
 & \stackrel{\mathrm{(i)}}{\leq}\frac{\left\{ 2|\bm{a}^{\top}\bm{u}_{l}^{\star}|+|\delta_{1}|\right\} \left\{ 2|\bm{a}^{\top}\bm{u}_{l}^{\star}|+|\delta_{2}|\right\} }{\left(\frac{1}{2}\right)^{2}-\frac{1}{8}}\nonumber\\
 & \leq72\left(\max\left\{ |\delta_{1}|,|\delta_{2}|\right\} \right)^{2}.
\label{eq:rank-r-u-hat-squared-bound-small}
\end{align}
Here, the inequality (i) makes use of the bound \eqref{eq:rank-r-l2-lower-bound}, Corollary~\ref{cor:eigengap-condition-bounds} (so that $1/2\leq |\lambda_l/\lambda_l^{\star}|\leq 2$), and the inequality $|\delta_3| \leq 1/8$. Therefore, we can take the triangle inequality to conclude that
\begin{align}
\min \left|\widehat{u}_{\bm{a},l} \pm \bm{a}^\top \bm{u}^\star_l \right|\leq \left|\widehat{u}_{\bm{a},l} \right| + \left|\bm{a}^\top \bm{u}^\star_l \right| \lesssim \max \left\{|\delta_1|, |\delta_2| \right\}, \label{eq:rank-r-u-hat-estimator-bound-small}
\end{align}
which combined with \eqref{eq:rank-r-delta12-bound} leads to the desired bound for this case.

\paragraph{Case 2: when \texorpdfstring{$|\bm{a}^{\top}\bm{u}_l^{\star}|$}{the truth}  is ``large''.} We now move on to the case where
\begin{align}
\label{eq:au-large-delta}
|\bm{a}^\top \bm{u}_l^\star| >  \max \left\{|\delta_1|, |\delta_2| \right\}.
\end{align}
By the same argument above in \eqref{eq:rank-r-u-hat-squared-bound-small}, we have
\begin{align}
\left|\frac{(\bm{a}^{\top}\bm{u}_{l})(\bm{a}^{\top}\bm{w}_{l})}{\bm{u}_{l}^{\top}\bm{w}_{l}} -(\bm{a}^{\top}\bm{u}_{l}^{\star})^{2}\right| & =\left|\frac{\left(\bm{a}^{\top}\bm{u}_{l}^{\star}\cdot\frac{\bm{u}_{l}^{\star\top}\bm{u}_{l}}{\lambda_{l}/\lambda_{l}^{\star}}+\delta_{1}\right)\left(\bm{a}^{\top}\bm{u}_{l}^{\star}\cdot\frac{\bm{u}_{l}^{\star\top}\bm{w}_{l}}{\lambda_{l}/\lambda_{l}^{\star}}+\delta_{2}\right)-(\bm{a}^{\top}\bm{u}_{l}^{\star})^{2}\left\{ \frac{\bm{u}_{l}^{\star\top}\bm{w}_{l}}{\lambda_{l}/\lambda_{l}^{\star}}\cdot\frac{\bm{u}_{l}^{\star\top}\bm{u}_{l}}{\lambda_{l}/\lambda_{l}^{\star}}+\delta_{3}\right\} }{\frac{\bm{u}_{l}^{\star\top}\bm{w}_{l}}{\lambda_{l}/\lambda_{l}^{\star}}\cdot\frac{\bm{u}_{l}^{\star\top}\bm{u}_{l}}{\lambda_{l}/\lambda_{l}^{\star}}+\delta_{3}}\right|\nonumber\\
 & =\left|\frac{\bm{a}^{\top}\bm{u}_{l}^{\star}\cdot\frac{\bm{u}_{l}^{\star\top}\bm{u}_{l}}{\lambda_{l}/\lambda_{l}^{\star}}\cdot\delta_{2}+\bm{a}^{\top}\bm{u}_{l}^{\star}\cdot\frac{\bm{u}_{l}^{\star\top}\bm{w}_{l}}{\lambda_{l}/\lambda_{l}^{\star}}\cdot\delta_{1}+\delta_{1}\delta_{2}-(\bm{a}^{\top}\bm{u}_{l}^{\star})^{2}\cdot\delta_{3}}{\frac{\bm{u}_{l}^{\star\top}\bm{w}_{l}}{\lambda_{l}/\lambda_{l}^{\star}}\cdot\frac{\bm{u}_{l}^{\star\top}\bm{u}_{l}}{\lambda_{l}/\lambda_{l}^{\star}}+\delta_{3}}\right|\nonumber\\
 & \lesssim\frac{\left|\bm{a}^{\top}\bm{u}_{l}^{\star}\right|\max\left\{ |\delta_{1}|,|\delta_{2}|\right\} +\left(\max\left\{ |\delta_{1}|,|\delta_{2}|\right\} \right)^{2}+\left(\bm{a}^{\top}\bm{u}_{l}^{\star}\right)^{2}|\delta_{3}|}{\left(\frac{1}{2}\right)^{2}-\frac{1}{8}}\nonumber\\
 & \lesssim\left|\bm{a}^{\top}\bm{u}_{l}^{\star}\right|\left(\frac{\sigma_{\mathrm{max}}\sqrt{\mu\kappa^{2}r\log n}}{\lambda_{\mathrm{min}}^{\star}}+\sigma_{\mathrm{max}}\sqrt{\mu\kappa^{4}r^{3}\log n}\max_{k\neq l}\frac{\big|\bm{a}^{\top}\bm{u}_{k}^{\star}\big|}{\big|\lambda_{l}^{\star}-\lambda_{k}^{\star}\big|}\right)\nonumber\\
 & \quad+\left(\bm{a}^{\top}\bm{u}_{l}^{\star}\right)^{2}\cdot\left(\frac{\sigma_{\mathrm{max}}^{2}\mu r^{2}\kappa^{4}\log n}{\left(\Delta_{l}^{\star}\right)^{2}}+\frac{\sigma_{\mathrm{max}}r^{2}\sqrt{\mu\kappa^{4}\log n}}{\lambda_{\mathrm{min}}^{\star}}\right),\label{eq:rank-r-u-hat-squared-bound-large}
\end{align}
where the last line relies on our upper bounds on $|\delta_1|,|\delta_2|$ and $|\delta_3|$ (see \eqref{eq:rank-r-delta12-bound} and \eqref{eq:rank-r-delta3-bound}) as well as the assumption \eqref{eq:au-large-delta}.

Recognizing the trivial fact (here, define $\max|a\pm b|=\max\{|a+b|,|a-b|\}$)
$$
\max |\widehat{u}_{\bm{a},l} \pm \bm{a}^\top \bm{u}_l^\star| \geq \left|\bm{a}^\top \bm{u}_l^\star\right|,
$$
we have
\begin{align}
\min\left|\widehat{u}_{\bm{a},l}\pm\bm{a}^{\top}\bm{u}_{l}^{\star}\right| & =\frac{\left|(\widehat{u}_{\bm{a},l})^{2}-(\bm{a}^{\top}\bm{u}_{l}^{\star})^{2}\right|}{\max\big|\widehat{u}_{\bm{a},l}\pm\bm{a}^{\top}\bm{u}_{l}^{\star}\big|}\leq\frac{\left|(\widehat{u}_{\bm{a},l})^{2}-(\bm{a}^{\top}\bm{u}_{l}^{\star})^{2}\right|}{\left|\bm{a}^{\top}\bm{u}_{l}^{\star}\right|} . \label{eq:ual-coarse-pm-ub1}
\end{align}
If $\frac{(\bm{a}^{\top}\bm{u}_{l})(\bm{a}^{\top}\bm{w}_{l})}{\bm{u}_{l}^{\top}\bm{w}_{l}}\geq0$,
then 
\[
\left|(\widehat{u}_{\bm{a},l})^{2}-(\bm{a}^{\top}\bm{u}_{l}^{\star})^{2}\right|=\left|\frac{(\bm{a}^{\top}\bm{u}_{l})(\bm{a}^{\top}\bm{w}_{l})}{\bm{u}_{l}^{\top}\bm{w}_{l}}-(\bm{a}^{\top}\bm{u}_{l}^{\star})^{2}\right|.
\]
If instead $\frac{(\bm{a}^{\top}\bm{u}_{l})(\bm{a}^{\top}\bm{w}_{l})}{\bm{u}_{l}^{\top}\bm{w}_{l}}<0$ holds,
then 
\[
\left|(\widehat{u}_{\bm{a},l})^{2}-(\bm{a}^{\top}\bm{u}_{l}^{\star})^{2}\right|\leq\left|(\widehat{u}_{\bm{a},l})^{2}+(\bm{a}^{\top}\bm{u}_{l}^{\star})^{2}\right|=\left|\frac{(\bm{a}^{\top}\bm{u}_{l})(\bm{a}^{\top}\bm{w}_{l})}{\bm{u}_{l}^{\top}\bm{w}_{l}}-(\bm{a}^{\top}\bm{u}_{l}^{\star})^{2}\right|.
\]
Therefore, putting these bounds together and invoking our bound \eqref{eq:rank-r-u-hat-squared-bound-large}, we deduce that
\begin{align}
 & \min\left|\widehat{u}_{\bm{a},l}\pm\bm{a}^{\top}\bm{u}_{l}^{\star}\right|
\leq\frac{\left|\frac{(\bm{a}^{\top}\bm{u}_{l})(\bm{a}^{\top}\bm{w}_{l})}{\bm{u}_{l}^{\top}\bm{w}_{l}}-(\bm{a}^{\top}\bm{u}_{l}^{\star})^{2}\right|}{\left|\bm{a}^{\top}\bm{u}_{l}^{\star}\right|}\nonumber\\
 & \lesssim\frac{\sigma_{\mathrm{max}}\sqrt{\mu\kappa^{2}r\log n}}{\lambda_{\mathrm{min}}^{\star}}+\sigma_{\mathrm{max}}\sqrt{\mu\kappa^{4}r^{3}\log n}\max_{k\neq l}\frac{\big|\bm{a}^{\top}\bm{u}_{k}^{\star}\big|}{\big|\lambda_{l}^{\star}-\lambda_{k}^{\star}\big|}+\left|\bm{a}^{\top}\bm{u}_{l}^{\star}\right|\cdot\left(\frac{\sigma_{\mathrm{max}}^{2}\mu r^{2}\kappa^{4}\log n}{\left(\Delta_{l}^{\star}\right)^{2}}+\frac{\sigma_{\mathrm{max}}r^{2}\sqrt{\mu\kappa^{4}\log n}}{\lambda_{\mathrm{min}}^{\star}}\right)\nonumber\\
 & \lesssim\frac{\sigma_{\mathrm{max}}r^{2}\sqrt{\mu\kappa^{4}\log n}}{\lambda_{\mathrm{min}}^{\star}}+\sigma_{\mathrm{max}}\sqrt{\mu\kappa^{4}r^{3}\log n}\max_{k\neq l}\frac{\big|\bm{a}^{\top}\bm{u}_{k}^{\star}\big|}{\big|\lambda_{l}^{\star}-\lambda_{k}^{\star}\big|}+\left|\bm{a}^{\top}\bm{u}_{l}^{\star}\right|\cdot\frac{\sigma_{\mathrm{max}}^{2}\mu r^{2}\kappa^{4}\log n}{\left(\Delta_{l}^{\star}\right)^{2}}.\label{eq:rank-r-u-hat-estimator-bound-large}
\end{align}
Taking collectively \eqref{eq:rank-r-u-hat-estimator-bound-small} and \eqref{eq:rank-r-u-hat-estimator-bound-large} establishes our estimation error bound for this case. 

\section{Proof for distributional guarantees (Theorem~\ref{thm:distribution-validity-rankr})}
\label{sec:proof-thm-distribution-validity-rankr}

\subsection{Distributional theory for linear forms of eigenvectors}

By virtue of Theorem \ref{thm:eigengap-condition}, we know that $\lambda_l$, $\bm{u}_l$ and $\bm{w}_l$ are all real-valued. 
Without loss of generality, assume that $\|\bm{a}\|_2 = 1$ and that $\bm{u}_l^{\star \top} \bm{u}_l>0$, which combined with Lemma \ref{lem:basic-uR-uL-rankr} and Notation~\ref{notation:leading-evectors} yield
\begin{equation}
\bm{u}_l^\top \bm{u}_l^\star = 1 - o(1), \quad \bm{w}_l^\top \bm{u}_l^\star = 1- o(1), \quad \bm{w}_l^\top \bm{u}_l = 1 - o(1).
\end{equation}
%


The main step lies in establishing the following claim 
\begin{align}
	\widehat{u}_{\bm{a},l}^{\mathsf{modified}} & = \bm{a}^{\top}\bm{u}^{\star}_l + \frac{(\bm{a}_{l}^{\perp})^{\top}\left(\bm{H}+\bm{H}^{\top}\right)\bm{u}_l^{\star}}{2\lambda_l^{\star}}
	+ o\left(\frac{\sigma_{\max}}{|\lambdastar_l|}\right) \label{eq:amazing-fact}\\
	& = \bm{a}^{\top}\bm{u}_l^{\star} + \frac{(\bm{a}_{l}^{\perp})^{\top}\left(\bm{H}+\bm{H}^{\top}\right)\bm{u}_l^{\star}}{2\lambda_l^{\star}}
	+ o\big(\sqrt{v^{\star}_{\bm{a},l}}\big)  \label{eq:amazing-fact2}
\end{align}
with $\bm{a}_l^{\perp}:=\bm{a}-(\bm{a}^{\top}\bm{u}_l^{\star}) \bm{u}_l^{\star}$, 
where the last line results from Lemma~\ref{lem:va-properties}. 
Let us take Claim \eqref{eq:amazing-fact} as given for now and come back to its proof in the sequel. To establish Theorem~\ref{thm:distribution-validity-rankr}, it suffices to pin down the distribution of $\frac{(\bm{a}_{l}^{\perp})^{\top}\left(\bm{H}+\bm{H}^{\top}\right)\bm{u}_l^{\star}}{2\lambda_l^{\star}}$ and show that it matches the distributional characterizations stated in Theorem~\ref{thm:distribution-validity-rankr}. This follows immediately from the classical Berry-Esseen theorem, which we defer to Lemma~\ref{LemPropertyInprod}.

The rest of the proof thus boils down to justifying the claim~\eqref{eq:amazing-fact}, for which we divide into two cases.

\subsubsection{The case when \texorpdfstring{$|\bm{a}^{\top}\bm{u}_l^{\star}|$}{the truth}  is not ``small''}

In this subsection, we focus on the scenario when $\sqrt{v^{\star}_{\bm{a},l}} \log n = o(|\bm{a}^\top \bm{u}_l^\star|) $ which, according to Lemma~\ref{Lem-var-control-rankr}, subsumes the case $  \sqrt{\widehat{v}_{\bm{a},l}} \log n = o(|\bm{a}^\top \bm{u}_l^\star|)$. Without loss of generality, we assume that $\bm{a}^\top \bm{u}_l^\star > 0$. By virtue of Lemma~\ref{lem:va-properties} and the assumption $\sigma_{\max}\asymp \sigma_{\min}$, we have
\[
v_{\bm{a},l}^{\star}\asymp\frac{\sigma_{\mathrm{max}}^{2}(1-|\bm{a}^{\top}\bm{u}_{l}^{\star}|^{2})}{|\lambda_{l}^{\star}|^{2}}\asymp\frac{\sigma_{\mathrm{max}}^{2}}{|\lambda_{l}^{\star}|^{2}},
\]
where the last line follows from our assumption that $|\bm{a}^{\top}\bm{u}_l^{\star}|$ is bounded away from 1. As a result, the regime considered in this subsection enjoys the following property:
\begin{align}
	\frac{\sigma_{\mathrm{max}}\log n}{|\lambda_{l}^{\star}|} = o(\bm{a}^{\top}\bm{u}_{l}^{\star}) .
	\label{eq:sigma-large-case-condition}
\end{align}


The key starting point is the following decomposition
\begin{align}
\big(\widehat{u}_{\bm{a},l}^{\mathsf{modified}}\big)^{2} & =\Big|\frac{(\bm{a}^{\top}\bm{u}_{l})(\bm{a}^{\top}\bm{w}_{l})}{\bm{u}_{l}^{\top}\bm{w}_{l}}\Big|\\
 & =\Bigg|\frac{\left\{ \frac{\bm{u}_{l}^{\star\top}\bm{u}_{l}}{\lambda_{l}/\lambda_{l}^{\star}}\cdot\left(\bm{a}^{\top}\bm{u}_{l}^{\star}+\frac{\bm{a}^{\top}\bm{H}\bm{u}_{l}^{\star}}{\lambda_{l}^{\star}}\right)+\tau_{1}\right\} \left\{ \frac{\bm{u}_{l}^{\star\top}\bm{w}_{l}}{\lambda_{l}/\lambda_{l}^{\star}}\cdot\left(\bm{a}^{\top}\bm{u}_{l}^{\star}+\frac{\bm{a}^{\top}\bm{H}^{\top}\bm{u}_{l}^{\star}}{\lambda_{l}^{\star}}\right)+\tau_{2}\right\} }{\frac{\bm{u}_{l}^{\star\top}\bm{w}_{l}}{\lambda_{l}/\lambda_{l}^{\star}}\cdot\frac{\bm{u}_{l}^{\star\top}\bm{u}_{l}}{\lambda_{l}/\lambda_{l}^{\star}}\left(1+\frac{\bm{u}_{l}^{\star\top}(\bm{H}+\bm{H}^{\top})\bm{u}_{l}^{\star}}{\lambda_{l}^{\star}}\right)+\tau_{3}}\Bigg|,\label{eq:estimator-squared-decompose-first-order-with-abs}
\end{align}
where 
\begin{equation}
\begin{cases}
\tau_{1} & :=\bm{a}^{\top}\bm{u}_{l}-\frac{\bm{u}_{l}^{\star\top}\bm{u}_{l}}{\lambda_{l}/\lambda_{l}^{\star}}\cdot\left(\bm{a}^{\top}\bm{u}_{l}^{\star}+\frac{\bm{a}^{\top}\bm{H}\bm{u}_{l}^{\star}}{\lambda_{l}^{\star}}\right),\\
\tau_{2} & :=\bm{a}^{\top}\bm{w}_{l}-\frac{\bm{u}_{l}^{\star\top}\bm{w}_{l}}{\lambda_{l}/\lambda_{l}^{\star}}\cdot\left(\bm{a}^{\top}\bm{u}_{l}^{\star}+\frac{\bm{a}^{\top}\bm{H}\bm{u}_{l}^{\star}}{\lambda_{l}^{\star}}\right),\\
\tau_{3} & :=\bm{u}_{l}^{\top}\bm{w}_{l}-\frac{\bm{u}_{l}^{\star\top}\bm{w}_{l}}{\lambda_{l}/\lambda_{l}^{\star}}\cdot\frac{\bm{u}_{l}^{\star\top}\bm{u}_{l}}{\lambda_{l}/\lambda_{l}^{\star}}\left(1+\frac{\bm{u}_{l}^{\star\top}(\bm{H}+\bm{H}^{\top})\bm{u}_{l}^{\star}}{\lambda_{l}^{\star}}\right).
\end{cases}\label{eq:defn-tau1-tau2-tau3}
\end{equation}
Here, $\tau_1,\tau_2,\tau_3$ encompass second- or higher-order terms in the Neumann series. As it turns out, these terms can be well-controlled, as stated in the following lemma. 
\begin{lemma} \label{lem:rank-r-inference-remainder-size}
	Instate the assumptions of Theorem \ref{thm:confidence-interval-validity-rankr}. With probability at least $1-O(n^{-6})$, we have
	\begin{align}
	|\tau_1|  = o\left(\frac{\sigma_{\mathrm{max}}}{|\lambda_{l}^{\star}|}\right), \quad 
	  |\tau_2|  = o\left(\frac{\sigma_{\mathrm{max}}}{|\lambda_{l}^{\star}|}\right), \quad 
	  |\tau_3| \lesssim \left(\frac{\sigma_{\mathrm{max}} \kappa^2 r \sqrt{\mu \log n}}{|\Delta_l^\star|}\right)^2 + o\left(\frac{\sigma_{\mathrm{max}}}{|\lambda_{l}^{\star}|}\right). 
	\end{align}
	Here, we recall that $\Delta_l^\star=\infty$ when $r = 1$, meaning that $|{\tau}_3| = o(\sigma_{\mathrm{max}}/|\lambda_{l}^{\star}|)$ when $r = 1$.
\end{lemma}

Recall from Theorem \ref{thm:eigengap-condition} and Lemma~\ref{lem:basic-uR-uL-rankr} that  $|\bm{u}_l^{\star \top} \bm{u}_l| = 1-o(1), |\bm{u}_l^{\star \top} \bm{w}_l| = 1-o(1)$ and $\lambda_l = (1+o(1)) \lambda_l^\star$. Therefore, Lemma \ref{lem:rank-r-inference-remainder-size} tells us that
\begin{align*}
|\tau_{1}|=o\left(\frac{\sigma_{\mathrm{max}}}{|\lambda_{l}^{\star}|}\right) & \leq o\left(|\bm{a}^{\top}\bm{u}_{l}^{\star}|\right)=o\left(\left|\frac{\bm{u}_{l}^{\star\top}\bm{u}_{l}}{\lambda_{l}/\lambda_{l}^{\star}}\cdot\bm{a}^{\top}\bm{u}_{l}^{\star}\right|\right);\\
|\tau_{2}|=o\left(\frac{\sigma_{\mathrm{max}}}{|\lambda_{l}^{\star}|}\right) & \leq o\left(|\bm{a}^{\top}\bm{u}_{l}^{\star}|\right)=o\left(\left|\frac{\bm{u}_{l}^{\star\top}\bm{w}_{l}}{\lambda_{l}/\lambda_{l}^{\star}}\cdot\bm{a}^{\top}\bm{u}_{l}^{\star}\right|\right).
\end{align*}
In addition, Lemma~\ref{lem:rank-r-inference-remainder-size} together with the assumption $\sigma_{\mathrm{max}} \kappa^2 r \sqrt{\mu \log n} = o\left(|\Delta_l^\star|\right)$ implies that $|{\tau}_3| = o(1)$.

With these bounds in place, we see that $\tau_1$ and $\frac{\bm{a}^{\top}\bm{H}\bm{u}_{l}^{\star}}{\lambda_{l}^{\star}}$ are indeed very small terms, meaning that the term in 
\eqref{eq:estimator-squared-decompose-first-order-with-abs} involving $\tau_1$ obeys
$\frac{\bm{u}_{l}^{\star\top}\bm{u}_{l}}{\lambda_{l}/\lambda_{l}^{\star}}\cdot\left(\bm{a}^{\top}\bm{u}_{l}^{\star}+\frac{\bm{a}^{\top}\bm{H}\bm{u}_{l}^{\star}}{\lambda_{l}^{\star}}\right)+{\tau}_{1}$ and $\frac{\bm{u}_{l}^{\star\top}\bm{u}_{l}}{\lambda_{l}/\lambda_{l}^{\star}}\bm{a}^{\top}\bm{u}_{l}^{\star}$ have the same signs. Similar conclusions hold as well for the terms involving $\tau_2$ and $\tau_3$. As a result, 
\begin{align*}
 & \mathsf{sign}\left\{ \frac{\left\{ \frac{\bm{u}_{l}^{\star\top}\bm{u}_{l}}{\lambda_{l}/\lambda_{l}^{\star}}\cdot\left(\bm{a}^{\top}\bm{u}_{l}^{\star}+\frac{\bm{a}^{\top}\bm{H}\bm{u}_{l}^{\star}}{\lambda_{l}^{\star}}\right)+\tau_{1}\right\} \left\{ \frac{\bm{u}_{l}^{\star\top}\bm{w}_{l}}{\lambda_{l}/\lambda_{l}^{\star}}\cdot\left(\bm{a}^{\top}\bm{u}_{l}^{\star}+\frac{\bm{a}^{\top}\bm{H}^{\top}\bm{u}_{l}^{\star}}{\lambda_{l}^{\star}}\right)+\tau_{2}\right\} }{\frac{\bm{u}_{l}^{\star\top}\bm{w}_{l}}{\lambda_{l}/\lambda_{l}^{\star}}\cdot\frac{\bm{u}_{l}^{\star\top}\bm{u}_{l}}{\lambda_{l}/\lambda_{l}^{\star}}\left(1+\frac{\bm{u}_{l}^{\star\top}(\bm{H}+\bm{H}^{\top})\bm{u}_{l}^{\star}}{\lambda_{l}^{\star}}\right)+\tau_{3}}\right\} \\
 & \quad=\mathsf{sign}\left\{ \frac{\left\{ \frac{\bm{u}_{l}^{\star\top}\bm{u}_{l}}{\lambda_{l}/\lambda_{l}^{\star}}\bm{a}^{\top}\bm{u}_{l}^{\star}\right\} \left\{ \frac{\bm{u}_{l}^{\star\top}\bm{w}_{l}}{\lambda_{l}/\lambda_{l}^{\star}}\bm{a}^{\top}\bm{u}_{l}^{\star}\right\} }{\frac{\bm{u}_{l}^{\star\top}\bm{w}_{l}}{\lambda_{l}/\lambda_{l}^{\star}}\cdot\frac{\bm{u}_{l}^{\star\top}\bm{u}_{l}}{\lambda_{l}/\lambda_{l}^{\star}}}\right\} =\mathsf{sign}\left\{ (\bm{a}^{\top}\bm{u}_{l}^{\star})^{2}\right\} =1.
\end{align*}
This indicates that we can safely remove the absolute value function in \eqref{eq:estimator-squared-decompose-first-order-with-abs} to obtain
\begin{equation}
\big(\widehat{u}_{\bm{a},l}^{\mathsf{modified}}\big)^{2} = \frac{\left\{\frac{\bm{u}_{l}^{\star\top}\bm{u}_{l}}{\lambda_{l}/\lambda_{l}^{\star}}  \left(\bm{a}^{\top}\bm{u}_{l}^{\star} + \frac{\bm{a}^\top \bm{H} \bm{u}_l^\star}{\lambda_l^\star}\right)+{\tau}_{1}\right\}\left\{\frac{\bm{u}_{l}^{\star\top}\bm{w}_{l}}{\lambda_{l}/\lambda_{l}^{\star}}  \left(\bm{a}^{\top}\bm{u}_{l}^{\star} + \frac{\bm{a}^\top \bm{H}^\top \bm{u}_l^\star}{\lambda_l^\star}\right)+{\tau}_{2}\right\}}{\frac{\bm{u}_{l}^{\star\top}\bm{w}_{l}}{\lambda_{l}/\lambda_{l}^{\star}} \frac{\bm{u}_{l}^{\star\top}\bm{u}_{l}}{\lambda_{l}/\lambda_{l}^{\star}} \left(1 + \frac{\bm{u}_l^{\star \top} (\bm{H} + \bm{H}^\top) \bm{u}_l^\star}{\lambda_l^\star}\right)+{\tau}_{3}}.
\label{eq:estimator-squared-decompose-first-order}
\end{equation}
Armed with this expression, the next lemma  develops the distributional characterization of $\widehat{u}_{\bm{a},l}^{\mathsf{modified}}$. 
\begin{lemma} \label{lem:rank-r-inference-squared-estimator-difference-size}
	Instate the assumptions of Theorem \ref{thm:confidence-interval-validity-rankr}.  With probability at least $1-O(n^{-6})$, it follows that
	\begin{align}
		&\left|\big(\widehat{u}_{\bm{a},l}^{\mathsf{modified}}\big)^{2} - \left(\bm{a}^{\top}\bm{u}^{\star}_l + \frac{1}{2\lambda^{\star}_l}(\bm{a}_l^{\perp })^{\top}\big(\bm{H}+\bm{H}^{\top}\big)\bm{u}^{\star}_l\right)^2 \right|  
	\nonumber \\
	 & \qquad \qquad 
	\lesssim | \bm{a}^\top \bm{u}_l^\star | \cdot (|{\tau_1}| + |{\tau_2}|) + (\bm{a}^\top \bm{u}_l^\star)^2 \cdot |{\tau_3}| + O\left(\frac{\sigma_{\mathrm{max}}^2 \log n}{\lambda_l^{\star 2}}\right).
	\end{align}
\end{lemma}
\noindent The proof of this result is given in Section~\ref{PfLemmaRankRInf}.

Given our assumption $\bm{a}^{\top} \bm{u}_l^{\star}>0$, the expression \eqref{eq:estimator-squared-decompose-first-order} together with the bounds on $\tau_1$, $\tau_2$ and $\tau_3$ immediately implies that
$$
\Bigg|\widehat{u}_{\bm{a},l}^{\mathsf{modified}} + \left( \bm{a}^{\top}\bm{u}^{\star}_l + \frac{1}{2\lambda^{\star}_l}\bm{a}_l^{\perp \top}\big(\bm{H}+\bm{H}^{\top}\big)\bm{u}^{\star}_l \right)\Bigg| \geq (1- o(1)) \cdot \bm{a}^{\top}\bm{u}^{\star}_l.
$$
This combined with Lemma \ref{lem:rank-r-inference-squared-estimator-difference-size} leads to
\begin{align*}
\Bigg|\widehat{u}_{\bm{a},l}^{\mathsf{modified}}-\left(\bm{a}^{\top}\bm{u}_{l}^{\star}+\frac{1}{2\lambda_{l}^{\star}}\bm{a}_{l}^{\perp\top}\big(\bm{H}+\bm{H}^{\top}\big)\bm{u}_{l}^{\star}\right)\Bigg| 
	& =\frac{\Big|(\widehat{u}_{\bm{a},l})^{2}-\left(\bm{a}^{\top}\bm{u}_{l}^{\star}+\frac{1}{2\lambda_{l}^{\star}}\bm{a}_{l}^{\perp\top}\big(\bm{H}+\bm{H}^{\top}\big)\bm{u}_{l}^{\star}\right)^{2}\Big|}{\Big|\widehat{u}_{\bm{a},l}+\left(\bm{a}^{\top}\bm{u}_{l}^{\star}+\frac{1}{2\lambda_{l}^{\star}}\bm{a}_{l}^{\perp\top}\big(\bm{H}+\bm{H}^{\top}\big)\bm{u}_{l}^{\star}\right)\Big|}\nonumber\\
 & \lesssim|\tau_{1}|+|\tau_{2}|+|\bm{a}^{\top}\bm{u}_{l}^{\star}|\cdot|\tau_{3}|+O\left(\frac{\sigma_{\mathrm{max}}^{2}\log n}{\lambda_{l}^{\star2}\cdot|\bm{a}^{\top}\bm{u}_{l}^{\star}|}\right).
\end{align*}
Taking this together with Lemma \ref{lem:rank-r-inference-remainder-size} as well as the condition \eqref{eq:sigma-large-case-condition}, we arrive at
\begin{align}
\Bigg|\widehat{u}_{\bm{a},l}^{\mathsf{modified}} - \left(\bm{a}^{\top}\bm{u}^{\star}_l + \frac{1}{2\lambda^{\star}_l}\bm{a}_l^{\perp \top}\big(\bm{H}+\bm{H}^{\top}\big)\bm{u}^{\star}_l \right)\Bigg| & \lesssim o \left(\frac{\sigma_{\mathrm{max}}}{|\lambda_l^\star|}\right) + | \bm{a}^\top \bm{u}_l^{\star} | \cdot \left(\frac{\sigma_{\mathrm{max}} \kappa^2 r \sqrt{\mu \log n}}{|\Delta_l^\star|}\right)^2 \nonumber \\
& = o \left(\frac{\sigma_{\mathrm{max}}}{|\lambda_l^\star|}\right),
\end{align}
where the last line results from our assumption $\bm{a}^\top \bm{u}_l^\star = o \left(\frac{\Delta_l^{\star 2} }{|\lambda_l^\star|  \sigma_{\mathrm{max}} \kappa^4 r^2 \mu \log n}\right)$. 

The proof of the claim~\eqref{eq:amazing-fact} is thus complete for this case.

\subsubsection{The case when \texorpdfstring{$|\bm{a}^{\top}\bm{u}_l^{\star}|$}{the truth}  is ``small''}
\label{sec:proof-au-small-inference}

We then move on to the scenario where $|\bm{a}^{\top}\bm{u}_l^{\star}|\lesssim \sqrt{v_{\bm{a},l}^{\star}}\log^{1.5} n$, which clearly subsumes the case with $|\bm{a}^{\top}\bm{u}_l^{\star}|\lesssim \sqrt{\widehat{v}_{\bm{a},l}}\log^{1.5} n$ (according to Lemma~\ref{Lem-var-control-rankr}). 
Once again, this restriction combined with Lemma \ref{lem:va-properties} and the assumption $\sigma_{\max}/\sigma_{\min}=O(1)$ requires that
\begin{align}
	|\bm{a}^{\top}\bm{u}_l^{\star}|\lesssim \frac{\sigma_{\max}\log^{1.5} n}{|\lambda_l^{\star}|}.
\end{align}

The proof is built upon the following ``first-order'' approximations of $\bm{a}^{\top}\bm{u}_l$ and $\bm{a}^{\top}\bm{w}_l$.  
\begin{lemma} 
\label{lem:au-distribution-rankr} 
Instate the assumptions of Theorem \ref{thm:confidence-interval-validity-rankr}. 
Fix any vector $\bm{a} \in \mathbb{R}^n$ with $\|\bm{a}\|_2=1$, and assume that 
\begin{align}
|\bm{a}^{\top}\bm{u}_{k}^{\star}|  &=o\left(\frac{|\lambda_{l}^{\star}-\lambda_{k}^{\star}|}{|\lambda_{l}^{\star}|\sqrt{\mu\kappa^{4}r^{3}\log n}}\right),
  \qquad \forall k\neq l.  
\end{align}
Then under Assumption $\ref{assumption:noise-size-rankr}$, with probability $1 - O(n^{-6})$ one has
\begin{align}
\label{eqn:coconut}
\left|\frac{\lambda_{l}}{\lambda_{l}^{\star}(\bm{u}_{l}^{\star\top}\bm{u}_{l})}\bm{a}^{\top}\bm{u}_{l}-\bm{a}^{\top}\bm{u}_{l}^{\star}-\frac{\bm{a}^{\top}\bm{H}\bm{u}_{l}^{\star}}{\lambda_{l}^{\star}}\right|=o\left(\frac{\sigma_{\max}}{|\lambda_{l}^{\star}|}\right) 
\end{align}	
%
%
If we further have $|\bm{a}^{\top}\bm{u}_{l}^{\star}|=o\left(1/\sqrt{\log n}\right)$, then with probability $1 - O(n^{-6})$, 
\begin{align}
\begin{cases}
\frac{1}{\bm{u}_l^{\star \top} \bm{u}_l} \bm{a}^\top \bm{u}_l = \bm{a}^\top \bm{u}_l^\star + \frac{(\bm{a}_l^{\perp})^{\top} \bm{H} \bm{u}_l^\star}{\lambda_l^\star} + o \left(\frac{\sigma_{\mathrm{max}}}{|\lambda_{l}^\star|} \right) ; \\
\frac{1}{\bm{u}_l^{\star \top} \bm{w}_l} \bm{a}^\top \bm{w}_l = \bm{a}^\top \bm{u}_l^\star + \frac{(\bm{a}_l^{\perp})^{\top} \bm{H}^\top \bm{u}_l^\star}{\lambda_l^\star} + o \left(\frac{\sigma_{\mathrm{max}}}{|\lambda_{l}^\star|} \right) .
\end{cases}
\end{align}
\end{lemma}
\begin{proof} See Appendix~\ref{sec:proof-lemma:au-distribution-rankr}. \end{proof}

In addition, we claim that the following relations hold
\begin{align}
\label{eq:sonata}
\begin{cases}
\frac{1}{\bm{u}_{l}^{\top}\bm{u}_{l}^{\star}}\bm{a}^{\top}\bm{u}_{l}=\bm{a}^{\top}\bm{u}_{l}+o\Big(\frac{\sigma_{\max}}{|\lambda_{l}^{\star}|}\Big),\\
\frac{1}{\bm{w}_{l}^{\top}\bm{u}_{l}^{\star}}\bm{a}^{\top}\bm{w}_{l}=\bm{a}^{\top}\bm{w}_{l}+o\Big(\frac{\sigma_{\max}}{|\lambda_{l}^{\star}|}\Big).
\end{cases}
\end{align}
If these claims were valid, then one would have
\begin{align*}
\widehat{u}_{\bm{a},l}^{\mathsf{modified}} & =\left(\bm{a}^{\top}\bm{u}_{l}+\bm{a}^{\top}\bm{w}_{l}\right)/2
  =\frac{1}{2\bm{u}_{l}^{\top}\bm{u}_{l}^{\star}}\bm{a}^{\top}\bm{u}_{l}^{\star}+\frac{1}{2\bm{w}_{l}^{\top}\bm{u}_{l}^{\star}}\bm{a}^{\top}\bm{u}_{l}^{\star}+o\Big(\frac{\sigma_{\max}}{|\lambda_{l}^{\star}|}\Big)\\
	& =\bm{a}^{\top}\bm{u}_{l}^{\star}+\frac{(\bm{a}_l^{\perp})^{\top}(\bm{H}+\bm{H}^{\top})\bm{u}_{l}^{\star}}{2\lambda_{l}^{\star}}+o\Big(\frac{\sigma_{\max}}{|\lambda_{l}^{\star}|}\Big),
\end{align*}
where the last inequality follows from Lemma~\ref{lem:au-distribution-rankr}. 
This validates the relation~\eqref{eq:amazing-fact} for this case, as long as the relations~\eqref{eq:sonata} hold true.

\paragraph*{Proof of the relations~\eqref{eq:sonata}:} As a direct consequence of Lemma~\ref{lem:au-distribution-rankr}, we can write 
\begin{align}
\frac{1}{\bm{u}_{l}^{\top}\bm{u}_{l}^{\star}}\bm{a}^{\top}\bm{u}_{l}-\bm{a}^{\top}\bm{u}_{l} & =\frac{1}{\bm{u}_{l}^{\top}\bm{u}_{l}^{\star}}\bm{a}^{\top}\bm{u}_{l}\cdot\big(1-\bm{u}^{\star\top}\bm{u}_{l}\big)\nonumber\\
 & =\left(\bm{a}^{\top}\bm{u}_{l}^{\star}+\frac{\bm{a}^{\top}\bm{H}\bm{u}_{l}^{\star}}{\lambda_{l}^{\star}}+o\Big(\frac{\sigma_{\max}}{|\lambda_{l}^{\star}|}\Big)\right)(1-\bm{u}^{\star\top}\bm{u}_{l}).
	\label{eq:adagio}
\end{align}
In view of Lemma \ref{lem:aHu-UB}, we have with probability at least $1 - O(n^{-10})$ that
\begin{align*}
	\left| \frac{\bm{a}^{\top} \bm{H} \ustar_l}{\lambdastar_l} \right|
	= O \left(\frac{\sigma_{\max}\sqrt{\log n}}{|\lambdastar_l|} 
	\right).
\end{align*}
In addition, Lemma~\ref{lem:basic-uR-uL-rankr} together with Assumption~\ref{assumption:noise-size-rankr}
guarantees that, with high probability, 
\begin{align*}
	| 1 - \bm{u}_l^{\star\top}\ur_l |  &= O\left(\frac{\kappa^{4}\sigma_{\max}^{2}n\log n}{(\lambda_{\mathrm{max}}^{\star})^{2}}+\frac{\mu\kappa^{4}r^{2}\sigma_{\max}^{2}\log n}{(\Delta_{l}^{\star})^{2}}\right) 
	=o\Big( \frac{1}{ {\log^{1.5} n} } \Big).
\end{align*}
Putting everything together and using the assumption $|\bm{a}^{\top}\ustar_l|\lesssim \frac{\sigma_{\max}{\log^{1.5} n}}{|\lambdastar_l|}$, we conclude that
\begin{align*}
	\frac{1}{\ur_l^{\top}\ustar}\bm{a}^{\top}\ur_l - \bm{a}^{\top}\ur = 
	 o\left(\frac{\sigmamax}{|\lambdastar_l|}\right)
\end{align*}
as claimed. The claim w.r.t.~$\bm{w}_l$ follows from exactly the same argument.

\subsection{Distributional theory for eigenvalues}
\label{sec:proof-thm:distribution-eigenvalues-rankr}

Take $\bm{a} = {\bm{u}^\star_l}$. By definition, we have $\bm{a}^\top \bm{u}_k^\star = 0$ for any $k \neq l$. 
The expression~\eqref{eqn:coconut} in Lemma~\ref{lem:au-distribution-rankr} thus indicates that, with probability $1 - O(n^{-6})$, we have 
\begin{align*}
\left|\frac{\lambda_{l}}{\lambda_{l}^{\star}(\bm{u}_{l}^{\star\top}\bm{u}_{l})}{\bm{u}^\star_l}^{\top}\bm{u}_{l}-{\bm{u}^\star_l}^{\top}\bm{u}_{l}^{\star}-\frac{{\bm{u}^\star_l}^{\top}\bm{H}\bm{u}_{l}^{\star}}{\lambda_{l}^{\star}}\right|=o\left(\frac{\sigma_{\max}}{|\lambda_{l}^{\star}|}\right) .
\end{align*}	
Rearranging terms and using the fact $\bm{u}^{\star\top}_l\bm{u}_{l}^{\star}=1$ further yield
\begin{align}
	|\lambda_{l} - \lambdastar_l - 
	{\bm{u}^\star_l}^{\top}\bm{H}\bm{u}_{l}^{\star}| = o\left(\sigma_{\max}\right).
	\label{eq:lambda-l-detailed-approximation}
\end{align}
Consequently, it is sufficient to characterize the distribution of ${\bm{u}^\star_l}^{\top}\bm{H}\bm{u}_{l}^{\star}.$ Setting $\bm{a} = {\bm{u}^\star_l}$ in Lemma~\ref{LemPropertyInprod}, we see that  $	W_{\lambda,l} = \frac{{\bm{u}^\star_l}^{\top}\bm{H}\bm{u}_{l}^{\star}}{\sqrt{v_{\lambda,l}^{\star}}} = \frac{{\bm{u}^\star_l}^{\top}\left(\bm{H}+\bm{H}^{\top}\right)\bm{u}_l^{\star}}{2\sqrt{v_{\lambda,l}^{\star}}}$ obeys
\begin{align}
	\sup_{z\in \real} \big| \, \mprob(W_{\lambda,l}\leq z) - \Phi(z) \,\big| \leq \frac{8}{\sqrt{\log n}} 
\end{align}
as claimed.

\subsection{Proof of Lemma \ref{lem:rank-r-inference-remainder-size}}

First of all, invoke the Neumann series (cf.~Lemma~\ref{lem:Neumann-expansion} and \eqref{eq:au-expansion-neumann})
and rearrange terms to obtain
\begin{align*}
\tau_{1} & =\bm{a}^{\top}\bm{u}_{l}-\frac{\bm{u}_{l}^{\star\top}\bm{u}_{l}}{\lambda_{l}/\lambda_{l}^{\star}}\left(\bm{a}^{\top}\bm{u}_{l}^{\star}+\frac{\bm{a}^{\top}\bm{H}\bm{u}_{l}^{\star}}{\lambda_{l}^{\star}}\right)\\
 & =\frac{\bm{u}_{l}^{\star\top}\bm{u}_{l}}{\lambda_{l}/\lambda_{l}^{\star}}\left[\bm{a}^{\top}\bm{u}_{l}^{\star}+\frac{1}{\lambda_{l}}\bm{a}^{\top}\bm{H}\bm{u}_{l}^{\star}+\sum_{s=2}^{\infty}\frac{1}{\lambda_{l}^{s}}\bm{a}^{\top}\bm{H}^{s}\bm{u}_{l}^{\star}+\sum_{k:k\neq l}\frac{\lambda_{k}^{\star}}{\lambda_{l}^{\star}}\frac{\bm{u}_{k}^{\star\top}\bm{u}_{l}}{\bm{u}_{l}^{\star\top}\bm{u}_{l}}\left\{ \sum_{s=0}^{\infty}\frac{1}{\lambda_{l}^{s}}\bm{a}^{\top}\bm{H}^{s}\bm{u}_{k}^{\star}\right\} \right]\\
 & \qquad-\frac{\bm{u}_{l}^{\star\top}\bm{u}_{l}}{\lambda_{l}/\lambda_{l}^{\star}}\cdot\left(\bm{a}^{\top}\bm{u}_{l}^{\star}+\frac{\bm{a}^{\top}\bm{H}\bm{u}_{l}^{\star}}{\lambda_{l}^{\star}}\right)\\
 & =\underset{=:\tau_{1,1}}{\underbrace{\frac{\bm{u}_{l}^{\star\top}\bm{u}_{l}}{\lambda_{l}/\lambda_{l}^{\star}}\left(\frac{1}{\lambda_{l}}-\frac{1}{\lambda_{l}^{\star}}\right)\bm{a}^{\top}\bm{H}\bm{u}_{l}^{\star}}}+\underset{=:\tau_{1,2}}{\underbrace{\sum_{k:k\neq l}\frac{\bm{u}_{k}^{\star\top}\bm{u}_{l}}{\lambda_{l}/\lambda_{k}^{\star}}\sum_{s=2}^{\infty}\frac{1}{\lambda_{l}^{s}}\bm{a}^{\top}\bm{H}^{s}\bm{u}_{k}^{\star}}}+\underset{=:\tau_{1,3}}{\underbrace{\sum_{k=1}^{r}\frac{\bm{u}_{k}^{\star\top}\bm{u}_{l}}{\lambda_{l}/\lambda_{k}^{\star}}\left\{ \sum_{s=2}^{\infty}\frac{1}{\lambda_{l}^{s}}\bm{a}^{\top}\bm{H}^{s}\bm{u}_{k}^{\star}\right\} }}
\end{align*}
and, similarly,
\begin{align*}
\tau_{3} & =\bm{u}_{l}^{\top}\bm{w}_{l}-\frac{\bm{u}_{l}^{\star\top}\bm{w}_{l}}{\lambda_{l}/\lambda_{l}^{\star}}\cdot\frac{\bm{u}_{l}^{\star\top}\bm{u}_{l}}{\lambda_{l}/\lambda_{l}^{\star}}\left(1+\frac{\bm{u}_{l}^{\star\top}(\bm{H}+\bm{H}^{\top})\bm{u}_{l}^{\star}}{\lambda_{l}^{\star}}\right)\\
 & =\left(\sum_{j=1}^{r}\frac{\lambda_{j}^{\star}}{\lambda_{l}}\big(\bm{u}_{j}^{\star\top}\bm{u}_{l}\big)\left\{ \sum_{s=0}^{\infty}\frac{1}{\lambda_{l}^{s}}\bm{H}^{s}\bm{u}_{j}^{\star}\right\} \right)\left(\sum_{j=1}^{r}\frac{\lambda_{j}^{\star}}{\lambda_{l}}\big(\bm{u}_{j}^{\star\top}\bm{w}_{l}\big)\left\{ \sum_{s=0}^{\infty}\frac{1}{\lambda_{l}^{s}}\bm{H}^{s}\bm{u}_{j}^{\star}\right\} \right)\\
 & \qquad-\frac{\bm{u}_{l}^{\star\top}\bm{w}_{l}}{\lambda_{l}/\lambda_{l}^{\star}}\cdot\frac{\bm{u}_{l}^{\star\top}\bm{u}_{l}}{\lambda_{l}/\lambda_{l}^{\star}}\left(1+\frac{\bm{u}_{l}^{\star\top}(\bm{H}+\bm{H}^{\top})\bm{u}_{l}^{\star}}{\lambda_{l}^{\star}}\right)\\
 & =\underset{=:\tau_{3,1}}{\underbrace{\frac{\bm{u}_{l}^{\star\top}\bm{w}_{l}}{\lambda_{l}/\lambda_{l}^{\star}}\cdot\frac{\bm{u}_{l}^{\star\top}\bm{u}_{l}}{\lambda_{l}/\lambda_{l}^{\star}}\cdot\bm{u}_{l}^{\star\top}(\bm{H}+\bm{H}^{\top})\bm{u}_{l}^{\star}\cdot\left(\frac{1}{\lambda_{l}}-\frac{1}{\lambda_{l}^{\star}}\right)}}+\underset{=:\tau_{3,2}}{\underbrace{\sum_{k:k\neq l}\frac{\bm{u}_{k}^{\star\top}\bm{w}_{l}}{\lambda_{l}/\lambda_{k}^{\star}}\cdot\frac{\bm{u}_{k}^{\star\top}\bm{u}_{l}}{\lambda_{l}/\lambda_{k}^{\star}}}}+\\
 & \quad\underset{=:\tau_{3,3}}{\underbrace{\sum_{(k_{1},k_{2})\neq(l,l)}\frac{\bm{u}_{k_{1}}^{\star\top}\bm{w}_{l}}{\lambda_{l}/\lambda_{k_{1}}^{\star}}\cdot\frac{\bm{u}_{k_{2}}^{\star\top}\bm{u}_{l}}{\lambda_{l}/\lambda_{k_{2}}^{\star}}\cdot\frac{2\bm{u}_{k_{1}}^{\star\top}\bm{H}\bm{u}_{k_{2}}^{\star}}{\lambda_{l}}}}+\underset{=:\tau_{3,4}}{\underbrace{\sum_{k_{1}=1}^{r}\sum_{k_{2}=1}^{r}\sum_{s_{1},s_{2}:s_{1}+s_{2}\geq2}\frac{\bm{u}_{k_{1}}^{\star\top}\bm{w}_{l}}{\lambda_{l}/\lambda_{k_{1}}^{\star}}\cdot\frac{\bm{u}_{k_{2}}^{\star\top}\bm{u}_{l}}{\lambda_{l}/\lambda_{k_{2}}^{\star}}\cdot\frac{\bm{u}_{k_{1}}^{\star\top}\bm{H}^{s_{1}+s_{2}}\bm{u}_{k_{2}}^{\star}}{\lambda_{l}^{s_{1}+s_{2}}}}}.
\end{align*}
In the special case where $r=1$, one has $\tau_{1,2}=\tau_{3,2}=\tau_{3,3}=0$.  In what follows, we develop bounds for these terms separately.

\paragraph{Controlling \texorpdfstring{${\tau}_{1,1}$}{} and \texorpdfstring{${\tau}_{3,1}$}{}.}
We have learned from Theorem~\ref{thm:eigengap-condition}, Lemma~\ref{lem:aHu-UB} and Lemma~\ref{lem:basic-uR-uL-rankr} that: under Assumption \ref{assumption:noise-size-rankr}, one has $\lambda_l = (1+o(1))\lambda_l^\star$, $\bm{u}_l^{\top} \bm{u}_l^\star=1-o(1)$, 
\begin{align}
	\left|\frac{\lambda_{l}-\lambda_{l}^{\star}}{\lambda_{l}}\right| & \lesssim \frac{\sigma_{\mathrm{max}}\sqrt{\mu\kappa^{2}r^{4}\log n}}{|\lambda_{l}^{\star}|}  \qquad \text{and} \label{eq:rank-r-inference-lambda-l-bound} \\
\max \left\{|\bm{a}^{\top}\bm{H}\bm{u}_{l}^{\star}|, |\bm{u}_{l}^{\star\top}\bm{H}\bm{u}_{l}^{\star}|\right\} & \lesssim\sigma_{\mathrm{max}}\sqrt{\log n}.
\end{align}
It then follows that
\begin{align*}
|\tau_{1,1}| & \lesssim\left|\bm{a}^{\top}\bm{H}\bm{u}_{l}^{\star}\right|\cdot\Big|\frac{\lambda_{l}-\lambda_{l}^{\star}}{(\lambda_{l}^{\star})^{2}}\Big|\lesssim\frac{\sigma_{\mathrm{max}}^{2}\kappa r^{2}\sqrt{\mu}\log n}{|\lambda_{l}^{\star}|^{2}}\leq\frac{\sigma_{\mathrm{max}}}{|\lambda_{l}^{\star}|}\cdot\frac{\sigma_{\mathrm{max}}\kappa r^{2}\sqrt{\mu}\log n}{|\lambda_{\mathrm{min}}^{\star}|}=o\Big(\frac{\sigma_{\mathrm{max}}}{|\lambda_{l}^{\star}|}\Big) \\
|\tau_{3,1}| & \lesssim\left|\bm{u}_{l}^{\star\top}\bm{H}\bm{u}_{l}^{\star}\right|\cdot\Big|\frac{\lambda_{l}-\lambda_{l}^{\star}}{(\lambda_{l}^{\star})^{2}}\Big|=o\Big(\frac{\sigma_{\mathrm{max}}}{|\lambda_{l}^{\star}|}\Big)
\end{align*}
with the proviso that $\sigma_{\max}\kappa r^{2}\sqrt{\mu}\log n=o\left(\lambda_{\min}^{\star}\right)$.

\paragraph{Controlling \texorpdfstring{${\tau}_{1,2}$}{} and \texorpdfstring{${\tau}_{3,3}$}{}.} We make the observation that for any $1 \leq k \leq r$,
\begin{align}
\left|\frac{\bm{a}^{\top}\bm{H}\bm{u}_{k}^{\star}}{\lambda_{l}}\right| & \overset{\mathrm{(i)}}{\lesssim}\frac{\sigma_{\mathrm{max}}\sqrt{\log n}}{|\lambda_{l}|}\asymp\frac{\sigma_{\mathrm{max}}\sqrt{\log n}}{|\lambda_{l}^{\star}|}\nonumber 
\overset{\mathrm{(ii)}}{=}o\left(\frac{\Delta_{l}^{\star}}{|\lambda_{l}^{\star}|\sqrt{\mu\kappa^{6}r^{3}\log^{2}n}}\sqrt{\log n}\right)\nonumber\\
& =o\left(\frac{|\lambda_{l}^{\star}-\lambda_{k}^{\star}|}{|\lambda_{l}^{\star}|\sqrt{\mu\kappa^{4}r^{3}\log n}}\right),
\end{align}
where (i) is a consequence of Lemma \ref{lem:aHu-UB}, and (ii) holds as long as $\sigma_{\max}\sqrt{\mu\kappa^{6}r^{3}\log^{2}n}=o(\Delta_{l}^{\star})$. Similarly, for any $1 \leq k_1, k_2 \leq r$,
\begin{align}
\left|\frac{\bm{u}_{k_1}^{\star \top} \bm{H} \bm{u}_{k_2}^\star }{\lambda_l}\right| = o\left(\frac{\Delta_{l}^{\star}}{|\lambda_{l}^{\star}|\sqrt{\mu\kappa^{6}r^{3}\log n}}\right). \label{eq:rank-r-inference-uk1-H-uk2-bound}
\end{align}
Additionally, one can invoke the inequality \eqref{eq:uk-ui-expansion-bound2} to deduce that for any $k \neq l$,
\begin{align}
\left|\bm{u}_{k}^{\star\top}\bm{u}_{l}\right| & \lesssim\left|1-\frac{\lambda_{k}^{\star}}{\lambda_{l}}\right|^{-1}\frac{\sigma_{\mathrm{max}}}{\lambda_{\mathrm{min}}^{\star}}\sqrt{\mu\kappa^{2}r\log n}
\lesssim\frac{\lambda_{l}^{\star}}{|\lambda_{l}^{\star}-\lambda_{k}^{\star}|}\cdot\frac{\sigma_{\mathrm{max}}}{\lambda_{\mathrm{min}}^{\star}}\sqrt{\mu\kappa^{2}r\log n}.	 \label{eq:rank-r-inference-inner-prod-bound}   \end{align}
Putting the above bounds together and using the conditions $|\bm{a}^{\top}\bm{u}_{k}^{\star}|  =o\Big(\frac{|\lambda_{l}^{\star}-\lambda_{k}^{\star}|}{|\lambda_{l}^{\star}|\sqrt{\mu\kappa^{4}r^{3}\log n}}\Big)$ and $|\bm{u}_l^{\star \top} \bm{u}_l| = 1 -o(1)$ give
\begin{align}
|{\tau}_{1,2}| & \lesssim \sum_{k\neq l,1\leq k\leq r}\left|\frac{\lambda_{k}^{\star}}{\lambda_{l}^{\star}}\right|\cdot\left|\frac{\bm{u}_{k}^{\star\top}\bm{u}_{l}}{\bm{u}_{l}^{\star\top}\bm{u}_{l}}\right|\left(\left|\frac{\bm{a}^{\top}\bm{H}\bm{u}_{k}^{\star}}{\lambda_{l}}\right| + |\bm{a}^{\top}\bm{u}_{k}^{\star}| \right) \nonumber \\
& =  o\left( \frac{\sigma_{\max}}{|\lambda_l^{\star}|} \right) + o\left(r\cdot\frac{\lambda_{\max}^{\star}}{|\lambda_{l}^{\star}-\lambda_{k}^{\star}|}\cdot\frac{\sigma_{\mathrm{max}}}{\lambda_{\mathrm{min}}^{\star}}\sqrt{\mu\kappa^{2}r\log n}\cdot\frac{|\lambda_{l}^{\star}-\lambda_{k}^{\star}|}{|\lambda_{l}^{\star}|\sqrt{\mu\kappa^{4}r^{3}\log n}}\right) =o\left(\frac{\sigma_{\mathrm{max}}}{|\lambda_{l}^{\star}|}\right).
\end{align}
In addition, for every pair $(k_1, k_2) \neq (l, l)$, Theorem \ref{thm:rankr-bounds} implies that
%
\begin{align}
\left(\bm{u}_{k_1}^{\star\top}\bm{u}_{l}\right)\left(\bm{u}_{k_2}^{\star\top}\bm{u}_{l}\right)=\begin{cases}
O\left(\frac{\sigma_{\mathrm{max}}}{\Delta_{l}^{\star}}\sqrt{\mu\kappa^{4}r\log n}\right), & \text{if }k_1=l,k_2\neq l\text{ or }k_1\neq l,k_2=l,\\
O\left(\frac{\sigma_{\mathrm{max}}^{2}}{(\Delta_{l}^{\star})^{2}}\mu\kappa^{4}r\log n\right), & \text{if }k_1\neq l,k_2\neq l, 
\end{cases}
\end{align}
which together with \eqref{eq:rank-r-inference-uk1-H-uk2-bound} and $\sigma_{\max}\sqrt{\mu\kappa^{6}r^{3}\log^{2}n}=o(\Delta_{l}^{\star})$ indicates that
\begin{align*}
I_1 := \left|\sum_{k_1\neq l}\sum_{k_2\neq l}\frac{\lambda_{k_1}^{\star}\lambda_{k_2}^{\star}}{\lambda_{l}^{2}}\left(\bm{u}_{k_1}^{\star\top}\bm{w}_{l}\right)\left(\bm{u}_{k_2}^{\star\top}\bm{u}_{l}\right)\cdot \frac{2\bm{u}_{k_1}^{\star\top}\bm{H}\bm{u}_{k_2}^{\star}}{\lambda_{l}}\right| & \lesssim r^{2}\kappa^{2}\max_{k\neq l,j\neq l}\big|\bm{u}_{k}^{\star\top}\bm{u}_{l}\big|\cdot\big|\bm{u}_{j}^{\star\top}\bm{w}_{l}\big| \cdot o\left(\frac{\Delta_{l}^{\star}}{|\lambda_{l}^{\star}|\sqrt{\mu\kappa^{6}r^{3}\log n}}\right)\\
& = o\left(\left( \frac{\sigma_{\mathrm{max}}}{\Delta_{l}^{\star}}\sqrt{\mu\kappa^{6}r^{3}\log n} \right)^2 \cdot \frac{\Delta_{l}^{\star}}{|\lambda_{l}^{\star}|\sqrt{\mu\kappa^{6}r^{3}\log n}}\right) \\
& = o \left(\frac{\sigma_{\mathrm{max}}}{|\lambda_l^\star|} \cdot \frac{\sigma_{\mathrm{max}}}{\Delta_{l}^{\star}}\sqrt{\mu\kappa^{6}r^{3}\log n} \right) \\
& = o \left(\frac{\sigma_{\mathrm{max}}}{|\lambda_l^\star|} \right); \\
I_2 := \left|\sum_{k_1\neq l}\frac{\lambda_{k_1}^{\star}\lambda_{l}^{\star}}{\lambda_{l}^{2}}\left(\bm{u}_{k_1}^{\star\top}\bm{w}_{l}\right)\left(\bm{u}_{l}^{\star\top}\bm{u}_{l}\right)\cdot \frac{2\bm{u}_{k_1}^{\star\top}\bm{H}\bm{u}_{l}^{\star}}{\lambda_{l}}\right| & \lesssim r\kappa\max_{{k_1}\neq l}\big|\bm{u}_{k_1}^{\star\top}\bm{w}_{l}\big|\cdot o\left(\frac{\Delta_{l}^{\star}}{|\lambda_{l}^{\star}|\sqrt{\mu\kappa^{6}r^{3}\log n}}\right)\\
& = o\left( \frac{\sigma_{\mathrm{max}}}{\Delta_{l}^{\star}}\sqrt{\mu\kappa^{6}r^{3}\log n} \cdot \frac{\Delta_{l}^{\star}}{|\lambda_{l}^{\star}|\sqrt{\mu\kappa^{6}r^{3}\log n}}\right) \\
& = o\left(\frac{\sigma_{\mathrm{max}}}{|\lambda_l^{\star}|} \right); \\
I_3 := \left|\sum_{k_2\neq l}\frac{\lambda_{l}^{\star}\lambda_{k_2}^{\star}}{\lambda_{l}^{2}}\left(\bm{u}_{l}^{\star\top}\bm{w}_{l}\right)\left(\bm{u}_{k_2}^{\star\top}\bm{u}_{l}\right)\cdot \frac{2\bm{u}_{l}^{\star\top}\bm{H}\bm{u}_{k_2}^{\star}}{\lambda_{l}}\right| & = o\left(\frac{\sigma_{\mathrm{max}}}{|\lambda_l^{\star}|} \right).
\end{align*}
Putting the preceding bounds together  yields
\begin{align}
	|{\tau}_{3,3}| \leq I_1 + I_2 + I_3 = o\left(\frac{\sigma_{\mathrm{max}}}{|\lambda_l^{\star}|} \right).
\end{align}

\paragraph{Controlling \texorpdfstring{${\tau}_{1,3}$}{} and \texorpdfstring{${\tau}_{3,4}$}{}.}
It can be deduced from \eqref{eq:rank-r-inference-lambda-l-bound} and  $|\bm{u}_l^{\star \top} \bm{u}_l| = 1 -o(1)$ that
\begin{align}
	|\tau_{1,3}| & \lesssim \sum_{k=1}^r \left|\frac{\lambda_{k}^{\star}}{\lambda_{l}^{\star}}\right|\cdot\left|\frac{\bm{u}_{k}^{\star\top}\bm{u}_{l}}{\bm{u}_{l}^{\star\top}\bm{u}_{l}}\right| \,  \sum_{s=2}^{\infty} \left|\frac{\bm{a}^\top \bm{H}^s \bm{u}_k^\star}{\lambda_l^s}\right| \lesssim \kappa \sum_{k=1}^{r}\sum_{s=2}^{\infty}\left|\frac{\bm{a}^{\top}\bm{H}^{s}\bm{u}_{k}^{\star}}{\lambda_{l}^{s}}\right|.
\end{align}
In addition, 
\begin{align}
\kappa\sum_{k=1}^r\sum_{s=2}^{\infty}\left|\frac{\bm{a}^{\top}\bm{H}^{s}\bm{u}_{k}^{\star}}{\lambda_{l}^{s}}\right| \lesssim
\frac{\kappa}{|\lambda_{l}^{\star}|}\sum_{k=1}^{r}\sum_{s=2}^{\infty}\left|\frac{2^{s}\bm{a}^{\top}\bm{H}^{s}\bm{u}_{k}^{\star}}{(\lambda_{\min}^{\star})^{s-1}}\right|,
\end{align}
where the last line follows since $|\lambda_l|\geq \lambda_{\min}^{\star} -\|\bm{H}\|\geq \lambda_{\min}^{\star} / 2$ (by the Bauer-Fike theorem and \eqref{eq:H-norm-simple-UB}). 	Applying Lemma \ref{lemma:perturbation-H} (or \eqref{eq:linear-form-Hs}) yields
\begin{align}
\frac{\kappa}{|\lambda_{l}^{\star}|}\sum_{k=1}^{r}\sum_{s=2}^{\infty}\left|\frac{2^{s}\bm{a}^{\top}\bm{H}^{s}\bm{u}_{k}^{\star}}{(\lambda_{\min}^{\star})^{s-1}}\right| & \leq\frac{\lambda_{\min}^{\star}\kappa}{|\lambda_{l}^{\star}|}\sqrt{\frac{\mu}{n}}\sum_{k=1}^{r}\sum_{s=2}^{\infty}\left(\frac{2c_{2}\sigma_{\max}\sqrt{n\log n}}{\lambda_{\min}^{\star}}\right)^{s}\leq\frac{\lambda_{\min}^{\star}\kappa r}{|\lambda_{l}^{\star}|}\sqrt{\frac{\mu}{n}}\frac{\left(\frac{2c_{2}\sigma_{\max}\sqrt{n\log n}}{\lambda_{\min}^{\star}}\right)^{2}}{1-\frac{2c_{2}\sigma_{\max}\sqrt{n\log n}}{\lambda_{\min}^{\star}}} \nonumber \\
& \asymp\frac{\lambda_{\min}^{\star}\kappa r}{|\lambda_{l}^{\star}|}\sqrt{\frac{\mu}{n}}\left(\frac{\sigma_{\max}\sqrt{n\log n}}{\lambda_{\min}^{\star}}\right)^{2}=o\left(\frac{\sigma_{\max}}{|\lambda_l^\star|}\right),
\end{align}
provided that $\sigma_{\max}\sqrt{\mu\kappa^{2}r^{2}n}\log n=o(\lambda_{\min}^{\star})$. The above bounds thus imply that
\begin{align}
|\tau_{1,3}|  = o\left(\frac{\sigma_{\max}}{|\lambda_l^\star|}\right).
\end{align}
Similarly, we can upper bound ${\tau}_{3,4}$ by
\begin{align}
|\tau_{3,4}| & \lesssim \kappa^2 \sum_{k_{1}=1}^{r}\sum_{k_{2}=1}^{r} \, \sum_{s_{1},s_{2}:s_{1}+s_{2}\geq 2}\left|\frac{\bm{u}_{k_{1}}^{\star\top}\bm{H}^{s_{1}+s_{2}}\bm{u}_{k_{2}}^{\star}}{\lambda_{l}^{s_{1}+s_{2}}}\right| \lesssim \frac{\kappa^2}{|\lambda_l^\star|}  \sum_{k_{1}=1}^{r}\sum_{k_{2}=1}^{r} \, \sum_{s_{1},s_{2}:s_{1}+s_{2}\geq 2} \left|\frac{2^{s_1+s_2} \bm{u}_{k_1}^{\star \top} \bm{H}^{s_1 + s_2} \bm{u}_{k_2^\star}}{(\lambda_{\mathrm{min}}^\star)^{s_1+s_2-1}}\right| \nonumber \\
& \leq \frac{\lambda_{\mathrm{min}}^\star \kappa^2 r^2}{|\lambda_l^\star|} \sqrt\frac{\mu}{n} \, \sum_{s_1,s_2:s_1 + s_2 \geq 2}\left(\frac{2c_{2}\sigma_{\max}\sqrt{n\log n}}{\lambda_{\min}^{\star}}\right)^{s_1 + s_2} \nonumber \\
& \stackrel{\mathrm{(i)}}{=} \frac{\lambda_{\mathrm{min}}^\star \kappa^2 r^2}{|\lambda_l^\star|} \sqrt\frac{\mu}{n} \left\{\left( \frac{1}{1-\frac{2c_{2}\sigma_{\max}\sqrt{n\log n}}{\lambda_{\min}^{\star}}} \right)^2 - 1 - 2\cdot \frac{2c_{2}\sigma_{\max}\sqrt{n\log n}}{\lambda_{\min}^{\star}} \right\} \nonumber \\
& \stackrel{\mathrm{(ii)}}{=} \frac{\lambda_{\mathrm{min}}^\star \kappa^2 r^2}{|\lambda_l^\star|} \sqrt\frac{\mu}{n} \cdot \left(\frac{2c_{2}\sigma_{\max}\sqrt{n\log n}}{\lambda_{\min}^{\star}}\right)^2 \cdot \frac{3 - \frac{4c_{2}\sigma_{\max}\sqrt{n\log n}}{\lambda_{\min}^{\star}}}{\left(1-\frac{2c_{2}\sigma_{\max}\sqrt{n\log n}}{\lambda_{\min}^{\star}}\right)^2} \nonumber \\
& \asymp\frac{\lambda_{\min}^{\star}\kappa^2 r^2}{|\lambda_{l}^{\star}|}\sqrt{\frac{\mu}{n}}\left(\frac{\sigma_{\max}\sqrt{n\log n}}{\lambda_{\min}^{\star}}\right)^{2}=o\left(\frac{\sigma_{\max}}{|\lambda_l^\star|}\right),
\end{align}
provided that $\sigma_{\max}\kappa^2 r^2\sqrt{n\mu}\log n=o(\lambda_{\min}^{\star}). $ Here, (i) and (ii)  make use of the elementary identity that $\sum_{s_1,s_2:s_1+s_2 \geq 2} a^{s_1+s_2} = \frac{1}{(1-a)^2} - 1 - 2a = \frac{3a^2(1-2a)}{(1-a)^2}$ for all $0 < a < 1$.

\paragraph{Controlling \texorpdfstring{$\tau_{3,2}$}{}.} Clearly,  one has $\tau_{3,2}=0$ if $r=1$. When $r\geq 2$, this term $\tau_{3,2}$ is the same term as $\mathcal{U}_1$ defined in \eqref{eq:estimator-squared-decompose-3-1}. Recalling our bound for inner product $|\bm{u}_k^{\star \top} \bm{w}_l|$ and $|\bm{u}_k^{\star \top} u_l|$ when $k \neq l$ from \eqref{eq:rank-r-inference-inner-prod-bound}, we deduce that
\begin{align}
|\tau_{3,2}| & \lesssim  \sum_{k: k \neq l} \left(\frac{\lambda_k^\star}{\lambda_l^\star}\right)^2 \cdot \left(\frac{\lambda_{l}^{\star}}{|\lambda_{l}^{\star}-\lambda_{k}^{\star}|}\cdot\frac{\sigma_{\mathrm{max}}}{\lambda_{\mathrm{min}}^{\star}}\sqrt{\mu\kappa^{2}r\log n}\right)^2 \leq r\left(\frac{\sigma_{\mathrm{max}} \kappa^2 \sqrt{\mu r \log n}}{|\Delta_l^\star|}\right)^2\nonumber \leq \left(\frac{\sigma_{\mathrm{max}} \kappa^2 r \sqrt{\mu \log n}}{|\Delta_l^\star|}\right)^2\nonumber
\end{align}

Putting all this together, we arrive at
\begin{align*}
|\tau_{1}| & \leq|\tau_{1,1}|+|\tau_{1,2}|+|\tau_{1,3}|=o\left(\frac{\sigma_{\mathrm{max}}}{|\lambda_{l}^{\star}|}\right)\\
|\tau_{3}| & \leq|\tau_{3,1}|+|\tau_{3,2}|+|\tau_{3,3}|+|\tau_{3,4}|=O\left(\frac{\sigma_{\mathrm{max}}\kappa^{2}r\sqrt{\mu\log n}}{|\Delta_{l}^{\star}|}\right)^{2}+o\left(\frac{\sigma_{\mathrm{max}}}{|\lambda_{l}^{\star}|}\right)
\end{align*}
as claimed. The bound on  $|{\tau}_2|$ can be established using exactly the same way as for $|\tau_{1}| $.

\subsection{Proof of Lemma \ref{lem:rank-r-inference-squared-estimator-difference-size}}
\label{PfLemmaRankRInf}
Direct calculation yields
\begin{align}
\big(\widehat{u}_{\bm{a},l}^{\mathsf{modified}}\big)^{2} - \left(\bm{a}^{\top}\bm{u}^{\star}_l + \frac{1}{2\lambda^{\star}_l}\bm{a}_l^{\perp \top}\big(\bm{H}+\bm{H}^{\top}\big)\bm{u}^{\star}_l\right)^2 & = \frac{\mathcal{R}}{\frac{\bm{u}_{l}^{\star\top}\bm{w}_{l}}{\lambda_{l}/\lambda_{l}^{\star}}\cdot\frac{\bm{u}_{l}^{\star\top}\bm{u}_{l}}{\lambda_{l}/\lambda_{l}^{\star}} \left(1 + \frac{\bm{u}_l^{\star \top} (\bm{H} + \bm{H}^\top) \bm{u}_l^\star}{\lambda_l^\star}\right)+\tau_{3}}, \label{eq:rank-r-inference-estimator-squared-difference-fraction}
\end{align}
where
\begin{align*}
\mathcal{R} & =\underbrace{\left\{ \frac{\bm{u}_{l}^{\star\top}\bm{u}_{l}}{\lambda_{l}/\lambda_{l}^{\star}}\cdot\left(\bm{a}^{\top}\bm{u}_{l}^{\star}+\frac{\bm{a}^{\top}\bm{H}\bm{u}_{l}^{\star}}{\lambda_{l}^{\star}}\right)+\tau_{1}\right\} \left\{ \frac{\bm{u}_{l}^{\star\top}\bm{w}_{l}}{\lambda_{l}/\lambda_{l}^{\star}}\cdot\left(\bm{a}^{\top}\bm{u}_{l}^{\star}+\frac{\bm{a}^{\top}\bm{H}^{\top}\bm{u}_{l}^{\star}}{\lambda_{l}^{\star}}\right)+\tau_{2}\right\} }_{=:\mathcal{R}_{1}}-\nonumber\\
 & \quad\underbrace{\left\{ \frac{\bm{u}_{l}^{\star\top}\bm{w}_{l}}{\lambda_{l}/\lambda_{l}^{\star}}\cdot\frac{\bm{u}_{l}^{\star\top}\bm{u}_{l}}{\lambda_{l}/\lambda_{l}^{\star}}\left(1+\frac{\bm{u}_{l}^{\star\top}(\bm{H}+\bm{H}^{\top})\bm{u}_{l}^{\star}}{\lambda_{l}^{\star}}\right)+\tau_{3}\right\} \cdot\left\{ \bm{a}^{\top}\bm{u}_{l}^{\star}+\frac{1}{2\lambda_{l}^{\star}}\bm{a}_{l}^{\perp\top}\big(\bm{H}+\bm{H}^{\top}\big)\bm{u}_{l}^{\star}\right\} ^{2}}_{=:\mathcal{R}_{2}}.
\end{align*}
Rearranging terms and using the above bounds, we can derive 
\begin{align*}
\mathcal{R}_{1} & =\frac{\bm{u}_{l}^{\star\top}\bm{w}_{l}}{\lambda_{l}/\lambda_{l}^{\star}}\cdot\frac{\bm{u}_{l}^{\star\top}\bm{u}_{l}}{\lambda_{l}/\lambda_{l}^{\star}}\cdot\left\{ \left(\bm{a}^{\top}\bm{u}_{l}^{\star}\right)^{2}+\left(\bm{a}^{\top}\bm{u}_{l}^{\star}\right)\cdot\frac{\bm{a}^{\top}(\bm{H}+\bm{H}^{\top})\bm{u}_{l}^{\star}}{\lambda_{l}^{\star}}+\frac{\bm{a}^{\top}\bm{H}\bm{u}_{l}^{\star}}{\lambda_{l}^{\star}}\cdot\frac{\bm{a}^{\top}\bm{H}^{\top}\bm{u}_{l}^{\star}}{\lambda_{l}^{\star}}\right\} \nonumber\\
 & \quad\quad+\bm{a}^{\top}\bm{u}_{l}^{\star}\cdot\left(\tau_{1}+\tau_{2}\right)+o\left(|\bm{a}^{\top}\bm{u}_{l}^{\star}|\cdot(|\tau_{1}|+|\tau_{2}))\right),\\
\mathcal{R}_{2} & =\frac{\bm{u}_{l}^{\star\top}\bm{w}_{l}}{\lambda_{l}/\lambda_{l}^{\star}}\cdot\frac{\bm{u}_{l}^{\star\top}\bm{u}_{l}}{\lambda_{l}/\lambda_{l}^{\star}}\cdot\Bigg\{(\bm{a}^{\top}\bm{u}_{l}^{\star})^{2}+\left(\bm{a}^{\top}\bm{u}_{l}^{\star}\right)\cdot\frac{\left(\bm{a}_{l}^{\perp}+(\bm{a}^{\top}\bm{u}_{l}^{\star})\bm{u}_{l}^{\star}\right)^{\top}(\bm{H}+\bm{H}^{\top})\bm{u}_{l}^{\star}}{\lambda_{l}^{\star}}\nonumber\\
 & \quad+\left(1+\frac{\bm{u}_{l}^{\star\top}(\bm{H}+\bm{H}^{\top})\bm{u}_{l}^{\star}}{\lambda_{l}^{\star}}\right)\cdot\frac{\left(\bm{a}_{l}^{\perp\top}(\bm{H}+\bm{H}^{\top})\bm{u}_{l}^{\star}\right)^{2}}{4(\lambda_{l}^{\star})^{2}}+\frac{\bm{u}_{l}^{\star\top}(\bm{H}+\bm{H}^{\top})\bm{u}_{l}^{\star}}{\lambda_{l}^{\star}}\cdot\frac{\bm{a}_{l}^{\perp\top}(\bm{H}+\bm{H}^{\top})\bm{u}_{l}^{\star}}{\lambda_{l}^{\star}}\Bigg\}\nonumber\\
 & \quad+(1+o(1))(\bm{a}^{\top}\bm{u}_{l}^{\star})^{2}\cdot\tau_{3}.
\end{align*}
By definition, we have $\bm{a}_l^{\perp} + (\bm{a}^\top \bm{u}_l^{\star}) \bm{u}_l^{\star} = \bm{a}$, 
which further gives
\begin{align*}
\left|\mathcal{R}\right| & =\left|\mathcal{R}_{1}-\mathcal{R}_{2}\right|\nonumber\\
 & =\Bigg|(1+o(1))\Bigg\{\frac{\bm{a}^{\top}\bm{H}\bm{u}_{l}^{\star}}{\lambda_{l}^{\star}}\cdot\frac{\bm{a}^{\top}\bm{H}^{\top}\bm{u}_{l}^{\star}}{\lambda_{l}^{\star}}-\left(1+\frac{\bm{u}_{l}^{\star\top}(\bm{H}+\bm{H}^{\top})\bm{u}_{l}^{\star}}{\lambda_{l}^{\star}}\right)\frac{\left(\bm{a}_{l}^{\perp\top}(\bm{H}+\bm{H}^{\top})\bm{u}_{l}^{\star}\right)^{2}}{4\lambda_{l}^{\star2}}-\nonumber\\
 & \quad\frac{\bm{u}_{l}^{\star\top}(\bm{H}+\bm{H}^{\top})\bm{u}_{l}^{\star}}{\lambda_{l}^{\star}}\cdot\frac{\bm{a}_{l}^{\perp\top}(\bm{H}+\bm{H}^{\top})\bm{u}_{l}^{\star}}{\lambda_{l}^{\star}}\Bigg\}+\bm{a}^{\top}\bm{u}_{l}^{\star}\left(\tau_{1}+\tau_{2}\right)+o\left(|\bm{a}^{\top}\bm{u}_{l}^{\star}|(|\tau_{1}|+|\tau_{2}))\right)-(1+o(1))(\bm{a}^{\top}\bm{u}_{l}^{\star})^{2}\tau_{3}\Bigg|\nonumber\\
 & \lesssim|\bm{a}^{\top}\bm{u}_{l}^{\star}|\cdot\left(|\tau_{1}|+|\tau_{2}|\right)+(\bm{a}^{\top}\bm{u}_{l}^{\star})^{2}\cdot|\tau_{3}|+ \frac{\sigma_{\mathrm{max}}^{2}\log n}{\lambda_{l}^{\star2}}
\end{align*}
with probability exceeding $1-O(n^{-6})$. 
Here, the last line uses Lemma \ref{lem:aHu-UB}. 
Substitution into \eqref{eq:rank-r-inference-estimator-squared-difference-fraction} yields
\begin{align*}
	& \left|\big(\widehat{u}_{\bm{a},l}^{\mathsf{modified}}\big)^{2}-\left(\bm{a}^{\top}\bm{u}_{l}^{\star}+\frac{1}{2\lambda_{l}^{\star}}\bm{a}_{l}^{\perp\top}\big(\bm{H}+\bm{H}^{\top}\big)\bm{u}_{l}^{\star}\right)^{2}\right| \leq\frac{|\mathcal{R}|}{1-o(1)} \\
	& \qquad \lesssim|\bm{a}^{\top}\bm{u}_{l}^{\star}|\cdot\left(|\tau_{1}|+|\tau_{2}|\right)+(\bm{a}^{\top}\bm{u}_{l}^{\star})^{2}\cdot|\tau_{3}|+\frac{\sigma_{\mathrm{max}}^{2}\log n}{\lambda_{l}^{\star2}}
\end{align*}
as claimed.


\subsection{Proof of Lemma \ref{lem:au-distribution-rankr}}
\label{sec:proof-lemma:au-distribution-rankr}

The first claim
	\begin{align}
		\left|\frac{\lambda_{l}}{\lambda_{l}^{\star}(\bm{u}_{l}^{\star\top}\bm{u}_{l})}\bm{a}^{\top}\bm{u}_{l}-\bm{a}^{\top}\bm{u}_{l}^{\star}-\frac{\bm{a}^{\top}\bm{H}\bm{u}_{l}^{\star}}{\lambda_{l}^{\star}}\right|=o\left(\frac{\sigma_{\max}}{|\lambda_{l}^{\star}|}\right) 
	\label{eq:UB103}
	\end{align}	
follows immediately from the bound on $|\tau_1|$ in Lemma~\ref{lem:rank-r-inference-remainder-size} and the facts $\lambda_l = (1+o(1)) \lambda_l^{\star}$ and $\bm{u}_l^{\top}\bm{u}_l^{\star}=1-o(1)$. As an immediate consequence, one has
\begin{align}
\left|\frac{\bm{a}^{\top}\bm{u}_{l}}{\bm{u}_{l}^{\top}\bm{u}_{l}^{\star}}\right| & \asymp\left|\frac{\lambda_{l}}{\lambda_{l}^{\star}}\frac{\bm{a}^{\top}\bm{u}_{l}}{\bm{u}_{l}^{\top}\bm{u}_{l}^{\star}}\right|\leq\left|\frac{\lambda_{l}}{\lambda_{l}^{\star}}\frac{\bm{a}^{\top}\bm{u}_{l}}{\bm{u}_{l}^{\top}\bm{u}_{l}^{\star}}-\bm{a}^{\top}\bm{u}_{l}^{\star}-\frac{\bm{a}^{\top}\bm{H}\bm{u}_{l}^{\star}}{\lambda_{l}^{\star}}\right|+\left|\bm{a}^{\top}\bm{u}_{l}^{\star}\right|+\left|\frac{\bm{a}^{\top}\bm{H}\bm{u}_{l}^{\star}}{\lambda_{l}^{\star}}\right| \nonumber\\
 & \lesssim o\Big(\frac{\sigma_{\max}}{|\lambda_{l}^{\star}|}\Big)+|\bm{a}^{\top}\bm{u}_{l}^{\star}|+\frac{\sigma_{\max}\sqrt{\log n}}{|\lambda_{l}^{\star}|}\asymp|\bm{a}^{\top}\bm{u}_{l}^{\star}|+\frac{\sigma_{\max}\sqrt{\log n}}{|\lambda_{l}^{\star}|}.
	\label{eq:aul-crude-bound123}
\end{align}

Next, the bound \eqref{eq:lambda-l-detailed-approximation}  tells us that
\begin{align*}
\frac{\lambda_{l}}{\lambda_{l}^{\star}}-1 & =\frac{\lambda_{l}-\lambda_{l}^{\star}}{\lambda_{l}^{\star}}=\frac{\bm{u}_{l}^{\star\top}\bm{H}\bm{u}_{l}^{\star}}{\lambda_{l}^{\star}}+o\Big(\frac{\sigma_{\max}}{|\lambda_{l}^{\star}|}\Big) .
\end{align*}
This in turn allows us to bound
\begin{align*}
\left|\left(\frac{\lambda_{l}}{\lambda_{l}^{\star}}-1\right)\frac{\bm{a}^{\top}\bm{u}_{l}}{\bm{u}_{l}^{\top}\bm{u}_{l}^{\star}}\right| & \lesssim\frac{\sigma_{\max}\sqrt{\log n}}{|\lambda_{l}^{\star}|}\left\{ |\bm{a}^{\top}\bm{u}_{l}^{\star}|+\frac{\sigma_{\max}\sqrt{\log n}}{|\lambda_{l}^{\star}|}\right\} =o\Big(\frac{\sigma_{\max}}{|\lambda_{l}^{\star}|}\Big),
\end{align*}
as long as $|\bm{a}^{\top}\bm{u}_{l}^{\star}|=o\left(1/\sqrt{\log n}\right)$.
Putting the above bounds together, we arrive at the advertised bound
\begin{align*}
 & \left|\frac{\bm{a}^{\top}\bm{u}_{l}}{\bm{u}_{l}^{\top}\bm{u}_{l}^{\star}}-\bm{a}^{\top}\bm{u}_{l}^{\star}-\frac{(\bm{a}_{l}^{\perp})^{\top}\bm{H}\bm{u}_{l}^{\star}}{\lambda_{l}^{\star}}\right|=\left|\frac{\lambda_{l}}{\lambda_{l}^{\star}}\frac{\bm{a}^{\top}\bm{u}_{l}}{\bm{u}_{l}^{\top}\bm{u}_{l}^{\star}}-\bm{a}^{\top}\bm{u}_{l}^{\star}-\frac{(\bm{a}_{l}^{\perp})^{\top}\bm{H}\bm{u}_{l}^{\star}}{\lambda_{l}^{\star}}-\left(\frac{\lambda_{l}}{\lambda_{l}^{\star}}-1\right)\frac{\bm{a}^{\top}\bm{u}_{l}}{\bm{u}_{l}^{\top}\bm{u}_{l}^{\star}}\right|\\
 & \quad\leq\left|\frac{\lambda_{l}}{\lambda_{l}^{\star}}\frac{\bm{a}^{\top}\bm{u}_{l}}{\bm{u}_{l}^{\top}\bm{u}_{l}^{\star}}-\bm{a}^{\top}\bm{u}_{l}^{\star}-\frac{(\bm{a}_{l}^{\perp})^{\top}\bm{H}\bm{u}_{l}^{\star}}{\lambda_{l}^{\star}}-\frac{\bm{u}_{l}^{\star\top}\bm{H}\bm{u}_{l}^{\star}}{\lambda_{l}^{\star}}\bm{a}^{\top}\bm{u}_{l}^{\star}\right|+O\left(\frac{\bm{u}_{l}^{\star\top}\bm{H}\bm{u}_{l}^{\star}}{\lambda_{l}^{\star}}\bm{a}^{\top}\bm{u}_{l}^{\star}\right)+o\Big(\frac{\sigma_{\max}}{|\lambda_{l}^{\star}|}\Big)\\
 & \quad\leq\left|\frac{\lambda_{l}}{\lambda_{l}^{\star}}\frac{\bm{a}^{\top}\bm{u}_{l}}{\bm{u}_{l}^{\top}\bm{u}_{l}^{\star}}-\bm{a}^{\top}\bm{u}_{l}^{\star}-\frac{\bm{a}^{\top}\bm{H}\bm{u}_{l}^{\star}}{\lambda_{l}^{\star}}\right|+O\left(\frac{\sigma_{\max}\sqrt{\log n}}{|\lambda_{l}^{\star}|}|\bm{a}^{\top}\bm{u}_{l}^{\star}|\right)+o\Big(\frac{\sigma_{\max}}{|\lambda_{l}^{\star}|}\Big)\\
 & \quad\leq o\Big(\frac{\sigma_{\max}}{|\lambda_{l}^{\star}|}\Big),
\end{align*}
where the penultimate inequality holds with the proviso that  $|\bm{a}^{\top}\bm{u}_{l}^{\star}|=o\left(1/\sqrt{\log n}\right)$.

The proof for the left eigenvector $\ul_l$ is essentially the same and is thus omitted for brevity.

\section{Proof for variance estimation (Theorem \ref{Lem-var-control-rankr})}
\label{sec:proof-lemma-var-control}

Without loss of generality, we assume $\|\bm{a}\|_2=1$ and $\lambda^{\star}_l=1$ in this section.

\subsection{Proof outline}
\paragraph{The estimation accuracy of $\vahat$.}
%
%
%
%
Before continuing, we note that under the assumption that $|\bm{a}^{\top}\bm{u}_{l}^{\star}| = o(\|\bm{a}\|_{2})=o(1)$, we have learned from Lemma~\ref{lem:va-properties} (or \eqref{eq:va-bounds-simple}) that
\begin{align}
	\tfrac{1}{2}\sigma_{\min}^{2}\|\bm{a}_{l}^{\perp}\|_{2}^{2}\leq v_{\bm{a},l}^{\star}\leq\sigma_{\max}^{2}\|\bm{a}_{l}^{\perp}\|_{2}\leq\sigma_{\max}^{2}\|\bm{a}\|_{2}^{2}.
\end{align}
In addition, the assumption $|\bm{a}^{\top}\bm{u}_l^{\star}|\leq (1-\epsilon) \|\bm{a}\|_2$ for some constant $0<\epsilon<1$ implies that 
\begin{align}
	1 = \|\bm{a}\|_2^2 \geq  \|\bm{a}_{l}^{\perp}\|_{2}^{2} = \|\bm{a}\|_{2}^{2}-(\bm{a}^{\top}\bm{u}_{l}^{\star})^{2} \geq 2\epsilon - \epsilon^2 > 0. 
\end{align}
Consequently, under the assumptions $\sigma_{\max}/\sigma_{\min} \asymp 1$ and $\epsilon \asymp 1$ one has
\begin{align}
	\label{eq:vstar-sigma-min2}
	 v_{\bm{a},l}^{\star} \asymp \sigma_{\max}^2 .
\end{align}

Recall the construction of our variance estimator
\begin{align}
\label{eq:brahms}
	\widehat{v}_{\bm{a},l}:=\frac{1}{4\lambda_{l}^{2}}\sum_{1\leq i,j\leq n}\big(\widehat{a}_{l,i}^{\perp}\widehat{u}_{l,j}+\widehat{a}_{l,j}^{\perp}\widehat{u}_{l,i}\big)^{2}\widehat{H}_{ij}^{2} ,
\end{align}
where $\widehat{\bm{a}}_l^{\perp} := \bm{a} - (\bm{a}^{\top}\widehat{\bm{u}}_l)\widehat{\bm{u}}_l$ and $\bm{a}_l^{\perp} := \bm{a} - (\bm{a}^{\top}\bm{u}_l^{\star}) \bm{u}_l^{\star} $. To control the size of $\vahat$,  we find it convenient to introduce a surrogate quantity
\begin{align}
\label{eq:beethoven}
\vatil \defn \frac{1}{4}\sum_{1\leq i,j\leq n}\left(a_{l,i}^{\perp}u_{l,j}^{\star}+a_{l,j}^{\perp}u_{l,i}^{\star}\right)^{2}\Hij^{2},
\end{align}
which turns out to be very close to our estimator $\vahat$. 
\begin{lemma}
\label{claim:spring-sonata}
	With probability exceeding $1-O(n^{-5})$, one has
	$\big|\vahat - \vatil \big| =  o(\vastar)$.
\end{lemma}
\begin{proof} See Appendix~\ref{sec:proof-claim-spring-sonata}. \end{proof}
\noindent This observation allows us to turn attention to bounding $\vatil$ instead. The result is this:
\begin{lemma}
\label{claim:shostakovich}
With probability exceeding $1-O(n^{-10})$, one has $|\vatil - \vastar|=  o(\vastar)$. 
\end{lemma}
\begin{proof} See Appendix~\ref{sec:proof-claim-shostakovich}. \end{proof}
\noindent Putting Lemmas~\ref{claim:spring-sonata}-\ref{claim:shostakovich} together and invoking the triangle inequality,
we establish the advertised bound
\begin{align}
\label{eq:weird-style}
	\vahat = \vastar + O( |\vahat - \vatil |+|\vatil - \vastar|) = (1+o(1)) \vastar .
\end{align}
The rest of this section is thus mainly dedicated to establishing Lemma~\ref{claim:spring-sonata} and Lemma~\ref{claim:shostakovich}.

\paragraph{The estimation accuracy of $\widehat{v}_{\lambda,l}$.} The proof for $\widehat{v}_{\lambda,l}=(1+o(1))v_{\lambda, l}^{\star}$ follows from a very similar argument, and is hence omitted.

\subsection{Proof of Lemma~\ref{claim:spring-sonata}}
\label{sec:proof-claim-spring-sonata}

Without loss of generality, assume $\langle{\bm{u}}_l, \bm{u}^{\star}_l\rangle >0$, which combined with Lemma~\ref{lem:basic-uR-uL-rankr} and Notation~\ref{sec:notation-eigen} gives
\begin{align}
	\langle{ \widehat{\bm{u}}}_l, \bm{u}^{\star}_l\rangle >0 \qquad \text{and} \qquad \langle{\bm{w}}_l, \bm{u}^{\star}_l\rangle >0. 
\end{align}
To prove this lemma, let us
first recall  that, under the assumptions stated in Theorem~\ref{thm:confidence-interval-validity-rankr} one has 
\begin{subequations}
\begin{align}
	\big\|\widehat{\bm{u}}_l-\bm{u}^{\star}_l\big\|_{\infty}  & \lesssim o\big(\tfrac{1}{\sqrt{\mu n\log n}}\big) , \qquad \|\widehat{\bm{u}}_{l}\|_{\infty} \leq 2\sqrt{\tfrac{\mu}{n}}, \label{eq:uhat-l-inf-norm-bound2}\\
	\big|\bm{a}^{\top}\widehat{\bm{u}}_{l}-\bm{a}^{\top}\bm{u}_{l}^{\star}\big| & \lesssim o\big(\tfrac{1}{\mu^2\log^2 n}\big) ,
	\label{eq:au-hat-bound3}
\end{align}
\end{subequations}
where \eqref{eq:uhat-l-inf-norm-bound2} comes from Lemma~\ref{lem:basic-uR-uL-rankr}, and \eqref{eq:au-hat-bound3} arises from \eqref{eq:rankr-linear-form-bound}  and \eqref{eq:au-sign-UB-125}.

We also need the following properties that control the difference between each term 
of $\vahat$ with that of $\vatil$. In particular, we introduce the following two lemmas.

\begin{lemma}
\label{lem:entrywise-error-SVD}
Suppose that Assumption~\ref{assumption:noise-size-rankr} holds and that $32\kappa^2 \sqrt{\mu r/n}\leq 1$. Then one has
\begin{align}
	\xi := \| \widehat{\bm{H}} - \bm{H}\|_{\infty}  = \| \bm{M}_{\mathsf{svd}}-\bm{M}^{\star}\|_{\infty} 
	\lesssim
	\sigma_{\max}\frac{\mu\kappa^{4}r\sqrt{\log n}}{\sqrt{n}}
	%
\end{align}
with probability exceeding $1-O(n^{-9})$. 
\end{lemma}
\begin{proof} See Appendix~\ref{sec:proof-lem:entrywise-error-SVD}. \end{proof}

\begin{lemma}
\label{LemA-diff}
With probability at least $1-O(n^{-8})$, for every $1\leq i,j \leq n$, one has 
\begin{align}
	|a_{l,i}^{\perp}| &\leq \big|a_{i}\big|+\big|u_{l,i}^{\star}\big| \leq  |a_i| + \sqrt{\mu/n} ; \\
	|\zeta_{ij}| &:= \left|a_{i}\widehat{u}_{l,j}+a_{j}\widehat{u}_{l,i} - (a_{i}u^{\star}_{l,j}+a_{j}u^{\star}_{l,i})\right| 
	\leq o\left(\frac{|a_{i}|+|a_{j}|+\big|u_{l,i}^{\star}\big|+\big|u_{l,j}^{\star}\big|}{\sqrt{\mu n\log n}}+\frac{1}{\mu n\log n}\right).\label{eq:a-diff}
\end{align} 	
\end{lemma}

\begin{proof}
Recalling that $\widehat{\bm{a}}_{l}^{\perp}=\bm{a}-(\bm{a}^{\top}\widehat{\bm{u}}_{l})\widehat{\bm{u}}_{l}$ and $\bm{a}_{l}^{\perp}=\bm{a}-(\bm{a}^{\top}\bm{u}_{l}^{\star})\bm{u}_{l}^{\star}$, we obtain as a consequence of inequalities~\eqref{eq:uhat-l-inf-norm-bound2} and \eqref{eq:au-hat-bound3} that 
\begin{align*}
\big|\widehat{a}_{l,i}^{\perp}-a_{l,i}^{\perp}\big| & =\big|\big[\bm{a}-(\bm{a}^{\top}\widehat{\bm{u}}_{l})\widehat{\bm{u}}_{l}-\bm{a}+(\bm{a}^{\top}\bm{u}_{l}^{\star})\bm{u}_{l}^{\star}\big]_{i}\big|=\big|(\bm{a}^{\top}\widehat{\bm{u}}_{l})\widehat{u}_{l,i}-(\bm{a}^{\top}\bm{u}_{l}^{\star})u_{l,i}^{\star}\big|\\
 & \leq\big|\bm{a}^{\top}\widehat{\bm{u}}_{l}\big|\cdot\big|\widehat{u}_{l,i}-u_{l,i}^{\star}\big|+\big|\bm{a}^{\top}\widehat{\bm{u}}_{l}-\bm{a}^{\top}\bm{u}_{l}^{\star}\big|\cdot\big|u_{l,i}^{\star}\big|\\
 & \leq o\Big(\tfrac{1}{\sqrt{\mu n\log n}}\Big)+o\Big(\tfrac{1}{\mu^{2}\log^{2}n}\Big)\cdot\sqrt{\tfrac{\mu}{n}}\\
 & =o\Big(\tfrac{1}{\sqrt{\mu n\log n}}\Big)
\end{align*}
\begin{align}
	\text{and}\qquad\big|a_{l,i}^{\perp}\big|\leq\big|a_{i}\big|+\big|\bm{a}^{\top}\bm{u}_{l}^{\star}\big|\cdot\big|u_{l,i}^{\star}\big|\leq\big|a_{i}\big|+\big|u_{l,i}^{\star}\big|\leq\big|a_{i}\big|+\sqrt{{\mu}/{n}}.	 
\end{align}
These in turn allow one to derive
\begin{align*}
 & \left|\widehat{a}_{l,i}^{\perp}\widehat{u}_{l,j}+\widehat{a}_{l,j}^{\perp}\widehat{u}_{l,i}-(a_{l,i}^{\perp}u_{j}^{\star}+a_{l,j}^{\perp}u_{i}^{\star})\right|\\
 & \qquad\leq\left|\widehat{a}_{l,i}^{\perp}\widehat{u}_{l,j}+\widehat{a}_{l,j}^{\perp}\widehat{u}_{l,i}-(a_{l,i}^{\perp}\widehat{u}_{l,j}+a_{l,j}^{\perp}\widehat{u}_{l,i})\right|+\left|a_{l,i}^{\perp}\widehat{u}_{l,j}+a_{l,j}^{\perp}\widehat{u}_{l,i}-(a_{l,i}^{\perp}u_{j}^{\star}+a_{l,j}^{\perp}u_{i}^{\star})\right|\\
 & \qquad\leq\left(\big|\widehat{a}_{l,i}^{\perp}-a_{l,i}^{\perp}\big|+\big|\widehat{a}_{l,j}^{\perp}-a_{l,j}^{\perp}\big|\right)\left(\big|u_{l,i}^{\star}\big|+\big|u_{l,j}^{\star}\big|+\|\widehat{\bm{u}}_{l}-\bm{u}_{l}^{\star}\|_{\infty}\right)+\left(\big|a_{l,i}^{\perp}\big|+\big|a_{l,j}^{\perp}\big|\right)\|\widehat{\bm{u}}_{l}-\bm{u}_{l}^{\star}\|_{\infty}\\
 & \qquad\leq o\left(\frac{1}{\sqrt{\mu n\log n}}\right)\left(\big|u_{l,i}^{\star}\big|+\big|u_{l,j}^{\star}\big|+o\left(\frac{1}{\sqrt{\mu n\log n}}\right)\right)+\left(\big|a_{i}\big|+\big|a_{j}\big|+\big|u_{l,i}^{\star}\big|+\big|u_{l,j}^{\star}\big|\right)\cdot o\left(\frac{1}{\sqrt{\mu n\log n}}\right)\\
 & \qquad=o\left(\frac{|a_{i}|+|a_{j}|+\big|u_{l,i}^{\star}\big|+\big|u_{l,j}^{\star}\big|}{\sqrt{\mu n\log n}}+\frac{1}{\mu n\log n}\right).
\end{align*}
%
\end{proof}
%

We are ready to establish the claim.
Invoking the above two lemmas yields 
\begin{align*}
\left|4\lambda_{l}^{2}\widehat{v}_{\bm{a},l}-4\widetilde{v}_{\bm{a},l}\right| 
& =\left|\sum\nolimits _{i,j}\big(a_{l,i}^{\perp}u_{l,j}^{\star}+a_{l,j}^{\perp}u_{l,i}^{\star}+\zeta_{ij}\big)^{2}\widehat{H}_{ij}^{2}-\sum\nolimits _{i,j}\big(a_{l,i}^{\perp}u_{l,j}^{\star}+a_{l,j}^{\perp}u_{l,i}^{\star}\big)^{2}H_{ij}^{2}\right|\\
& \leq \sum\nolimits _{i,j} 2 \big|\zeta_{ij}\big|\left(\big|a_{l,i}^{\perp}u_{l,j}^{\star}+a_{l,j}^{\perp}u_{l,i}^{\star}\big| + \big|\zeta_{ij}\big|\right)\left(|H_{ij}|+\xi\right)^{2}\\
 & \quad + \sum\nolimits _{i,j}2 \left(\big|a_{l,i}^{\perp}u_{l,j}^{\star}+a_{l,j}^{\perp}u_{l,i}^{\star}\big|+\big|\zeta_{ij}\big|\right)^{2}\xi\left(|H_{ij}|+\xi\right)\\
 & \overset{(\mathrm{i})}{\lesssim}\frac{1}{n\sqrt{\log n}}\sum\nolimits _{i,j}\left(a_{i}^{2}+a_{j}^{2}+\big|u_{l,i}^{\star}\big|^{2}+\big|u_{l,j}^{\star}\big|^{2}+\frac{1}{\mu n\log n}\right)\left(|H_{ij}|+\xi\right)^{2}\\
 & \quad+\frac{\mu}{n}\sum\nolimits _{i,j}\xi\left(a_{i}^{2}+a_{j}^{2}+\big|u_{l,i}^{\star}\big|^{2}+\big|u_{l,j}^{\star}\big|^{2}+\frac{1}{\mu n\log n}\right)\left(|H_{ij}|+\xi\right)\\
 & \lesssim\,\frac{1}{n\sqrt{\log n}}\Bigg\{\underset{=:\alpha_{1}}{\underbrace{\sum\nolimits _{i,j}\left(a_{i}^{2}+a_{j}^{2}+\big|u_{l,i}^{\star}\big|^{2}+\big|u_{l,j}^{\star}\big|^{2}+\frac{1}{\mu n\log n}\right)H_{ij}^{2}}}\Bigg\}\\
 & \qquad+\frac{1}{n\sqrt{\log n}}\underset{=:\alpha_{2}}{\underbrace{\sum\nolimits _{i,j}\left(a_{i}^{2}+a_{j}^{2}+\big|u_{l,i}^{\star}\big|^{2}+\big|u_{l,j}^{\star}\big|^{2}+\frac{1}{\mu n\log n}\right)\xi^{2}}}\\
 & \qquad+\frac{\mu}{n}\underset{=:\alpha_{3}}{\underbrace{\sum\nolimits _{i,j}\left(a_{i}^{2}+a_{j}^{2}+\big|u_{l,i}^{\star}\big|^{2}+\big|u_{l,j}^{\star}\big|^{2}+\frac{1}{\mu n\log n}\right)\xi\left(|H_{ij}|+\xi\right)}},
\end{align*}
where (i) follows since (according to \eqref{eq:a-diff} and the incoherence condition)
\begin{align*}
\begin{cases}
\big|a_{l,i}^{\perp}u_{j}^{\star}+a_{l,j}^{\perp}u_{i}^{\star}\big|+\big|\zeta_{ij}\big| & \lesssim\left(\big|a_{i}\big|+\big|a_{j}\big|+\big|u_{l,i}^{\star}\big|+\big|u_{l,j}^{\star}\big|+\sqrt{\frac{1}{\mu n\log n}}\right)\left(\big|u_{l,i}^{\star}\big|+\big|u_{l,j}^{\star}\big|+\sqrt{\frac{1}{\mu n\log n}}\right)\\
 & \lesssim\sqrt{a_{i}^{2}+a_{j}^{2}+\big|u_{l,i}^{\star}\big|^{2}+\big|u_{l,j}^{\star}\big|^{2}+\frac{1}{\mu n\log n}}\cdot\sqrt{\frac{\mu}{n}}\\
 & \lesssim\sqrt{\frac{\mu}{n}\left(a_{i}^{2}+a_{j}^{2}+\big|u_{l,i}^{\star}\big|^{2}+\big|u_{l,j}^{\star}\big|^{2}+\frac{1}{\mu n\log n}\right)};\\
\big|\zeta_{ij}\big|\left(\big|a_{l,i}^{\perp}u_{j}^{\star}+a_{l,j}^{\perp}u_{i}^{\star}\big|+\big|\zeta_{ij}\big|\right) & \lesssim\left(\big|a_{i}\big|+\big|a_{j}\big|+\big|u_{l,i}^{\star}\big|+\big|u_{l,j}^{\star}\big|+\sqrt{\frac{1}{\mu n\log n}}\right)^{2}\sqrt{\frac{1}{\mu n\log n}}\cdot\sqrt{\frac{\mu}{n}}\\
 & \lesssim\left(a_{i}^{2}+a_{j}^{2}+\big|u_{l,i}^{\star}\big|^{2}+\big|u_{l,j}^{\star}\big|^{2}+\frac{1}{\mu n\log n}\right)\frac{1}{n\sqrt{\log n}}.
\end{cases}
\end{align*}
%
We then control $\alpha_1$, $\alpha_2$ and $\alpha_3$ separately.
\begin{itemize}
\item
Regarding the term $\alpha_{1}$, it is first seen (using the assumption $\|\bm{a}\|_2^2=1$) that
\[
\mathbb{E}[\alpha_{1}]\lesssim\sum\nolimits _{i,j}\left(a_{i}^{2}+a_{j}^{2}+\big|u_{l,i}^{\star}\big|^{2}+\big|u_{l,j}^{\star}\big|^{2}+\frac{1}{\mu n\log n}\right)\sigma_{\max}^{2}\asymp n\sigma_{\max}^{2}.
\]
In addition,  Bernstein's inequality tells us that with probability exceeding $1-O(n^{-10})$, 
\begin{align*}
\left|\alpha_{1}-\mathbb{E}[\alpha_{1}]\right| & \lesssim\sqrt{\sum_{i,j}\mathsf{Var}\left[\left(a_{i}^{2}+a_{j}^{2}+\frac{\mu}{n}\right)H_{ij}^{2}\right]\log n}+\max_{i,j}\left\{ \left(a_{i}^{2}+a_{j}^{2}+\frac{\mu}{n}\right)H_{ij}^{2}\right\} \log n\\
 & \lesssim\sqrt{\sum_{i,j}\left(a_{i}^{2}+a_{j}^{2}+\frac{\mu}{n}\right)^{2}B^{2}\sigma_{\max}^{2}\log n}+\max_{i,j}\left(a_{i}^{2}+a_{j}^{2}+\frac{\mu}{n}\right)B^{2}\log n\\
 & \lesssim\sqrt{\sum_{i,j}\left\{ \left(a_{i}^{2}+a_{j}^{2}\right)^{2}+\left(\frac{\mu}{n}\right)^{2}\right\} B^{2}\sigma_{\max}^{2}\log n}+B^{2}\log n\\
 & \lesssim\sqrt{\sum_{i,j}\left\{ a_{i}^{2}+a_{j}^{2}+\left(\frac{\mu}{n}\right)^{2}\right\} B^{2}\sigma_{\max}^{2}\log n}+B^{2}\log n\\
 & \lesssim B\sigma_{\max}\sqrt{n\log n}+B^{2}\log n\lesssim\sigma_{\max}^{2}n,
\end{align*}
where the fourth line follows since $a_i^2+a_j^2 \leq\|\bm{a}\|_{2}^2=1$,
and the last inequality comes from Assumption \ref{assumption:noise-size-rankr}.
Consequently, 
\begin{align}
\label{eq:alpha-UB5}
\alpha_{1}\leq\left|\alpha_{1}-\mathbb{E}[\alpha_{1}]\right|+\mathbb{E}[\alpha_{1}]\lesssim\sigma_{\max}^{2}n.
\end{align}

\item By virtue of Lemma \ref{lem:entrywise-error-SVD} and the assumption
$\|\bm{a}\|_{2}=1$, the second term $\alpha_{2}$ can be upper bounded
by 
\begin{align*}
\alpha_{2} & =\sum\nolimits _{i,j}\left(a_{i}^{2}+a_{j}^{2}+\big|u_{l,i}^{\star}\big|^{2}+\big|u_{l,j}^{\star}\big|^{2}+\frac{1}{\mu n\log n}\right)\xi^{2}=2n\left\{ \|\bm{a}\|_{2}^{2}+\|\bm{u}_{l}^{\star}\|_{2}^{2}+\frac{1}{\mu\log n}\right\} \xi^{2}\\
 & \asymp n \xi^{2} \lesssim\sigma_{\max}^{2}\mu^{2}\kappa^{8}r^{2}\log n.
\end{align*}

\item When it comes to $\alpha_{3}$, we first make the observation that
\begin{align*}
\xi\left(|H_{ij}|+\xi\right) & \lesssim\begin{cases}
\mu\xi^{2}\log n & \text{if }|H_{ij}|\lesssim\xi\mu\log n\\
\frac{1}{\mu\log n}\left(|H_{ij}|+\xi\right)^{2} & \text{if }|H_{ij}|\gtrsim\xi\mu\log n
\end{cases}\\
 & \lesssim\mu\xi^{2}\log n+\frac{1}{\mu\log n}\left(|H_{ij}|+\xi\right)^{2}\ind\{|H_{ij}|\gtrsim\xi\mu\log n\}\\
 & \asymp \mu\xi^{2}\log n+\frac{1}{\mu\log n}H_{ij}^{2}.
\end{align*}
As a consequence, one can derive
\begin{align*}
\alpha_{3} & \lesssim\mu\log n\underset{=\alpha_{2}}{\underbrace{\sum\nolimits _{i,j}\left(a_{i}^{2}+a_{j}^{2}+\big|u_{l,i}^{\star}\big|^{2}+\big|u_{l,j}^{\star}\big|^{2}+\frac{1}{\mu n\log n}\right)\xi^{2}}}  \\
 & \qquad + \frac{1}{\mu\log n}\underset{=\alpha_{1}}{\underbrace{\sum\nolimits _{i,j}\left(a_{i}^{2}+a_{j}^{2}+\big|u_{l,i}^{\star}\big|^{2}+\big|u_{l,j}^{\star}\big|^{2}+\frac{1}{\mu n\log n}\right)H_{ij}^{2}}}\\
 & =\alpha_{2}\mu\log n+\frac{1}{\mu\log n}\alpha_{1}.
\end{align*}
%
\end{itemize}
Putting the above bounds together, we conclude that
\begin{align*}
\left|4\lambda^{2}_l \widehat{v}_{\bm{a},l}-4\lambda^{\star2}_l \widetilde{v}_{\bm{a},l}\right| & \lesssim\frac{1}{n\sqrt{\log n}}\alpha_{1}+\frac{1}{n\sqrt{\log n}}\alpha_{2}+\frac{\mu}{n}\left(\alpha_{2}\mu\log n+\frac{1}{\mu\log n}\alpha_{1}\right)\\
 & \asymp\frac{1}{n\sqrt{\log n}}\alpha_{1}+\alpha_{2}\frac{\mu^{2}\log n}{n}
	 \lesssim \left( \frac{1}{\sqrt{\log n}}+ \frac{ \mu^4\kappa^{8}r^2 \log^2 n }{n} \right) \sigma_{\max}^2   \\
	& = o(\sigma_{\max}^{2}),
\end{align*}
provided that $\mu^{4} \kappa^8 r^2 \log^{2}n=o(n)$. This combined with the fact that $\lambda_l = (1+o(1)) \lambda^{\star}_l$ (cf.~Corollary~\ref{cor:eigengap-condition-bounds}) and the inequality~\eqref{eq:vstar-sigma-min2} gives
\[
	\widehat{v}_{\bm{a},l} - \widetilde{v}_{\bm{a},l} = o(\vastar). 
\]

\subsection{Proof of Lemma~\ref{claim:shostakovich}}
\label{sec:proof-claim-shostakovich}

The randomness in the quantity $\vatil$ purely comes from $\Hij$. Given that $\vatil$ is the sum of independent random variables, it is easily seen that $\Exs [\vatil] = \vastar$.
Invoke Bernstein's inequality \cite{tropp2015introduction} to guarantee that with probability at least $1- O(n^{-10})$,  
\begin{align*}
\big|4(\vatil-\vastar)\big| & =\big|4(\vatil-\Exs[\vatil])\big|=\Big|\sum_{1\leq i,j\leq n}\left(a_{l,i}^{\perp}u_{l,j}^{\star}+a_{l,j}^{\perp}u_{l,i}^{\star}\right)^{2}(\Hij^{2}-\Exs[\Hij^{2}])\Big|\\
 & \lesssim\sqrt{\sum_{i,j}\var\left[(a_{l,i}^{\perp}u_{l,j}^{\star}+a_{l,j}^{\perp}u_{l,i}^{\star})^{2}H_{ij}^{2}\right]\log n}+\max_{i,j}\left\{ \left|(a_{l,i}^{\perp}u_{l,j}^{\star}+a_{l,j}^{\perp}u_{l,i}^{\star})^{2}B^{2}\right|\right\} \log n\\
 & \stackrel{(\mathrm{i})}{\lesssim}\sqrt{\sum_{i,j}(a_{l,i}^{\perp})^{4}(u_{l,j}^{\star})^{4}\sigma_{\max}^{2}B^{2}\log n}+\max_{i,j}(a_{l,i}^{\perp}u_{l,j}^{\star})^{2}B^{2}\log n\\
 & \lesssim\sqrt{\Big\{\max_{i,j}(a_{l,i}^{\perp})^{2}(u_{l,j}^{\star})^{2}\Big\}\sum_{i,j}(a_{l,i}^{\perp}u_{l,j}^{\star})^{2}\sigma_{\max}^{2}B^{2}\log n}+\max_{i,j}(a_{l,i}^{\perp}u_{l,j}^{\star})^{2}B^{2}\log n\\
 & \stackrel{(\mathrm{ii})}{\lesssim}B\sigma_{\max}\sqrt{\frac{\mu\log n}{n}}+\frac{\mu B^{2}\log n}{n}\stackrel{(\mathrm{iii})}{=}o\big(\sigma_{\max}^{2}\big)=o(\vastar),
\end{align*}
where the inequality (i) uses the fact that $\var[\Hij^2] \leq \Exs[\Hij^4] \leq B^2 \Exs[\Hij^2] \leq B^2 \sigma_{\max}^2$,  (ii) follows since $\max_{i,j} (a_{l,i}^{\perp}u^\star_{j})^2 \leq \linf{\ustar}^2\leq \mu/n$ and $\|\bm{a}_l^{\perp}\|_2^2\leq \|\bm{a}\|_2^2=\|\ustar\|_2^2=1$,  (iii) arises from Assumption \ref{assumption:noise-size-rankr}, and the last identity follows from \eqref{eq:vstar-sigma-min2}. 
Combining the above pieces, we obtain $|\vatil - \vastar|=  o(\vastar)$ as claimed.


\subsection{Proof of Lemma~\ref{lem:entrywise-error-SVD}}
\label{sec:proof-lem:entrywise-error-SVD}

For notational convenience, we denote by $\bm{M}_{\mathsf{svd}}=\bm{U}_{\mathsf{svd}} \bm{\Sigma}_{\mathsf{svd}} \bm{V}_{\mathsf{svd}}^{\top}$ the compact SVD of $\bm{M}_{\mathsf{svd}}$ (or equivalently, the rank-$r$ SVD of $\bm{M}$). We shall also define the rotation matrix
\begin{equation}
	\bm{Q}:= \underset{\bm{R}\in\mathcal{O}^{r\times r}}{\arg\min} ~\Bigg\|\left[\begin{array}{c}
\bm{U}_{\mathsf{svd}}\\
\bm{V}_{\mathsf{svd}}
\end{array}\right]\bm{R}-\left[\begin{array}{c}
\bm{U}^{\star}\\
\bm{U}^{\star}
\end{array}\right]\Bigg\|_{\mathrm{F}},\label{eq:defn-Q-Usvd}
\end{equation}
where $\mathcal{O}^{r\times r}$ denotes the set of $r\times r$ orthonormal
matrices.

We first record a perturbation bound regarding the singular subspace of $\bm{M}_{\mathsf{svd}}$. 
\begin{lemma}
\label{lem:Usvd-Q-2inf}
Instate the assumptions of Lemma~\ref{lem:entrywise-error-SVD}. With probability exceeding $1-O(n^{-9})$ one has 
\begin{align}
\label{eq:Usvd-Vsvd-2bound2}
\max\big\{ \|\bm{U}_{\mathsf{svd}}\bm{Q}-\bm{U}^{\star}\|_{2,\infty}, \|\bm{V}_{\mathsf{svd}}\bm{Q}-\bm{U}^{\star}\|_{2,\infty} \big\} 
 \lesssim \kappa^2 \gamma  \sqrt{\frac{\mu r}{n}},
\end{align}
where $\gamma \asymp \frac{\sigma_{\max}\sqrt{n\log n}}{\lambda_{\min}^{\star}}$. 
\end{lemma}
\noindent As immediate consequences of Lemma~\ref{lem:Usvd-Q-2inf} and Assumption~\ref{assumption:noise-size-rankr}, we have
\begin{align}
\label{eq:Usvd-Vsvd-2bound}
\begin{cases}
\|\bm{U}_{\mathsf{svd}}\|_{2,\infty} & \leq\|\bm{U}_{\mathsf{svd}}\bm{Q}-\bm{U}^{\star}\|_{2,\infty}+\|\bm{U}^{\star}\|_{2,\infty}\lesssim \kappa \sqrt{\frac{\mu r}{n}},\\
\|\bm{V}_{\mathsf{svd}}\|_{2,\infty} & \leq\|\bm{V}_{\mathsf{svd}}\bm{Q}-\bm{V}^{\star}\|_{2,\infty}+\|\bm{V}^{\star}\|_{2,\infty}\lesssim \kappa \sqrt{\frac{\mu r}{n}},
\end{cases}
\end{align}
provided that $\kappa^2 \sqrt{\mu r/n}\leq 1$.

We are now ready to control the entrywise error of $\bm{M}_{\mathsf{svd}}$. Towards this, we make note of the following bound
\begin{align}
\|\ensuremath{\bm{M}_{\mathsf{svd}}}-\bm{M}^{\star}\|_{\infty} & \lesssim\|\bm{Q}^{\top}\bm{\Sigma}_{\mathsf{svd}}\bm{Q}-\bm{\Lambda}^{\star}\|\cdot\|\bm{U}_{\mathsf{svd}}\|_{2,\infty}\|\bm{V}_{\mathsf{svd}}\|_{2,\infty}\nonumber \\
 & \quad+\|\bm{\Lambda}^{\star}\|\left(\|\bm{U}_{\mathsf{svd}}\|_{2,\infty}+\|\bm{V}^{\star}\|_{2,\infty}\right)\left(\|\bm{U}_{\mathsf{svd}}\bm{Q}-\bm{U}^{\star}\|_{2,\infty}+\|\bm{U}_{\mathsf{svd}}\bm{Q}-\bm{U}^{\star}\|_{2,\infty}\right),
\label{eq:Msvd-decomp-1}
\end{align}
which can be obtained by combining the inequalities (C.17)-(C.18) in \cite{abbe2017entrywise}. In addition, following the argument in \cite[Appendix~C.3.3]{abbe2017entrywise}, we obtain
\begin{equation}
\label{eq:Q-Sigma-bound}
\|\bm{Q}^{\top}\bm{\Sigma}_{\mathsf{svd}}\bm{Q}-\bm{\Lambda}^{\star}\|\lesssim\gamma\lambda_{\min}^{\star},
\end{equation}
where $\gamma$ is defined in \eqref{eq:defn-gamma}. Substituting \eqref{eq:Usvd-Vsvd-2bound2}, \eqref{eq:Usvd-Vsvd-2bound} and \eqref{eq:Q-Sigma-bound} into \eqref{eq:Msvd-decomp-1} and combining terms, we reach
\begin{align}
\|\ensuremath{\bm{M}_{\mathsf{svd}}}-\bm{M}^{\star}\|_{\infty}
\lesssim & ~\gamma\cdot\frac{\mu\kappa^{2}r}{n}\lambda_{\min}^{\star}+\gamma\cdot\frac{\mu\kappa^{4}r}{n}\lambda_{\min}^{\star}\asymp\gamma\cdot\frac{\mu\kappa^{4}r}{n}\lambda_{\min}^{\star} \nonumber\\
\asymp & ~\sigma_{\max}\sqrt{\frac{\mu^{2}\kappa^{8}r^2\log n}{n}}.
\label{eq:M-Msvd-bound}
\end{align}
Finally, it follows from our construction that
$\widehat{\bm{H}}-\bm{H}= \bm{M}-\bm{M}_{\mathsf{svd}}-\bm{H} = \bm{M}^{\star}- \bm{M}_{\mathsf{svd}}$, which combined with \eqref{eq:M-Msvd-bound} establishes this lemma.

%

%
%

\begin{proof}[Proof of Lemma \ref{lem:Usvd-Q-2inf}]
In order to invoke \cite[Theorem 2.1]{abbe2017entrywise} and the symmetric dilation trick in \cite[Section 3.3]{abbe2017entrywise}, we need to first verify the assumptions required in \cite[Section 2.1]{abbe2017entrywise}. To this end, we introduce the following auxiliary quantity and function
\begin{subequations}
\begin{align}
\label{eq:defn-gamma}
\gamma &:= c_{\gamma}\frac{\sigma_{\max}\sqrt{n\log n}}{\lambda_{\min}^{\star}} \\
%
\label{eq:defn-varphi}
\varphi(x) &:=\begin{cases}
c_{\varphi}\left(\frac{\sigma_{\max}\sqrt{n\log n}}{\lambda_{\min}^{\star}}x+\frac{B\log n}{\lambda_{\min}^{\star}}\right),\quad & \text{if }\frac{1}{\sqrt{n}}\leq x\leq1\\
c_{\varphi}\left(\frac{\sigma_{\max}\sqrt{n\log n}}{\lambda_{\min}^{\star}}+\frac{B\sqrt{n}\log n}{\lambda_{\min}^{\star}}\right)x, & \text{if }0\leq x<\frac{1}{\sqrt{n}}
\end{cases}
\end{align}
\end{subequations}
for some sufficiently large constants $c_{\gamma},c_{\varphi}>0$. 
We then make the following observations:
\begin{itemize}
\item To begin with, the two-to-infinity norm of the truth can be bounded by
\begin{align}
\label{eq:Mstar-2-inf-UB}
\|\bm{M}^{\star}\|_{2,\infty}=\|\bm{U}^{\star}\bm{\Lambda}^{\star}\bm{U}^{\star\top}\|_{2,\infty}\leq\|\bm{U}^{\star}\|_{2,\infty}\cdot\|\bm{\Lambda}^{\star}\|\cdot\|\bm{U}^{\star}\|\leq\lambda_{\max}^{\star}\sqrt{\frac{\mu r}{n}}
= \lambda_{\min}^{\star}\sqrt{\frac{\mu \kappa^2 r}{n}} .
\end{align}
\item In view of the matrix Bernstein inequality \cite{tropp2015introduction},  with probability exceeding $1-O(n^{-6})$ one has
\begin{equation}
\|\bm{M}-\bm{M}^{\star}\|=\|\bm{H}\|\lesssim\sigma_{\max}\sqrt{n\log n}=\frac{\sigma_{\max}\sqrt{n\log n}}{\lambda_{\min}^{\star}}\lambda_{\min}^{\star} \leq \gamma \lambda_{\min}^{\star}.
\end{equation}
\item
Denote by $\bm{H}_{i,\cdot}$ the $i$th row of $\bm{H}$. For any
fixed $\bm{W}\in\mathbb{R}^{n\times r}$, the matrix Bernstein inequality yields
\begin{align}
\|\bm{H}_{i,\cdot}\bm{W}\|_{2} & \lesssim\sqrt{\sum_{j=1}^{n}\sigma_{i,j}^{2}\|\bm{W}_{j,\cdot}\|_{2}^{2}\log n}+B\|\bm{W}\|_{2,\infty}\log n\lesssim\sigma_{\max}\|\bm{W}\|_{\mathrm{F}}\sqrt{\log n}+B\|\bm{W}\|_{2,\infty}\log n \nonumber\\
 & =\lambda_{\min}^{\star}\|\bm{W}\|_{2,\infty}\left\{ \frac{\|\bm{W}\|_{\mathrm{F}}}{\sqrt{n}\|\bm{W}\|_{2,\infty}}\frac{\sigma_{\max}\sqrt{n\log n}}{\lambda_{\min}^{\star}}+\frac{B\log n}{\lambda_{\min}^{\star}}\right\} 
\label{eq:H-W-2-norm}
\end{align}
with probability exceeding $1-O(n^{-10})$. This combined with \eqref{eq:defn-varphi} and $\frac{1}{\sqrt{n}}\leq \frac{\|\bm{W}\|_{\mathrm{F}}}{\sqrt{n}\|\bm{W}\|_{2,\infty}} \leq 1$ gives
\begin{align}
\mathbb{P}\left\{ \|\bm{H}_{i,\cdot}\bm{W}\|_{2}\geq\lambda_{\min}^{\star}\|\bm{W}\|_{2,\infty}\varphi\left(\frac{\|\bm{W}\|_{\mathrm{F}}}{\sqrt{n}\|\bm{W}\|_{2,\infty}}\right)\right\} \geq1-O\big(n^{-10}\big).
\end{align}
\end{itemize}

With the above observations in place,  \cite[Theorem 2.1]{abbe2017entrywise} together with the dilation trick in \cite[Section 3.3]{abbe2017entrywise} implies that:  with probability at least $1-O(n^{-5})$, one has
\begin{align}
\|\bm{U}_{\mathsf{svd}}\bm{Q}-\bm{U}^{\star}-\bm{H}\bm{U}^{\star}(\bm{\Lambda}^{\star})^{-1}\|_{2,\infty} & =\|\bm{U}_{\mathsf{svd}}\bm{Q}-\bm{M}\bm{U}^{\star}(\bm{\Lambda}^{\star})^{-1}\|_{2,\infty}\lesssim\kappa^{2}\gamma\|\bm{U}^{\star}\|_{2,\infty}+\gamma\|\bm{M}^{\star}\|_{2,\infty}/\lambda_{\min}^{\star}\nonumber\\
 & \lesssim\gamma\kappa^{2}\sqrt{\frac{\mu r}{n}}+\gamma\sqrt{\frac{\mu\kappa^{2}r}{n}}\asymp\gamma\sqrt{\frac{\mu\kappa^{4}r}{n}},\label{eq:Usvd-10}
\end{align}
%
%
where the last line arises from the incoherence condition as well as the inequality \eqref{eq:Mstar-2-inf-UB}. 
In addition, 
\begin{align*}
\|\bm{U}^{\star}(\bm{\Lambda}^{\star})^{-1}\|_{2,\infty} & \leq\|\bm{U}^{\star}\|_{2,\infty}\|(\bm{\Lambda}^{\star})^{-1}\|\leq\frac{1}{\lambda_{\min}^{\star}}\sqrt{\frac{\mu r}{n}},\\
\|\bm{U}^{\star}(\bm{\Lambda}^{\star})^{-1}\|_{\mathrm{F}} & \leq\|\bm{U}^{\star}\|_{\mathrm{F}}\|(\bm{\Lambda}^{\star})^{-1}\|\leq\frac{\sqrt{r}}{\lambda_{\min}^{\star}},
\end{align*}
which taken collectively with \eqref{eq:H-W-2-norm} (by setting $\bm{W}=\bm{U}^{\star}(\bm{\Lambda}^{\star})^{-1}$) demonstrate that
\begin{align}
\|\bm{H} \bm{U}^{\star}(\bm{\Lambda}^{\star})^{-1}\|_{2,\infty} & \lesssim\frac{\sigma_{\max}\sqrt{r\log n}}{\lambda_{\min}^{\star}}+\frac{B}{\lambda_{\min}^{\star}}\sqrt{\frac{\mu r\log^{2}n}{n}}
\lesssim \frac{\sigma_{\max}\sqrt{\mu r\log n}}{\lambda_{\min}^{\star}} \label{eq:H-W2-UB}
\end{align}
with probability $1-O(n^{-9})$. 
Taking it together with \eqref{eq:Usvd-10} and using the triangle inequality immediately yield
\begin{align*}
\|\bm{U}_{\mathsf{svd}}\bm{Q}-\bm{U}^{\star}\|_{2,\infty} & \leq\|\bm{U}_{\mathsf{svd}}\bm{Q}-\bm{U}^{\star}-\bm{H}\bm{U}^{\star}(\bm{\Lambda}^{\star})^{-1}\|_{2,\infty}+\|\bm{H}\bm{U}^{\star}(\bm{\Lambda}^{\star})^{-1}\|_{2,\infty}\\
 & \lesssim\gamma\sqrt{\frac{\mu\kappa^{4}r}{n}}+\frac{\sigma_{\max}\sqrt{\mu r\log n}}{\lambda_{\min}^{\star}}\asymp\gamma\sqrt{\frac{\mu\kappa^{4}r}{n}}
\end{align*}
as claimed. The part concerning $\bm{V}_{\mathsf{svd}}$ follows from the same argument. 
\end{proof}
\begin{remark}
	The careful reader would note that the results in \cite{abbe2017entrywise} require another relation between the incoherence condition and $\gamma$ (namely, Assumption A1). This relation cannot be satisfied with the current choice of $\gamma$ given in \eqref{eq:defn-gamma}, unless we  increase it to
	\begin{align}
		\gamma = c_{\gamma}\frac{\sigma_{\max}\sqrt{n\log n}}{\lambda_{\min}^{\star}} + \sqrt{\frac{\mu\kappa^2 r}{n}}.
		\label{eq:defn-gamma-large}
	\end{align}
	Fortunately, the second term of \eqref{eq:defn-gamma-large} can be removed with slightly more refined analysis, provided that Assumption~\ref{assumption:noise-size-rankr} holds. In short, the analysis of \cite{abbe2017entrywise} is built upon a sequence of auxiliary leave-one-out estimates $\bm{M}^{(m)}$ ($1\leq m\leq n$) --- obtained by zeroing out the $m$th row/column of  $\bm{M}$ --- which allows us to approximate $\bm{M}$ while being statistically independent of the data in the $m$th row/column. However, the approximation error $\bm{M}^{(m)}-\bm{M}$ does not decrease to zero even if $\bm{H} = \bm{0}$, leading to a bias term that is non-vanishing as $\sigma_{\max}\rightarrow 0$.  To address this issue, it suffices to replace $\bm{M}^{(m)}$ by $\widetilde{\bm{M}}^{(m)}$, where $\widetilde{\bm{M}}^{(m)}$ is obtained by replacing all entries in the $m$th row/column of $\bm{M}$ by their expected values.  This allows us to ensure that   $\widetilde{\bm{M}}^{(m)}-\bm{M}\rightarrow \bm{0}$ as $\sigma_{\max}\rightarrow 0$, which in turn leads to the removal of the second term of  \eqref{eq:defn-gamma-large}.  The refined proof is nearly identical to the original proof in \cite{abbe2017entrywise}, and is hence omitted here for the sake of brevity. 
\end{remark}

\section{Proof for minimax lower bounds (Theorem~\ref{thm:minimax-l2-au-eigengap}) }
\label{sec:proof-minimax-lower-bounds}

In the following, we prove each lower bound repectively. 

\subsection{Proof of Eqs.~(\ref{eq:minimax-L2-u}) and (\ref{eq:minimax-au-eigen-gap})}

Without loss of generality, consider any $1\leq l\leq r$ and any
$k\neq l$. Consider the following hypotheses regarding the eigen-decomposition
of $\bm{M}^{\star}\in\mathbb{R}^{n\times n}$: 
\begin{align*}
\mathcal{H}_{0}: & \quad\bm{M}=\bm{M}^{\star}+\bm{H}=\sum_{j=1}^{r}\lambda_{j}^{\star}\bm{u}_{j}^{\star}\bm{u}_{j}^{\star\top}+\bm{H};\\
\mathcal{H}_{k}: & \quad\bm{M}=\bm{M}_{k}+\bm{H}=\lambda_{l}^{\star}\widetilde{\bm{u}}_{l}\widetilde{\bm{u}}_{l}^{\top}+\lambda_{k}^{\star}\widetilde{\bm{u}}_{k}\widetilde{\bm{u}}_{k}^{\top}+\sum_{j:j\neq k,j\neq l}\lambda_{j}^{\star}\bm{u}_{j}^{\star}\bm{u}_{j}^{\star\top}+\bm{H}.
\end{align*}
In words, $\mathcal{H}_{k}$ is obtained by perturbing the $l$th
and the $k$th eigenvectors under $\mathcal{H}_{0}$, with the remaining
eigenvectors unaltered. In particular, for any $k\neq l$, we shall
pick $\widetilde{\bm{u}}_{l}$ and $\widetilde{\bm{u}}_{k}$ such
that they are equivalent to $\bm{u}_{l}^{\star}$ and $\bm{u}_{k}^{\star}$
modulo global rotation, namely,
\begin{equation}
\bm{u}_{l}^{\star}\bm{u}_{l}^{\star\top}+\bm{u}_{k}^{\star}\bm{u}_{k}^{\star\top}=\widetilde{\bm{u}}_{l}\widetilde{\bm{u}}_{l}^{\top}+\widetilde{\bm{u}}_{k}\widetilde{\bm{u}}_{k}^{\top};\label{eq:constraint-u1-u2-hypotheses}
\end{equation}
as we shall see, this rotational invariance constraint (\ref{eq:constraint-u1-u2-hypotheses})
plays a pivotal role in understanding the effect of the eigen-gap
upon estimation accuracy. In what follows, we let $\mathbb{P}_{0}$
(resp.~$\mathbb{P}_{k}$) denote the distribution of $\bm{M}$ under
$\mathcal{H}_{0}$ (resp.~$\mathcal{H}_{k}$), and let $\mathbb{P}_{0,i,j}$
(resp.~$\mathbb{P}_{k,i,j}$) stand for the distribution of $M_{ij}$
under $\mathcal{H}_{0}$ (resp.~$\mathcal{H}_{k}$).

\paragraph{Step 1: calculation of KL divergence.} We first calculate
the KL divergence of $\mathbb{P}_{0}$ from $\mathbb{P}_{k}$ as follows
\begin{equation}
\mathsf{KL}\left(\mathbb{P}_{k}\,\|\,\mathbb{P}_{0}\right)\overset{(\mathrm{i})}{=}\sum_{1\leq i,j\leq n}\mathsf{KL}\left(\mathbb{P}_{k,i,j}\,\|\,\mathbb{P}_{0,i,j}\right)\overset{(\mathrm{ii})}{=}\sum_{1\leq i,j\leq n}\frac{\left(M_{ij}^{\star}-(\bm{M}_{k})_{ij}\right)^{2}}{2\sigma_{ij}^{2}}\leq\frac{\|\bm{M}^{\star}-\bm{M}_{k}\|_{\mathrm{F}}^{2}}{2\sigma_{\min}^{2}}.\label{eq:KL-UB-1}
\end{equation}
Here, (i) holds since KL divergence is additive for independent distributions
\cite[Chapter 2.4]{tsybakov2009introduction}, and (ii) follows since
\[
\mathsf{KL}\left(\mathcal{N}(\mu_{1},\sigma^{2})\,\|\,\mathcal{N}(\mu_{2},\sigma^{2})\,\right)=\frac{(\mu_{1}-\mu_{2})^{2}}{2\sigma^{2}}.
\]
In addition, a little algebra reveals that
\begin{align*}
\bm{M}^{\star}-\bm{M}_{k} & =\lambda_{k}^{\star}\left(\bm{u}_{l}^{\star}\bm{u}_{l}^{\star\top}+\bm{u}_{k}^{\star}\bm{u}_{k}^{\star\top}\right)+\left(\lambda_{l}^{\star}-\lambda_{k}^{\star}\right)\bm{u}_{l}^{\star}\bm{u}_{l}^{\star\top}-\left\{ \lambda_{k}^{\star}\left(\widetilde{\bm{u}}_{l}\widetilde{\bm{u}}_{l}^{\top}+\widetilde{\bm{u}}_{k}\widetilde{\bm{u}}_{k}^{\top}\right)+\left(\lambda_{l}^{\star}-\lambda_{k}^{\star}\right)\widetilde{\bm{u}}_{l}\widetilde{\bm{u}}_{l}^{\top}\right\} \\
 & =\left(\lambda_{l}^{\star}-\lambda_{k}^{\star}\right)\left(\bm{u}_{l}^{\star}\bm{u}_{l}^{\star\top}-\widetilde{\bm{u}}_{l}\widetilde{\bm{u}}_{l}^{\top}\right),
\end{align*}
where the last relation makes use of the condition (\ref{eq:constraint-u1-u2-hypotheses}).
Therefore, we continue the bound (\ref{eq:KL-UB-1}) to reach
\begin{align}
\mathsf{KL}\left(\mathbb{P}_{k}\,\|\,\mathbb{P}_{0}\right)\leq\frac{\|\bm{M}^{\star}-\bm{M}_{k}\|_{\mathrm{F}}^{2}}{2\sigma_{\min}^{2}} & =\frac{\left(\lambda_{l}^{\star}-\lambda_{k}^{\star}\right)^{2}}{2\sigma_{\min}^{2}}\left\Vert \bm{u}_{l}^{\star}\bm{u}_{l}^{\star\top}-\widetilde{\bm{u}}_{l}\widetilde{\bm{u}}_{l}^{\top}\right\Vert _{\mathrm{F}}^{2}\leq\frac{\left(\lambda_{l}^{\star}-\lambda_{k}^{\star}\right)^{2}}{\sigma_{\min}^{2}}\left\Vert \bm{u}_{l}^{\star}-\widetilde{\bm{u}}_{l}\right\Vert _{2}^{2},\label{eq:KL-UB-2}
\end{align}
where the last line holds due to the following inequality that holds
for any unit vectors $\bm{u}$ and $\bm{v}$:
\begin{align}
\|\bm{u}\bm{u}^{\top}-\bm{v}\bm{v}^{\top}\|_{\mathrm{F}}^{2} & =\|\bm{u}\|_{2}^{4}+\|\bm{v}\|_{2}^{4}-2\langle\bm{u}\bm{u}^{\top},\bm{v}\bm{v}^{\top}\rangle=2-2(\bm{u}^{\top}\bm{v})^{2}\nonumber \\
 & =(2-2\bm{u}^{\top}\bm{v})(1+\bm{u}^{\top}\bm{v})=\tfrac{1}{2}\|\bm{u}-\bm{v}\|_{2}^{2}\cdot\|\bm{u}+\bm{v}\|_{2}^{2}\nonumber \\
 & \leq2\|\bm{u}-\bm{v}\|_{2}^{2}.\label{eq:uv-f-bound}
\end{align}

\paragraph{Step 2: bounding minimax probability of error.} Define
the minimax probability of error as follows
\begin{align}
	p_{\mathrm{e},k}:=\inf_{\psi}\max\Big\{\mathbb{P}\left\{ \psi\text{ rejects }\mathcal{H}_{0}\mid\mathcal{H}_{0}\right\} ,\,\mathbb{P}\left\{ \psi\text{ rejects }\mathcal{H}_{k}\mid\mathcal{H}_{k}\right\} \Big\},	
\end{align}
where the infimum is over all tests. Standard minimax lower bounds
\cite[Theorem 2]{tsybakov2009introduction} tell us that: if
\[
\mathsf{KL}\left(\mathbb{P}_{k}\,\|\,\mathbb{P}_{0}\right)\leq1/16,
\]
then one necessarily has $p_{\mathrm{e},k}\geq1/5$. This taken collectively
with the upper bound (\ref{eq:KL-UB-2}) implies that: if $\frac{\left(\lambda_{l}^{\star}-\lambda_{k}^{\star}\right)^{2}}{\sigma_{\min}^{2}}\left\Vert \bm{u}_{l}^{\star}-\widetilde{\bm{u}}_{l}\right\Vert _{2}^{2}\leq1/16$,
or equivalently, if 
\[
\|\bm{u}_{l}^{\star}-\widetilde{\bm{u}}_{l}\|_{2}\leq\frac{\sigma_{\min}}{4\big|\lambda_{l}^{\star}-\lambda_{k}^{\star}\big|},
\]
then the minimax probability of error $p_{\mathrm{e},k}$ is lower
bounded by $1/5$.

\paragraph{Step 3(a): establishing minimax $\ell_{2}$ lower bounds.}
The above minimax probability of testing error has direct implications
on $\ell_{2}$ eigenvector estimation. Suppose that $\|\bm{u}_{l}^{\star}-\widetilde{\bm{u}}_{l}\|_{2}=\frac{\sigma_{\min}}{4|\lambda_{l}^{\star}-\lambda_{k}^{\star}|}$,
then the calculation in (\ref{eq:KL-UB-2}) indicates that
\[
\|\bm{M}^{\star}-\bm{M}_{k}\|_{\mathrm{F}}\leq\sqrt{2}\big|\lambda_{l}^{\star}-\lambda_{k}^{\star}\big|\cdot\|\bm{u}_{l}^{\star}-\widetilde{\bm{u}}_{l}\|_{2}<\tfrac{\sigma_{\min}}{2},
\]
and hence $\bm{M}_{k}\in\mathcal{M}(\bm{M}^{\star})$. Moreover, when
$\sigma_{\min}\leq\big|\lambda_{l}^{\star}-\lambda_{k}^{\star}\big|$,
one has
\begin{align*}
\|\bm{u}_{l}^{\star}+\widetilde{\bm{u}}_{l}\|_{2} & =\|2\bm{u}_{l}^{\star}-(\bm{u}_{l}^{\star}-\widetilde{\bm{u}}_{l})\|_{2}\geq\|2\bm{u}_{l}^{\star}\|_{2}-\|\bm{u}_{l}^{\star}-\widetilde{\bm{u}}_{l}\|_{2}\geq2-\frac{\sigma_{\min}}{4\big|\lambda_{l}^{\star}-\lambda_{k}^{\star}\big|}\\
 & >\frac{\sigma_{\min}}{4\big|\lambda_{l}^{\star}-\lambda_{k}^{\star}\big|}=\|\bm{u}_{l}^{\star}-\widetilde{\bm{u}}_{l}\|_{2},
\end{align*}
meaning that $\|\bm{u}_{l}^{\star}-\widetilde{\bm{u}}_{l}\|_{2}=\min\|\bm{u}_{l}^{\star}\pm\widetilde{\bm{u}}_{l}\|_{2}$.
Thus, the standard reduction scheme described in \cite[Chapter 2.2]{tsybakov2009introduction}
leads to
\[
\inf_{\widehat{\bm{u}}_{l}}\sup_{\bm{A}\in\mathcal{M}_{0}(\bm{M}^{\star})}\mathbb{E}\Big[\min\|\widehat{\bm{u}}_{l}\pm\bm{u}_{l}(\bm{A})\|_{2}\Big]\gtrsim p_{\mathrm{e},k}\|\bm{u}_{l}^{\star}-\widetilde{\bm{u}}_{l}\|_{2}\asymp\frac{\sigma_{\min}}{\big|\lambda_{l}^{\star}-\lambda_{k}^{\star}\big|},
\]
where the infimum is taken over all estimator for $\bm{u}_{l}(\bm{A})$.
Since the preceding bound holds for all $k\neq l$ , we conclude that
\[
\inf_{\widehat{\bm{u}}_{l}}\sup_{\bm{A}\in\mathcal{M}_{0}(\bm{M}^{\star})}\mathbb{E}\Big[\min\|\widehat{\bm{u}}_{l}\pm\bm{u}_{l}(\bm{A})\|_{2}\Big]\gtrsim\max_{k:k\neq l}\frac{\sigma_{\min}}{\big|\lambda_{l}^{\star}-\lambda_{k}^{\star}\big|}=\frac{\sigma_{\min}}{\Delta_{l}^{\star}}.
\]

\paragraph{Step 3(b): establishing minimax lower bounds on estimating
linear functionals of eigenvectors.} The preceding minimax probability
of error also has direct implications on estimating linear functionals
of eigenvectors. In order to satisfy the rotational invariance constraint
(\ref{eq:constraint-u1-u2-hypotheses}), we set
\[
\left[\widetilde{\bm{u}}_{l},\widetilde{\bm{u}}_{k}\right]=\left[\bm{u}_{l}^{\star},\bm{u}_{k}^{\star}\right]\small\left[\begin{array}{cc}
\cos\theta_{k} & \sin\theta_{k}\\
-\sin\theta_{k} & \cos\theta_{k}
\end{array}\right]
\]
for some $\theta_{k}\in[-\pi/2,\pi/2]$. Before continuing, we shall
also make precise the connection between $\|\bm{u}_{l}^{\star}-\widetilde{\bm{u}}_{l}\|_{2}$
and $\theta_{k}$. Specifically, the above relation allows one to
derive
\[
\bm{u}_{l}^{\star}-\widetilde{\bm{u}}_{l}=\bm{u}_{l}^{\star}\left(1-\cos\theta_{k}\right)+\bm{u}_{k}^{\star}\sin\theta_{k}=2\bm{u}_{l}^{\star}\sin^{2}\tfrac{\theta_{k}}{2}+\bm{u}_{k}^{\star}\sin\theta_{k}
\]
and, as a consequence,
\begin{equation}
\begin{cases}
\|\bm{u}_{l}^{\star}-\widetilde{\bm{u}}_{l}\|_{2}\leq\left(1-\cos\theta_{k}\right)+|\sin\theta_{k}|=2\sin^{2}\tfrac{\theta_{k}}{2}+|\sin\theta_{k}|\leq3|\theta_{k}|,\\
\|\bm{u}_{l}^{\star}-\widetilde{\bm{u}}_{l}\|_{2}\geq|\sin\theta_{k}|\geq\tfrac{2}{\pi}|\theta_{k}|.
\end{cases}\label{eq:ul-UB-theta}
\end{equation}

In what follows, we shall take $\mathsf{sign}(\theta_{k})=\mathsf{sign}\big(\frac{\bm{a}^{\top}\bm{u}_{k}^{\star}}{\bm{a}^{\top}\bm{u}_{l}^{\star}}\big)$,
and generate the magnitude $|\theta_{k}|$ as follows
\begin{equation}
|\theta_{k}|\sim\mathsf{Uniform}\Big(\Big[\tfrac{\sigma_{\min}}{120\big|\lambda_{l}^{\star}-\lambda_{k}^{\star}\big|},\,\tfrac{\sigma_{\min}}{12\big|\lambda_{l}^{\star}-\lambda_{k}^{\star}\big|}\Big]\Big),\label{eq:theta-uniform-bound}
\end{equation}
which combined with (\ref{eq:ul-UB-theta}) guarantees that
\begin{equation}
\frac{\sigma_{\min}}{60\pi\big|\lambda_{l}^{\star}-\lambda_{k}^{\star}\big|}\leq\|\bm{u}_{l}^{\star}-\widetilde{\bm{u}}_{l}\|_{2}\leq\frac{\sigma_{\min}}{4\big|\lambda_{l}^{\star}-\lambda_{k}^{\star}\big|}.\label{eq:ul-ul-L2-risk}
\end{equation}

We aim to translate the difficulty in distinguishing $\mathcal{H}_{0}$
and $\mathcal{H}_{k}$ into a fundamental lower bound on estimating
the linear functional. Towards this, we are in need of computing the
difference of the linear functional under these two hypotheses (namely,
$\bm{a}^{\top}\bm{u}_{l}^{\star}-\bm{a}^{\top}\widetilde{\bm{u}}_{l}$).
The above expressions give
\begin{align}
\left|\bm{a}^{\top}\bm{u}_{l}^{\star}-\bm{a}^{\top}\widetilde{\bm{u}}_{l}\right| & =\left|2\big(\bm{a}^{\top}\bm{u}_{l}^{\star}\big)\sin^{2}\tfrac{\theta_{k}}{2}+\big(\bm{a}^{\top}\bm{u}_{k}^{\star}\big)\sin\theta_{k}\right|\overset{(\mathrm{i})}{=}2\left|\bm{a}^{\top}\bm{u}_{l}^{\star}\right|\sin^{2}\tfrac{\theta_{k}}{2}+\left|\bm{a}^{\top}\bm{u}_{k}^{\star}\right|\cdot|\sin\theta_{k}|\label{eq:au-bound1}\\
 & \asymp\theta_{k}^{2}\left|\bm{a}^{\top}\bm{u}_{l}^{\star}\right|+|\theta_{k}|\cdot\left|\bm{a}^{\top}\bm{u}_{k}^{\star}\right|,\nonumber 
\end{align}
where the identity (i) results from the condition $\mathsf{sign}(\theta_{k})=\mathsf{sign}\big(\frac{\bm{a}^{\top}\bm{u}_{k}^{\star}}{\bm{a}^{\top}\bm{u}_{l}^{\star}}\big)$.

We shall also control $\left|\bm{a}^{\top}\bm{u}_{l}^{\star}+\bm{a}^{\top}\widetilde{\bm{u}}_{l}\right|$,
for which we divide into two cases:
\begin{itemize}
\item If $|\bm{a}^{\top}\bm{u}_{l}^{\star}|\geq\frac{\sigma_{\min}}{|\lambda_{l}^{\star}-\lambda_{k}^{\star}|}|\bm{a}^{\top}\bm{u}_{k}^{\star}|\geq12|\theta_{k}|\cdot|\bm{a}^{\top}\bm{u}_{k}^{\star}|$,
then the identity $\bm{u}_{l}^{\star}+\widetilde{\bm{u}}_{l}=2\bm{u}_{l}^{\star}-\big(\bm{u}_{l}^{\star}-\widetilde{\bm{u}}_{l}\big)$
together with (\ref{eq:au-bound1}) yields
\begin{align*}
\left|\bm{a}^{\top}\bm{u}_{l}^{\star}+\bm{a}^{\top}\widetilde{\bm{u}}_{l}\right| & =\left|2|\bm{a}^{\top}\bm{u}_{l}^{\star}|-2|\bm{a}^{\top}\bm{u}_{l}^{\star}|\sin^{2}\tfrac{\theta_{k}}{2}-|\bm{a}^{\top}\bm{u}_{k}^{\star}|\cdot|\sin\theta_{k}|\right|\\
 & \geq2\left|\bm{a}^{\top}\bm{u}_{l}^{\star}\right|\cos^{2}\tfrac{\theta_{k}}{2}-|\bm{a}^{\top}\bm{u}_{k}^{\star}|\cdot|\theta_{k}|\\
 & \geq\left|\bm{a}^{\top}\bm{u}_{l}^{\star}\right|-|\bm{a}^{\top}\bm{u}_{k}^{\star}|\cdot|\theta_{k}|\geq\tfrac{1}{2}|\bm{a}^{\top}\bm{u}_{l}^{\star}|\\
 & \gtrsim\theta_{k}^{2}\left|\bm{a}^{\top}\bm{u}_{l}^{\star}\right|+|\theta_{k}|\cdot\left|\bm{a}^{\top}\bm{u}_{k}^{\star}\right|,
\end{align*}
where we have used the fact that $\sigma_{\min}\leq|\lambda_{l}^{\star}-\lambda_{k}^{\star}|$
(so that $|\theta_{k}|\leq\frac{\sigma_{\min}}{12|\lambda_{l}^{\star}-\lambda_{k}^{\star}|}\leq\frac{1}{12}$
and hence $\cos^{2}\tfrac{\theta_{k}}{2}\geq\frac{1}{2}$).
\item Suppose now that $|\bm{a}^{\top}\bm{u}_{l}^{\star}|<\frac{\sigma_{\min}}{|\lambda_{l}^{\star}-\lambda_{k}^{\star}|}|\bm{a}^{\top}\bm{u}_{k}^{\star}|$.
As one can easily verify (which we omit for brevity), the scheme (\ref{eq:theta-uniform-bound})
guarantees that with probability exceeding $1/2$, one has
\begin{align*}
\left|\bm{a}^{\top}\bm{u}_{l}^{\star}+\bm{a}^{\top}\widetilde{\bm{u}}_{l}\right| & =\left|2|\bm{a}^{\top}\bm{u}_{l}^{\star}|-2|\bm{a}^{\top}\bm{u}_{l}^{\star}|\sin^{2}\tfrac{\theta_{k}}{2}-|\bm{a}^{\top}\bm{u}_{k}^{\star}|\cdot|\sin\theta_{k}|\right|\\
 & \gtrsim\max\left\{ |\bm{a}^{\top}\bm{u}_{l}^{\star}|,\,|\bm{a}^{\top}\bm{u}_{k}^{\star}|\cdot|\theta_{k}|\right\} \\
 & \gtrsim\theta_{k}^{2}\left|\bm{a}^{\top}\bm{u}_{l}^{\star}\right|+|\theta_{k}|\cdot\left|\bm{a}^{\top}\bm{u}_{k}^{\star}\right|.
\end{align*}
\end{itemize}
Putting the above cases together reveals that
\begin{align*}
\left|\bm{a}^{\top}\bm{u}_{l}^{\star}+\bm{a}^{\top}\widetilde{\bm{u}}_{l}\right| & \gtrsim\theta_{k}^{2}\left|\bm{a}^{\top}\bm{u}_{l}^{\star}\right|+|\theta_{k}|\cdot\left|\bm{a}^{\top}\bm{u}_{k}^{\star}\right|
\end{align*}
with probability exceeding $1/2$. Consequently, one can find $|\theta_{k}|\in\Big[\tfrac{\sigma_{\min}}{120\big|\lambda_{l}^{\star}-\lambda_{k}^{\star}\big|},\,\tfrac{\sigma_{\min}}{12\big|\lambda_{l}^{\star}-\lambda_{k}^{\star}\big|}\Big]$
such that
\[
\min\left|\bm{a}^{\top}\bm{u}_{l}^{\star}\pm\bm{a}^{\top}\widetilde{\bm{u}}_{l}\right|\gtrsim\theta_{k}^{2}\left|\bm{a}^{\top}\bm{u}_{l}^{\star}\right|+|\theta_{k}|\cdot\left|\bm{a}^{\top}\bm{u}_{k}^{\star}\right|\asymp\left|\bm{a}^{\top}\bm{u}_{l}^{\star}\right|\frac{\sigma_{\min}^{2}}{\big|\lambda_{l}^{\star}-\lambda_{k}^{\star}\big|^{2}}+\left|\bm{a}^{\top}\bm{u}_{k}^{\star}\right|\frac{\sigma_{\min}}{\big|\lambda_{l}^{\star}-\lambda_{k}^{\star}\big|},
\]
where we recall that $\min|a\pm b|:=\min\{|a-b|,|a+b|\}$.

Suppose for the moment that $\bm{M}_{k}\in\mathcal{M}_0(\bm{M}^{\star})$.
Applying the standard reduction scheme \cite[Chapter 2.2]{tsybakov2009introduction}
once again yields
\begin{align*}
\inf_{\widehat{u}_{\bm{a},l}}\sup_{\bm{A}\in\mathcal{M}_{0}(\bm{M}^{\star})}\mathbb{E}\Big[\min\big|\widehat{u}_{\bm{a},l}\pm\bm{a}^{\top}\bm{u}_{l}(\bm{A})\big|\Big] & \gtrsim p_{\mathrm{e},k}\min\left|\bm{a}^{\top}\bm{u}_{l}^{\star}\pm\bm{a}^{\top}\widetilde{\bm{u}}_{l}\right|\gtrsim\frac{\sigma_{\min}^{2}\left|\bm{a}^{\top}\bm{u}_{l}^{\star}\right|}{\big|\lambda_{l}^{\star}-\lambda_{k}^{\star}\big|^{2}}+\frac{\sigma_{\min}\left|\bm{a}^{\top}\bm{u}_{k}^{\star}\right|}{\big|\lambda_{l}^{\star}-\lambda_{k}^{\star}\big|},
\end{align*}
where the infimum is taken over all estimators for $\bm{a}^{\top}\bm{u}_{l}(\bm{A})$.
Recognizing that the preceding inequality holds for all $k\neq l$,
we immediately arrive at the advertised claim
\begin{align*}
\inf_{\widehat{u}_{\bm{a},l}}\sup_{\bm{A}\in\mathcal{M}_{0}(\bm{M}^{\star})}\mathbb{E}\Big[\min\big|\widehat{u}_{\bm{a},l}\pm\bm{a}^{\top}\bm{u}_{l}(\bm{A})\big|\Big] & \gtrsim\sigma_{\min}^{2}\max_{k:k\neq l}\frac{\left|\bm{a}^{\top}\bm{u}_{l}^{\star}\right|}{\big|\lambda_{l}^{\star}-\lambda_{k}^{\star}\big|^{2}}+\sigma_{\min}\max_{k:k\neq l}\frac{\left|\bm{a}^{\top}\bm{u}_{k}^{\star}\right|}{\big|\lambda_{l}^{\star}-\lambda_{k}^{\star}\big|}\\
 & =\left|\bm{a}^{\top}\bm{u}_{l}^{\star}\right|\frac{\sigma_{\min}^{2}}{\Delta_{l}^{\star2}}+\sigma_{\min}\max_{k:k\neq l}\frac{\left|\bm{a}^{\top}\bm{u}_{k}^{\star}\right|}{\big|\lambda_{l}^{\star}-\lambda_{k}^{\star}\big|}.
\end{align*}

Finally, it remains to justify that $\bm{M}_{k}\in\mathcal{M}_{0}(\bm{M}^{\star})$
for all $k\neq l$. When $|\theta_{k}|\leq\frac{\sigma_{\min}}{12|\lambda_{l}^{\star}-\lambda_{k}^{\star}|},$
invoking (\ref{eq:KL-UB-2}) and (\ref{eq:ul-UB-theta}) yields
\[
\left\Vert \bm{M}^{\star}-\bm{M}_{k}\right\Vert _{\mathrm{F}}^{2}\leq2\left(\lambda_{l}^{\star}-\lambda_{k}^{\star}\right)^{2}\left\Vert \bm{u}_{l}^{\star}-\widetilde{\bm{u}}_{l}\right\Vert _{2}^{2}\leq18\left(\lambda_{l}^{\star}-\lambda_{k}^{\star}\right)^{2}|\theta_{k}|^{2}M<\sigma_{\min}^{2}/4.
\]
This means that $\bm{M}_{k}\in\mathcal{M}_{0}(\bm{M}^{\star})$, thus
concluding the proof.

\subsection{Proof of Eq.~(\ref{eq:minimax-au-sigma})}

Consider the following two hypotheses
\begin{align*}
\mathcal{H}_{0}: & \quad\bm{M}=\bm{M}_{0}+\bm{H}:=\lambda_{l}^{\star}\bm{u}_{l}^{\star}\bm{u}_{l}^{\star\top}+\sum_{j:j\neq l}\lambda_{j}^{\star}\bm{v}_{j}\bm{v}_{j}^{\top}+\bm{H},\\
\mathcal{H}_{\bm{a}}: & \quad\bm{M}=\bm{M}_{\bm{a}}+\bm{H}:=\lambda_{l}^{\star}\widetilde{\bm{u}}_{\bm{a}}\widetilde{\bm{u}}_{\bm{a}}^{\top}+\sum_{j:j\neq l}\lambda_{j}^{\star}\bm{v}_{j}\bm{v}_{j}^{\top}+\bm{H},
\end{align*}
where the $\bm{v}_{j}$'s are orthonormal vectors obeying $\langle\bm{v}_{j},\bm{a}\rangle=\langle\bm{v}_{j},\bm{u}_{l}^{\star}\rangle=0$
for any $j\neq l$, and $\widetilde{\bm{u}}_{\bm{a}}$ is defined
as 
\[
\widetilde{\bm{u}}_{\bm{a}}:=\frac{1}{\|\bm{u}_{l}^{\star}+\delta\bm{a}_{\perp}\|_{2}}(\bm{u}_{l}^{\star}+\delta\bm{a}_{\perp}),\qquad\text{with }\bm{a}_{\perp}=\bm{a}-(\bm{a}^{\top}\bm{u}_{l}^{\star})\bm{u}_{l}^{\star}\,\text{ and }\,\delta=\frac{\sigma_{\min}}{12|\lambda_{l}^{\star}|\cdot\|\bm{a}_{\perp}\|_{2}}.
\]
Recognizing the simple fact $\langle\bm{a}_{\perp},\bm{u}_{l}^{\star}\rangle=0$
, we can derive
\[
\|\bm{u}_{l}^{\star}+\delta\bm{a}_{\perp}\|_{2}=\sqrt{1+\delta^{2}\|\bm{a}_{\perp}\|_{2}^{2}}\qquad\text{and}\qquad\widetilde{\bm{u}}_{\bm{a}}:=\frac{1}{\sqrt{1+\delta^{2}\|\bm{a}_{\perp}\|_{2}^{2}}}(\bm{u}_{l}^{\star}+\delta\bm{a}_{\perp}).
\]
Our proof proceeds as follows. Without loss of generality, we shall
assume that $\bm{a}^{\top}\bm{u}^{\star}\geq0$. 
\begin{itemize}
\item Let $\mathbb{P}_{0}$ (resp.~$\mathbb{P}_{\bm{a}}$) denote the distribution
of $\bm{M}$ under $\mathcal{H}_{0}$ (resp.~$\mathcal{H}_{\bm{a}}$).
Repeating the derivation in (\ref{eq:KL-UB-1}) gives
\begin{equation}
\mathsf{KL}\left(\mathbb{P}_{\bm{a}}\,\|\,\mathbb{P}_{0}\right)\leq\frac{\|\bm{M}_{0}-\bm{M}_{\bm{a}}\|_{\mathrm{F}}^{2}}{2\sigma_{\min}^{2}}=\frac{(\lambda_{l}^{\star})^{2}}{2\sigma_{\min}^{2}}\|\bm{u}_{l}^{\star}\bm{u}_{l}^{\star\top}-\widetilde{\bm{u}}_{\bm{a}}\widetilde{\bm{u}}_{\bm{a}}^{\top}\|_{\mathrm{F}}^{2}\leq\frac{(\lambda_{l}^{\star})^{2}}{\sigma_{\min}^{2}}\|\bm{u}_{l}^{\star}-\widetilde{\bm{u}}_{\bm{a}}\|_{2}^{2},\label{eq:KL-UB-1-1}
\end{equation}
where the last inequality comes from (\ref{eq:uv-f-bound}). In addition,
\begin{align}
\|\bm{u}_{l}^{\star}-\widetilde{\bm{u}}_{\bm{a}}\|_{2} & =\Big\|\bm{u}_{l}^{\star}-\frac{1}{\sqrt{1+\delta^{2}\|\bm{a}_{\perp}\|_{2}^{2}}}(\bm{u}_{l}^{\star}+\delta\bm{a}_{\perp})\Big\|_{2}\nonumber \\
 & \leq\big\|\bm{u}_{l}^{\star}-(\bm{u}_{l}^{\star}+\delta\bm{a}_{\perp})\big\|_{2}+\Big(1-\frac{1}{\sqrt{1+\delta^{2}\|\bm{a}_{\perp}\|_{2}^{2}}}\Big)\cdot\big\|\bm{u}_{l}^{\star}+\delta\bm{a}_{\perp}\big\|_{2}\nonumber \\
 & \leq\delta\|\bm{a}_{\perp}\|_{2}+\frac{\sqrt{1+\delta^{2}\|\bm{a}_{\perp}\|_{2}^{2}}-1}{\sqrt{1+\delta^{2}\|\bm{a}_{\perp}\|_{2}^{2}}}\cdot\left(1+\delta\|\bm{a}_{\perp}\|_{2}\right)\nonumber \\
 & \leq\delta\|\bm{a}_{\perp}\|_{2}+\left(1+\delta\|\bm{a}_{\perp}\|_{2}\right)\delta\|\bm{a}_{\perp}\|_{2}\leq3\delta\|\bm{a}_{\perp}\|_{2},\label{eq:ul-ua-L2bound}
\end{align}
where we have made use of the fact that $\delta\|\bm{a}_{\perp}\|_{2}=\frac{\sigma_{\min}}{12|\lambda_{l}^{\star}|}\leq\frac{1}{12}$.
Combining this with (\ref{eq:KL-UB-1-1}) and our choice of $\delta$,
we arrive at
\begin{equation}
\mathsf{KL}\left(\mathbb{P}_{\bm{a}}\,\|\,\mathbb{P}_{0}\right)\leq\frac{9\delta^{2}(\lambda_{l}^{\star})^{2}\|\bm{a}_{\perp}\|_{2}^{2}}{\sigma_{\min}^{2}}=\frac{1}{16}.\label{eq:KL-UB-1-1-1}
\end{equation}
\item Define the minimax probability of error as follows
\[
p_{\mathrm{e},\bm{a}}:=\inf_{\psi}\max\Big\{\mathbb{P}\left\{ \psi\text{ rejects }\mathcal{H}_{0}\mid\mathcal{H}_{0}\right\} ,\,\mathbb{P}\left\{ \psi\text{ rejects }\mathcal{H}_{\bm{a}}\mid\mathcal{H}_{\bm{a}}\right\} \Big\},
\]
where the infimum is over all tests. It then follows from \cite[Theorem 2]{tsybakov2009introduction}
that: if 
\[
\mathsf{KL}\left(\mathbb{P}_{k}\,\|\,\mathbb{P}_{0}\right)\leq1/16,
\]
one necessarily has $p_{\mathrm{e},\bm{a}}\geq1/5$. This taken collectively
with the upper bound (\ref{eq:KL-UB-1-1-1}) implies that: the minimax
probability of error $p_{\mathrm{e},\bm{a}}$ in distinguishing $\mathcal{H}_{0}$
and $\mathcal{H}_{\bm{a}}$ is indeed lower bounded by $1/5$.
\item Combining (\ref{eq:ul-ua-L2bound}) with the fact $\delta\|\bm{a}_{\perp}\|_{2}=\frac{\sigma_{\min}}{12|\lambda_{l}^{\star}|}$
gives
\[
\|\widetilde{\bm{u}}_{\bm{a}}-\bm{u}_{l}^{\star}\|_{2}\leq\frac{\sigma_{\min}}{4|\lambda_{l}^{\star}|},
\]
thus indicating that $\bm{M}_{\bm{a}}\in\mathcal{M}_{1}(\bm{M}^{\star})$.
Apply the standard reduction scheme \cite[Chapter 2.2]{tsybakov2009introduction}
to yield
\begin{equation}
\inf_{\widehat{u}_{\bm{a},l}}\sup_{\bm{A}\in\mathcal{M}_{1}(\bm{M}^{\star})}\mathbb{E}\Big[\min\big|\widehat{u}_{\bm{a},l}\pm\bm{a}^{\top}\bm{u}_{l}(\bm{A})\big|\Big]\gtrsim\min\big|\bm{a}^{\top}(\bm{u}_{l}^{\star}\pm\widetilde{\bm{u}}_{\bm{a}})\big|.\label{eq:LB-au-sigma-1}
\end{equation}
Everything then boils down to lower bounding $\min\big|\bm{a}^{\top}(\bm{u}_{l}^{\star}\pm\widetilde{\bm{u}}_{\bm{a}})\big|$.
On the one hand, it is seen that
\begin{align*}
\big|\bm{a}^{\top}(\bm{u}_{l}^{\star}-\widetilde{\bm{u}}_{\bm{a}})\big| & =\left|\bm{a}^{\top}\bm{u}_{l}^{\star}-\frac{1}{\sqrt{1+\delta^{2}\|\bm{a}_{\perp}\|_{2}^{2}}}(\bm{a}^{\top}\bm{u}_{l}^{\star}+\delta\bm{a}^{\top}\bm{a}_{\perp})\right|\\
 & \geq\frac{\delta\big|\bm{a}^{\top}\bm{a}_{\perp}\big|}{\sqrt{1+\delta^{2}\|\bm{a}_{\perp}\|_{2}^{2}}}-\big|\bm{a}^{\top}\bm{u}_{l}^{\star}\big|\Bigg(1-\frac{1}{\sqrt{1+\delta^{2}\|\bm{a}_{\perp}\|_{2}^{2}}}\Bigg)\\
 & \geq\tfrac{1}{2}\delta\|\bm{a}_{\perp}\|_{2}^{2}-\big|\bm{a}^{\top}\bm{u}_{l}^{\star}\big|\cdot\delta^{2}\|\bm{a}_{\perp}\|_{2}^{2}\\
 & \geq\tfrac{1}{4}\delta\|\bm{a}_{\perp}\|_{2}^{2}.
\end{align*}
On the other hand, one can employ the assumption $\bm{a}^{\top}\bm{u}^{\star}\geq0$
to derive
\begin{align*}
\big|\bm{a}^{\top}(\bm{u}_{l}^{\star}+\widetilde{\bm{u}}_{\bm{a}})\big| & =\left|\bm{a}^{\top}\bm{u}_{l}^{\star}+\frac{1}{\sqrt{1+\delta^{2}\|\bm{a}_{\perp}\|_{2}^{2}}}(\bm{a}^{\top}\bm{u}_{l}^{\star}+\delta\bm{a}^{\top}\bm{a}_{\perp})\right|\geq\bm{a}^{\top}\bm{u}_{l}^{\star}+\frac{1}{\sqrt{1+\delta^{2}\|\bm{a}_{\perp}\|_{2}^{2}}}\delta\|\bm{a}_{\perp}\|_{2}^{2}\\
 & \geq\tfrac{1}{2}\delta\|\bm{a}_{\perp}\|_{2}^{2}.
\end{align*}
Putting all this together, we conclude that 
\[
\min\big|\bm{a}^{\top}(\bm{u}_{l}^{\star}\pm\widetilde{\bm{u}}_{\bm{a}})\big|\geq\tfrac{1}{4}\delta\|\bm{a}_{\perp}\|_{2}^{2},
\]
which taken collectively with (\ref{eq:LB-au-sigma-1}) and our choice
$\delta=\frac{\sigma_{\min}}{12|\lambda_{l}^{\star}|\cdot\|\bm{a}_{\perp}\|_{2}}$
yields
\[
\inf_{\widehat{u}_{\bm{a},l}}\sup_{\bm{A}\in\mathcal{M}_{1}(\bm{M}^{\star})}\mathbb{E}\Big[\min\big|\widehat{u}_{\bm{a},l}\pm\bm{a}^{\top}\bm{u}_{l}(\bm{A})\big|\Big]\gtrsim\min\big|\bm{a}^{\top}(\bm{u}_{l}^{\star}\pm\widetilde{\bm{u}}_{\bm{a}})\big|\gtrsim\frac{\sigma_{\min}}{|\lambda_{l}^{\star}|}\|\bm{a}_{\perp}\|_{2}.
\]
\end{itemize}

\subsection{Proof of Eq.~(\ref{eq:minimax-L2-sigma})}
The proof of this part uses Fano's inequality. 
We start by constructing a proper packing set of the space: $N$ unit vectors $\left\{ \bm{v}_{i}\right\} _{1\leq i\leq N}$
within the subspace perpendicular to $\{\bm{u}_{i}^\star\}_{i\neq l}$ obeying
\begin{itemize}
\item $\langle\bm{v}_{i},\bm{u}^\star_{j}\rangle = 0$ for all $1\leq i\leq N$, $1\leq j\leq r$
and $j\neq l$;
\item there exists some sufficiently small constant $c_{3}>0$ such that
\[
\|\bm{v}_{i}-\bm{v}_{j}\|_{2}\geq c_{3}\frac{\sigma_{\min}\sqrt{n}}{|\lambda_{l}^{\star}|},\qquad1\leq i\neq j\leq N;
\]
\item there exists some sufficiently large constant $c_{4}>0$ such that
\[
\|\bm{v}_{i}-\bm{u}_{l}^{\star}\|_{2}\leq c_{4}\frac{\sigma_{\min}\sqrt{n}}{|\lambda_{l}^{\star}|},\qquad1\leq i\leq N.
\]
\end{itemize}
Standard packing arguments \cite[Chapter 5.1]{wainwright2019high}
imply that $N$ can be as large as 
\begin{equation}
N=\exp\Big(n\log\frac{c_{4}}{2c_{3}}\Big).\label{eq:packing-bound}
\end{equation}
In addition, when $4c_{4}\frac{\sigma_{\min}\sqrt{n}}{|\lambda_{l}^{\star}|}<1$,
it follows that
\begin{align}
\|\bm{v}_{i}-\bm{v}_{j}\|_{2} & \leq\|\bm{v}_{i}-\bm{u}_{l}^{\star}\|_{2}+\|\bm{v}_{j}-\bm{u}_{l}^{\star}\|_{2}\leq2c_{4}\frac{\sigma_{\min}\sqrt{n}}{|\lambda_{l}^{\star}|},\label{eq:vi-vj-UB}\\
\|\bm{v}_{i}+\bm{v}_{j}\|_{2} & \geq2\|\bm{v}_{i}\|_{2}-\|\bm{v}_{i}-\bm{v}_{j}\|_{2}\geq2-2c_{4}\frac{\sigma_{\min}\sqrt{n}}{|\lambda_{l}^{\star}|}>2c_{4}\frac{\sigma_{\min}\sqrt{n}}{|\lambda_{l}^{\star}|},
\end{align}
thus indicating that
\begin{equation}
\min\|\bm{v}_{i}\pm\bm{v}_{j}\|_{2}\geq\min\left\{ c_{3},2c_{4}\right\} \frac{\sigma_{\min}\sqrt{n}}{|\lambda_{l}^{\star}|},\qquad1\leq i\neq j\leq N.\label{eq:vij-min-LB}
\end{equation}

The next step is to associate each vector $\bm{v}_{i}$ $(1\leq i\leq N$)
with a hypothesis as follows
\begin{align*}
\mathcal{H}_{i}: & \quad\bm{M}=\bm{M}_{i}+\bm{H}:=\lambda_{l}^{\star}\bm{v}_{i}\bm{v}_{i}^{\top}+\sum_{j:j\neq l}\lambda_{j}^{\star}\bm{u}_{j}^{\star}\bm{u}_{j}^{\star\top}+\bm{H},\qquad1\leq i\leq N.
\end{align*}
If we denote by $\mathbb{P}_{i}$ the distribution of $\bm{M}$ under
$\mathcal{H}_{i}$ $(1\leq i\leq N)$, then we can invoke the argument
in (\ref{eq:KL-UB-1}) and (\ref{eq:uv-f-bound}) to upper bound 
\begin{align}
\mathsf{KL}\left(\mathbb{P}_{i}\,\|\,\mathbb{P}_{j}\right) & \leq\frac{\|\bm{M}_{i}-\bm{M}_{j}\|_{\mathrm{F}}^{2}}{2\sigma_{\min}^{2}}=\frac{(\lambda_{l}^{\star})^{2}}{2\sigma_{\min}^{2}}\|\bm{v}_{i}\bm{v}_{i}^{\top}-\bm{v}_{j}\bm{v}_{j}^{\top}\|_{\mathrm{F}}^{2}\leq\frac{(\lambda_{l}^{\star})^{2}}{\sigma_{\min}^{2}}\|\bm{v}_{i}-\bm{v}_{j}\|_{2}^{2}\nonumber \\
 & \leq4c_{4}^{2}n\label{eq:KL-pi-pj}
\end{align}
for any $i\neq j$, where the last inequality follows from (\ref{eq:vi-vj-UB}).
This upper bound on KL divergence plays an important role in invoking
Fano's inequality. More specifically, recall that Fano's inequality
\cite[Corollary 2.6]{tsybakov2009introduction} asserts that if 
\[
\frac{1}{N}\sum_{i=2}^{N}\mathsf{KL}\left(\mathbb{P}_{i}\,\|\,\mathbb{P}_{1}\right)\leq\frac{1}{4}\log N,
\]
then the minimax probability of testing error must satisfy
\[
p_{\mathrm{e},N}:=\inf_{\psi}\max_{1\leq i\leq N}\mathbb{P}\left\{ \psi\neq i\mid\mathcal{H}_{i}\right\} \geq1/4,
\]
where the infimum is over all tests. Combining this with the upper
bound (\ref{eq:KL-pi-pj}) and the packing number (\ref{eq:packing-bound}),
we see that if 
\begin{equation}
4c_{4}^{2}n\leq\frac{N}{4(N-1)}\log N=\frac{N}{4(N-1)}\cdot n\log\frac{c_{4}}{2c_{3}},\label{eq:KL-requirement}
\end{equation}
then one would have $p_{\mathrm{e},N}\geq1/4$. Clearly, the condition
(\ref{eq:KL-requirement}) would hold as long as $c_{4}/c_{3}$ is
sufficiently large. 

To finish up, it suffices to invoke the standard reduction scheme
\cite[Chapter 2.2]{tsybakov2009introduction} to obtain
\begin{align*}
\inf_{\widehat{\bm{u}}_{l}}\sup_{\bm{A}\in\mathcal{\mathcal{M}}_{2}}\mathbb{E}\Big[\min\|\widehat{\bm{u}}_{l}\pm\bm{u}_{l}(\bm{A})\|_{2}\Big] & \gtrsim\min_{1\leq i\neq j\leq N}\left\{ \min\|\bm{v}_{i}\pm\bm{v}_{j}\|_{2}\right\} \gtrsim\frac{\sigma_{\min}\sqrt{n}}{\big|\lambda_{l}^{\star}\big|},
\end{align*}
where the last inequality is a consequence of the condition (\ref{eq:vij-min-LB}).
This concludes the proof.

\section{A few more auxiliary lemmas}
\label{sec:auxiliary-lemmas}

\begin{lemma}
\label{lem:va-properties}
	Suppose that $\bm{H}$ satisfies Assumptions~\ref{assumption-H}. Fix any $\bm{a}\in \mathbb{R}^n$. Then one has
\begin{align}
\label{eq:expression-va}
	\mathsf{Var}\left[\frac{1}{2\lambda_l^{\star}}\bm{a}^{\top}(\bm{H}+\bm{H}^{\top})\bm{u}_l^{\star}\right]=\frac{1}{4(\lambda_l^{\star})^{2}}\sum_{1\leq i,j\leq n}(a_{i}u_{l,j}^{\star}+a_{j}u_{l,i}^{\star})^{2}\sigma_{ij}^{2} =: v_{\bm{a},l}^{\star}, 
\end{align}
which satisfies 
\begin{align}
\label{eq:va-bounds}
\frac{1}{2}\left( \|\bm{a}\|_2^2  + \big(\bm{a}^{\top}\bm{u}_l^{\star}\big)^{2}\right)\frac{\sigma_{\min}^{2}}{(\lambda_l^{\star})^{2}} 
\leq v_{\bm{a},l}^{\star}
\leq\frac{1}{2}\left( \|\bm{a}\|_2^2  +\big(\bm{a}^{\top}\bm{u}_l^{\star}\big)^{2}\right)\frac{\sigma_{\max}^{2}}{(\lambda_l^{\star})^{2}}.
\end{align}
\end{lemma}
\begin{proof} See Appendix \ref{sec:proof-lem-va-properties}. \end{proof}
	If we further have $\bm{a}^{\top}\bm{u}_l^{\star} = o(\|\bm{a}\|_2)$, then one has
\begin{align}
\label{eq:va-bounds-simple}
	\frac{1+o(1)}{2} \frac{\sigma_{\min}^{2} \|\bm{a}\|_2^2}{(\lambda_l^{\star})^{2}} 
\leq v_{\bm{a},l}^{\star}
	\leq
	\frac{1+o(1)}{2} \frac{\sigma_{\max}^{2}\|\bm{a}\|_2^2}{(\lambda_l^{\star})^{2}}.
\end{align}

\begin{lemma}
\label{lem:aHu-UB}
Suppose that $\bm{H}$ satisfies Assumptions~\ref{assumption-H} and obeys $B\leq\sigma_{\max}\sqrt{\frac{n}{\mu\log n}}$. 
Fix any $\bm{a}\in \mathbb{R}^n$. Then with probability at least $1 - O(n^{-10})$, we have 
	\begin{align}
	\label{eq:property-one-inprod}
		\max\left\{ \left| \bm{a}^{\top} \bm{H} \bm{u}^{\star}_l \right|, \left| \bm{a}^{\top} \bm{H}^{\top} \bm{u}^{\star}_l \right|, \left| \bm{a}^{\top} (\bm{H}+\bm{H}^{\top}) \bm{u}^{\star}_l \right| \right\} \lesssim \sigmamax \ltwo{\bm{a}} \sqrt{\log n} 
	\end{align}
for all $1\leq l\leq r$.  
\end{lemma}
\begin{proof} See Appendix \ref{sec:proof-lemma-PropertyInprod}.\end{proof}

\begin{lemma}
\label{LemPropertyInprod}
Fix any unit vector $\bm{a}\in \mathbb{R}^n$, and consider another fixed unit vector $\bm{u}^{\star}$ obeying $\|\bm{u}^{\star}\|_{\infty}\leq \sqrt{\mu/n}$. 
	Suppose that $\bm{H}$ satisfies Assumption~\ref{assumption-H} and obeys $B\leq \sigma_{\min}\sqrt{n/(\mu\log n)}$. 
Then the distribution of $W \defn \frac{\bm{a}^{\top}\left(\bm{H}+\bm{H}^{\top}\right)\bm{u}^{\star}}{2\sqrt{v^{\star}}}$
satisfies
	\begin{align}
	\label{eq:property-two-inprod}
	\sup_{z\in \real} \big| \, \mprob(W\leq z) - \Phi(z) \,\big| \leq \frac{8}{\sqrt{\log n}},
	\end{align}
	where $\Phi(\cdot)$ is the CDF of a standard Gaussian, and $v^{\star}:=\mathsf{Var}\big(\frac{1}{2}\bm{a}^{\top}\left(\bm{H}+\bm{H}^{\top}\right)\bm{u}^{\star}\big)$. 
\end{lemma}
\begin{proof} See Appendix \ref{sec:proof-lemma-PropertyInprod}.\end{proof}

\begin{lemma}
\label{lem:basic-uR-uL-rankr}
Suppose that $\bm{H}$ satisfies Assumptions~\ref{assumption-H} and \ref{assumption:noise-size-rankr-revised}. Let $\ur_l$ (resp.~$\ul_l$) be the $l$th right (resp.~left) eigenvector of $\bm{M}$ obeying $\ur^{\top}_l\ustar>0$ and $\ur^{\top}_l\ul_l>0$, and set $\widehat{\bm{u}}_l := \frac{1}{\|\ur_l + \ul_l\|_2}(\ur_l + \ul_l) $. Then with probability exceeding $1-O(n^{-10})$, one has
\begin{equation}
\begin{cases}
\bm{u}_{l}^{\star\top}\bm{u}_{l} & =1-O\left(\frac{\kappa^{4}\sigma_{\max}^{2}n\log n}{(\lambda_{\mathrm{max}}^{\star})^{2}}+\frac{\mu\kappa^{4}r^{2}\sigma_{\max}^{2}\log n}{(\Delta_{l}^{\star})^{2}}\right) \geq 49/50,\\
\bm{u}_{l}^{\star\top}\bm{w}_{l} & =1-O\left(\frac{\kappa^{4}\sigma_{\max}^{2}n\log n}{(\lambda_{\mathrm{max}}^{\star})^{2}}+\frac{\mu\kappa^{4}r^{2}\sigma_{\max}^{2}\log n}{(\Delta_{l}^{\star})^{2}}\right) \geq 49/50,\\
\widehat{\bm{u}}^{\top}_l\bm{u}^{\star}_l & =1-O\left(\frac{\kappa^{4}\sigma_{\max}^{2}n\log n}{(\lambda_{\mathrm{max}}^{\star})^{2}}+\frac{\mu\kappa^{4}r^{2}\sigma_{\max}^{2}\log n}{(\Delta_{l}^{\star})^{2}}\right) \geq 49/50,
%
\end{cases}
\label{eq:ur-ul-uhat-rankr}
\end{equation}
\begin{align}
	\text{and}\qquad \|\bm{u}_{l}+\bm{w}_{l}\|_{2}=2+O\left(\frac{\kappa^{4}\sigma_{\max}^{2}n\log n}{(\lambda_{\mathrm{max}}^{\star})^{2}}+\frac{\mu\kappa^{4}r^{2}\sigma_{\max}^{2}\log n}{(\Delta_{l}^{\star})^{2}}\right) \geq 9/5.
\end{align}
In addition, if Assumption~\ref{assumption:noise-size-rankr} holds, then with probability at least $1-O(n^{-10})$, 
\begin{align}
	& \max\left\{ \|\ur_l - \bm{u}^{\star}_l\|_{\infty},\|\ul_l-\bm{u}^{\star}_l\|_{\infty},\|\widehat{\bm{u}}_l-\bm{u}^{\star}_l\|_{\infty}\right\} \nonumber\\
	& \qquad \lesssim \frac{\sigma_{\mathrm{max}}}{\lambda_{\mathrm{min}}^{\star}}\sqrt{\mu\kappa^{4}r\log n}+\frac{\sigma_{\mathrm{max}}}{\Delta_{l}^{\star}}\sqrt{\frac{\mu^{2}\kappa^{4}r^{3}\log n}{n}}
	\lesssim o\left(\frac{1}{\sqrt{\mu n\log n}}\right)
\label{eq:ur-ul-uhat-inf-rankr} .
\end{align}
\end{lemma}
\begin{proof} See Appendix \ref{sec:proof-lem-basic-uR-uL}.  \end{proof}

\subsection{Proof of Lemma \ref{lem:va-properties}}
\label{sec:proof-lem-va-properties}
Without loss of generality, assume that $\lambda_l^{\star}=1$ and that $\|\bm{a}\|_2=1$. The first claim \eqref{eq:expression-va} follows from direct calculations. Regarding the second claim \eqref{eq:va-bounds}, we first establish the lower bound 
\begin{align*}
v_{\bm{a},l}^{\star} & \geq\frac{1}{4}\sum_{1\leq i,j\leq n}(a_{i}u_{l,j}^{\star}+a_{j}u_{l,i}^{\star})^{2}\sigma_{\min}^{2}\\
 & =\frac{1}{4}\sum_{1\leq i,j\leq n}\Big(a_{i}^{2}(u_{l,j}^{\star})^{2}+a_{j}^{2}(u_{l,i}^{\star})^{2}+2a_{i}u_{l,i}^{\star}a_{j}u_{l,j}^{\star}\Big)\sigma_{\min}^{2}\\
 & =\frac{1}{2}\left( 1 +\big(\bm{a}^{\top}\bm{u}_l^{\star}\big)^{2}\right)\sigma_{\min}^{2},
\end{align*}
where we have used the fact that $\sum_{i,j}a_{i}^{2}(u_{l,j}^{\star})^{2}=\sum_i a_{i}^{2}\sum_j(u_{l,j}^{\star})^{2} = 1$. 
A similar argument  leads to the advertised upper bound $v_{\bm{a},l}^{\star}\leq\frac{1}{2}(1+(\bm{a}^{\top}\bm{u}_l^{\star})^{2})\sigma_{\max}^{2}$.


\subsection{Proofs of Lemma \ref{lem:aHu-UB} and Lemma~\ref{LemPropertyInprod}}
\label{sec:proof-lemma-PropertyInprod}

Without loss of generality, we shall assume $\lambda^{\star}_l=1$ throughout the proof. 
We shall also assume that $|H_{ij}|\leq B$ for all $1\leq i,j\leq n$.\footnote{Otherwise, we can first look at the truncated version $\widetilde{H}_{ij}:= H_{ij} \ind_{\{|H_{ij}|\leq B\}}$ (which is also zero-mean with variance obeying $\var[\widetilde{H}_{ij}]=(1+o(1))\sigma_{ij}^2$ under our assumptions), and then argue that $H_{ij}=\widetilde{H}_{ij}$ ($\forall i,j$) with high probability. This argument is fairly standard and is hence omitted for brevity.}   
By direct calculations, the following quantity of interest can be expressed as
\begin{align}
\label{eq:take-a-bow}
	\frac{\bm{a}^{\top}(\bm{H}+\bm{H}^\top) \ustar}{2}
	= \frac{1}{2} \sum_{i,j} (a_{i}u^\star_{j} + a_{j}u^\star_{l,i}) H_{ij} 
	=: \Tone. 
\end{align}
The above quantity is the sum of $n^2$ independent zero-mean random variables, each obeying
\begin{align}
	\var \left[(a_{i}u^\star_{l,j} + a_{j}u^\star_{l,i}) H_{ij}\right]
	&= (a_{i}u^\star_{l,j} + a_{j}u^\star_{l,i})^2 \sigma^2_{ij}; \label{eq:var-single} \\
	|(a_{i}u^\star_{l,j} + a_{j}u^\star_{l,i}) H_{ij}| &\leq |a_{i}u^\star_{l,j} + a_{j}u^\star_{l,i}|\, B.
\end{align}

\paragraph*{Proof of Lemma \ref{lem:aHu-UB}:} We shall only establish the lemma for the quantity $\bm{a}^{\top}(\bm{H}+\bm{H}^{\top})\bm{u}_l^{\star}$; the bounds on $\bm{a}^{\top}\bm{H}\bm{u}_l^{\star}$ and $\bm{a}^{\top}\bm{H}^{\top}\bm{u}_l^{\star}$ follow similarly. 
Invoking Bernstein's inequality for bounded random variables \cite{tropp2015introduction}, we guarantee that with probability at least $1- O(n^{-11})$,  
\begin{align*}
\Big|\sum_{i,j}(a_{i}u_{l,j}^{\star}+a_{j}u_{l,i}^{\star})H_{ij}\Big| & \lesssim\sqrt{\sum_{i,j}\var\left[(a_{i}u_{l,j}^{\star}+a_{i}u_{l,j}^{\star})H_{ij}\right]\log n}+\max_{i,j}\left\{ \left|(a_{i}u_{l,j}^{\star}+a_{j}u_{l,i}^{\star})B\right|\right\} \log n\\
 & \leq\sqrt{\sum_{i,j}(a_{i}u_{l,j}^{\star}+a_{j}u_{l,i}^{\star})^{2}\sigmamax^{2}\log n}+\max_{i,j}\left\{ \left|(a_{i}u_{l,j}^{\star}+a_{j}u_{l,i}^{\star})B\right|\right\} \log n\\
 & \lesssim\sqrt{\big(\sigmamax^{2}\log n\big)\sum_{i,j}\big\{(a_{i}u_{l,j}^{\star})^{2}+(a_{j}u_{l,i}^{\star})^{2}\big\}}+B\|\bm{u}_{l}^{\star}\|_{\infty}\|\bm{a}\|_{2}\log n\\
 & \asymp\sqrt{\big(\sigmamax^{2}\log n\big)\|\bm{a}\|_{2}^{2}\|\bm{u}_{l}^{\star}\|_{2}^{2}}+B\|\bm{u}_{l}^{\star}\|_{\infty}\|\bm{a}\|_{2}\log n.
\end{align*}
%
Now, as an immediate consequence of the incoherence condition and the identity $\ltwo{\bm{a}} = \ltwo{\bm{u}_l^{\star}} = 1$, we have
\begin{align*}
	\Big|\sum_{i,j} (a_{i}u^\star_{l,j} + a_{j}u^\star_{l,i}) H_{ij}\Big| 
	\lesssim
	\sigmamax\sqrt{\log n} 
	+ 
	B \log n \sqrt{\frac{\mu}{n}}.
\end{align*}
Combining this with the expression~\eqref{eq:take-a-bow} and the assumption $B\leq\sigma_{\max}\sqrt{\frac{n}{\mu\log n}}$, we complete the proof of the inequality~\eqref{eq:property-one-inprod} via the union bound.

\paragraph*{Proof of Lemma~\ref{LemPropertyInprod}:}
To establish this inequality, the key is to make use of the following Berry-Esseen type inequality \cite[Theorem 3.7]{chen2010normal}
.
\begin{theorem}
\label{thm:Berry-Esseen}
Let $\xi_1, \xi_2, \ldots, \xi_n$ be independent zero-mean random variables
satisfying $\sum_{i=1}^n \var[\xi_i]=1$. Then 
\begin{align}
\label{eq:awesome-berry-esseen}	
\sup_{z\in \real} \Big|\,\mprob\Big(\sum_{i=1}^n \xi_i \leq z\Big) - \Phi(z) \,\Big| \leq 10\gamma,
\qquad
	\text{where } \gamma \defn \sum_{i=1}^n \Exs [|\xi_i|^3].
\end{align}
\end{theorem}
In order to apply Theorem \ref{thm:Berry-Esseen}, let us define
\begin{align*}
	\xi_{ij} \defn \frac{1}{\sqrt{\var[\Tone]}} \frac{1}{2} (a_{i}u^\star_{j} + a_{i}u^\star_{j}) H_{ij} \quad \text{and} \quad W=\sum_{i,j}\xi_{ij} ,
\end{align*}
where it follows from the equality~\eqref{eq:var-single} that
\begin{align}
	v^{\star} := \var[\Tone] = \frac{1}{4} \sum_{i,j} (a_{i}u^\star_{j} + a_{j}u^\star_{i})^2 \sigma^2_{ij}.
\end{align}
By definition, it is easily seen that $\sum_{i,j} \var[\xi_{ij}]=1.$
Therefore, the property~\eqref{eq:awesome-berry-esseen} follows. In order to establish inequality~\eqref{eq:property-two-inprod}, it suffices to upper bound $\gamma$ in Theorem~\ref{thm:Berry-Esseen}.

To this end, we first make the observation that 
\begin{align*}
	\gamma = \sum_{i,j} \Exs [ |\xi_{ij}|^3 ]
	&= \frac{1}{\var[\Tone]^{3/2}}  
	\frac{1}{8} 
	\sum_{i,j} \Exs \left[ |(a_{i}u^\star_{j} + a_{i}u^\star_{j}) H_{ij}|^3 \right] \\
	&\leq \frac{1}{8 (v^{\star})^{3/2}}  
	\sum_{i,j} \Exs \left[\max_{i,j} |(a_{i}u^\star_{j} + a_{i}u^\star_{j}) H_{ij}|
	\cdot |(a_{i}u^\star_{j} + a_{i}u^\star_{j}) H_{ij}|^2 \right] \\
	&\leq \frac{1}{8 (v^{\star})^{3/2}}  
	\sum_{i,j} \Exs \left[\ltwo{\bm{a}} \|\bm{u}^\star\|_{\infty}  B
	\cdot |(a_{i}u^\star_{j} + a_{i}u^\star_{j}) H_{ij}|^2 \right] \\
	&\leq \frac{1}{8 (v^{\star})^{3/2}}  
	\ltwo{\bm{a}} \|\bm{u}^\star\|_{\infty}  B (4{v^{\star}}) = \frac{1}{2(v^{\star})^{1/2}}  
	\ltwo{\bm{a}} \|\bm{u}^\star\|_{\infty}  B .
\end{align*}
%
Re-organizing terms and using $\ltwo{\bm{a}} = \ltwo{\bm{u}^{\star}}= 1$ as well as the incoherence condition, 
we have 
\begin{align}
	\gamma \leq  B \sqrt{\frac{\mu}{4v^{\star} n}}. \label{eq:gamma-UB1}
\end{align}
It thus remains to lower bound $v^{\star}$. 
Again, use Lemma \ref{lem:va-properties} and the fact $\|\bm{a}\|_2=1$ to arrive at $4 v^{\star} \geq  2 \sigmamin^2$. 
Substitution into \eqref{eq:gamma-UB1} allows us to further control $\gamma$ as 
\begin{align*}
	\gamma \leq \frac{B}{\sigmamin}\sqrt{\frac{\mu}{ 2n}} \leq \frac{1}{\sqrt{2\log n}},
\end{align*}
where the last inequality arises since $B\leq \sigma_{\min}\sqrt{n/(\mu\log n)}$. Applying Theorem \ref{thm:Berry-Esseen} concludes the proof.


\subsection{Proof of Lemma \ref{lem:basic-uR-uL-rankr}}
\label{sec:proof-lem-basic-uR-uL}

To begin with, Theorem~\ref{thm:eigengap-condition} tells that that $\bm{u}_l$ and $\bm{w}_l$ are both real-valued  under Assumption~\ref{assumption:noise-size-rankr-revised}.  
It has also been shown in Theorem~\ref{thm:rankr-bounds} that 
\begin{align*}
\begin{cases}
\bm{u}_{l}^{\star\top}\bm{u}_{l} & =1-O\left(\frac{\kappa^{4}\sigma_{\max}^{2}n\log n}{(\lambda_{\mathrm{max}}^{\star})^{2}}+\frac{\mu\kappa^{4}r^{2}\sigma_{\max}^{2}\log n}{(\Delta_{l}^{\star})^{2}}\right) \geq 49/50 \\
\big|\bm{u}_{l}^{\star\top}\bm{w}_{l}\big| & =1-O\left(\frac{\kappa^{4}\sigma_{\max}^{2}n\log n}{(\lambda_{\mathrm{max}}^{\star})^{2}}+\frac{\mu\kappa^{4}r^{2}\sigma_{\max}^{2}\log n}{(\Delta_{l}^{\star})^{2}}\right) \geq 49/50
\end{cases}
\end{align*}
provided that $\bm{u}_{l}^{\star\top}\bm{u}_{l}>0$ and that Assumption \ref{assumption:noise-size-rankr-revised} holds.
This immediately indicates that (i) $\|\bm{u}_l-\ustar_l\|_2 = \sqrt{2-2\langle \bm{u}_l, \ustar_l \rangle} \leq 1/5$ and  $\|\bm{u}_l+\ustar_l\|_2=2\|\ustar_l\|_2 - \|\bm{u}_l-\bm{u}_l^{\star}\|_2 \geq 9/5$, and (ii) one either has $\ul^{\top}_l\bm{u}_l^{\star}\geq 49/50$
or $\ul^{\top}_l\bm{u}_l^{\star}\leq -49/50$. If $\ul^{\top}_l\bm{u}^{\star}_l \leq 49/50$,
then one necessarily has
\begin{align*}
\ul^{\top}_l\ur_l & =\ul^{\top}_l\bm{u}^{\star}_l+\ul^{\top}_l(\ur_l-\bm{u}_l^{\star}) \leq -49/50 + \|\ul_l\|_{2}\|\ur_l-\bm{u}^{\star}_l\|_{2}\\
 & \leq -49/50 + 1/5 <0,
\end{align*}
thus contradicting the assumption that $\ul_l^{\top}\ur_l>0$. As a result, one concludes that 
\begin{align}
 \bm{u}_{l}^{\star\top}\bm{w}_{l}   =1-O\left(\frac{\kappa^{4}\sigma_{\max}^{2}n\log n}{(\lambda_{\mathrm{max}}^{\star})^{2}}+\frac{\mu\kappa^{4}r^{2}\sigma_{\max}^{2}\log n}{(\Delta_{l}^{\star})^{2}}\right) \geq 49/50.
\end{align}
%

Further, it follows from Theorem~\ref{thm:rankr-bounds} that 
\begin{align*}
\min\big\{\left\Vert \bm{u}_{l}-\bm{u}_{l}^{\star}\right\Vert _{\infty},\left\Vert \bm{u}_{l}+\bm{u}_{l}^{\star}\right\Vert _{\infty}\big\} & \lesssim\frac{\sigma_{\mathrm{max}}}{\lambda_{\mathrm{min}}^{\star}}\sqrt{\mu\kappa^{4}r\log n}+\frac{\sigma_{\mathrm{max}}}{\Delta_{l}^{\star}}\sqrt{\frac{\mu^{2}\kappa^{4}r^{3}\log n}{n}}\\
 & \leq o\Big(\frac{1}{\sqrt{\mu n\log n}}\Big),
\end{align*}
where the last line holds under Assumption~\ref{assumption:noise-size-rankr}. 
Suppose that $\|\ur_l+\bm{u}_l^{\star}\|_{\infty}\leq\|\ur_l-\bm{u}_l^{\star}\|_{\infty}$.
Then the above inequality implies that
\[
\|\ur_l+\bm{u}_l^{\star}\|_{2}\leq\sqrt{n}\|\ur_l+\bm{u}_l^{\star}\|_{\infty} 
= \sqrt{n} \min\left\{ \|\ur_l-\bm{u}_l^{\star}\|_{\infty},\|\ur_l+\bm{u}_l^{\star}\|_{\infty}\right\} \ll 1,
\]
which is contradictory to the relation $\|\ur_l+\bm{u}_l^{\star}\|_{2}\geq 9/5$
shown above. Therefore, we must have
\begin{equation}
\|\ur_l-\bm{u}_l^{\star}\|_{\infty}=\min\left\{ \|\ur_l-\bm{u}_l^{\star}\|_{\infty},\|\ur_l+\bm{u}_l^{\star}\|_{\infty}\right\} \lesssim O\left(\frac{\sigma_{\max}\sqrt{\mu\log n}}{\big|\lambda^{\star}\big|}\right)
= o\Big(\frac{1}{\sqrt{\mu n\log n}}\Big) 
\end{equation}
as claimed. Similarly, we shall also have $\|\ur_l-\bm{u}_l^{\star}\|_{2} =\min\{\|\ur_l\pm\bm{u}_l^{\star}\|_{2}\}$ under Assumption~\ref{assumption:noise-size-rankr-revised}, which combined with  Theorem~\ref{thm:rankr-bounds} gives
\begin{align}
	\label{eq:L2-bound-ul-wl-sign}
\max\left\{ \|\bm{u}_{l}-\bm{u}_{l}^{\star}\|_{2},\|\bm{w}_{l}-\bm{u}_{l}^{\star}\|_{2}\right\} \lesssim O\left(\frac{\kappa^{4}\sigma_{\max}^{2}n\log n}{(\lambda_{\mathrm{max}}^{\star})^{2}}+\frac{\mu\kappa^{4}r^{2}\sigma_{\max}^{2}\log n}{(\Delta_{l}^{\star})^{2}}\right).
\end{align}

In addition, the above results further imply that
\begin{align}
\|\bm{u}_{l}+\bm{w}_{l}\|_{2}-2 & =\frac{\|\bm{u}_{l}+\bm{w}_{l}\|_{2}^{2}-4}{\|\bm{u}_{l}+\bm{w}_{l}\|_{2}+2}\overset{(\mathrm{i})}{\asymp}\|\bm{u}_{l}+\bm{w}_{l}\|_{2}^{2}-4 \nonumber\\
 & =\|\bm{u}_{l}\|_{2}^{2}+\|\bm{w}_{l}\|_{2}^{2}-2+2\{\langle\bm{u}_{l},\bm{w}_{l}\rangle-1\} \nonumber\\
 & \overset{(\mathrm{ii})}{=}2\{\langle\bm{u}_{l},\bm{w}_{l}\rangle-\langle\bm{u}_{l}^{\star},\bm{u}_{l}^{\star}\rangle\}=2\{\langle\bm{u}_{l}-\bm{u}_{l}^{\star},\bm{w}_{l}\rangle+\langle\bm{u}_{l}^{\star},\bm{w}_{l}-\bm{u}_{l}^{\star}\rangle\} \nonumber\\
 & \overset{(\mathrm{iii})}{=}O\left(\frac{\kappa^{4}\sigma_{\max}^{2}n\log n}{(\lambda_{\mathrm{max}}^{\star})^{2}}+\frac{\mu\kappa^{4}r^{2}\sigma_{\max}^{2}\log n}{(\Delta_{l}^{\star})^{2}}\right),
	\label{eq:ul-plus-wl-norm}
\end{align}
where (i) holds since $2\leq\|\bm{u}_{l}+\bm{w}_{l}\|_{2}+2\leq4$,
(ii) follows since $\|\bm{u}_{l}\|_{2}=\|\bm{w}_{l}\|_{2}=1$, and
(iii) comes from Cauchy-Schwarz and the $\ell_{2}$ bound \eqref{eq:L2-bound-ul-wl-sign}.

To finish up, it remains to show that the claimed bounds hold when $\bm{u}_l$ is replaced by $\widehat{\bm{u}}_l$. Regarding the bound on $\widehat{\bm{u}}_{l}^{\top}\bm{u}_{l}^{\star}$, we make the observation that
\begin{align*}
\widehat{\bm{u}}_{l}^{\top}\bm{u}_{l}^{\star} & =\frac{1}{2}(\bm{u}_{l}^{\top}\bm{u}_{l}^{\star}+\bm{w}_{l}^{\top}\bm{u}_{l}^{\star})+\widehat{\bm{u}}_{l}^{\top}\bm{u}_{l}^{\star}-\frac{1}{2}(\bm{u}_{l}+\bm{w}_{l})^{\top}\bm{u}_{l}^{\star}\\
 & =\frac{1}{2}(\bm{u}_{l}^{\top}\bm{u}_{l}^{\star}+\bm{w}_{l}^{\top}\bm{u}_{l}^{\star})+\left(\frac{1}{\|\bm{u}_{l}+\bm{w}_{l}\|_{2}}-\frac{1}{2}\right)(\bm{u}_{l}+\bm{w}_{l})^{\top}\bm{u}_{l}^{\star}\\
 & =1-O\left(\frac{\kappa^{4}\sigma_{\max}^{2}n\log n}{(\lambda_{\mathrm{max}}^{\star})^{2}}+\frac{\mu\kappa^{4}r^{2}\sigma_{\max}^{2}\log n}{(\Delta_{l}^{\star})^{2}}\right)+\frac{2-\|\bm{u}_{l}+\bm{w}_{l}\|_{2}}{2\|\bm{u}_{l}+\bm{w}_{l}\|_{2}}\cdot(\bm{u}_{l}+\bm{w}_{l})^{\top}\bm{u}_{l}^{\star}\\
 & =1-O\left(\frac{\kappa^{4}\sigma_{\max}^{2}n\log n}{(\lambda_{\mathrm{max}}^{\star})^{2}}+\frac{\mu\kappa^{4}r^{2}\sigma_{\max}^{2}\log n}{(\Delta_{l}^{\star})^{2}}\right),
\end{align*}
where the last line comes from \eqref{eq:ul-plus-wl-norm}. The advertised $\ell_{\infty}$ norm bounds for $\widehat{\bm{u}}_l$ follow from the same arguments; we omit it for brevity.

\section{Example: a symmetric case with homoscedastic Gaussian noise}
\label{sec:example-homoscedastic-Gaussian}

In this section, we isolate a simple example to illustrate the potential applicability
of our main results in the presence of certain symmetric noise matrices.
Specifically, consider the symmetric and homoscedastic case with Gaussian
noise, namely, 
\begin{subequations}
\label{eq:homoscedastic-Gaussian-model}
\begin{align}
H_{ij} & =H_{ji}\ \overset{\text{i.i.d.}}{\sim}\ \mathcal{N}(0,\sigma^{2}),\qquad1\leq i\leq j\leq n; \\ \label{eq:Gaussian-symmetric-noise}
\bm{M}^{\star} & =\sum_{l=1}^{r}\lambda_{l}^{\star}\bm{u}_{l}^{\star}\bm{u}_{l}^{\star\top}, \qquad \bm{M}=\bm{M}^{\star} + \bm{H}.  
\end{align}
\end{subequations}
Clearly, this setting differs from the model in Assumption~\ref{assumption-H} in that both $\bm{H}$
and $\bm{M}$ are now symmetric matrices. To invoke our theorems, we need to first asymmetrize the data matrix. 
Towards this, we propose a procedure as follows, motivated by a simple asymmetrization trick pointed out by \cite{chen2018asymmetry}. 
\begin{enumerate}
\item Generate $\widetilde{\bm{M}} = \bm{M} + \widetilde{\bm{H}}$, where $\widetilde{\bm{H}}$ is a skew-symmetric matrix obeying
\begin{align}
	\widetilde{H}_{ij} = - \widetilde{H}_{ji} \ \overset{\text{i.i.d.}}{\sim}\ \mathcal{N}(0,\sigma^{2}),\qquad1\leq i<j\leq n. 
\end{align}
\item Run Algorithm~\ref{alg:CI-au} with input data matrix $\widetilde{\bm{M}}$. 
\end{enumerate}
In short, the above scheme injects additional noise to the data matrix before implementing the inferential procedures. The theoretical support is this:   
\begin{corollary}
	Consider the model \eqref{eq:homoscedastic-Gaussian-model}. Then Theorems~\ref{thm:confidence-interval-validity-rankr}, \ref{thm:distribution-validity-rankr} and \ref{Lem-var-control-rankr} continue to hold, if Algorithms~\ref{alg:CI-au} and \ref{alg:CI-lambda} take $\widetilde{\bm{M}}$ as the input matrix.
\end{corollary}
\begin{proof}
	The key observation is that: under our construction, the quantities $H_{ij}+\widetilde{H}_{ij}$ and $H_{ij}-\widetilde{H}_{ij}$ are uncorrelated and hence independent. Therefore, the effective noise matrix $\bm{H}_{\mathsf{eff}}:=\bm{H}+\widetilde{\bm{H}}$ is in general asymmetric and contains i.i.d.~entries drawn from $\mathcal{N}(0,2\sigma^2)$. Given that
\begin{align}
	\bm{H}_{\mathsf{eff}}+\bm{H}_{\mathsf{eff}}^{\top} = \bm{H}+\bm{H}^{\top} + \widetilde{\bm{H}} + \widetilde{\bm{H}}^{\top} = \bm{H}+\bm{H}^{\top},
\end{align}
applying Theorems~\ref{thm:confidence-interval-validity-rankr}, \ref{thm:distribution-validity-rankr} and \ref{Lem-var-control-rankr} immediately concludes the proof. 
\end{proof}
Interestingly, the above findings indicate that: a simple asymmetrization trick via proper noise injection allows one to construct confidence intervals under symmetric and homoscedastic Gaussian noise. We caution that the above procedure requires prior knowledge on the noise level $\sigma$, which can often be reliably estimated in the homoscedastic case.

While noise injection does not affect the first-order uncertainty term $\bm{a}^{\top} (\bm{H}+\bm{H}^{\top}) \bm{u}_l^{\star}$, it does inflate the higher-order residual terms. For practical purposes, we recommend running the above procedure independently for multiple times and returning the ``average'' of them. This is summarized as follows, which might help improve practical performance.
\begin{enumerate}

\item Generate $K$ matrices $\bm{M}^{(k)}=\bm{M}+\bm{H}^{(k)}$ $(1\leq k\leq K$),
where $\{\bm{H}^{(k)}\}$ are independent skew-symmetric matrices
obeying
\begin{equation}
H_{ij}^{(k)}=-H_{ji}^{(k)}\ \overset{\text{i.i.d.}}{\sim}\ \mathcal{N}(0,\sigma^{2}),\qquad1\leq i<j\leq n.\label{eq:noise-skew-symmetric}
\end{equation}

\item For each $1\leq k\leq K$, run Algorithm~\ref{alg:CI-au} with input data matrix $\bm{M}^{(k)}$. We denote by $\lambda_l^{(k)}$ the $l$th eigenvalue of $\bm{M}^{(k)}$, $\bm{u}_l^{(k)}$ the associated right eigenvector, $\widehat{u}_{\bm{a},l}^{(k)}$ the resulting estimator for $\bm{a}^{\top}\bm{u}_l^{\star}$, 
respectively. 
Here, we calibrate the global signs of $\{\bm{u}_{l}^{(k)}\}$  by ensuring $\langle \bm{u}_{l}^{(k)}, \bm{u}_{l}^{(k+1)} \rangle \geq 0$ for all $1\leq k<K$. 
We also let $\widehat{v}_{\bm{a},l}^{(k)}$ and $\widehat{v}_{\lambda,l}^{(k)}$  represent the resulting variance estimators.  

\item Average the estimators
\begin{equation}
\widehat{u}_{\bm{a},l}=\frac{1}{K}\sum_{k=1}^{K}\widehat{u}_{\bm{a},l}^{(k)}\qquad\text{and}\qquad\lambda_{l}=\frac{1}{K}\sum_{k=1}^{K}\lambda_{l}^{(k)} .
\label{eq:ua-lambda-symmetric}
\end{equation}

\item For any prescribed coverage level $1-\alpha$, return the confidence intervals for $\bm{a}^{\top}\bm{u}_l^{\star}$ and $\lambda_l^{\star}$ as follows 
\begin{subequations}
\label{eq:CI-au-lambda-symmetric}
\begin{align}
	\mathsf{CI}_{1-\alpha}^{\bm{a}}  &:= \Big[\, \widehat{u}_{\bm{a},l}\pm\Phi^{-1}(1-\alpha/2)\sqrt{\widehat{v}_{\bm{a},l}^{(1)}} \,\Big]; \\
	\mathsf{CI}_{1-\alpha}^{\lambda} &:= \Big[\, \lambda_{l} \pm \Phi^{-1}(1-\alpha/2)\sqrt{\widehat{v}_{\lambda,l}^{(1)}}\,\Big]. 
\end{align}
\end{subequations}

\end{enumerate}

\section{Symmetrize or not? Some high-level interpretation}
\label{sec:interpretation-symmetrization}

The numerical experiments reported in Section~\ref{SecNumericals} reveal the potential benefit of $\mathsf{Spectral}\text{-}\mathsf{asym}$ compared to $\mathsf{Spectral}\text{-}\mathsf{sym}$. Here, we provide some informal interpretation. 

If we hope the eigenvector $\bm{u}_2$ to behave as a reliable estimate of $\bm{u}_2^{\star}$ (meaning that $|\bm{u}_2^{\top}\bm{u}_2^{\star}|=1-o(1)$), then we would necessarily require $|\bm{u}_2^{\top}\bm{u}_1^{\star}|=o(1)$.  To develop an understanding about the size of $\bm{u}_2^{\top}\bm{u}_1^{\star}$, we recall from Neumann's series that (see \cite{chen2018asymmetry} or Appendix~\ref{sec:preliminaries})
\begin{align}
\bm{u}_{1}^{\star\top}\bm{u}_{2}=\frac{\lambda_{1}^{\star}\left(\bm{u}_{1}^{\star\top}\bm{u}_{2}\right)}{\lambda_{2}}+\frac{\lambda_{2}^{\star}\left(\bm{u}_{2}^{\star\top}\bm{u}_{2}\right)}{\lambda_{2}}\sum_{s=1}^{\infty}\frac{\bm{u}_{1}^{\star\top}\bm{H}^{s}\bm{u}_{2}^{\star}}{\lambda_{2}^{s}}+\frac{\lambda_{1}^{\star}\left(\bm{u}_{1}^{\star\top}\bm{u}_{2}\right)}{\lambda_{2}}\sum_{s=1}^{\infty}\frac{\bm{u}_{1}^{\star\top}\bm{H}^{s}\bm{u}_{1}^{\star}}{\lambda_{2}^{s}}, 
\end{align}	
or equivalently,
\begin{align}
\bm{u}_{1}^{\star\top}\bm{u}_{2}=\frac{\lambda_{2}^{\star}\left(\bm{u}_{2}^{\star\top}\bm{u}_{2}\right)}{\lambda_{2}-\lambda_{1}^{\star}-\lambda_{1}^{\star}\sum_{s=1}^{\infty}\frac{\bm{u}_{1}^{\star\top}\bm{H}^{s}\bm{u}_{1}^{\star}}{\lambda_{2}^{s}}}\sum_{s=1}^{\infty}\frac{\bm{u}_{1}^{\star\top}\bm{H}^{s}\bm{u}_{2}^{\star}}{\lambda_{2}^{s}}. 
\label{eq:u1-u2-prod-expression}
\end{align}
Consequently, an important condition that allows one to control quantity $\bm{u}_{1}^{\star\top}\bm{u}_{2}$ is to ensure that each summand in expression \eqref{eq:u1-u2-prod-expression} is as small as possible. 

Towards this end, let us single out the second-order term $\bm{u}_{1}^{\star\top}\bm{H}^{2}\bm{u}_{2}^{\star}$, which suffices for us to develop some intuition (here, we assume that $\lambda_2 \approx \lambda_2^{\star}$ so that the denominator $\lambda_2^{2}$ does not affect the order of this term). 
\begin{itemize}
\item If $\bm{H}$ is asymmetric, then it has been shown in \cite{chen2018asymmetry} (see also Appendix~\ref{sec:preliminaries}) that 
\begin{align}
\label{eq:u1-u2-prod-example-asym}
\big|\bm{u}_{1}^{\star\top}\bm{H}^{2}\bm{u}_{2}^{\star}\big|\lesssim\sigma_{\max}^{2}\sqrt{\mu n}\log n	
\end{align}
with high probability, which is exceedingly small given the assumption that $\sigma_{\max} \sqrt{n\log n}\ll 1$. 
\item If $\bm{H}$ is replaced by $\widetilde{\bm{H}}:=\frac{1}{2}(\bm{H}+\bm{H}^{\top})$, then in general there is no guarantee that $\bm{u}_{1}^{\star\top}\widetilde{\bm{H}}^{2}\bm{u}_{2}^{\star}$ can be equally well-controlled.  Take the case \eqref{eq:u1-u2-example-rank2} and \eqref{eq:noise-variance-example-rank2} for example: straightforward calculations reveal
\begin{equation}
	\mathbb{E}\left[\bm{u}_{1}^{\star\top}\widetilde{\bm{H}}^{2}\bm{u}_{2}^{\star}\right]=\bm{u}_{1}^{\star\top}\mathsf{Var}\big(\tfrac{1}{2}(\bm{H}+\bm{H}^{\top})\big)\bm{u}_{2}^{\star}=\tfrac{1}{4} \sigma_{1}^{2}n \asymp \sigma_{\max}^2 n,
\end{equation}
which far exceeds the order \eqref{eq:u1-u2-prod-example-asym} with asymmetric data matrices. If this is the case, then one would have to require a better control of the multiplicative term $\frac{\lambda_{2}^{\star}}{\lambda_{2}-\lambda_{1}^{\star}-\lambda_{1}^{\star}\sum_{s=1}^{\infty}\frac{\bm{u}_{1}^{\star\top}\bm{H}^{s}\bm{u}_{1}^{\star}}{\lambda_{2}^{s}}}$ in  \eqref{eq:u1-u2-prod-expression}, which is often equivalent to imposing a more stringent separation condition on the eigenvalue pair $(\lambda_1^{\star},\lambda_2^{\star})$. 
\end{itemize}
The take-home message is this: $\mathsf{Spectral}\text{-}\mathsf{sym}$ might sometimes be suboptimal when dealing with heteroscedastic noise, particularly when the variance structure of the noise matrix is somewhat ``aligned'' with the true eigen-structure.  As a result, $\mathsf{Spectral}\text{-}\mathsf{sym}$ might not be as robust as $\mathsf{Spectral}\text{-}\mathsf{asym}$ when performing eigenvector estimation.

\vspace*{1cm}

\bibliographystyle{alpha}
\bibliography{bibfile_asymmetric}

\end{document}